\documentclass[a4paper, 11pt]{article}
\usepackage{header2}
\usepackage{permutation}
\usepackage{multirow}

\usepackage[font={small}]{caption}
\begin{document}

\title
{\makebox[0pt][c]{A Combinatorial Approach to Rauzy-type Dynamics II:}\\
\makebox[0pt][c]{the Labelling Method and a Second Proof}\\ of the KZB Classification Theorem}

\author{Quentin De Mourgues\\
%$^{\dagger, \ddagger}$\\
  \multicolumn{1}{p{.7\textwidth}}{\centering\emph{\small 
\rule{0pt}{14pt}%
LIPN, Universit\'e Paris 13\\
99, av.~J.-B.~Cl\'ement, 93430 Villetaneuse, France\\
{\tt quentin.demourgues@lipn.fr}}}
}

\maketitle

% !!!! latex abstract is buggy and breaks a lot of formatting in
% following text !!!!
% http://latex-community.org/forum/viewtopic.php?t=3996

\begin{center}
\begin{minipage}{.9\textwidth}
\noindent
\small
{\bf Abstract.}
Rauzy-type dynamics are group actions on a collection of combinatorial
objects. The first and best known example (the Rauzy dynamics) concerns an action on
permutations, associated to interval exchange transformations (IET)
for the Poincar\'e map on compact orientable translation surfaces. The
equivalence classes on the objects induced by the group action 
have been classified by Kontsevich and Zorich, and by Boissy through methods involving both combinatorics algebraic geometry, topology and dynamical systems.
Our precedent paper \cite{DS17} as well as the one of Fickenscher \cite{Fic16} proposed an ad hoc combinatorial proof of this classification. 

\qquad However, unlike those two previous combinatorial proofs, we develop in this paper a general method, called the labelling method, which allows one to classify Rauzy-type dynamics in a much more systematic way. We apply the method to the Rauzy dynamics and obtain a third combinatorial proof of the classification. The method is versatile and will be used to classify three other Rauzy-type dynamics in follow-up articles.

\qquad Another feature of this paper is to introduce an algorithmic method to work with the sign invariant of the Rauzy dynamics. With this method we can prove most of the identities appearing in the literature so far (\cite{KZ03},\cite{Del13} \cite{Boi13} \cite{DS17}...) in an automatic way. 
\end{minipage}
\end{center}

\newpage

\tableofcontents

\newpage

%%%%%%%%%%%%%%%%%%%%%%%%%%%%%%%%%%%%%%%%%%%%%%%%%%%%%%%

%-------------------------------------------------------
\section{Permutational diagram monoids and groups}\label{sec.algsett}
%Definitions}
%-------------------------------------------------------
The labelling method that we introduce in section 2 is a method to classify Rauzy-type dynamics. In this section we define what we mean by Rauzy-type dynamics.

Let $X_I$ be a set of labelled combinatorial objects, with elements
labelled from the set $I$. We shall identify it with the set
$[n]=\{1,2,\ldots,n\}$ (and use the shortcut
$X_n \equiv X_{[n]}$).\footnote{The use of $I$ instead of just $[n]$
  is a useful notation when considering substructures: if $x \in X_n$,
  and $x' \subseteq x$ w.r.t.\ some notion of inclusion, it may be
  conventient to say that $x' \in X_I$ for $I$ a suitable subset of
  $[n]$, instead that the canonical one.}
The symmetric group $\kS_n$ acts naturally on $X_n$, by producing the
object with permuted labels.

Vertex-labeled graphs (or digraphs, or hypergraphs) are a typical
example.  
Set partitions are a special case of hypergraphs (all vertices have
degree 1).
% , and each label being repeated exactly once).  
Matchings are a special case of partitions, in which all blocks have
size 2.
% (they can be seen as perfect matchings on the complete graph $\cK_{2n}$).  
Permutations $\s$ are a special case of matchings, in which each block
$\{i,j\}$ has $i \in \{1,2,\ldots,n\}$ and $j = \s(i)+n \in
\{n+1,n+2,\ldots,2n\}$.
% (they can be seen as perfect matchings
% on the complete bipartite graph $\cK_{n,n}$).

We will consider dynamics over spaces of this type, generated by
operators of a special form that we now introduce:
\begin{definition}[Monoid and group operators]
\label{def.monoidoperator}
We say that $A$ is a \emph{monoid operator} on set $X_n$, if, for a map $\aa: X_n \to \kS_n$, it consists of the map on $X_n$ defined by
\be
A(x) = \aa_{x} x
\ef,
\ee
where the action $\aa_x x$ is in the sense of the symmetric-group
action over $X_n$.  We say that $A$ is a \emph{group operator} if,
furthermore, $\aa_{A(x)}=\aa_x$.
\end{definition}
\noindent
Said informally, the function $\aa$ ``poses a question'' to the
structure $x$. The possible answers are different permutations, by which we act on $x$.
Note that is all our applications the set $Y_n=\{\aa_x| x\in X_n\}$ of all possible permutations has a 
much smaller cardinality than $X_n$ and $\kS_n$, i.e.\ very few
`answers' are possible. In the Rauzy case, $|X_n|=|\kS_n|=n!$
while $|Y_n|=n$. The asymptotic behaviour is similar ($|X_n|$ is at
least exponential in $n$, while $|Y_n|$ is linear) in all of our
applications.
\begin{figure}
\begin{center}\includegraphics[scale=.75]{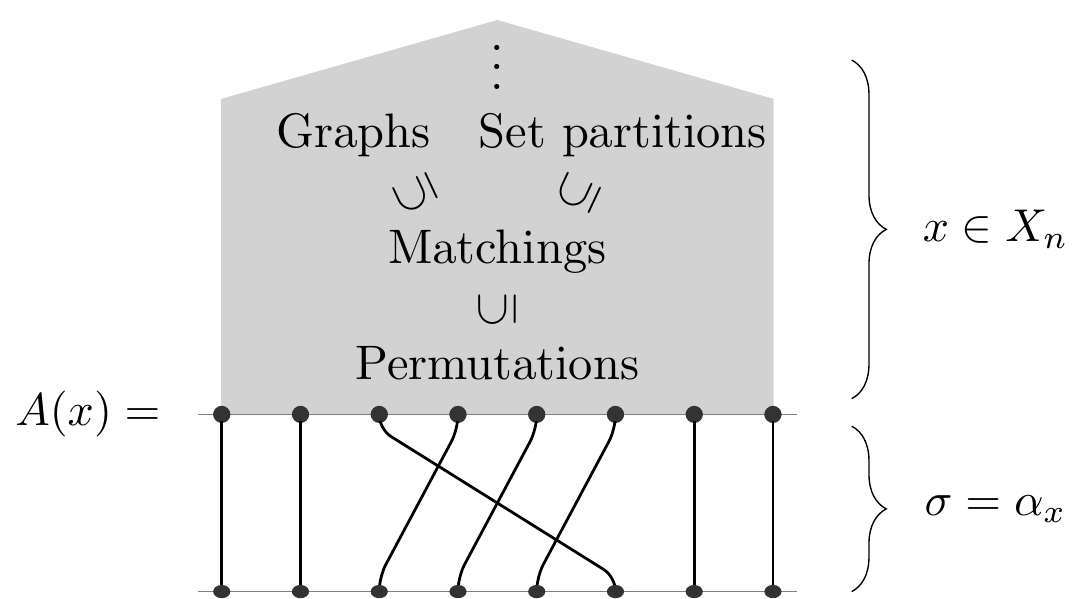}\end{center}
\caption{Let $x\in X_n$ be a combinatorial structure on $n$ vertices. $A_x$ choose a permutation $\s$ depending on $x$ and permutes the vertices of $x$ according to it. We represent this symmetric-group action as a diagrammatic action from below.}
\end{figure}
Clearly we have:

%% W.r.t.\ the generalised diagram operators outlined in the
%% introduction, this is clearly a special case, that may be called
%% \emph{permutational}.  

\begin{proposition}
\label{prop.Groupoperatorsareinvertible}
Group operators are invertible.
\end{proposition}

\pf For a given value of $n$, let $A$ be a group operator on the set
$X_n$.  The property $\aa_{A(x)}=\aa_x$ implies that, for all $k \in \bN$,
$A^k(x) = (\aa_{x})^k x$. Thus, for all $x$ there exists an integer
$d_A(x) \in \bN^+$ such that $A^{d_A(x)}(x)=x$. More precisely,
$d_A(x)$ is (a divisor of) the l.c.m.\ of the cycle-lengths of
$\aa_{x}$. Call
$d_A = \lcm_{x \in X_n} d_A(x)$ (i.e., more shortly, the l.c.m.\ over
$y \in Y_n$ of the cycle-lengths of $\alpha_y$). Then $d_A$ is a
finite integer, and we can pose $A^{-1} = A^{d_A-1}$. The reasonings
above show that $A$ is a bijection on $X_n$, and $A^{-1}$ is its
inverse.  \qed

\begin{definition}[monoid and group dynamics]
We call a \emph{monoid dynamics} the datum of a family of spaces
$\{X_n\}_{n \in \bN}$ as above, and a finite collection $\cA=\{A_{i}\}$ of
monoid operators.
We call a \emph{group dynamics} the analogous structure, in which all
$A_i$'s are group operators.

For a monoid dynamics on the datum $S_n = (X_n, \cA)$, we say that 
$x, x' \in X_n$ are \emph{strongly connected}, $x \sim x'$, if there
exist words $w, w' \in \cA^*$ such that $w x = x'$ and $w' x' = x$.

For a group dynamics on the datum $S_n = (X_n, \cA)$, we say that $x,
x' \in X_n$ are \emph{connected}, $x \sim x'$, if there exists a word
$w \in \cA^*$ such that $w x = x'$.
% In the two cases, $\sim$ is clearly an equivalence relation.
\end{definition}
\noindent
Here the action $wx$ is in the sense of monoid action.  Being
connected is clearly an equivalence relation, and coincides with the
relation of being graph-connected on the Cayley Graph associated to
the dynamics, i.e.\ the digraph with vertices in $X_n$, and edges $x
\longleftrightarrow_{i} x'$ if $A_i^{\pm 1} x = x'$. An analogous
statement holds for strong-connectivity, and the associated Cayley
Digraph.

We call such dynamics Rauzy-type dynamics since the original Rauzy dynamics that inspired this definition is also of this type. 

%% \footnote{I.e.\ the graph with vertex-set identified
%%   with $X_n$, and edges which are labeled from $\cA \cup \cA^{-1}$ and
%%   oriented: we have the oriented edge $(x,x')$ with label $A$ if
%%   $A(x)=x^{-1}$, and in this case we also have the oriented edge
%%   $(x',x)$ with label $A^{-1}$.}

\begin{definition}[classes of configurations]
Given a dynamics as above, and $x \in X_n$, we define 
$C(x) \subseteq X_n$, the \emph{class of $x$}, as the set of
configurations connected to $x$, $C(x)=\{ x' \,:\, x\sim x'\}$.
\end{definition}

We will call those classes the Rauzy classes of the dynamics.

\begin{definition}[Invariant of classes]
Given a dynamics as above, we say that $f: X_n \to G_n$ is an invariant of the dynamics if $f(x)=f(x')$ for every pair $(x,x')$ of connected combinatorial structures.
\end{definition}

Thus given a dynamics $(X_n, \cA)$ our goal is to find a family of invariant $(f_i)$ such that two combinatorial structures are in the same Rauzy class if and only if they have the same invariants.

\section{The labelling method}

The labelling method propose a sequence of steps to classify a dynamics $(X_n,(A_k))$ for which we already know the set of invariants. That is to say, we wish to prove that the invariants we identified completely characterize the Rauzy classes.
\label{sec_labelling_method}

\subsection{Definition and application}

We place ourself in the framework of section \ref{sec.algsett}. Let
$X_n$ be a combinatorial set and $(A_k)$ be the set of operators of
the dynamics.  Before introducing the labelling method, let us lay
down a number of necessary definitions.

\begin{definition}[$(k,r)$-coloring and reduction]\label{def_coloring}
Let $x\in X_n$ be a combinatorial structure, and $k$, $r$ integers
with $k+r=n$.  A \emph{$(k,r)$-coloring} $c$ of $x$ is a coloring of
the $n$ vertices of $x$ into a black set of $k$ vertices and a gray
set of $r$ vertices, such that, calling $y$ the restriction of $x$ to
the black vertices, $y$ is also a combinatorial structure (thus $y \in
X_k$). We call $y$ the \emph{reduction} of $x$.
\end{definition}

\noindent
For example, for a matching $m\in \cM_n=X_{2n}$, in principle we shall
color in black and gray the $2n$ points. However, in order for $y$ to
be a matching, in a valid $(k,r)$-coloring we need that the endpoints
of an edge of the matching are either both black or both gray (in
particular $r$ and $k$ are even).

Our first important notion is that of a boosted sequence. 

\begin{definition}[Boosted sequence and boosted dynamics]
Let $x\in X_n$ and let $y$ be the reduction of a $(k,r)$-coloring $c$
of $x$. Let $S$ be a sequence of operators in $X_k$ such that $S(y)=
y'$.

If a sequence $S'$ in $X_n$ is such that $y'$ is the reduction of
$(x',c')=S'(x,c)$, we say that $S'$ is a \emph{boosted sequence of $S$
  for $(x,c)$}.  See figure~\ref{fig_Boosted_sequence}.

\begin{figure}[tb!]
\begin{center}
%% $\begin{tikzcd}
%%   (x,c) \arrow{d}{red} \arrow{r}{S'} &  (x',c') \arrow{d}{red}   \\
%%  y \arrow{r}{S} &    y' \\
%% \end{tikzcd}
%% \quad\quad
%% \begin{tikzcd}
%%   (x,c) \arrow{d}{red} \arrow{r}{S''} &  (x'',c'') \arrow{d}{red}   \\
%%  y \arrow{r}{S} &    y' \\
%% \end{tikzcd}$
%% \\
$
\begin{tikzcd}
&& (x'',c'')
\arrow[ldd, bend left,  "\textrm{red}"]
\\
  (x,c) \arrow{d}{\textrm{red}} 
\arrow[urr, bend left, "S''"]
\arrow{r}{S'} &  (x',c') \arrow{d}{\textrm{red}} &  \\
 y \arrow{r}{S} &    y' & \\
\end{tikzcd}$
\end{center}
\caption{\label{fig_Boosted_sequence} $S'$ is a boosted sequence of
  $S$ for $(x,c)$, because $\textrm{red} \circ S'$ gives the same
  result as $S \circ \textrm{red}$. Note that $S'$ may be not unique:
  in the diagram we show a second boosted sequence of $S$ for $(x,c)$,
  namely $S''$. It is not necessarily the case that $(x',c') :=
  S'(x,c)$ coincides with $(x'',c'') := S''(x,c)$, but merely that
  their reductions coincide (they must be both $y'$).}
\end{figure}
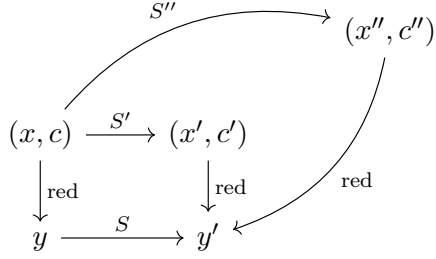

A dynamics $((X_n)_n,(A_k)_k)$ \emph{has a boosted dynamics} if for
every $n$, $x\in X_n$, $(k,r)$-coloring $c$ of $x$ and
sequence $S$ as above,
% $n$, for every $x\in X_n$, for every $(k,r)$-coloring $c$ of $x$ and
% for every sequence of operators $S$ in $X_k$ 
there exists a boosted sequence $S'$.
% of $S$ for $(x,c)$.
\end{definition}

\noindent Thus in the following diagram: $\begin{tikzcd}[column sep=15pt] 
(x,c) \arrow{d}[swap]{red} \\
y \arrow{r}{S}& y'
\end{tikzcd}$, a boosted sequence $S'$ is a lift of $S$ such that the square becomes commutative: $\begin{tikzcd}[column sep=15pt] 
(x,c) \arrow{d}[swap]{red} \arrow{r}{S'} & (x',c') \arrow{d}[swap]{red} \\
y \arrow{r}{S}& y'
\end{tikzcd}$.

Note that, for a given $(n,x,c,S)$ as above, 
% in the definition of the boosted dynamics, 
the boosted sequence $S'$ is not necessarily unique. This is rather
obvious when the dynamics is a group dynamics, as in this case, even
without boosting, the set of sequences $S: x_1 \to x_2$ is a coset
(w.r.t.\ sequence concatenation) within the group of all possible
sequences, i.e., if $S x_1 = x_2$, every sequence $S': x_2 \to x_2$ is
such that $(S')^k S x_1 = x_2$ for all $k \in bZ$.
% This was actually already the case in the Rauzy dynamics $\cS_n$. 
Nonetheless, with abuse of notation, in our definition of the boosted
dynamics we will often refer to `the' boosted sequence $S'$, by this
referring to the most natural boosted sequence, w.r.t.\ some notion
changing from case to case and depending on the invariants.

As we said in the introduction, we represent combinatorial structures
in $X_n$ as $n$ vertices on the real line (indicized from left to
right) with the combinatorial structure (be it matching, set
partition, graph or hypergraph) placed above in the upper half
plane. This is done in order to allow the dynamics operator to act
graphically from below (they will be represented as permutations,
contained within a horizontal strip), as is customarily the case for
diagram algebras (like the partition algebra, the braid group, and so
on)

This graphical representation is also useful in order to produce a
certain construction on configurations, which we now describe.  Place
two new vertices, with label $0$ and $n+1$, to the left and right of
the existing vertices $1, \ldots, n$, respectively. Then, consider the intervals between the pairs of adjacent vertices $(i,i+1)$, for $i=0,\ldots,n$. This modification of the structure $x
\in X_n$ will be called a \emph{combinatorial structure with intervals}.

If we have a $(k,r)$ coloring of $x$, in the construction of the
combinatorial structure with intervals we will set the two extra points,
$0$ and $n+1$, as black.

\begin{figure}[t]
\begin{center}
\includegraphics[scale=.5]{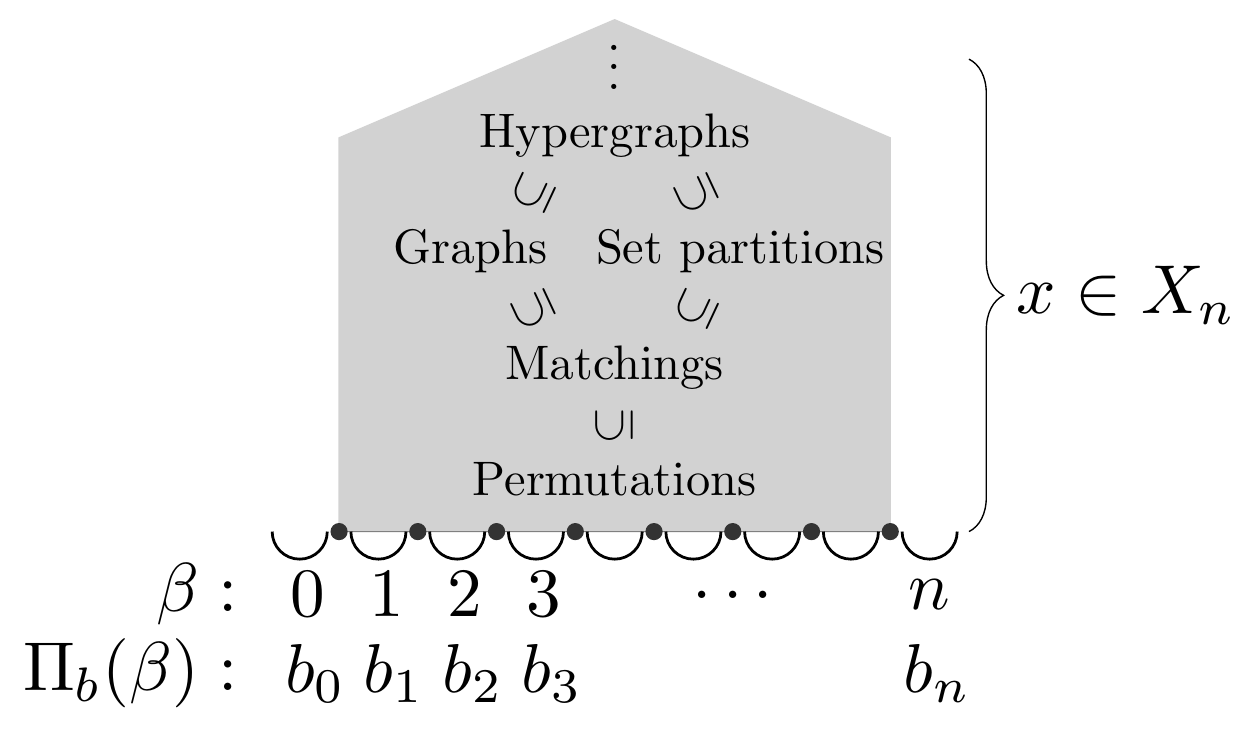}
\end{center}
\caption{\label{fig_arc_comb_structure}A representation of a combinatorial structure $x\in X_n$ with a labelling $\Pi_b.$ The intervals are represented by the bottoms arcs.}
\end{figure}

Let $\Sigma=\{b_0 \ldots,b_n\}$ be an alphabet of distinct symbols. We
define a \emph{labelling of $x\in X_n$} to be a bijection
\[\begin{array}{ccccc}
\Pi_b & : & \{0,\ldots,n\} & \to & \Sigma
\end{array}\]
i.e., a labelling of the intervals between the vertex of $x$ in the construction above,
with the symbols from $\Sigma$. See figure
\ref{fig_arc_comb_structure}.

In the following treatment, once $\Pi_b$ is given, it will be
convenient to designate intervals either by their positions $\beta \in
\{0,\ldots,n\}$ or by their labels $b\in \Sigma$, depending from the
situation. As a convention, and in order to avoid confusion, we will
use greek and latin letters as above in the two cases.

% Obviously given a label we obtain the position by applying $\Pi_b^{-1}$ and inversely.

We need a definition for a consistent action of the dynamics on the
labelling. Recall that in section \ref{sec.algsett} we have defined an
operator $A$ of the dynamics as the action of an application $\aa: X_n \to \kS_n$, so that 
$A x = \aa_{x} x$ is the symmetric-group action of $\aa_x$ on $x$.

\begin{definition}[Labelling and dynamics]
For $x \in X_n$, let $(x,\Pi_b)$ be a combinatorial structure with
labelling. The operator $\hat{A}$ is the \emph{labelled extension} of
$A$ if it acts on $(x,\Pi_b)$ as follows:
% We extend the action of the dynamics to the labelled structure by defining the action of each operator $A$ on $(x,\Pi_b)$: 
\[
\hat{A}(x,\Pi_b)=(\aa_{x}(x),\Pi_b \circ \alpha'_{x})
\]
where, in analogy with the Definition \ref{def.monoidoperator}
$\aa: X_n \to \kS_n$ $\aa': X_{n} \to \kS_{n+1}$. 
\end{definition}

Thus, in this case, not only $A$ chooses a permutation $\aa_x$,
depending on $x$, by which it acts on the $n$ vertices of $x$, but it
also chooses a permutation $\aa'_{x}$, again depending only on $x$,
by which it acts on the set of $n+1$ labels, in both cases in the sense
of the symmetric-group action.

When there will be no confusion, we will often write $\hat{A}$ for $A$.

Such a definition is not unique, the choice will be made so as to verify the following property:

\begin{definition}[Labelling, vertices, and boosted dynamics]
Let $x\in X_n$, $c$ a $(k,r)$-coloring of $x$, $y\in X_k$ the
corresponding reduction, and $\Pi_b$ a labelling of $y$.

We say that a gray vertex $\ell$ of $x$ is \emph{within an interval}
$\beta$ of $y$, and write $\ell \in b$, if the black vertices of $x$,
corresponding to the vertices $\beta$ and $\beta+1$ of $y$, are the
first black vertices to the left and to the right of $\ell$ (note that
$\beta$ and $\beta+1$ may be $0$ and $k+1$, i.e.\ the external
vertices added in the interval construction).
% i.e. they respectively maximize $\{i<\ell\,|\, i \text{ black}\}$ and minimize $\{i>\ell\,|\,i \text{ black}\}$.  
We say that the labelling is \emph{compatible with the boosted
  dynamics} if the following statement holds:

Let $x$, $c$, $y$ and $\Pi_b$ as above, $S$ a sequence in $X_k$, and
$S'$ a boosted sequence of $S$ for $(x,c)$.  If $v$ is a gray vertex
of $x$ within the interval $b\in\Sigma$ of $(y,\Pi_b)$, then the image of
$v$ in $(x',c')=S'(x,c)$ must be within the image of the interval $b$ in
$(y',\Pi'_b)=S(y,\Pi_b)$. In other word $v$ is within the interval with position $\Pi'^{-1}_b(b)$ in $(y',\Pi'_b)$.
\end{definition}

\noindent
In other words, if the labelling is compatible with the boosted
dynamics, we can keep track of the positions of the gray vertices by
understanding how the labelling of a reduced combinatorial structure
$y \in X_k$ evolves in the extended labelled dynamics. 

\begin{definition}\label{def_monomprob}
Given a combinatorial structure $x\in X_n$ with a labelling $\Pi_b$,
let $L(x,\Pi_b)$ be defined as
\be
L(x,\Pi_b)=\{\Pi'_b\ |\ \exists S \textrm{~with~}
S(x,\Pi_b)=(x,\Pi'_b)\}
\ef,
\ee
that is, the set of all the labellings that are reachable from
$(x,\Pi_b)$ by a loop $S$ of the dynamics.

We say that the \emph{$r$-point monodromy of $x$ is known} if the
following problem is solved: for every ordered $r$-uple of intervals
$(b_1,\ldots,b_r)$ in the structure with intervals $(x,\Pi_b)$, and every
$r$-uple $(\beta_1,\ldots,\beta_r)$ of disctinct integers in
$\{0,\ldots,n\}$, we know whether there exists a $\Pi'_b \in
L(x,\Pi_b)$, and thus an $S$ such that $S(x,\Pi_b)=(x,\Pi'_b)$, with
the property that $S(b_a)=\beta_a$ for all $a=1,\ldots,r$.
\end{definition}

\noindent
The definition above is somewhat redundant. In fact it is clear that,
for $\pi \in \kS_{n+1}$,
\be
L(x,\Pi_b \circ \pi)
=
L(x,\Pi_b) \, \pi
\ee
(and in particular, by choosing $\pi=\Pi_b^{-1}$, it suffices to
describe the set $L$ for the canonical left-to-right labelling). Also,
if we have a group dynamics, the dependence from $x$ is only through
the class of $x$, more precisely, if $S_0$ is a sequence in the labelled
dynamics from $(x,\Pi_b)$ to $(x',\Pi'_b)$, then
\be
L(x,\Pi_b)
=
L(x',\Pi'_b)
\ee
because, for each $\Pi''_b \in L(x',\Pi'_b)$ realised with the sequence
$S$, then $\Pi''_b \in L(x,\Pi_b)$ is realised with the sequence
$S_0^{-1} S S_0$, and similarly with $x \leftrightarrow x'$ and $S_0
\leftrightarrow S_0^{-1}$.

The interest in these properties is highlighted in the crucial
definition below, of the `labelling method' as a whole:

\begin{definition}\label{def_amenable_lb}
 A dynamics $((X_n)_n,(A_k)_k)$ is \emph{amenable to the labelling method} if:
\begin{enumerate}
\item It has a boosted dynamics.
\item There is a labelling compatible with the boosted dynamics.
\item For a suitable value of $r>0$, the $r$-point monodromy of any class is known.
%% Given a combinatorial structure $x\in X_n$ with a labelling $\Pi_b$, the set 
%% \[
%% L(x,\Pi_b)=\{\Pi'_b\ |\ \exists S \textrm{~with~} S(x,\Pi_b)=(x,\Pi'_b)\}
%% \] 
%% of all the labellings that are reachable from $(x,\Pi_b)$ by a loop
%% $S$, is known (i.e., we have a formula for the boolean function
%% $\Pi'_b \stackrel{?}{\in} L(x,\Pi_b)$).
%% % understandable (i.e. we can prove a quantifying theorem).
\end{enumerate}
\end{definition}
\noindent
When proving a classification theorem for a group
dynamics,\footnote{In case of a monoid dynamics the goal needs to be
  modified accordingly, and there is no general recipe. If the Cayley
  digraphs of the classes has a strongly connected component, we can
  restrict the study to this portion of the classes, with small
  difference w.r.t.\ the group case. Another possible sensible approach
  consists in modifying the monoid dynamics into a group dynamics, by
  adding new or modifying existant operators, so that the weak-connectedness in the original
  Cayley digraph of the classes turns into ordinary connectedness in
  the new Cayley graph.}
we shall usually proceed in two steps: first, guessing the `right' set
$\mathrm{Inv}_n$ of invariants for the dynamics on $X_n$, and then,
proving that the set is right. In order to do the latter, i.e.\ to
prove that there exists exactly one class per invariant, we need to
show two things:
\begin{description}
\item[existence] For every $\mathrm{inv}\in \mathrm{Inv}_n$ there
  exists a combinatorial structure $x\in X_n$ with
  invariant~$\mathrm{inv}$.
\item[completeness] The invariant discriminates the classes, i.e.\ 
for every pair $x_1,x_2\in X_n$ with the same invariant
$\mathrm{inv}\in \mathrm{Inv}_n$ there exists a sequence $S$ such that
$S(x_1)=x_2$.
\end{description} 
The first item of the list is often achieved relatively easily, by
constructing a candidate $x\in X_n$ either directly or by
induction. The difficulty of the construction depends mainly on how
tractable and explicit the invariant set is.

The second part, which is sensibly more complicated, is where the
labelling method comes forth.

The method is aimed at constructing an inductive step, so we assume
that the invariant is complete up to sizes $n'<n$, and that the
$r$-point monodromy problem (with an appropriate value of $r$
specified below) is solved for all configurations up to sizes $n'<n$,
and we analyse how we can prove that the invariant is complete at size
$n$. At a later moment, we shall solve the $r$-point monodromy problem (with the same value of $r$) at
size $n$, in light of the $r$-point monodromy problem at smaller size,
and the completeness at size up to $n$. The proof strategy will change from dynamics to dynamics but we will indicate some elements of proof that remain mostly identical in the next section.

So, we have two arbitrary $x_1$, $x_2\in X_n$, combinatorial structure
with the same invariant $\mathrm{inv}\in \mathrm{Inv}_n$, and we want
to show that they are connected by the dynamics.

%
%The first ingredient that we shall use is the notion of
%$(k,r)$-colorings and reductions. In certain circumstances (and surely
%when the $r$-point monodromy problem is solved) we can interpret
%$(k,r)$-colorings and reductions the other way around, i.e.\ as ways
%of embedding a class $C_k$ at size $k$ (or a portion of it) into a
%class $C_n$ at size $n$. Different colourings may correspond to
%different embeddings, so the same (portion of) class $C_k$ will appear
%inside $C_n$ with several copies (and it may appear, again with
%several copies, within other classes $C'_n$, while also other classes
%$C'_k$ may appear in $C_n$, and some portions of $C_n$ may have no
%intelligible reductions to classes at smaller size).

\noindent
\rule{0pt}{14pt}% 
Our first goal is {\bf find $x'_1$ and $x'_2$ in $X_n$,} connected to
$x_1$ and $x_2$ respectively, with the following property: there
exists two $(k,r)$-coloring\footnote{We shall adopt $r$ to be the
  smallest value such that $X_k$ is a valid combinatorial set, and in
  most of our applications it is just $r=1$ or $r=2$.}
$c_1$ and $c_2$ of $x'_1$ and $x'_2$ respectively, such that
$y_1:=\textrm{red}(x'_1,c_1)$ 
and $y_2:=\textrm{red}(x'_2,c_2)$ 
% $y_2$ the reductions of $(x'_1,c_1)$ and $(x'_2,c_2)$ 
have the same invariant $\mathrm{inv}'\in \mathrm{Inv}_k$. Which, by
induction hypothesis, implies that they are in the same class.
% That is, $x'_1$ and $x'_2$ are in (possibly different) embeddings of the
%same class.  
Call $S_1$ and $S_2$ the sequences such that $x'_1 = S_1
x_1$ and $x'_2 = S_2 x_2$.  \\
We call $x_1'$ and $x_2'$ normal forms as their structure will often be rather specific in order to prove the required property. Thus we can rename our first goal as 'finding normal forms for every pair $(x_1,x_2)$' or 'normalizing $x_1$ and $x_2$'.

\begin{figure}[t]
\begin{center}
$\begin{tikzcd}
\textrm{size $n$} 
&
x_1 \arrow{r}{S_1} 
& x'_1 \arrow[d,dashrightarrow]
\arrow{rr}{S'_* S'} 
&& x'_2  \arrow[d,dashrightarrow]
& x_2 \arrow[swap]{l}{S_2} 
\\
\makebox[0pt][c]{\raisebox{6.5pt}{\textrm{size $n$,}}}%
\makebox[0pt][c]{\raisebox{-6.5pt}{\textrm{$(k,r)$-colored}}}
% \makebox[0pt][c]{\raisebox{-6.5pt}{\textrm{$(n-r,r)$-colored}}}
&
& (x'_1,c_1) \arrow{d}{\textrm{red}} \arrow{r}{S'} 
& (x_3,c_3)  \arrow{d}{\textrm{red}} \arrow{r}{S'_*} 
& (x'_2,c_2) \arrow{d}{\textrm{red}} 
% & x_2  \arrow[swap]{l}{S_2} 
\\
\textrm{size $n-r$}
&
& y_1 \arrow{r}{S} 
\arrow[d,dashrightarrow]
& y_2 \arrow{r}{S_*} 
\arrow[d,dashrightarrow]
& y_2 
\arrow[d,dashrightarrow]
\\
\makebox[0pt][c]{\raisebox{6.5pt}{\textrm{size $n-r$,}}}%
\makebox[0pt][c]{\raisebox{-6.5pt}{\textrm{labelled}}}
&
& (y_1,\Pi_b) \arrow{r}{\hat{S}} 
& (y_2,\Pi'_b) \arrow{r}{\hat{S}_*}&
(y_2,\Pi''_b)
\\
\end{tikzcd}$
\end{center}
\caption{\label{fig_labelling_method} Outline of the proof of
  connectivity between $x_1$ and $x_2$, using the labelling
  method. The sequence $S$ sends $y_1$ to $y_2$, however the intervals
  containing the gray vertices of $x'_1$ may not be at their correct
  place, in order to match with those of $x'_2$. The sequence $S_*$
  corrects for this. Thus the boosted sequence $S'_* S'$ sends $x'_1$
  to $x'_2$, and $x_1$ and $x_2$ are connected.}
\end{figure}
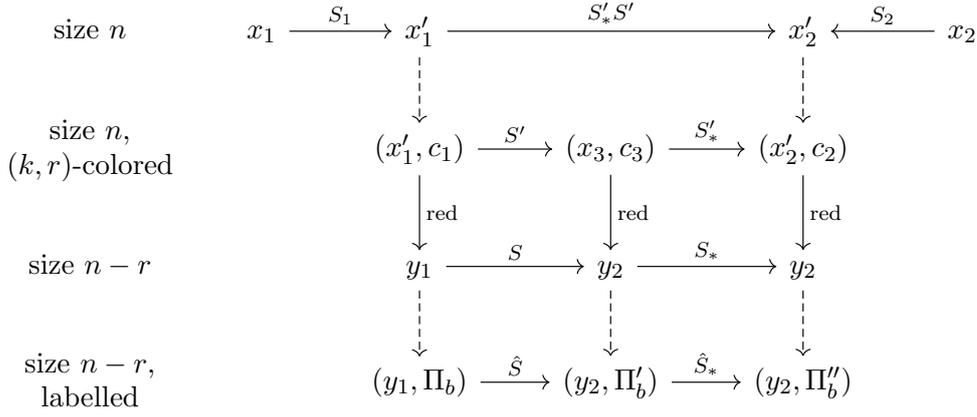

\noindent
\rule{0pt}{14pt}% 
Then, use induction to {\bf find a sequence $S$ such that
  $S(y_1)=y_2$}.
% , since $y_1$ and $y_2$ have the same invariant and thus are in the same class.
\\
\noindent
\rule{0pt}{14pt}% 
Since the dynamics is amenable to the labelling method, there exists a
boosted dynamics and a labelling compatible with it. {\bf Let us choose
such a labelling.}
\\
\noindent
\rule{0pt}{14pt}% 
Let $\Pi_b$ be a labelling of $y_1$. Say that the $r$ gray vertices of
$(x'_1,c_1)$ are within the intervals with labels $b^1,\ldots,b^r$ of
$(y_1,\Pi_b)$, in this order.  In addition, the $r$ gray vertices of
$(x'_2,c_2)$ are within the intervals with position, say, $\beta^1,
\ldots,\beta^r$ of $y_2$, in this order. {\bf We will be concerned
  with these two $r$-uples.}
% labels $b^{\prime 1}$ and $b^{\prime 2}$ of $(y_2,\Pi'_b)$. Let $\alpha_1=\Pi_b^{\prime-1}(b^{\prime 1})$ and $\alpha_2=\Pi_b^{\prime-1}(b^{\prime 2})$ be the positions of this two arcs on $y_2.$
\\
\noindent
\rule{0pt}{14pt}% 
By the boosted dynamics, there exists a sequence $S'$ which is {\bf a
  boosted sequence of $S$} for $(x'_1,c_1)$. Let
$(x_3,c_3)=S'(x'_1,c_1)$. If it were $(x_3,c_3)=(x'_2,c_2)$, we would
be done, but, this needs not to be the case, 
%and $(x_3,c_3)$ may be contained in yet another embedded
%copy
%  (the third one besides the ones for $x'_1$ and for $x'_2$) 
%of the class $\textrm{inv}'$ at size $k$.
  So, we know that $(x_3,c_3)$ and $(x'_2,c_2)$ are not equal in general, nonetheless the
reduction of both of them is $y_2$.  \\
\noindent
\rule{0pt}{14pt}% 
By compatibility of the labelling with the boosted dynamics, {\bf we
  know the location of the $r$ gray vertices of $(x_3,c_3)$}. Namely,
by virtue of our notation (in which the alphabet $\Pi_b$ changes
alongside the dynamics), we have that the $r$ gray vertices are within
the intervals with labels $b^1$, \ldots, $b^r$ of
$(y_2,\Pi'_b)=S(y_1,\Pi_b)$.  Thus $x_3$ and $x'_2$ possibly differ by
the position of the $r$ vertices which are gray in the colouring, that
are within the intervals with position $(\Pi'_b)^{-1}(b^a)$ and $\beta^a$
of $y_2$, for $x_3$ and $x'_2$, respectively, for $1 \leq a \leq r$.
\\
\noindent
\rule{0pt}{14pt}% 
Since the $r$-point monodromy problem is solved for $y_2$, {\bf we know
whether there exists a sequence $S_*$ that sends the label $b^a$ to the positions $\beta^a$}. More precisely, a sequence $S_* \in L(y_2,\Pi'_b)$ such that the
labeling $\Pi''_b=S_*(\Pi'_b)$ has $(\Pi''_b)^{-1}(b^a)=\beta_a$ for
all $1 \leq a \leq r$.
\\
\noindent
\rule{0pt}{14pt}% 
{\bf If $S_*$ exists,} let $S'_*$ be the boosted sequence of $S_*$ for $(x_3,c_3)$. Then,
again by compatibility of the labelling with the boosted dynamics, we have $S'_*(x_3,c_3)=(x'_2,c_2)$. 
Thus we have shown that $x_1$ and $x_2$ are connected, namely
\be
x_2 = (S_2)^{-1} S'_* S' S_1 x_1
\ee
See figure \ref{fig_labelling_method}
\\
\noindent
\rule{0pt}{14pt}% 
Conversely, {\bf if $S^*$ does not exist, the labelling method has failed.}

\subsection{Difficulties and implementation\label{sec_diff_lab}}

Let us discuss the possible difficulties one can encounter when trying
to employ the labelling method.  We can identify two main tasks:
proving that the dynamics is amenable to the method, and then applying
it to the `proof that $x_1 \sim x_2$' outlined above. 
More precisely, those two tasks are divided into fives items:
\begin{itemize}
\item Proving the existence of a boosted dynamics.
\item Finding a labelling compatible with it. 
\item Proving the $r$-point monodromy for some $r$.
\item Initiating the labelling method by finding normal forms.
\item And finally carrying out the rest of the labelling method. 
\end{itemize}
 Each of which has  its own difficulties.

In the first part: 
\begin{itemize}
\item Proving the existence of a boosted dynamics is often easy. The
  dynamics we study (at least as far as our investigations have led us) are `regular', in the
  sense that the operators in $X_{n+1}$ and in $X_n$ do not differ
  much from each other.  Thus, given an operator
  $S=A_i$ on $X_{n-r}$ acting on $y$, the corresponding operator $A_i$ on $X_n$ will 
act on $(x,c)$ in much the same way and we only need to correct the sequence a 'little' to obtain a valid boosted sequence. 

For example in three dynamics that we have studied the boosted sequences of an operator $A_i$ has been respectively $A_i$, $A_i^j$ for some $j$ and $A_j A_i A_j^{-1}$ for some $A_j$ (with $A_j$ possibly the identity).

  This heuristic may not necessarily work in the case of a monoid
  dynamics, since, for a given operator $A_{i}$
  on $X_{n-r}$, it may be not the case that $A_i$ is well-defined on $X_{n}$. Thus when working with monoid dynamics the verification that the boosted dynamics exists is a crucial step that can very well fail (which in turn dashes any hope of applying the method).
\item Proving that there is a labelling compatible with the boosted
  dynamics is normally straigthforward. The objective is to extend the actions of the operators $(A_i)_i$ on the labelling so that the property that the labelling on $X_{n-r}$ keeps track of the gray vertices on $X_n$ is verified. Since, at this point, we have already defined the boosted dynamics we know how the gray vertices moves when applying the boosted sequence on $X_n$ corresponding to an operator on $X_{n-r}$. Thus we just define the action of $A_i$ on the labelling so that the labels of the intervals containing the gray vertices are moved along the gray vertices.

Moreover, given the regularity of the operators, $\alpha'(\cdot)$ is rather similar to $\alpha(\cdot)$.

\item Solving the $r$-point monodromy problem
is the most difficult part of the amenability conditions. We need to investigate the structure of the set 
$L(x,\Pi_b)$ which depends heavily on the invariant set.

Let us recall that our goal is to find a set of loops $L_g(x)$ on $x\in X_n$ such that the group generated by $L_g(x)$ is either all of $L(x,\Pi_b)$ (in which case we have all the information we can ever obtain from the monodromy group) or a subgroup containing at least enough information to answer the $r$-point monodromy problem for some $r$.

Let us note that we can choose (at least in the case of a group dynamics) a specific $x_0$ since all the $L(x,\Pi_b)$ are isomorphic.\footnote{In case of a monoid dynamics the situation can be a bit more difficult but as long as the dynamics is regular enough, a similar argument should be useable.} Thus we can choose a combinatorial structure $x_0$ for which loops that generate the monodromy are easy to obtain (either by directly finding such sequences or by choosing a combinatorial structure that works well with an induction).

Thus the proof is performed by induction at size $n$, \textbf{after} the
induction step for the classification theorem has been established at
size $n$. I.e.\ we suppose both statements true at size $n'<n$ then we
apply the labelling method to prove the classification theorem at size
$n$, and finally we work out the monodromy problem at size $n$, in
this very order for all of our applications. 

Indeed, in this way we will only need to construct a combinatorial structure $x_0$ with the required structure for every valid invariant $inv$ and then the classification theorem at size $n$ will implies that for any starting combinatorial structure $x$ with invariant $inv$ there is a path to $x_0$.

Finally let us conclude with a caveat that will be of importance for the fifth task. Recall that the $r$-point monodromy problem is stated as follows: Given $(x,\Pi_b)$, labels $b_{i_1},\ldots,b_{i_r}\in \Sigma$ and positions $\beta_1,\ldots,\beta_r$ is there a loop that sends, for every $1 \leq k \leq r$, the label $b_{i_k}$ to the position $\beta_k$ ?\\
However the answer to this question can be somewhat implicit, that is, the answer is yes if and only if $P((x,\Pi_b), b_{i_1},\ldots,b_{i_r},\beta_1,\ldots,\beta_r)$ is true for a given property $P$ depending on the structure of the invariant. 

\end{itemize}

For the application of the labelling method there are two steps which can necessitate an
important amount of non-automatised work. 

\begin{itemize}
\item At the very beginning, we shall certify that we can reach
configurations $x'_1$ and $x'_2$ which are of some normal forms.

They are two approaches to this problem. The first one, used in \cite{Rau79} for the standardisation procedure (see also lemma 3.2 in \cite{DS17}), \cite{DS17} for the T-structure and \cite{Fic16} for the piece-wise order reversing permutation of a given type, consists in directly finding a sequence connecting $x$ to a certain $x'$ in normal form. We refer to the three previously cited article to see implementations of this method. 

The second approach is a very nice trick to make this kind of problem easier. Let $T_A$ be the set of normal forms and let $T_B\subseteq T_A$. We will want to choose $T_B$ such that if $y\in T_B$ and $(x,c)$ has reduction $y$ then $x\in T_A$ or a similar property.

 We want to prove that for every $x$ there exists $x'\in T_A$ such that $x\sim x'$. We proceed in two steps:
 
First we construct for every invariant $inv \in Inv_n$ a combinatorial structure $x_0\in T_B$. This construction is often done by induction and can replace the \textbf{existence} step of the labelling method (indeed we are exactly constucting a candidate $x\in X_n$ of a very special form for every $inv\in Inv_n$). Note that since combinatorial structure in $T_B$ have a very constrained structure, the proof of the existence of an $x_0\in T_B$ for every invariant can be long. In this paper a whole section (Sections \ref{sec_IIXshaped}) has to be dedicated to it.

Then the proof is performed by induction at size $n$, \textbf{before} the induction step for the classification theorem has been established at size $n$.\\
 More precisely, since the classification theorem is proven at size $n'<n$, we know that every class (at this size) contains a $x_0\in T_B$. We start from our $x$ at size $n$, we choose a $(n-r,r)$-coloring $c$. Then the reduction $y$ is connected to a $y_0\in T_B$ by a sequence $S$. Consider $x'=S(x,c)$ either $x'\in T_A$ due to our choice of $T_B$ (that is the best case scenario) or at least it should be close to a $x''\in T_A$ with only the gray vertices not at the correct positions. Then using the $r$-point monodromy (as we did in the labelling method) one should be able to prove that $x' \sim x''$. 

Clearly this is only a heuristic since we need to construct the correct $T_B$ for this to work. If $T_A$ is a very large class then we may just have $T_B=T_A$ i.e. a local modification (by adding $r$ vertices) of a combinatorial structure $x$ in $T_A$ remains in $T_A$. 

We will employ this strategy for the labelling proof of the KZB classification in this article. In this case $T_A$ is the set of shift-irreducible standard families and $T_B$ is the set of $I_2X$ permutations.  

Exceptional classes also cause problems, as they often break the ``one
class per invariant'' landscape (which is modified into ``one
non-exceptional class per invariant''). As a result, even with the
normals forms at hand, we may have that $y_1$ and $y_2$ share
the same invariant but are not in the same class, because one of the
two is in an exceptional class.  When we had to face this issue, we
could solve it by choosing $x'_1$ and $x'_2$ appropriately so that
$y_1$ and $y_2$ are easily checked to not be in an exceptional
class (thus constraining even more the normal forms). However, in the cases we have studied so far, certifying this has required the complete characterisation of the structure of the
exceptional classes, a knowledge which may remain elusive in more
complicated dynamics.
\item The step concerning the sequence $S_*$ also needs a certain amount of
care. Indeed, as we have anticipated, the solution to the $r$-point
monodromy problem involve a property $P$ that could be difficult to exploit due to its potentially implicit
form.

At this point of the proof we have the following: $(x_3,c_3)$ and $(x'_2,c2)$ with the same invariant, their reduction $y_2$ with a labelling $\Pi'_b$ and $r$ labels $b_{i_1},\ldots,b_{i_r}$ and $r$ positions $\beta_1,\ldots \beta_r$ such the gray vertices of $(x_3,c_3)$ are within the interval with labels $(b_{i_k})$ and the gray vertices of $(x'_2,c_2)$ are within the interval with positions $(\beta_{k})$. The question is whether this information is enough to deduce that $P((y,\Pi'_b), b_{i_1},\ldots,b_{i_r},\beta_1,\ldots,\beta_r)$ is true.

In this article $P$ is rather explicit so we have enough information and the step is straightforward. However let us quiclky describe an example for which it is non-trivial. 
The involution dynamics (studied in \cite{DS18}) has a unique invariant which is a graph. Moreover   
the labelling of the intervals of a combinatorial structure $x$ correspond bijectively to the labelling of the edges of the graph. Then the property associated to the 1-point monodromy is the following:\\
$P((x,\Pi_b),b, \beta)$ is true if and only if the edges corresponding to $b$ and $\beta$ are isomorphic.
Thus to prove that $P$ is true we must construct an automorphism $g$ of the graph $G$ of $(y_2,\Pi'_b)$ that sends the edge $b$ to the edge $\beta$. However the only information we dispose of is that the graphs $G_1$ and $G_2$ associated to $(x_3,c_3)$ and $(x'_2,c_2)$ are isomorphic by some $f$ and that they are obtained from $G$ by a local modification involving the edges $b$ and $\beta$ respectively. \\ It is the case that without further information on the structure of $G$ and $G_1$, $f$ is not enough to recover an automorphism $g$ of $G$ sending $b$ to $\beta$. In that case our solution is to consider a subset of classed for which the graph invariants have enough structure to recover $g$ from $f$. Then we prove that the classification for this subset of classes yields a classification for all the classes.

For the rest, and assuming that the sequence $S_*$ with the
appropriate properties always exists (otherwise, as we said, we would
be in trouble), all the other steps are elementary, and follow a
non-ambiguous roadmap, with no need of case-to-case inventions.
\end{itemize}

There is one last thing to point out. In the labelling proof we implicitly assumed that
no gray vertices were adjacent since ortherwise there would be two or more gray vertices
within the same interval and thus the number of intervals of interest would be strictly less
than $r$ the number of gray vertices. This assumption does not cost much: it suffices to have
normal forms for which the $(k, r)$-coloring does not produce adjacent gray vertices. In our
application $r$ is small so such property is easy to satisfy.

Nevertheless, we could consider the case where more than one gray vertices are within
the same interval. The same reasoning as above occurs with the following difference: within
an interval the order of the vertices might get permuted by the dynamics thus the data of
interest are the labels $b_1 , . . . , b_{r_1}$ of the $r_1 < r$ intervals containing the gray vertices plus a
permutation $\pi_1 , . . . ,\pi_{r_1}$ indicating the order of the vertices within the same interval.

The permutations $\pi_i$ depend on the dynamics: for the Rauzy dynamics the order remains
fixed so we have $\pi_i = id$. For the involution dynamics \cite{DS18} we have $\pi_i = id$ or $\pi_i = \omega$ the order reversing permutation. We will make use of this slight generalisation for the involution dynamics so we refer to this paper for an example.

We have spend some time detailling how to handle the different problems that can arise in the labelling method, now we summarize this discussion by presenting the organisation of such proof. 

We organise the proof in two steps: the preparation of the induction and its execution. 

\begin{itemize}
\item{Preparation of the induction:}
\begin{itemize}
\item{Part 1:} Definition of the invariants and study of it. The study answers the question what happens to the invariants when adding/removing vertices to a combinatorial structure.
\item{Part 2:} Definition of a boosted dynamics and definition of a compatible labelling and statement of the $r$-point monodromy.
\item{Part 3:} Construction of a $x\in T_b$ for every invariant $inv$.
\end{itemize}
\item{Execution of the induction:}
\begin{itemize}
	\item{Part 4:} Proof that every $(x_1,x_2)$ are connected to a $(x'_1,x'_2)$ in normal forms.
	\item{Part 5:} Proof of the classification theorem by the labelling method
	\item{Part 6:} Proof of the r-point monodromy.
 \end{itemize} 
  \end{itemize} 

\noindent A final remark. The labelling method can only be applied if the set of invariants is known and such set changes completely from case to case. However computing the set $L(x,\Pi_b)$ is of use for find the invariants as the structure of the monodromy group gives indirect information on the underlying invariant. 

For example in this article, the cycle invariant can really be found thanks to the monodromy. Indeed, by studying the monodromy one would notice that the labels are partitionned into $k$ subset of a given length $\lam_1,\ldots,\lam_k$ and that the labels of a given subset can shift in a cyclic way (i.e. one loop can shift all the labels of a given subset by 1 or 2 in a cyclic way and leave the others in place). This information is quite enough to reconstruct the full structure of the cycle invariant.

\section{Definition of the Rauzy dynamics}
We will now apply the labelling method defined in the previous section to give an original proof of the classification of the Rauzy dynamics. Other proofs have been achieved in \cite{Boi12} (Boissy uses the classification proof for the extended Rauzy dynamics appearing in \cite{KZ03}) using geometric methods and in \cite{Fic16} and \cite{DS17} using combinatorial methods.

The extended Rauzy classes (classified in \cite{KZ03}) are of interest in the translation surface field as they are in one-to-one correspondance with the connected components of the strata of the moduli space of abelian differentials (see \cite{Vee82}). As for the non-extended Rauzy classes, they are in one-to-one correspondance with the connected components of the strata of the moduli space of abelian differentials with a marked zero as shown in \cite{Boi12}.
%-------------------------------------------------------
\subsection{The Rauzy dynamics}
% : dynamics and extended dynamics on
%   permutations, and dynamics on matchings}
\label{ssec.3families}
%-------------------------------------------------------

\begin{figure}[b!]
\[
\begin{array}{cc}
m=\big( (16)(24)(37)(58) \big) \in \matchs_8
&
\s=[41583627] \in \kS_8 \subseteq \matchs_{16}
\\
&
\textrm{matching diagram representation}
\\
\raisebox{7pt}{
\setlength{\unitlength}{10pt}
\begin{picture}(10,4)(0,-1)
\put(0,0.2){\includegraphics[scale=2.5]{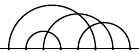}}
\put(.5,-.25){\rotatebox{30}{\makebox[0pt][r]{\scriptsize{$\bm{1}$}}}}
\put(6.8,-.25){\rotatebox{30}{\makebox[0pt][r]{\scriptsize{$m(1)=\bm{6}$}}}}
\end{picture}
}
&
\rule{0pt}{72pt}%
\raisebox{7pt}{
\setlength{\unitlength}{10pt}
\begin{picture}(20,4)(0,-1)
\put(0,0){\includegraphics[scale=2.5]{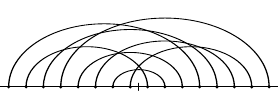}}
\put(.5,-.25){\rotatebox{30}{\makebox[0pt][r]{\scriptsize{$\bm{1}$}}}}
\put(14.5,-.25){\rotatebox{30}{\makebox[0pt][r]{\scriptsize{$n+\s(1)=\bm{12}$}}}}
\end{picture}
}
% \raisebox{2pt}{\includegraphics[scale=0.5]{FigFol/Figure1_fig_match_perm_ex.pdf}}
\\
\rule{0pt}{22pt}%
\s=[41583627] \in \kS_8
&
\s=[41583627] \in \kS_8
\\
\textrm{diagram representation}
&
\textrm{matrix representation}
\\
\raisebox{12.5pt}{%
\setlength{\unitlength}{10pt}
\begin{picture}(11,3)
\put(0,0){\includegraphics[scale=2.5]{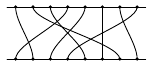}}
\put(.9,-.15){\rotatebox{-30}{\scriptsize{$\bm{1}$}}}
\put(3.6,5.4){\rotatebox{-30}{\makebox[0pt][c]{\scriptsize{$\s(1)=\bm{4}$}}}}
\end{picture}
}
&
\rule{0pt}{90pt}\raisebox{0pt}{%
\setlength{\unitlength}{10pt}
\begin{picture}(8,8)
\put(0,0){\includegraphics[scale=1]{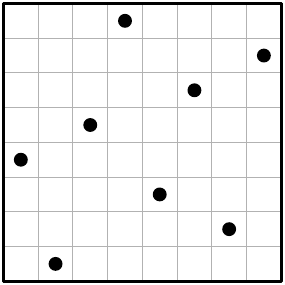}}
\put(0.4,-.6){\scriptsize{$\bm{1}$}}
\put(-.3,3.4){\makebox[0pt][r]{\scriptsize{$\s(1)=\bm{4}$}}}
\end{picture}
}
\end{array}
\]
\caption{\label{fig_match_representation} Diagram representations of
  matchings and permutations, and matrix representation of permutations.}
\end{figure}
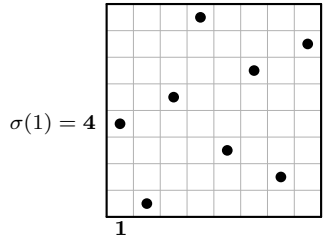

Let $\kS_n$ denote the set of permutations of size $n$, and $\kM_n$
the set of matchings over $[2n]$, thus with $n$ arcs.  Let us call
$\coxe$ the permutation $\coxe(i)=n+1-i$.

A permutation $\s \in \kS_n$ can be seen as a special case of a
matching over $[2n]$, in which the first $n$ elements are paired to
the last $n$ ones, i.e.\ the matching $m_{\s}\in \kM_n$ associated to
$\s$ is $m_{\s}=\{ (i,\s(i)+n) \,|\, i\in [n]\}.$

We say that $\s \in \kS_n$ is \emph{irreducible} if $\coxe \s$ doesn't
leave stable any interval $\{1,\ldots,k\}$, for $1 \leq k < n$,
i.e.\ if $\{\s(1), \ldots, \s(k)\} \neq \{n-k+1,\ldots,n\}$ for any
$k=1,\ldots,n-1$.
% non-trivial interval of consecutive indices,
% within $[n]$, 
We also say that $m \in \kM_n$ is \emph{irreducible} if
% the corresponding permutation of type $2^n$ in $\kS_{2n}$ is (or,
% equivalently, if 
it does not match an interval
$\{1,\ldots,k\}$ to an interval $\{2n-k+1,\ldots,2n\}$.
%doesn't leave
% stable any non-trivial interval of consecutive indices within $[n]$.
Let us call $\kSirr_n$ and $\kMirr_n$
% $\perms_n$ and $\matchs_n$ 
the corresponding sets of irreducible configurations.
%, permutations and matchings respectively.

%% Call $\gamma_{L,i,n}$ the permutation with cycles
%% $(1\,2\,\cdots\,i-1)(i)(i+1)\cdots(n)$, and
%% $\gamma_{R,i,n}$ the permutation with cycles
%% $(n\,n-1\,\cdots\,i+1)(i)(i-1)\cdots(1)$.

We represent matchings over $[2n]$ as arcs in the upper half plane,
connecting pairwise $2n$ points on the real line (see figure
\ref{fig_match_representation}, top left). Permutations, being a
special case of matching, can also be represented in this way (see
figure \ref{fig_match_representation}, top right), however, in order
to save space and improve readibility, we rather represent them as
arcs in a horizontal strip, connecting $n$ points at the bottom
boundary to $n$ points on the top boundary (as in
Figure~\ref{fig_match_representation}, bottom left).  Both sets of
points are indicised from left to right. We use the name of
\emph{diagram representation} for such representations.

We will also often represent configurations as grids filled with one
bullet per row and per column (and call this \emph{matrix
  representation} of a permutation). We choose here to conform to the
customary notation in the field of Permutation Patterns, by adopting
the algebraically weird notation, of putting a bullet at the
\emph{Cartesian} coordinate $(i,j)$ if $\s(i)=j$, so that the identity
is a grid filled with bullets on the \emph{anti-diagonal}, instead
that on the diagonal. An example is given in figure
\ref{fig_match_representation}, bottom right.

Let us define a special set of permutations (in cycle notation)
\begin{subequations}
\label{eqs.opecycdef}
\begin{align}
\gamma_{L,n}(i)
&=
(i-1\;i-2\;\cdots\;1)(i)(i+1)\cdots(n)
\ef;
\\
\gamma_{R,n}(i)
&=
(1)(2)\cdots(i)(i+1\;i+2\;\cdots\;n)
\ef;
\end{align}
\end{subequations}
i.e., in a picture in which the action is diagrammatic, and acting on
structures $x \in X_n$ from below,
\begin{align*}
\gamma_{L,n}(i):
&
\quad
\setlength{\unitlength}{8.75pt}
\begin{picture}(8,2)(-.1,.8)
\put(-.3,0.3){\includegraphics[scale=1.75]{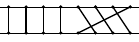}}
\put(0,0){$\scriptstyle{1}$}
\put(3,0){$\scriptstyle{i}$}
\put(7,0){$\scriptstyle{n}$}
\end{picture}
&
\gamma_{R,n}(i):
&
\quad
\setlength{\unitlength}{8.75pt}
\begin{picture}(8,2)(-.1,.8)
\put(-.3,0.3){\includegraphics[scale=1.75]{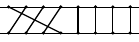}}
\put(0,0){$\scriptstyle{1}$}
\put(4,0){$\scriptstyle{i}$}
\put(7,0){$\scriptstyle{n}$}
\end{picture}
\end{align*}
Of course, $\coxe \gamma_{L,n}(i) \coxe = \gamma_{R,n}(n+1-i)$.

\begin{figure}[tb!]
\begin{align*}
L\;
\Big(
\;
\raisebox{-15pt}{\reflectbox{\includegraphics[scale=1.5]{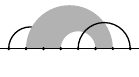}}}
\;
\Big)
&=
\raisebox{-15pt}{\reflectbox{\includegraphics[scale=1.5]{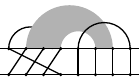}}}
&
R\;
\Big(
\;
\raisebox{-15pt}{\includegraphics[scale=1.5]{FigFol/Figure1_fig_Mlr2bul.pdf}}
\;
\Big)
&=
\raisebox{-15pt}{\includegraphics[scale=1.5]{FigFol/Figure1_fig_Mlr1bul.pdf}}
\\
L\;\Big(
\;
\raisebox{14.2pt}{\includegraphics[scale=1.75, angle=180]{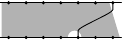}}
% \raisebox{-9pt}{\includegraphics[scale=1.75]{FigFol/Figure1_fig_PPlr3.pdf}}
\;
\Big)
&=
\raisebox{27.2pt}{\includegraphics[scale=1.75, angle=180]{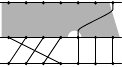}}
% \raisebox{-9pt}{\includegraphics[scale=1.75]{FigFol/Figure1_fig_PPlr4.pdf}}
\rule{0pt}{38pt}
&
R\;\Big(
\;
\raisebox{-9pt}{\includegraphics[scale=1.75]{FigFol/Figure1_fig_PPlr1_simp.pdf}}
\;
\Big)
&=
\raisebox{-24pt}{\includegraphics[scale=1.75]{FigFol/Figure1_fig_PPlr2_simp.pdf}}
\end{align*}
\caption{\label{fig.defDyn3}Our two main examples of dynamics, in
  diagram representation. Top: the $\matchs_n$ case. top: the
  $\perms_n$ case.}
\end{figure}

\begin{figure}[b!!]
\begin{align*}
L\;\left(
\;
\raisebox{-25pt}{\includegraphics[scale=.499]{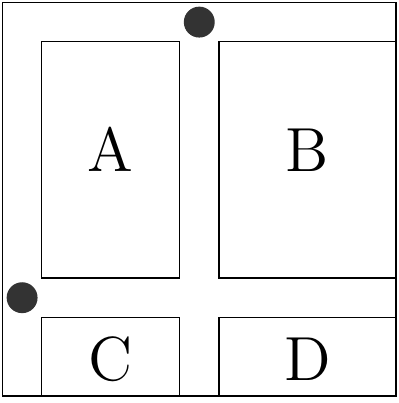}}
\;
\right)
&=
\raisebox{-25pt}{\includegraphics[scale=.499]{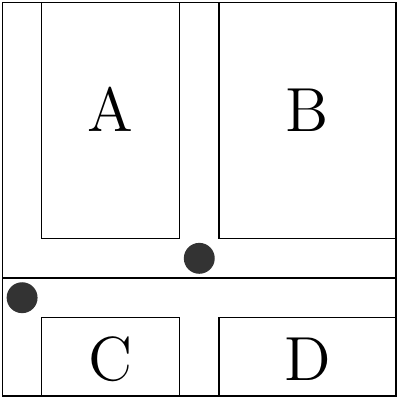}}
\rule{0pt}{38pt}
&
R\;\left(
\;
\raisebox{-25pt}{\includegraphics[scale=.499]{FigFol/Figure1_fig_matrix_dyna_1.pdf}}
\;
\right)
&=
\raisebox{-25pt}{\includegraphics[scale=.499]{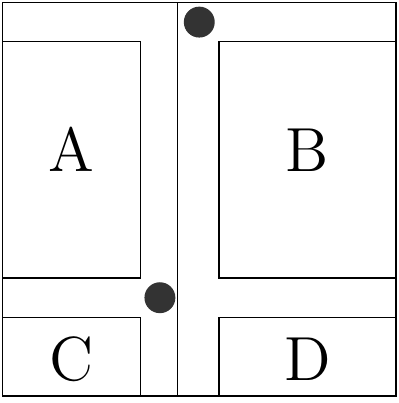}}
\end{align*}
\caption{\label{fig.defDyn3mat} The Rauzy dynamics concerning
  permutations, in matrix representation.
  }
\end{figure}

% !!!! HAKING FOR THE SPECIFIC PAGING OF FIGURES !!!!
% \newpage

The Rauzy dynamics $\perms_n$ that we study in this article is defined as a special case of the $\matchs_n$ dynamics:
\begin{description}
\item[$\matchs_n$\;:]\ The space of configuration is $\kMirr_n$,
  irreducible $n$-arc matchings. There are two generators, $L$ and
  $R$,  with $\a_L(m)=\gamma_{L,2n}(m(1))$ and $\a_R(m)=\gamma_{R,2n}(m(2n))$ ($\a$ is as in Definition
  \ref{def.monoidoperator}).  See Figure~\ref{fig.defDyn3}, top.
\item[\phantom{$\matchs_n$}\gostrR{$\perms_n$}\;:]\ The space of
  configuration is $\kSirr_n$, irreducible permutations of size
  $n$. Again, there are two generators, $L$ and $R$. If permutations
  are seen as matchings such that indices in $\{1,\ldots,n\}$ are
  paired to indices in $\{n+1,\ldots,2n\}$, the dynamics coincide with
  the one given above. See Figure~\ref{fig.defDyn3}, bottom.
%% \end{description}
%% ppppppppppppppppppp
%% \begin{description}
\end{description}

\begin{figure}
\begin{align*}
&
\includegraphics{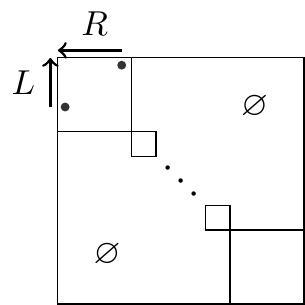}
&&
\includegraphics{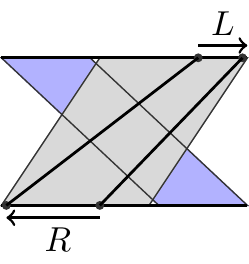}
\end{align*}
\caption{Left: A reducible permutation in matrix representation; the
  $\perms$ dynamic acts on the first block with $L$ and $R$. Right: A reducible permutation in diagram
  representation; the $\perms$ dynamics acts on the gray part with $L$
  and $R$, while leaving the blue part unchanged.
\label{fig.ex_reducibility}}
\end{figure}

The motivation for restricting to irreducible permutations and
matchings shall be clear at this point: a non-irreducible permutation
is a grid with a non-trivial block-decomposition. The operators $L$
and $R$ only act on the first block (say, of size $k$), so that the study of the dynamics trivially
reduces to the study of the $\perms_k$ dynamics on these blocks (see figure
\ref{fig.ex_reducibility}).

This simple observation, however, comes with a disclaimer: To apply the labelling method we must guarantee that the outcome of our manipulations on irreducible configurations is still irreducible during the induction steps in the classification theorem. We explain this problem in the overview section. (The solution will be to choose normal forms that take this problem into account).

%-------------------------------------------------------
\subsection{Definition of the invariants}
\label{ssec.invardefs}
%-------------------------------------------------------

The main purpose of this article, is to
characterise the classes appearing in the so called Rauzy dynamics $\perms_n$ using the labelling method. This gives a third combinatorial proof of this result, the two preceding ones can be found in \cite{DS17} and \cite{Fic16}.

In this section we recall the definition of the
invariants, the proof of their invariance can be found in \cite{DS17}.

% -------------------------------------------------------
\subsubsection{Cycle invariant}
\label{sssec.cycinv}
%-------------------------------------------------------

\noindent
Let $\s$ be a permutation, identified with its diagram.  An edge of
$\s$ is a pair $(i^-,j^+)$, for $j=\s(i)$, where $-$ and $+$ denote
positioning at the bottom and top boundary of the diagram.  Perform
the following manipulations on the diagram: (1)~replace each edge with
a pair of crossing edges; more precisely, replace each edge endpoint,
say $i^-$, by a black and a white endpoint, $i_b^-$ and $i_w^-$ (the
black on the left), then introduce the edges $(i_b^-,j_w^+)$ and
$(i_w^-, j_b^+)$.  (2)~connect by an arc the points $i_w^\pm$ and
$(i+1)_b^\pm$, for $i=1,\ldots,n-1$, both on the bottom and the top of
the diagram; (3)~connect by an arc the top-right and bottom-left
endpoints, $n_w^+$ and $1_b^-$. Call this arc the ``$-1$ mark''.

\begin{figure}[tb!]
\begin{center}
\begin{tabular}{cp{0cm}c}
\includegraphics[scale=5]{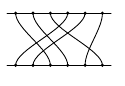}
&&
\includegraphics[scale=5]{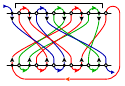}
\end{tabular}
\end{center}
\caption{\label{fig.exLamMatch}Left: an irreducible permutation,
  $\s=[\,451263\,]$. Right: the construction of the cycle
  structure. Different cycles are in different colour. The length of a
  cycle or path, defined as the number of top (or bottom) arcs, is
  thus 2 for red and violet, and 1 for blue.  The green arrows are the
  endpoints of the rank path, which is in blue.  As a result, in this
  example $\lam(\s)=(2,2)$ (for the cycles of color red and violet),
  $r(\s)=1$ (corresponding to the rank path of color blue), and
  $\ell(\s)=2$.}
\end{figure}

The resulting structure is composed of a number of closed cycles, and
one open path connecting the top-left and bottom-right endpoints, that
we call the \emph{rank path}. If it is a cycle that goes through the
$-1$ mark (and not the rank path), we call it the \emph{principal
  cycle}.  Define the length of an (open or closed) path as the number
of top (or bottom) arcs (connecting a white endpoint to a black
endpoint) in the path. These numbers are always positive integers (for
$n>1$ and irreducible permutations).  The length $r$ of the rank path
will be called the \emph{rank} of $\s$, and $\lam = \{ \lam_i \}$, the
collection of lengths of the cycles, will be called the \emph{cycle
  structure} of $\s$. Define $\ell(\s)$ as the number of cycles in
$\s$ (this does not include the rank).  See
Figure~\ref{fig.exLamMatch}, for an example.

Note that this quantity does \emph{not} coincide with the ordinary
path-length of the corresponding paths. The path-length of a cycle of
length $k$ is $2k$, unless it goes through the $-1$ mark, in which
case it is $2k+1$. Analogously, if the rank is $r$, the path-length of
the rank path is $2r+1$, unless it goes through the $-1$ mark, in
which case it is $2r+2$. (This somewhat justifies the name of ``$-1$
mark'' for the corresponding arc in the construction of the cycle
invariant.)

In the interpretation within the geometry of translation surfaces, the
cycle invariant is exactly the collection of conical singularities in
the surface (we have a singularity of $2k\pi$ on the surface, for
every cycle of length $\lam_i=k$ in the cycle invariant, and the rank
corresponds to the `marked singularity' of \cite{Boi12}, see
article \cite{DS17}).

It is easily seen that
\be
r + \sum_i \lam_i = n-1
\ef,
\label{eq.size_inv_cycle}
\ee
this formula is called the 
\emph{dimension formula}. 
% \emph{Gauss-Bonnet formula}. 
Moreover, in
the list $\{r,\lam_1,\ldots,\lam_{\ell}\}$, there is an even number of
even entries. This is part of theorem \ref{thm_arf_value} stated next section and proven during the induction step.

We have
\begin{proposition}
\label{prop.cycinv1}
The pair $(\lam,r)$ is invariant in the $\perms$ dynamics.
\end{proposition}
\noindent
For a proof, see \cite{DS17} section 3.1.

We have also shown in \cite{DS17} appendix B that cycles of length 1 have an especially simple behaviour and can thus be omitted from the classification theorem. Thus \textbf{all the classes we consider in this article have a cycle invariant $\lam$ with no parts of length 1.} 

%-------------------------------------------------------
\subsubsection{Sign invariant}
\label{ssec.arf_inv_intro}
%-------------------------------------------------------

\begin{figure}[tb!]
\[
\begin{array}{cp{4mm}c}
\includegraphics[scale=4]{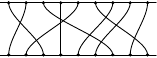}
&&
\includegraphics[scale=4]{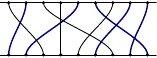}
%% \\
%% \includegraphics[scale=1.6]{FigFol/Figure1_fig_corde_ex1.pdf}
%% &&
%% \includegraphics[scale=1.6]{FigFol/Figure1_fig_corde_ex2.pdf}
\end{array}
\]
\caption{\label{fig.exarf}Left: an example of permutation,
  $\s=[\,251478396\,]$. Right: an example of subset $I=\{1,2,6,8,9\}$
  (labels are for the bottom endpoints, edges in $I$ are in
  blue). There are two crossings, out of the maximal number
  $\binom{|I|}{2}=10$, thus $\chi_I=8$ in this case, and this set
  contributes $(-1)^{|I|+\chi_I}=(-1)^{5+8}=-1$ to $A(\s)$.}
 %% On the
 %%  bottom part of the figure, the same configurations seen as chord
 %%  diagrams.}
\end{figure}

\noindent
For $\s$ a permutation, let $[n]$ be identified to the set of edges
(e.g., by labeling the edges w.r.t.\ the bottom endpoints, left to
right). For $I \subseteq [n]$ a set of edges, define $\chi(I)$ as the
number of pairs $\{i',i''\} \subseteq I$ of non-crossing edges.  Call
\be
\Abar(\s) := \sum_{I \subseteq [n]} (-1)^{|I| + \chi(I)}
\label{eq.arf1stDef}
\ee
the 
\emph{Arf invariant} of $\s$ (see Figure \ref{fig.exarf} for an
example). Call $s(\s)= \mathrm{Sign}(\Abar(\s)) \in \{-1,0,+1\}$ the
\emph{sign} of $\s$.

Both the quantity $\Abar(\cdot)$ and $s(\cdot)$ are invariant for the dynamics $\perms$. The proof can be found in \cite{DS17} section 4.2. However, in order to illustrate how to automatically compute the arf identities, the proof will be given once more in equation (\ref{prop.signinvdyn}).

 There exists an important relationship between the arf invariant and the cycle invariant that we describe below. The proof of this theorem will be done during the induction (see the proof overview section for more details).

\begin{theorem}\label{thm_arf_value}
Let $\s$ be a permutation with cycle invariant $(\lambda,r)$ and let $\ell$ be the number of cycles (not including the rank) of $\s$ i.e. $\ell=|\lam|$.
\begin{itemize}
\item The list $\lambda\bigcup\{r\}$ has an even number of even parts. 
\item $\Abar(\s)=\begin{cases}\pm2^{\frac{n+\ell}{2}}& 
% \text{ if the number of even parts in the list $\lambda\bigcup\{r\}$
%   is 0.}\\
\text{ if there are no even parts in the list $\lambda\bigcup\{r\}$}\\
0 &\text{ otherwise.} \end{cases}$
\end{itemize}
\end{theorem}

as a consequence we have:
\begin{proposition}\label{pro_sign_inv}
The sign of $\s$ can be written as
$s(\s)=2^{-\frac{n+\ell}{2}}\Abar(\s)$, where $\ell$ is the number of
cycles of $\s$.
\end{proposition}

\noindent

% -------------------------------------------------------
\subsection{Exceptional classes}
\label{sssec.ididp}
% -------------------------------------------------------

the invariants described above allow to characterise all classes
for the dynamics on irreducible configurations, \emph{with two
  exceptions}. These two exceptional classes are called $\Id_n$ and $\tree_n$.

$\Id_n$ is called the \emph{`hyperelliptic class'},
because the Riemann surface associated to $\Id_n$ is hyperelliptic. Similarly $\Id'_n$ is often referred to as the \emph{`hyperelliptic class' with a marked point}).

Those two classes were studied in details in Appendix C of \cite{DS17}. In this article, and as outlined in the labelling method section, we will need to certify that the reduced permutations $y_1$ and $y_2$ we obtain from $(x'_1,c_1)$ and $(x'_2,c_2)$ do not fall into exceptional classes. This will not be very difficult but we will need a short lemma that we included into appendix A. Moreover, to make sense of this lemma we reproduce in this appendix a few results of the appendix C.1 of the paper \cite{DS17}.

\begin{table}[t!!]
\[
\begin{array}{|c|c|c|}
\cline{2-3}
\multicolumn{1}{c|}{}
&
\textrm{$n$ even}
&
\textrm{$n$ odd}
\\
\hline
(\lam,r) \text{ of } \Id_n & 
(\emptyset,n-1)&
\raisebox{-5pt}{\rule{0pt}{15pt}}%
(\{\frac{n-1}{2}\}, \frac{n-1}{2})
\\ 
\hline
(\lam,r) \text{ of } \tree_n &
\raisebox{-5pt}{\rule{0pt}{15pt}}%
(\{ \frac{n-2}{2}, \frac{n-2}{2} \},1)& 
(\{n-2\},1)
\\ \hline
\end{array}
\]
\[
\begin{array}{|c|cccccccc|}
\cline{2-9}
\multicolumn{1}{c|}{n \mod 8}
&
0&1&2&3&4&5&6&7
\\ \hline
s \text{ of } \Id_n & 
+&0&-&-&-&0&+&+
\\ 
\hline
s \text{ of } \tree_n &
+&+&0&-&-&-&0&+\\
\hline
\end{array}
\]
\caption{Cycle, rank and sign invariants of the exceptional classes. The
  sign \hbox{$s \in \{-1,0,+1\}$} is shortened into $\{-,0,+\}$.
\label{table_invariant_id_tree}}
\end{table}

The cycle and sign invariants of these classes depend from their size
mod~4, and are described in Table~\ref{table_invariant_id_tree}.  

%-------------------------------------------------------
\subsection{The classification theorems}
\label{ssec.theointro}
%-------------------------------------------------------

For the case of the $\perms$ dynamics, we have a classification
involving the cycle structure $\lambda(\s)$, the rank $r(\s)$ and sign
$s(\s)$ described in Section \ref{ssec.invardefs}

\begin{theorem}\label{thm.Main_theorem}
Besides the exceptional classes $\Id$ and $\tree$, which have cycle
and sign invariants described in Table~\ref{table_invariant_id_tree},
the number of classes with cycle invariant $(\lam,r)$ (no
$\lambda_i=1$) depends on the number of even elements in the list
$\{\lambda_i\} \cup \{r\}$, and is, for $n \geq 9$,
\begin{description}
\item[zero,] if there is an odd number of even elements;
\item[\phantom{zero,}\gostrR{one,}] if there is a positive even number of even elements; the class then has sign 0.
\item[\phantom{zero,}\gostrR{two,}] if there are no even elements at all. The two classes then have
non-zero opposite sign invariant.
\end{description}
%% if both $r$ and all $\lam_i$ are odd, and $\lam_i>1$ for
%% all $i$, and they have
%% opposite sign invariant, and there exists exactly one primitive class 
%% with cycle invariant
%% $(\lam,r)$, if both $r$ and all $\lam_i$ are odd, and $\lam_i>1$
For $n \leq 8$ the number of classes with given cycle
invariant may be smaller than the one given above, and the list in
Table \ref{tab.smallsizeThm1} gives a complete account.

As a consequence, two permutations $\s$ and $\s'$, not of
$\Id$ or $\tree$ type, are in the same class iff they have the same
cycle and sign invariant.
\end{theorem}

% !!!! POSITION HACKED !!!!
\begin{table}[t!!]
\[
\begin{array}{r||c|c|cccccc}
n & \Id & \Id' & \multicolumn{6}{c}{\textrm{non-exceptional classes}}
\\
\hline
4 & \emptyset|3- &  \\
5 & 2|2 & 3|1- \\
6 & \emptyset|5+ & 22|1 & \emptyset|5- \\
7 & 3|3+ & 5|1+ & 2|4 & 4|2 & 3|3- & 5|1- \\
8 & \emptyset|7+ & 33|1+ & 
22|3 & 32|2 & 42|1 & 33|1- & \emptyset|7+ & \emptyset|7-
\end{array}
\]
\caption{\label{tab.smallsizeThm1}List of invariants $(\lam,r,s)$ for
  $n \leq 8$, for which the corresponding class in the $\perms_n$
  dynamics exists. We shorten $s$ to $\{-,+\}$ if valued $\{-1,+1\}$,
  and omit it if valued~0.}
\end{table}

\noindent
The original theorem from Kontsevich and Zorich \cite{KZ03} classifies the extended Rauzy dynamics $\permsex$, however we have shown in section 6.5 of \cite{DS17} that this theorem is a simple corollary of Theorem \ref{thm.Main_theorem}.

%%%%%%%%%%%%%%%%%%%%%%%%%%%%%%%%%%%%%%%%%%%%%%%%%%%%%%%
\section{Proof overview}
\label{ssec_pf_overview}
\label{ch.rauzy2}

In this section we present a proof of the classifcation of the
Rauzy classes of the dynamics $\cS_n$ by applying the labelling
method.

A proof of classification via the labelling method, always proceeds in two parts. First we prove the
amenability of the dynamics $\cS_n$ in section \ref{sec_arc_label}
(although, as usually occurs in this family of proofs, the proof of
the $r$-point monodromy must be deferred to the main induction). Then
we carry out the main induction in section \ref{sec_induction}, in
which we apply the labelling method to prove the classification
theorem and the $r$-point monodromy theorem.

In an instance where the labelling method can be applied `smoothly'
(as is the case with the involution dynamics in \cite{DS18}, though it does have some interesting specific difficulties), the main induction is
constituted of three statements: the \emph{existence} and
\emph{completeness} part of the classifcation theorem, and then the
$r$-point monodromy of the set $L(x,\Pi_b)$. Moreover, nothing
prevents from presenting the main induction immediately after that the
amenability has been established, i.e. no complicated normal forms have to be defined.

The dynamics $\cS_n$ is not such a smooth instance. Therefore the main
induction, that we now detail, will include a certain number of
unavoidable technical statements, and several sections (i.e.\ Sections
\ref{sec_edge_addition} to \ref{sec_IIXshaped}) are needed in
preparation of it. In particular, we will construct our normal forms using the trick introduced in the section \ref{sec_diff_lab} of the labelling method. The set $T_A$ will be the set of shift-irreducible permutations and the set $T_B\subseteq T_A$ that of the $I_2X$-permutations.

\bigskip
\noindent
\textbf{Main induction:} We demonstrate by induction on $n$ the seven
following statements:

\noindent
\begin{tabular}{rl}
% \imagetop{\rule{0pt}{3pt}
\raisebox{-\height+\baselineskip}{1.} &
% \imagetop{
\raisebox{-\height+\baselineskip}{\begin{minipage}{.9\textwidth}
\begin{proposition}\label{pro_shift_irreducible} Every class contains a shift-irreducible standard family (cf definition \ref{def_shift_irre}).
\end{proposition}
\end{minipage}}
\end{tabular}

\noindent
\begin{tabular}{rl}
% \imagetop{\rule{0pt}{3pt}
\raisebox{-\height+\baselineskip}{2.} &
% \imagetop{
\raisebox{-\height+\baselineskip}{\begin{minipage}{.9\textwidth}
\textbf{Theorem \ref{thm_arf_value}} on the relationship between the arf invariant and the cycle invariant.
\end{minipage}}
\end{tabular}

\noindent
\begin{tabular}{rl}
% \imagetop{\rule{0pt}{3pt}
\raisebox{-\height+\baselineskip}{3.} &
% \imagetop{
\raisebox{-\height+\baselineskip}{\begin{minipage}{.9\textwidth}
\begin{proposition}[Existence part of the classification theorem]\label{pro_existence}
For every valid invariant $(\lam,r,s)$ (i.e., every invariant in the
list of Thm.~\ref{thm.Main_theorem}) there exists a permutation with
invariant $(\lam,r,s)$.
\end{proposition}
\end{minipage}}
\end{tabular}

\noindent
\begin{tabular}{rl}
% \imagetop{\rule{0pt}{3pt}
\raisebox{-\height+\baselineskip}{4.} &
% \imagetop{
\raisebox{-\height+\baselineskip}{\begin{minipage}{.9\textwidth}
\begin{proposition}[First step of the labelling method: the normal forms]\label{pro_labelling_first_step}
Let $\s_1,\s_2\in \kS_n$ be two irreducible permutation with invariant
$(\lam,r,s)$. There exist $\s'_1$ and $\s'_2$, connected to $\s_1$ and
$\s_2$ respectively, with the following property:\\ let $c_1$ be the
$(2n-2,2)$-coloring of $\s'_1$ where the edge $e'_1=(\s_1^{\prime-1}(1),1)$ is
grayed, and call $\tau_1$ the reduction of $(\s'_1,c_1)$; define the
analogous quantities for $\s_2'$ (i.e., the coloring $c_2'$, the edge
$e_2'=(\s_2^{\prime-1}(1),1)$, and the configuration $\tau_2$); then $\tau_1$
and $\tau_2$ are irreducible, have the same invariant $(\lam',r',s')$,
and none of them is in an exceptional class.
\end{proposition}
\end{minipage}}
\end{tabular}

\noindent
\begin{tabular}{rl}
% \imagetop{\rule{0pt}{3pt}
\raisebox{-\height+\baselineskip}{5.} &
% \imagetop{
\raisebox{-\height+\baselineskip}{\begin{minipage}{.9\textwidth}
\begin{proposition}[Completeness part of the classification theorem]\label{pro_completeness}
Every pair of permutations $\s$ and $\s'$ with the same invariant $(\lam,r,s)$ are connected. 
\end{proposition}
\end{minipage}}
\end{tabular}

\noindent
\begin{tabular}{rl}
% \imagetop{\rule{0pt}{3pt}
\raisebox{-\height+\baselineskip}{6.} &
% \imagetop{
\raisebox{-\height+\baselineskip}{\begin{minipage}{.9\textwidth}
\begin{proposition}\label{pro_i2x}
Every non exceptional class contains a $I_2X$-permutation (Defined above corollary \ref{cor_chara_shift}).
%  and Corollary \ref{cor_chara_shift}).
\end{proposition}
\end{minipage}}
\end{tabular}

\noindent
\begin{tabular}{rl}
% \imagetop{\rule{0pt}{3pt}
\raisebox{-\height+\baselineskip}{7.} &
% \imagetop{
\raisebox{-\height+\baselineskip}{\begin{minipage}{.9\textwidth}
\textbf{Theorem \ref{thm_2monodromy}} on the 2-point monodromy problem
for this dynamics.
\end{minipage}}
\end{tabular}

\bigskip
\noindent The three major steps of the present proof are to establish
Propositions \ref{pro_existence} and \ref{pro_completeness} and
Theorem~\ref{thm_2monodromy}.

As we commented in the Section \ref{sec_diff_lab}, the
most difficult step in the  organisation of the abstract
labelling method, when specialised to the $\cS_n$ dynamics, is finding the normal forms. In this case, this step takes the form of Proposition \ref{pro_labelling_first_step} within the main induction using both Propositions \ref{pro_shift_irreducible} and \ref{pro_i2x} as support.

This proposition is difficult to establish because it requires that
$\tau_1$ and $\tau_2$ satisfy simultaneously four unrelated properties:
being irreducible, having the same cycle invariant $(\lam',r')$,
having the same sign invariant $s'$, and not being in an exceptional
class.  Let us outline the steps of the proof of this proposition:
\paragraph*{Cycle invariant:} In Section \ref{sec_shift_irreducible}
(Proposition \ref{pro_std_d_cycle_inv}), we observe that if $\s_1'$ and
$\s_2'$ are standard and are both of type $X(r,i)$ or of type
$H(r_1,r_2)$, then their reductions have the same cycle invariant
$(\lam',r')$.

\paragraph*{Irreducibility:} We introduce a special subset of
standard families that we call the \emph{shift-irreducible standard
  families} (see Definition \ref{def_shift_irre}). By definition, if
$\s$ is in a shift-irreducible standard family, then the reduction
$\tau$ of $(\s,c)$ is (almost always and in our case always) irreducible, where $c$ is the
$(2n-2,2)$-coloring of $\s$ where edge $e=(\s^{-1}(1),1)$ is
grayed. 

\paragraph*{Sign invariant:} In Section \ref{sec_arf} we prove
Proposition \ref{lem_q_i_sign}. Applied to our case, it shows that if
$\s'_1$ and $\s_2'$ have invariant $(\lam,r,s)$ with $s\neq 0$ and are
shift-irreducible standard permutation of type $X(r,i)$ then $\tau_1$
and $\tau_2$ also have sign invariant $s$. There are two other cases
concerning the permutations of type $H$, for which the underlying mechanism is
too technical to be described concisely. Let us only note here that they make use of Theorem \ref{thm_arf_value} and
Proposition \ref{cor.signsurjq2x} (also established in Section
\ref{sec_arf}), respectively.

\paragraph*{Exceptional classes:} Finally we must certify that neither
$\tau_1$ nor $\tau_2$ are in an exceptional class. This is achieved
by showing that, in a standard family $(\s_i)_{1\leq i\leq n-1}$, at most one of the reductions of the permutation $(\s_i,c_i)$ where $c_i$ is the $(2n-2,2)$-coloring in which the edge ($\s^{-1}_i(1),1$) is gray, is in a exceptional class (a lemma established in the appendix discussing the structure of exceptional classes). As a result, we can always avoid this case.\\

Thus, by choosing $\s_1'$ and $\s_2'$ in shift-irreducible standard
families, we can take care simultaneously of those four properties and prove proposition \ref{pro_labelling_first_step}.

However, contrarily to the existence of standard families, proving that every irreducible permutation is
connected to a shift-irreducible standard family by directly finding a sequence is too hard. Thus we make use of the trick introduced in section \ref{sec_diff_lab} and we find a special subset $T_B \subseteq T_A$ of standard permutations, named \emph{$I_2X$ permutations}, for which the associated standard family is
shift-irreducible (see corollary \ref{cor_chara_shift}). We dedicate the whole Section \ref{sec_IIXshaped} to the construction of one such permutation for every $(\lam,r,s)$, which, on top of this, are not in an exceptional class (this result is the content of Theorems \ref{thm_I_2X} and~\ref{thm_I_2X_odd}).
%  a special type of shift-irreducible standard permutation. 

Then, as we outlined in section \ref{sec_diff_lab}, in the main induction we use the stronger induction hypothesis of Proposition \ref{pro_i2x} (statement 6) in order to prove
Proposition \ref{pro_shift_irreducible} (statement 1). Finally
Proposition \ref{pro_i2x} (statement 6) is demonstrated by applying
the completeness part of the classification theorem (Proposition
\ref{pro_completeness}, statement 5) to the existence of a $I_2X$
permutation for every $(\lam,r,s)$ (proven in Section
\ref{sec_IIXshaped}, Theorems \ref{thm_I_2X} and~\ref{thm_I_2X_odd}). 

Thus we can organise this article as follows 
\begin{itemize}
\item{Section 5:} We define the boosted dynamics, a labelling compatible with it and state the 2-point monodromy theorem.
\item{Sections 6-7:} We study the cycle invariant when adding or removing edges on permutations and define both the shift-irreducible family and the $I_2X$-permutations.
\item{Section 8:} We study the arf invariant when adding or removing edges on permutations and some other local modification.
\item{Section 9:} Construct a $I_2X$ permutation for every possible $(\lam,r,s)$.
\item{Section 10:} Proceed with the induction, proving the 7 statements presented in this overview.
\end{itemize}

We note that the structure is very close to the organisation presented at the very end of the labelling section: Part 1 corresponds to sections 6-7-8, part 2 to section 5, part 3 to section 9 and parts 4,5,6 to section 10.

This organisation will reappear in all proofs using the labelling method with some minor differences on a case-by-case basis.

%-------------------------------------------------------
\section{Amenability of $\cS_n$ to the labelling method}
\label{sec_arc_label}

In this section we introduce the notions necessary to show that
$\cS_n$ is amenable to the labelling method. Thus our task is to define a boosted dynamics,
determine a labelling $\Pi_b$ compatible with it as well as to state the
$2$-point monodromy theorem for set $L(\s,\Pi_b)$ for any permutation
$\s$.\\
Let us overview the content of the section: 

We begin the section by introducing the boosted dynamics.

Then the second subsection introduces the definition of a special set of
labellings, that we call \emph{consistent labellings}, and proves a
number of properties. In short, consistent labellings are labellings on an alphabet
that respects the structure of the cycle invariant.

The third subsection extends the dynamics to the labelled case,
proves its compatibility with the boosted dynamics (cf.\ Theorem
\ref{thm_keep_track_label}), and introduces the theorem pertinent to
the 2-point monodromy problem (cf.\ Corollary~\ref{cor_2monodromy}).

%-------------------------------------------------------
\subsection{Boosted dynamics}
\label{ssec.reddyn}
%-------------------------------------------------------

Call \emph{pivots} of $\s$ the two edges $(1,\s(1))$ and
$(\s^{-1}(n),n)$. For a pair $(\s,c)$, we say that $\s$ is \emph{proper} if no gray edge
of $\s$ is a pivot.  In this case, the dynamics on $\tau$ (the reduction of $\s$) extends to
the \emph{boosted dynamics} on $\s$, as follows:

for every operator $H$ (i.e., $H\in\{L,R\}$), we define
$\alpha_H(\s,c)$ as the smallest positive integer such that
$H^{\alpha_H(\s,c)}(\s)$ is proper, and, for a sequence 
$S=H_k \cdots H_2 H_1$ acting on $\tau$, the sequence $B(S)$, the
\emph{boosted sequence} of $S$, acting on $\s$ is
$B(S)=H_k^{\alpha_k} \cdots H_2^{\alpha_2} H_1^{\alpha_1}$ for the
appropriate set of $\alpha_j$'s.

A simple verification (by induction on the number of gray edges) shows that $B(S)$ makes the following square commutes.
\begin{center}
$\begin{tikzcd}
  (x,c) \arrow{d}{\textrm{red}} 
\arrow{r}{B(S)} &  (x',c') \arrow{d}{\textrm{red}} &  \\
 y \arrow{r}{S} &    y' & \\
\end{tikzcd}$
\end{center} 
Thus $B(S)$ is a well-defined boosted sequence.

%-------------------------------------------------------
\subsection{Preliminaries: Definition of a \emph{consistent labelling}}

In the labelling method we labelled the intervals of a combinatorial structure with intervals. In this particular case, the intervals corresponds to the bottom and top arcs introduced in the construction of the cycle invariant. Thus in the following we will rather label the arcs (top and bottom), and instead of saying \textbf{'the gray vertex is within the interval with label $b$'} we will say \textbf{'the gray edge is within the arcs with labels $t$ and $b$ (for the top arc and bottom arc respectively)'}.   

Also note the change of terminology from graying vertices to graying edges. Indeed, recall definition \ref{def_coloring}: when choosing a coloring $c$ for a permutation we must garantee that the reduction is still a permutation and thus we will always gray edges (or equivalently the two vertices connected by the edge) rather than just vertices.)\\

%Our notion of arcs departs slightly from the definition inherited from
%the combinatorial structure with arcs as we define both bottom and top
%arcs in our diagram representation of permutations. 

% \begin{definition}
Let $\s$ be a permutation of size $n$. The procedure to construct the
cycle invariant $(\lambda,r)$ (as described in Section
\ref{sssec.cycinv}) involves the introduction of $n-1$ top and bottom
arcs connecting adjacent top and bottom vertices. As our proof
requires it, we also add one top arc, to the left of the other top
arcs (i.e to the left of the edge $(\s^{-1}(1),1)$) and one bottom
arc, to the right of the bottom arcs (i.e to the right of the edge
$(n,\s(n))$). Those two arcs are clearly added at the two endpoints of
the rank path (in red in the figure below).

We number the top and bottom arcs from left to right, and refer to
them by their position: the bottom arc $\beta\in \{1,\ldots,n\}$ is
the $\beta$-th bottom arc, counting from the left.
% \begin{center}
\[
\begin{array}{cc}
\raisebox{10.9pt}{ \includegraphics[scale=.5]{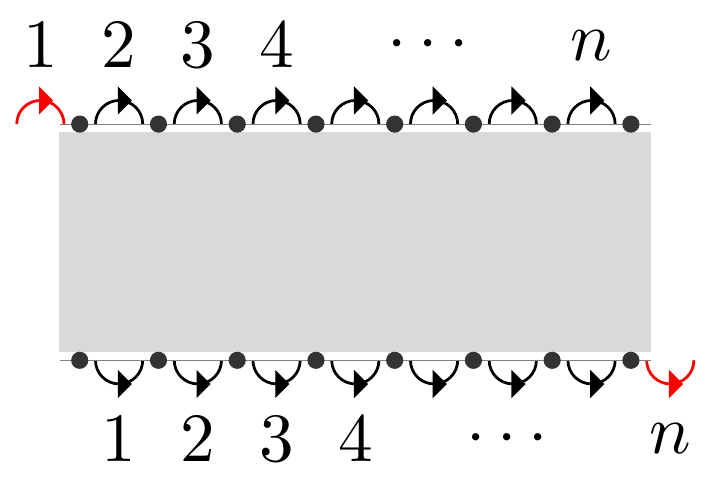}}
&
\raisebox{11pt}{\includegraphics[scale=.5]{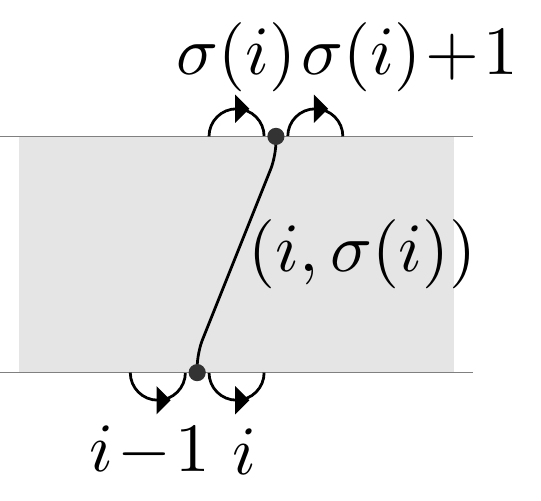}}
\end{array}
\]
% \end{center}

By convention, the variables used to name the positions of the top
(bottom) arcs will be $\alpha$ (respectively $\beta$), in order not to
make confusion with other parts of the diagram (for which we will use
$i,j,\ldots$ or $x,y,\ldots$).
% \end{definition}

\begin{definition}\label{def_consecutive}
We say that two (bottom) arcs $\beta,\beta'$ are \emph{consecutive}
(in a cycle) if they are inside the same cycle (or the rank path), 
and they are consecutive in the cyclic order induced by the cycle (or
the total order induced by the path). This occurs in one of the three
graphical patterns:
\begin{center}\begin{tabular}{ccc}
    \includegraphics[scale=.5]{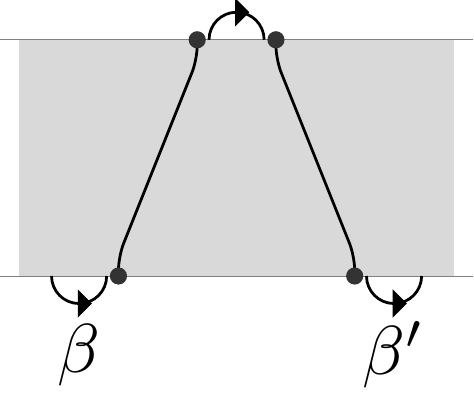} &
    \includegraphics[scale=.5]{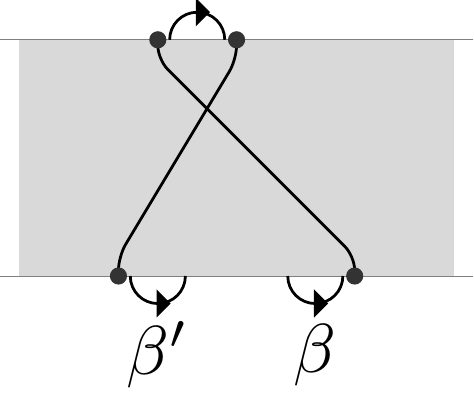} &
    \includegraphics[scale=.5]{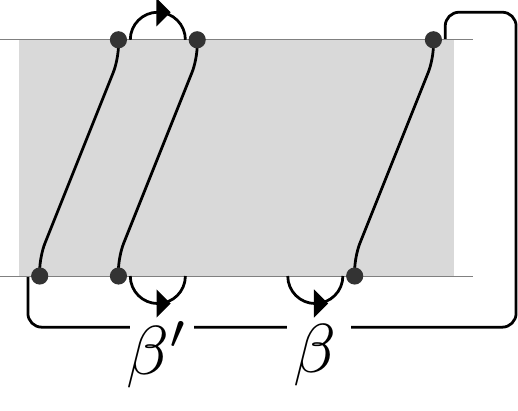} 
\end{tabular}
\end{center}
In formulas:
\[
\beta'=\s^{-1}(\s(\beta+1)+1) \text{ if } \s(\beta+1) <n \text{ and }
\beta'=\s^{-1}(\s(1)+1) \text{ if } \s(\beta+1) =n
\]
We define consecutive arcs for top arcs similarly.
\end{definition}

\begin{remark}
As we have seen above, when representing graphically the consecutive
arcs, we need three figures depending on the different cases (edges
crossing or not, and edges ending at the north-est corner of the
diagram). However, these cases are treated in a very similar way, and,
in the graphical explanation of our following properties, we shall
mostly draw consecutive arcs by representing the case of non-crossing
and non-corner edges, i.e.\ the left-most of the drawings above. It is
intended that the underlying reasonings remain valid for the other
cases.
\end{remark}

Next we define suitable alphabets used to label the top and bottom
arcs of a permutation.

\begin{notation}
For all $j$, let $\Sigma_{i,j}=\{b_{0,i,j},\ldots,b_{i-1,i,j}\}$ and
$\Sigma'_{i,j}=\{t_{0,i,j},\ldots,t_{i-1,i,j} \}$ be a pair of
alphabets which label the bottom arcs and the top arcs respectively of
a cycle of length $i$, and let
$\Sigma_{r}=\{b^{\rm rk}_{0},\ldots,b^{\rm rk}_{r}\}$ and
$\Sigma'_{r}=\{t^{\rm rk}_{0},\ldots,t^{\rm rk}_{r}\}$ be the pair of
alphabets used to label the bottom arcs and the top arcs respectively
of the rank path.
\end{notation}

Note that the labels of the rank range from $0$ to $r$, instead of
form 0 to $r-1$, since we added a left-most top arc and a right-most
bottom arc, which are now part of the rank, according to the
construction of the arcs of a permutation outlined above.

\begin{figure}[b]
\begin{center}\begin{tabular}{c} \includegraphics[scale=.7]{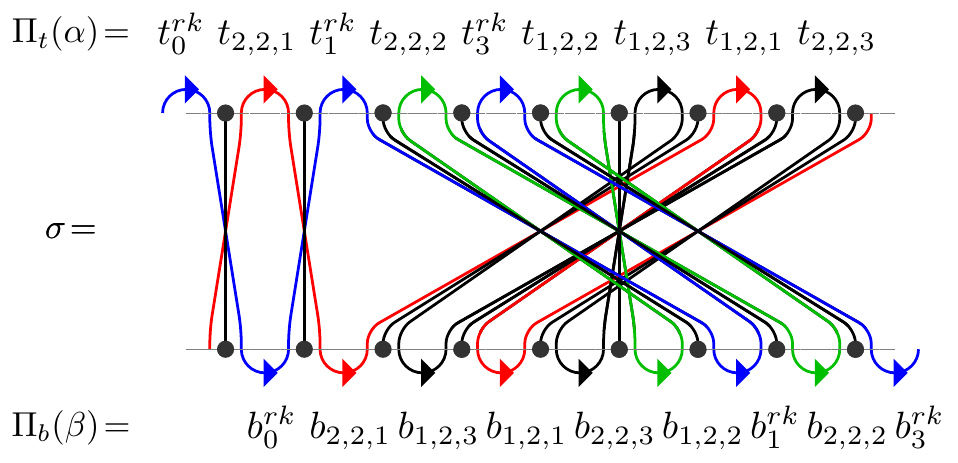}  \includegraphics[scale=.7]{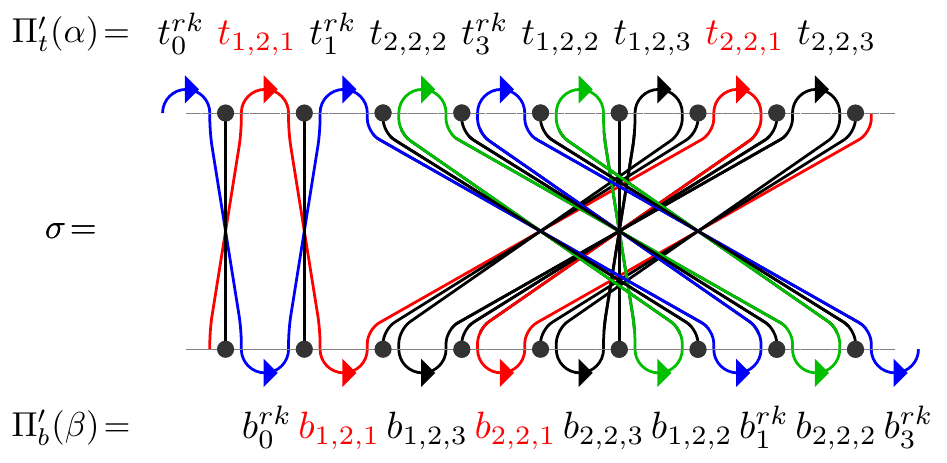} \end{tabular}\end{center}
\caption{\label{fig_example_labelling} Two consistent labellings $(\Pi_b,\Pi_t)$ and $(\Pi'_b,\Pi'_t)$ of a permutation $\s$ with cycle invariant $(\{2,2,2\},2)$. Following definition \ref{def_shift_exch} we have $(\Pi'_b,\Pi'_t)=\mathrm{Sh}^1_{2,1}((\Pi_b,\Pi_t))$ }
\end{figure}

Finally, we can introduce our notion of consistent labelling.
% In essence, a consistent labelling labels the arcs of the
% permutation in a way  that respects the cycle invariant. 
\begin{definition}[Consistent labelling]
\label{def_consistent_lab}
Let $\s$ be a permutation with invariant
$(\lam=\{\lam_1^{m_1},\ldots,\lam_k^{m_k}\},r)$ and define a
\emph{consistent labelling} to be a pair $(\Pi_b,\Pi_t)$ of
bijections:
\[
\begin{array}{ccccc}
\Pi_b & : & \{1,\ldots,n\} & \to & \Sigma_b= 
\Sigma_r \cup \Big[ \bigcup_{i=1}^k \Big( \bigcup_{j=1}^{m_i}
  \Sigma_{\lam_i,j} \Big) \Big]
\\
\Pi_t & : & \{1,\ldots,n\} &\to &\Sigma_t= 
\Sigma'_r \cup \Big[ \bigcup_{i=1}^k \Big( \bigcup_{j=1}^{m_i}
  \Sigma'_{\lam_i,j} \Big) \Big]
\end{array}
\]
such that
\begin{enumerate}
\item Two arcs within the same cycle have labels within the same
  alphabet. Thus if $S_b=\{(\beta_k)_{1\leq k \leq \lam_i}\}$ and
  $S_t=\{(\alpha_k)_{1\leq k \leq \lam_i}\}$ are the sets of bottom
  (respectively top) arcs of a cycle of length  $\lam_i$, then
  $\Pi_b(S_b)=\Sigma_{\lam_i,j}$ and $\Pi_t(S_t)=\Sigma'_{\lam_i,j}$
  for some $1 \leq j\leq m_i$.
\item Two consecutive arcs of a cycle of length $\lam_i$ have labels
  with consecutive indices: if $\beta$ and $\beta'$ are consecutive,
  then $\Pi_b(\beta)=b_{k,\lam_i,j}$ for some $k<\lam_i$ and $j\leq
  m_i$ and $\Pi_b(\beta')=b_{k+1,\lam_i,j}$, where $k+1$ is intended
  modulo $\lam_i$. Likewise for top arcs.\\[-7mm]
\begin{center}
\begin{tabular}{cc}
  \includegraphics[scale=.5]{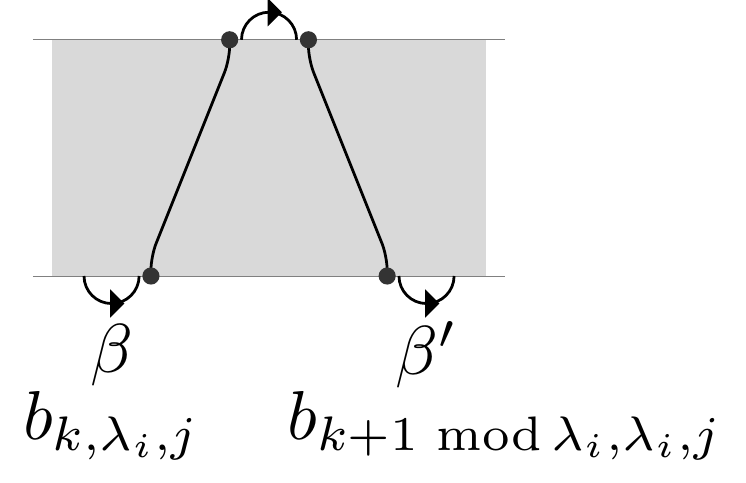} &
  \raisebox{24pt}{\includegraphics[scale=.5]{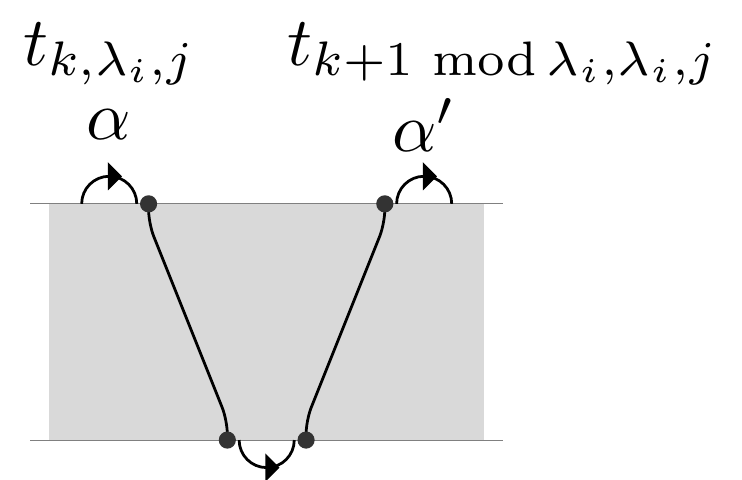}}
\end{tabular}
\end{center}
~\\[-12mm]
\item The bottom right arc $\beta$ and the top left arc $\alpha$ of an
  edge $i$ are labeled by the same indices:
\[
\text{ if } \beta=i,\alpha=\s(i), \text{ then }\, 
\left\{
\begin{array}{c}
\Pi_t(\alpha)=t_{k,\lam_{\ell},j} \ \Leftrightarrow \ 
\Pi_b(\beta)=b_{k,\lam_\ell,j} 
\\
\rule{0pt}{12pt}
\Pi_t(\alpha)=t^{\rm rk}_{k} \ \Leftrightarrow \ 
\Pi_b(\beta)=b^{\rm rk}_{k}
\end{array}
\right.
\] \\[-12mm]
\begin{center}
\includegraphics[scale=.5]{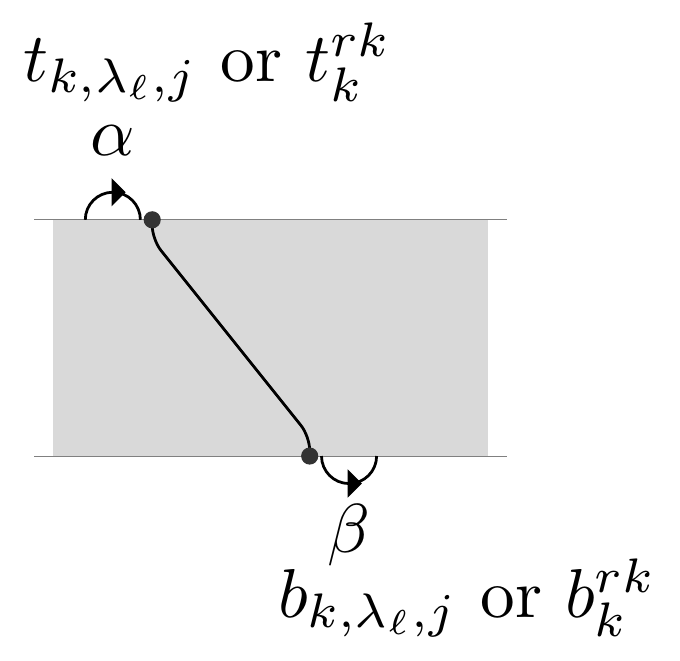}\end{center}~\\[-14mm]
\item Let $(\alpha_i)_{0\leq i\leq r}$ and $(\beta_i)_{0\leq i\leq r}$ be the top (respectively bottom) arcs of the rank ordered along the path (thus $\alpha_0=1$, $\beta_0=\s^{-1}(1)$, $\alpha_1=\s(\beta_0+1)+1$ $\ldots$ $\alpha_r=\s(n)$, $\beta_r=n$). Then \[\forall \ 0\leq i\leq r, \ \Pi_t(\alpha_i)=t^{\rm rk}_{i} \text{ and } \Pi_b(\beta_i)=b^{\rm rk}_{i} \] \\[-10mm]
\begin{center} \includegraphics[scale=.5]{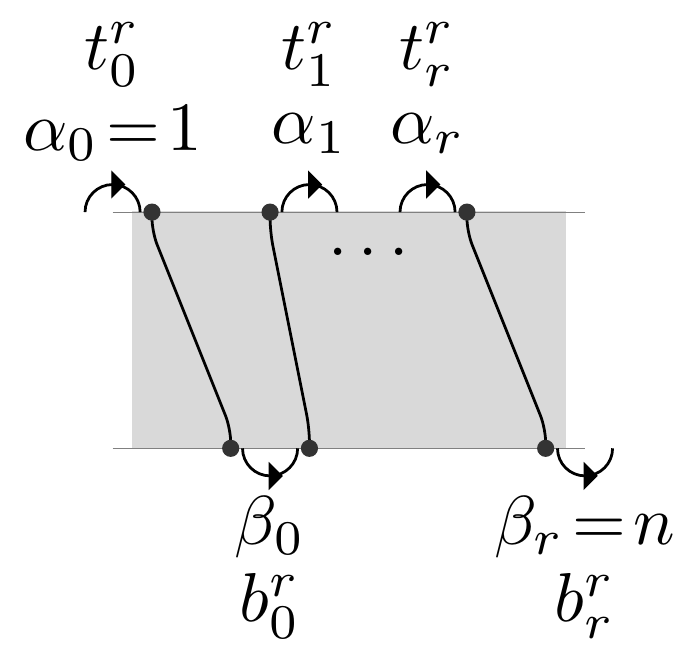}\end{center}~\\[-14mm]
\end{enumerate}
\end{definition}
\noindent
Figure~\ref{fig_example_labelling} provides two examples of consistent
labellings.

\begin{lemma}\label{rk_top_from_bot}
Let $\s$ be a permutation and 
$\Pi_b \,:\,\{1,\ldots,n\} \to \Sigma_b$ a labelling of bottom arcs
verifying property 1, 2 and 4 of Definition
\ref{def_consistent_lab}. Then there exists a unique 
$\Pi_t \,:\,\{1,\ldots,n\} \to \Sigma_t$ such that $(\Pi_b,\Pi_t)$ is
a consistent labelling.
\end{lemma}

\begin{pf}
Let  $(\Pi_b,\Pi_t)$ a be a consistent labelling. Then by property 3
we must have 
$\Pi_t(\alpha)= \begin{cases} t_{i,\lam_\ell,j} & \text{if }
  \Pi_b(\s^{-1}(\alpha))=b_{i,\lam_\ell,j}\\ t^{\rm rk}_{i} & \text{if
  } \Pi_b(\s^{-1}(\alpha))=b^{\rm rk}_{i}\\
\end{cases}.$
This uniquely defines $\Pi_t$.  
\qed
\end{pf}

\noindent
This lemma implies that the data $(\s,\Pi_b)$ or $(\s,\Pi_t)$ are
sufficient to reconstruct $(\s,(\Pi_b,\Pi_t))$. Thus, occasionally, we
will consider just $(\s,\Pi_b)$ rather than $(\s,(\Pi_b,\Pi_t))$.

\begin{lemma}\label{lem_verif_consistent}
Let $(\s,(\Pi_b,\Pi_t))$ be a permutation with a labelling. We have
that $(\Pi_b,\Pi_t)$ is a consistent labelling (i.e.\ it verifies
properties 1 to 4) if property 2 is true for bottom arcs, property 3
is true and $\Pi_t(1)=t^{\rm rk}_0$.
\end{lemma}

\begin{pf}
Clearly property 1 is implied by property 2 and 3, and property 4 is
implied by property 2 and 3 applied to the rank, plus the fact that
$\Pi_t(1)=t^{\rm rk}_0$, which, by property 2, fixes the order of the
other labels of the rank.  \qed
\end{pf}

It is clear from the definition that two consistent labellings of a
permutation can only differ by the two following operations: a cyclic
shift of the labels within a cycle (due to property 2), and the
permutation of the `$j$' labels of two cycles of same size (due to
property 1). In particular, the labels of the rank, and their order,
coincide in all consistent labelings of a given permutation
(by property 4).

We state this property more formally in the following definition and
proposition.

\begin{definition}\label{def_shift_exch}
Let $\s$ be a permutation with invariant
$(\lam=\{\lam_1^{m_1},\ldots,\lam_k^{m_k}\},r)$, and define operators
on consistent labelling $\Pi_b \,:\, \{1,\ldots,n\} \to \Sigma$
\begin{itemize}
\item The \emph{shift operator} $\mathrm{Sh}^m_{\lambda_i,j}$, that shifts by
  $m$ the labels of the $j$th cycle of length $\lambda_i$:
\[
\forall m\geq1, \ \forall \, 1\leq i \leq k, \ \forall \, 1\leq j\leq
m_i
\] 
\[\mathrm{Sh}^m_{\lambda_i,j}(\Pi_b)(\beta)=
\begin{cases} 
b_{\ell+m,\lam_i,j} & \text{if } \exists \ell\, / \Pi_b(\beta)=b_{\ell,\lam_i,j}\\
\Pi_b(\beta) & \text{otherwise. } 
\end{cases} 
\]
\item The \emph{exchange operator} $\mathrm{Ex}_{\lambda_i,j_1,j_2}$,
  that exchanges the labels of the $j_1$th and $j_2$th cycles of
  length $\lam_i$:
\[
\forall \, 1\leq i \leq k, \ \forall \, 1\leq j_1,j_2 \leq m_i
\]
\[
\mathrm{Ex}_{\lambda_i,j_1,j_2}(\Pi_b)(\beta)= \begin{cases} 
b_{\ell,\lam_i,j_2} & \text{if }\exists \ell\, / \Pi_b(\beta)=b_{\ell,\lam_i,j_1}  \\
b_{\ell,\lam_i,j_1} & \text{if }\exists \ell\, / \Pi_b(\beta)=b_{\ell,\lam_i,j_2}  \\
\Pi_b(\beta) & \text{otherwise} \\
\end{cases}
\]
\end{itemize}
\end{definition}

\begin{proposition}[Set of consistent labellings]
\label{pro_set_lab}
Let $\s$ be a permutation with invariant
$(\lam=\{\lam_1^{m_1},\ldots,\lam_k^{m_k}\},r)$. Two consistent
labellings are obtained from one another by a sequence of shift and
exchange operators.
\end{proposition}

\begin{lemma}
\label{pro_size_consistent}
Let $\s$ be a permutation with invariant
$(\lam=\{\lam_1^{m_1},\ldots,\lam_k^{m_k}\},r)$ and let $\Pi$ be the
set of consistent labellings. Then
\[
|\Pi|=\prod_{i=1}^k ( m_i! \lam_i^{m_i} )
\ef.
\]
\end{lemma}

%-------------------------------------------------------
\subsection{The amenability}

Let us define two families of permutations (in cycle notation):
\begin{align*}
\gamma_{t,n}(i) &=
(1)(2)\cdots(i)(i+1\ n\ n-1 \ \cdots \ i+2)
\\
\gamma_{b,n}(i) &=
(2\ 3 \ \cdots \  i-1\ 1)(i)(i+1)\cdots(n)
\end{align*}
Note that, with this notation, 
$L(\s)=\gamma^{-1}_{t,n}(\s(1))\circ \s$ 
and $R(\s)=\s \circ \gamma_{b,n}(\s^{-1}(n))$.

\begin{figure}[tb]
\begin{center}
\begin{tabular}{cc}
 \raisebox{10pt}{\includegraphics[scale=.5]{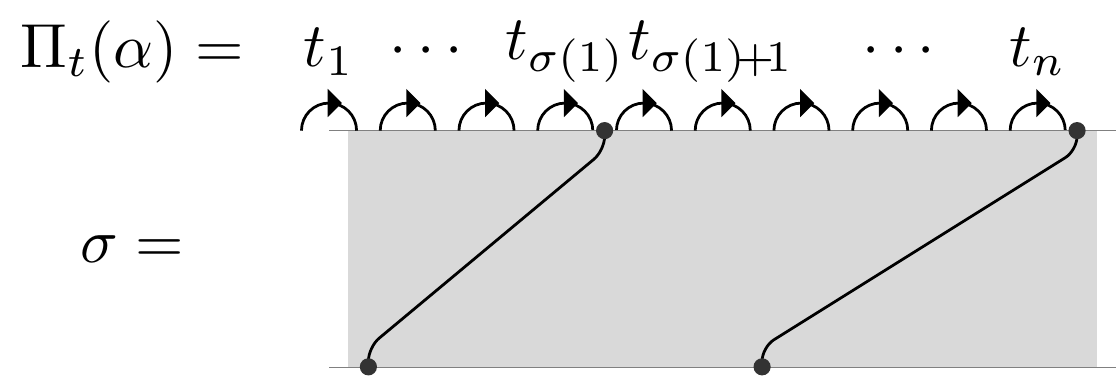}}& \multirow{-5}{*}{\vspace{0in} \includegraphics[scale=.5]{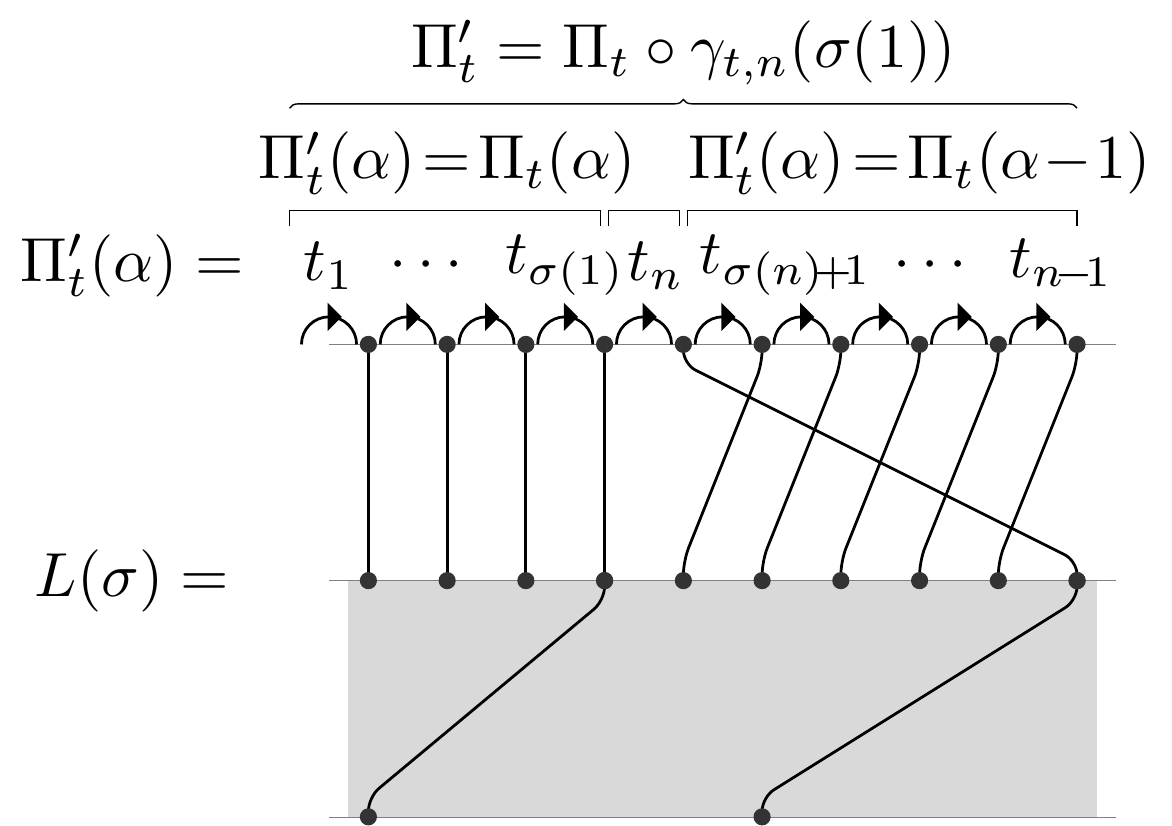}}  \\
  \raisebox{0pt}{\includegraphics[scale=.5]{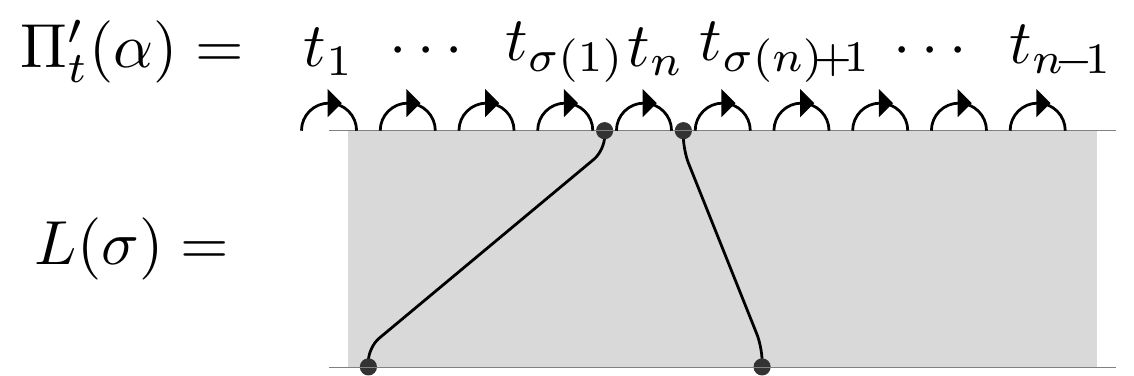}} & \\
\end{tabular}
\end{center}
\caption{\label{fig_definition_piprime} The action on the dynamics
  (with an $L$ operator) on a consistent labelling.}
\end{figure}

The dynamics $\cS_n$ can be naturally extended in order to act also on
the labelling.
\begin{definition}[Action of the dynamics $\cS_n$ on the labelling]
\label{def.actSnLabe}
Let $(\s, (\Pi_b,\Pi_t))$ be a permutation with a consistent
labelling. Then $L(\s, (\Pi_b,\Pi_t))= (L(\s),(\Pi_b,\Pi_t\circ
\gamma_{t,n}(\s(1)))$ 
and  $R(\s, (\Pi_b,\Pi_t))= (R(\s),(\Pi_b \circ \gamma_{b,n}(\s^{-1}(n)),\Pi_t))$.
\end{definition}
\noindent
Refer to Figure \ref{fig_definition_piprime} for an illustration.

There are two good reasons for this to be `the correct way' of
extending the dynamics to the labelled case.  First, it is compatible
with the structure of the cycle invariant, as the image by the
dynamics of a consistent labelling is another consistent labelling
(this is the content of Theorem~\ref{thm_consistent_lab}). Moreover,
as is crucially needed by our methods, it is compatible with the
boosted dynamics (see Theorem~\ref{thm_keep_track_label}).

\begin{theorem}\label{thm_consistent_lab}
Let $(\s,(\Pi_b,\Pi_t))$ be a permutation with a consistent labelling.
% and let $X\in \{L,R\}$ 
Then $L(\s,(\Pi_b,\Pi_t))$ and $R(\s,(\Pi_b,\Pi_t))$ 
are permutations with a consistent labelling.
\end{theorem}

\begin{pf}
By symmetry, we can consider just the action of the operator $L$.
Let us show that 
$( L(\s),(\Pi'_b=\Pi_b,\Pi'_t=\Pi_t\circ \gamma_{t,n}(\s(1)) )$ is a
permutation with a consistent labelling. First note that
\[
L(\s)(i)=\begin{cases}
\s(i) & \text{if } \s(i)\leq \s(1),\\
 \s(i)+1 & \text{if } \s(1)+1\leq \s(i) \leq n-1\\
 \s(1)+1 & \text{if } \s(i)=n,\\
\end{cases}
\]
By Lemma \ref{lem_verif_consistent}, we need to check that $\Pi'_b$
verifies property 2, $(\Pi'_b,\Pi'_t)$ verifies property 3, and
$\Pi'_t(1)=t^{\rm rk}_0$.

Let $\alpha>1$ be a top arc. There are three possibilities:
$\alpha\leq \s(1)$, $\alpha=\s(1)+1$ and $\alpha\geq \s(1)+2$.
\begin{itemize}
\item If $\alpha\leq \s(1)$ then $\Pi'_t(\alpha)=\Pi_t(\alpha)$ by
  definition of $\Pi'_t$ (cf.\ Figure \ref{fig_definition_piprime}
  right). Let $\beta=L(\s)^{-1}(\alpha-1)$ and
  $\beta'=L(\s)^{-1}(\alpha)$ be the two consecutive bottom arcs
  associated to $\alpha$ in $L(\s)$. We also have
  $\beta=\s^{-1}(\alpha-1)$ and $\beta'=\s^{-1}(\alpha)$ by definition
  of $L(\s)$, since $\alpha\leq \s(1)$.  Thus
  $\Pi'_b(\beta)=b_{i,\lam_\ell,j}$ (or $b^{\rm rk}_i$) and
  $\Pi'_b(\beta')=b_{i+1\mod \lam_\ell,\lam_\ell,j}$ (or 
  $b^{\rm rk}_{i+1}$), by property 2 for $(\s,\Pi_b)$, since
  $\Pi'_b=\Pi_b$ and $\beta$ and $\beta'$ are consecutive in $\s$.

  Likewise, $\Pi'_t(\alpha)=t_{i+1\mod \lam_\ell,\lam_\ell,j}$ (or
  $t^{\rm rk}_{i+1}$) and
  $\Pi'_b(\beta')=b_{i+1\lam_\ell,\lam_\ell,j}$ (or 
  $b^{\rm rk}_{i+1}$) by property 3 for $(\s,(\Pi_b,\Pi_t))$, since
  $\Pi'_b=\Pi_b$ and $\Pi'_t(\alpha)=\Pi_t(\alpha)$.

  Thus in this case $\Pi'_b=\Pi_b$ verifies property 2 and
  $(\Pi'_b,\Pi'_t)$ verifies property 3 (see also figure
  \ref{fig_ex_pf_constitent_1}).

  Note that the same reasoning applies for $\alpha=1$, from which we
  deduce that $\Pi'_t(1)=t^{\rm rk}_0$ and
  $\Pi'_b(L(\s)^{-1}(1))=b^{\rm rk}_0$.
\begin{figure}[h!]
{\centering
\begin{tabular}{cc}
\includegraphics[scale=.5]{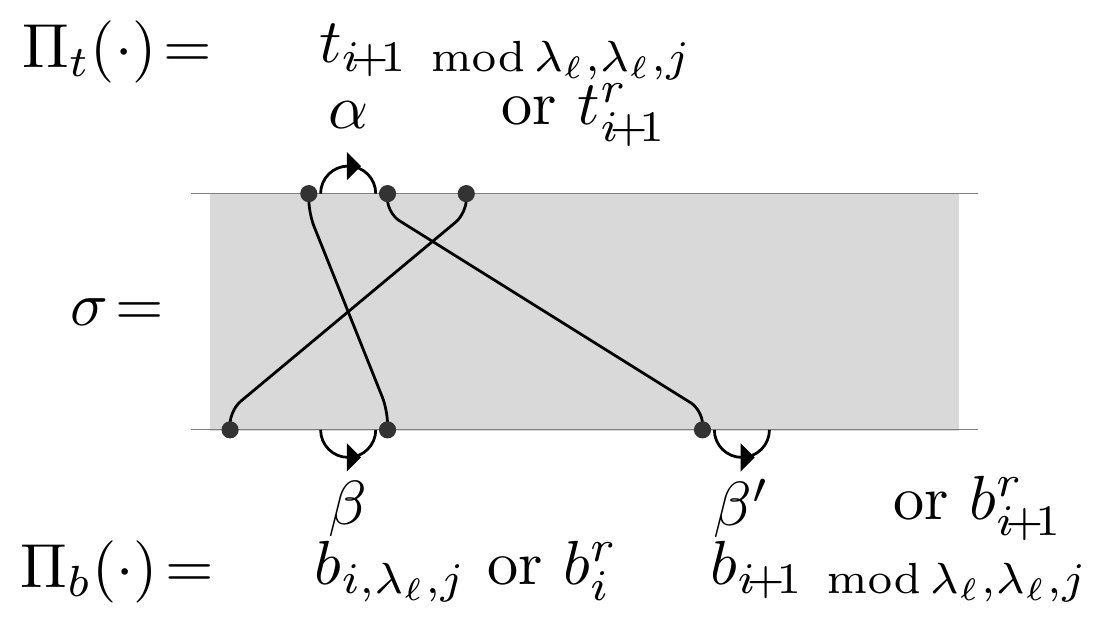}& \includegraphics[scale=.5]{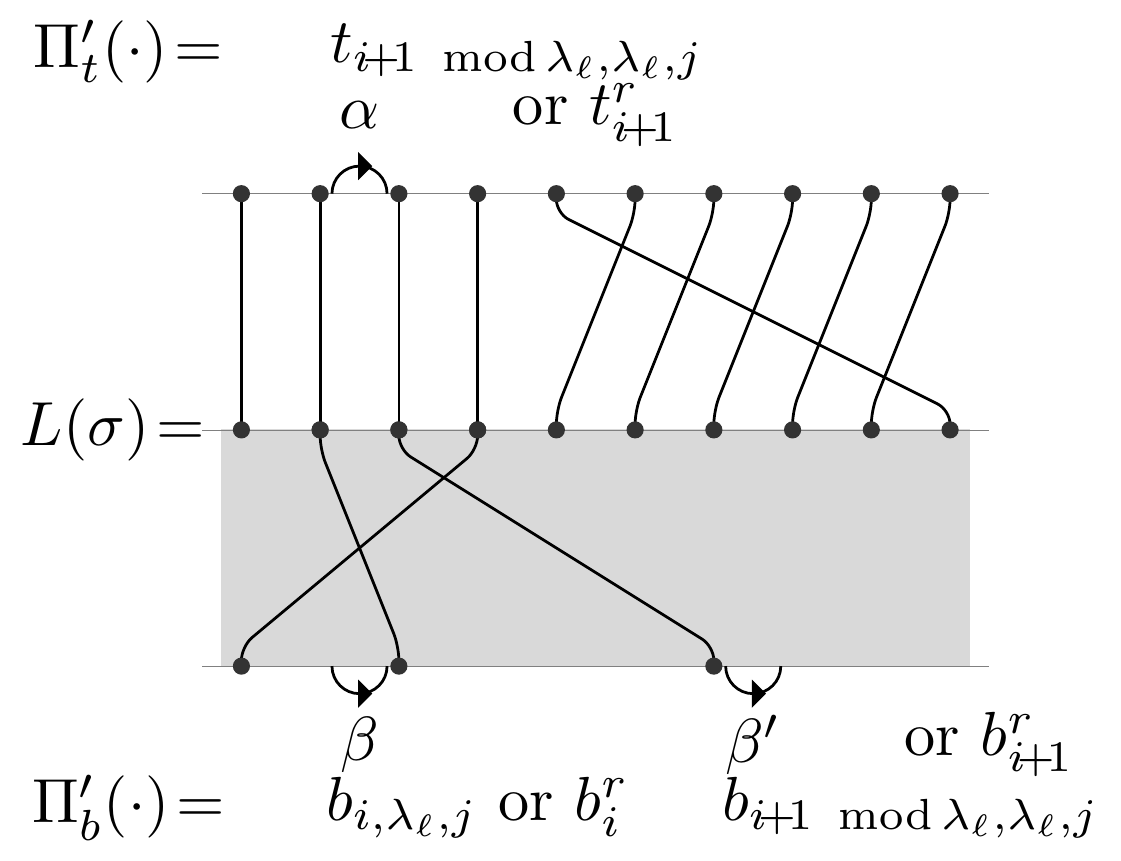}
\end{tabular}
\caption{\label{fig_ex_pf_constitent_1} The proof that $\Pi'_b$
  verifies property 2 and $(\Pi'_b,\Pi'_t)$ verifies property 3 for
  the case $\alpha\leq \s(1)$.}
}
\end{figure}
\item If $\alpha= \s(1)+1$, let $\alpha'=n$. Then
  $\Pi'_t(\alpha)=\Pi_t(\alpha')$ by definition of $\Pi'_t$
  (cf. Figure \ref{fig_definition_piprime}). Let
  $\beta'=L(\s)^{-1}(n)$ and $\beta'=L(\s)^{-1}(\alpha)$ be the two
  consecutive bottom arcs associated to $\alpha$ in $L(\s)$ (this is
  the special case of consecutive arcs, cf.\ the third figure in
  Definition \ref{def_consecutive}). We have $\beta=\s^{-1}(n-1)$ and
  $\beta'=\s^{-1}(n)$ by definition of $L(\s)$, thus $\beta$ and
  $\beta'$ are also the two consecutive bottom arcs associated to
  $\alpha$ in $\s$.

  Thus $\Pi'_b(\beta)=b_{i,\lam_\ell,j}$ (or $b^{\rm rk}_i$) and
  $\Pi'_b(\beta')=b_{i+1\mod \lam_\ell,\lam_\ell,j}$ (or 
  $b^{\rm rk}_{i+1}$), by property 2 for $(\s,\Pi_b)$, since
  $\Pi'_b=\Pi_b$ and $\beta$ and $\beta'$ are consecutive in $\s$.

  Likewise, $\Pi'_t(\alpha)=t_{i+1\mod \lam_\ell,\lam_\ell,j}$ (or
  $t^{\rm rk}_{i+1}$) and
  $\Pi'_b(\beta')=b_{i+1\lam_\ell,\lam_\ell,j}$ (or 
  $b^{\rm rk}_{i+1}$), by property 3 for $(\s,(\Pi_b,\Pi_t))$, since
  $\Pi'_b=\Pi_b$ and $\Pi'_t(\alpha=\s(1)+1)=\Pi_t(\alpha'=n)$.

  Thus, also in this case $\Pi'_b=\Pi_b$ verifies property 2 and
  $(\Pi'_b,\Pi'_t)$ verifies property 3 (see also Figure
  \ref{fig_ex_pf_constitent_2}).
\begin{figure}[h!]
{\centering
\begin{tabular}{cc}
\includegraphics[scale=.5]{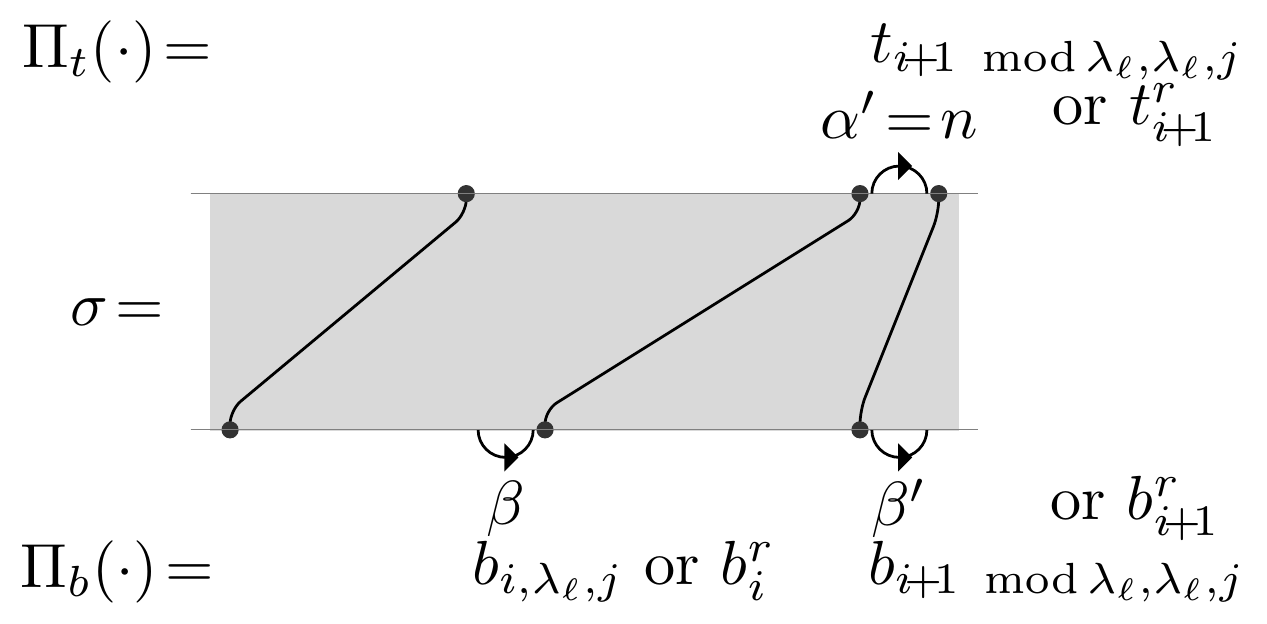}& \includegraphics[scale=.5]{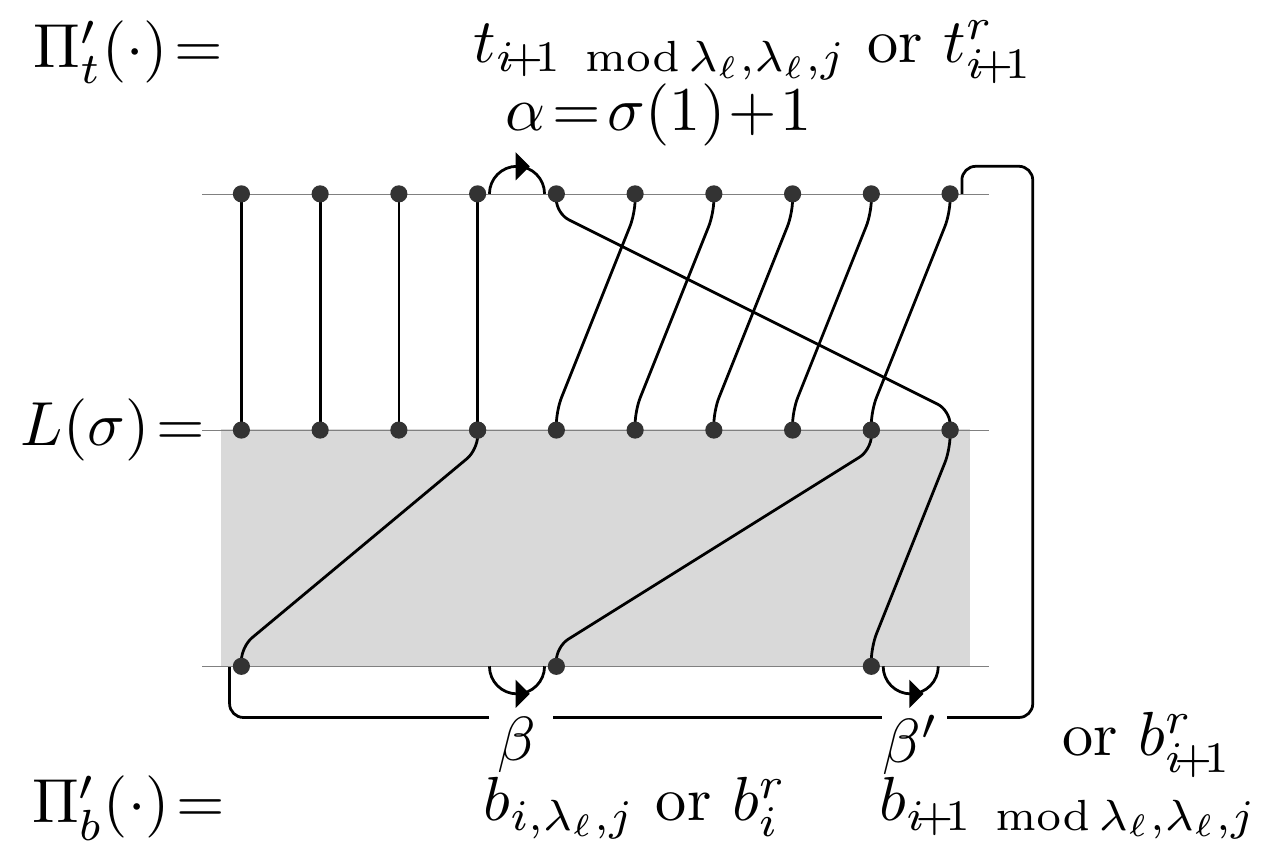}
\end{tabular}
\caption{\label{fig_ex_pf_constitent_2} The proof that $\Pi'_b$
  verifies property 2 and $(\Pi'_b,\Pi'_t)$ verifies property 3 for
  the case $\alpha= \s(1)+1$.}
}
\end{figure}

\item If $\alpha\geq \s(1)+2$, let $\alpha'=\alpha-1$. Then
  $\Pi'_t(\alpha)=\Pi_t(\alpha')$ by definition of $\Pi'_t$ (cf. Figure
  \ref{fig_definition_piprime}, right). Let
  $\beta=L(\s)^{-1}(\alpha-1)$ and $\beta'=L(\s)^{-1}(\alpha)$ be the
  two consecutive bottom arcs associated to $\alpha$. We have
  $\beta=\s^{-1}(\alpha')$ and $\beta'=\s^{-1}(\alpha'-1)$ by
  definition of $L(\s)$, thus $\beta$ and $\beta'$ are also the two
  consecutive bottom arcs associated to $\alpha$ in $\s$.

  Thus $\Pi'_b(\beta)=b_{i,\lam_\ell,j}$ (or $b^{\rm rk}_i$) and
  $\Pi'_b(\beta')=b_{i+1\mod \lam_\ell,\lam_\ell,j}$ (or 
  $b^{\rm rk}_{i+1}$), by property 2 for $(\s,\Pi_b)$, since
  $\Pi'_b=\Pi_b$ and $\beta$ and $\beta'$ are consecutive in $\s$.

  Likewise $\Pi'_t(\alpha)=t_{i+1\mod \lam_\ell,\lam_\ell,j}$ (or
  $t^{\rm rk}_{i+1}$) and
  $\Pi'_b(\beta')=b_{i+1\lam_\ell,\lam_\ell,j}$ (or 
  $b^{\rm rk}_{i+1}$) by property 3 for $(\s,(\Pi_b,\Pi_t))$, since
  $\Pi'_b=\Pi_b$ and $\Pi'_t(\alpha)=\Pi_t(\alpha')$.
%\begin{figure}[h!]
%{\centering
%\begin{tabular}{cc}
%\includegraphics[scale=.5]{figure/fig_L_rauzy_lab_pf_2_1}& \includegraphics[scale=.5]{figure/fig_L_rauzy_lab_pf_2_2}
%\end{tabular}
%\caption{}
%}
%\end{figure}
\end{itemize}
\qed
\end{pf}

\begin{corollary}\label{cor_consistent_shift_exchange}
Let $(\s,\Pi_b)$ be a permutation with a consistent labelling, and let
$S$ be a `loop' (i.e.\ a sequence in $\{L,R\}^*$ such that
$S(\s)=\s$). Then $S(\Pi_b)$ is a consistent labelling on $\s$, and is
obtained from $\Pi_b$ by a sequence of shift and exchange operators.
\end{corollary}
\begin{pf}
This is immediate from Theorem \ref{thm_consistent_lab} and
Proposition \ref{pro_set_lab}. \qed
\end{pf}

Let us now prove the compatibility of the labelling with the boosted dynamics.

\begin{figure}[tb]
\begin{center}
$\begin{tikzcd}[%
    ,row sep = 3ex
    ,/tikz/column 1/.append style={anchor=base east}
    ,/tikz/column 2/.append style={anchor=base east}
    ]
 \raisebox{-15pt}{\includegraphics[scale=.5]{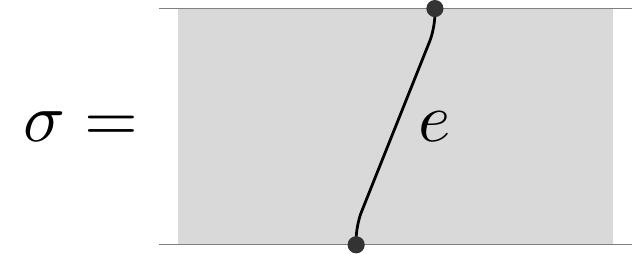}}\arrow{d}{red} \arrow{r}{B(S)} &  \raisebox{-15pt}{\includegraphics[scale=.5]{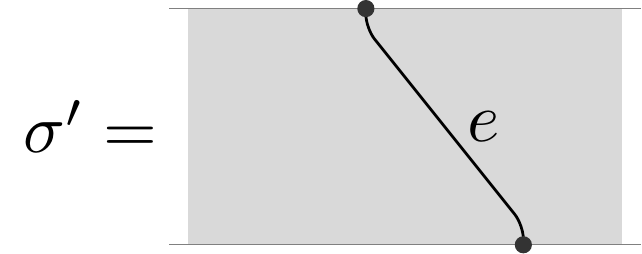}}\arrow{d}{red}   \\
  \raisebox{-32pt}{\includegraphics[scale=.5]{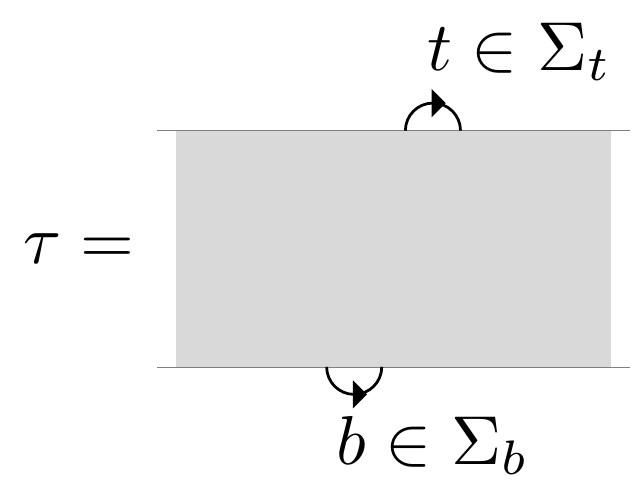}}\arrow{r}{S} &    \raisebox{-32pt}{\includegraphics[scale=.5]{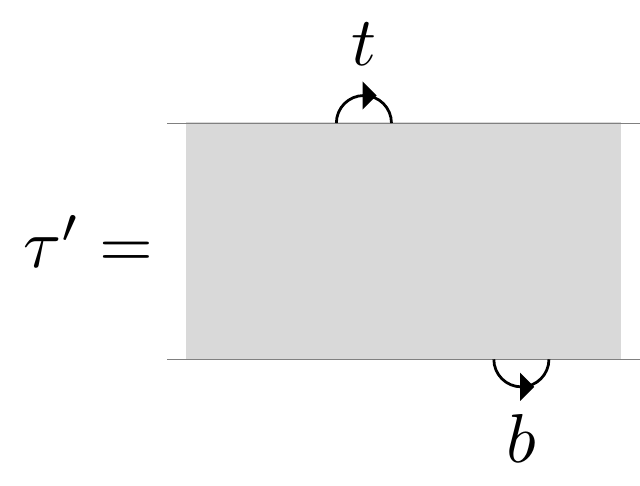}} \\[-.5cm]
\end{tikzcd}$
\end{center}
\caption{\label{fig_keep_track_label} The gray edges of $\s$ are
  transported along the labels of the arcs when applying the boosted
  dynamics. For instance, $e\in(b,t)$ in $\tau,\Pi$.}
\end{figure}

\begin{theorem}[The labelling is compatible with the boosted dynamics]
\label{thm_keep_track_label}
Let $\s\in \kS_n$ be a permutation, and let $\tau\in\kS_k$ be the
reduction of a $(2k,2r)$-coloring $c$ of $\s$. 
% Hence the gray edges of $\s$ are inserted within different bottom
% and top arcs of $\tau$.
Let $(\Pi_b,\Pi_t)$ be a consistent labelling of $\tau$, and let $S$
be a sequence of operators for the dynamics $\kS_k$, acting on
$(\tau,(\Pi_b,\Pi_t))$ as described in Definition \ref{def.actSnLabe}.
Let $e$ be a gray edge of $\s$ with endpoints within the arcs with
labels $t=\Pi_t(\alpha)$ and $b=\Pi_b(\beta)$ of $\tau$, for some
$\alpha$ and $\beta$.

Then, in $(\s',c')=B(S)(\s,c)$ the gray edge $e$ has its endpoints
within the arcs with labels $t$ and $b$ of
$(\tau',(\Pi'_b,\Pi'_t))=S(\tau,(\Pi_b,\Pi_t))$. See Figure
\ref{fig_keep_track_label}.
\end{theorem}

\begin{pf}
By induction on the length of the sequence, and by symmetry, we can
consider just the case $S=L$.

Let us consider $\s\in \kS_n$, $c$ a $(2k,2r)$-coloring of $\s$,
$\tau\in \kS_k$ the corresponding reduction and $(\Pi_b,\Pi_t)$ a
consistent labelling of $\tau$, as in the statement of the theorem.

Let $e \in (t,b)$ be a gray edge of $(\s,c)$, which is inserted within
the top arc with label $t$ and bottom arc with label $b$ of 
$(\tau,(\pi_b,\Pi_t))$.

We must show that the gray edge $e$ of $(\s',c')=B(L)(\s,c)$ is still
within the arcs with labels $t$ and $b$ of
$(\tau',\Pi_b',\Pi_t')=L(\tau,\Pi_b,\Pi_t)$.

First, let $k_2$ be the size of the block of gray edges immediately to
the left of the right pivot (that is, the edge $(\s^{-1}(n),n)$).  By
the mechanisms of the boosted dynamics, in $(\s,c)$ we have
$B(L)=L^{k_2}$ (see Figure \ref{fig_keep_track_label_for_L}).  We
shall consider three cases, depending on the position of $\alpha$, the
arc labelled by $t$:
\begin{itemize}
\item If $\tau(1)<\alpha=\Pi^{-1}_t(t)<k$, then for
  $(\tau',(\Pi'_b,\Pi'_t))=L(\tau,(\Pi_b,\Pi_t))$ we have
  $\Pi'^{-1}_t(t)=\alpha+1$. Since $e=(i,\sigma(i))$ is within the
  arcs $t$ and $b$ of $(\tau,(\Pi_b,\Pi_t))$, we have
  $\s(1)<\s(i)<n-k_2$, thus $L^{k_2}(e)=(i,\s(i)+k_2)$ in
  $(\s',c')=L^{k_2}(\s,c)$. Therefore in $(\s',c')$ the edge $e$ is
  inserted within the arcs with labels $t$ and $b$ of
  $(\tau',(\Pi'_b,\Pi'_t))$. See
  Figure~\ref{fig_keep_track_label_for_L} with the edge $e_1$ inserted
  between $t_1$ and $b_1$.

\item If $\alpha=\Pi^{-1}_t(t)<\tau(1)$, we are in a situation which
  is almost identical to the previous one, and we omit to discuss it.

\item If $\alpha=\Pi^{-1}_t(t)=n$, then for
  $(\tau',(\Pi'_b,\Pi'_t))=L(\tau,(\Pi_b,\Pi_t))$ we have
  $\Pi'^{-1}_t(t)=\tau(1)+1$. Since $e=(i,\sigma(i))$ is inserted
  within $t$ and $b$ in $(\tau,(\Pi_b,\Pi_t))$, we have
  $n-k_2\leq\s(i)<n$, thus $L^{k_2}(e)=(i,j)$ with $j$ such that
  $s(1)<j \leq \s(1)+k_2$ in $(\s',c')=L^{k_2}(\s,c)$. Thus in
  $(\s',c')$ the gray edge $e$ is inserted within the arcs with labels
  $t$ and $b$ of $\tau'$. See Figure~\ref{fig_keep_track_label_for_L}
  with the edge $e_2$ inserted within $t_2$ and $b_2$.
\end{itemize}

\begin{figure}[h]
\begin{center}
$\begin{tikzcd}[%
    ,row sep = 3ex
    ,/tikz/column 1/.append style={anchor=base east}
    ,/tikz/column 2/.append style={anchor=base east}
    ]
 \raisebox{-15pt}{\includegraphics[scale=.5]{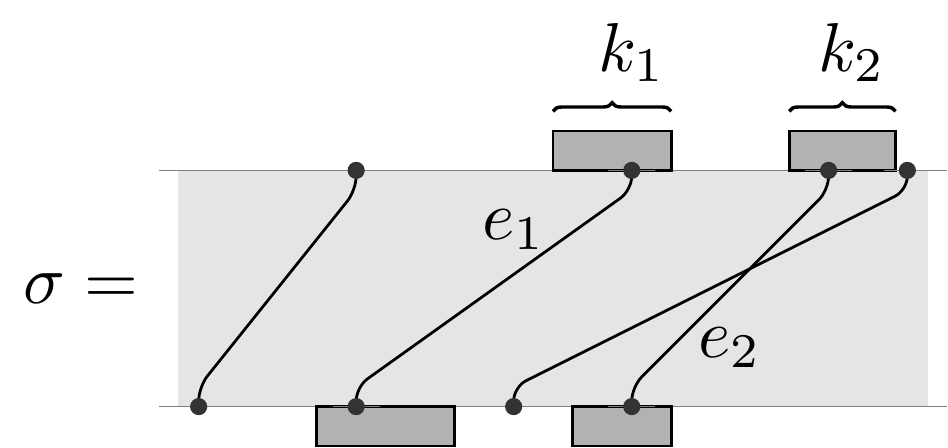}}\arrow{d}{red} \arrow{r}{B(L)=L^{k_2}} &  \raisebox{-15pt}{\includegraphics[scale=.5]{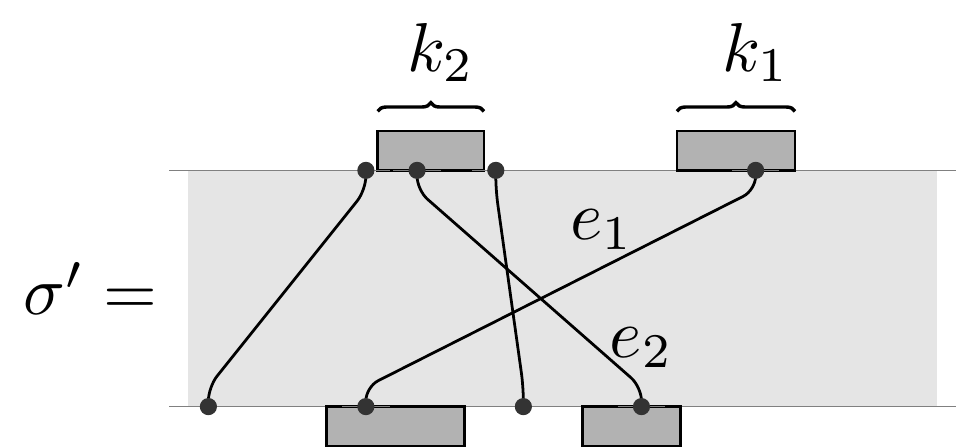}}\arrow{d}{red}   \\
  \raisebox{-32pt}{\includegraphics[scale=.5]{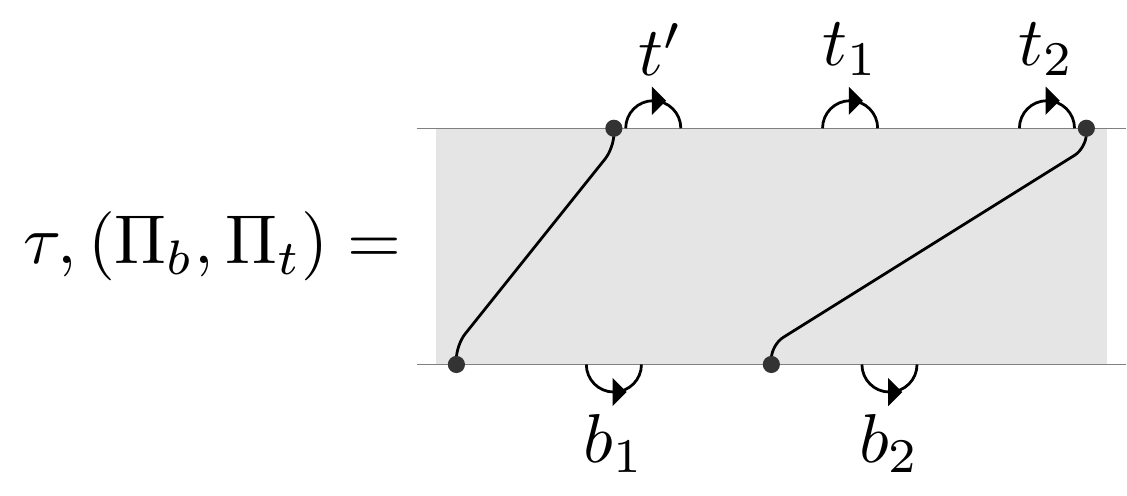}}\arrow{r}{L} &    \raisebox{-32pt}{\includegraphics[scale=.5]{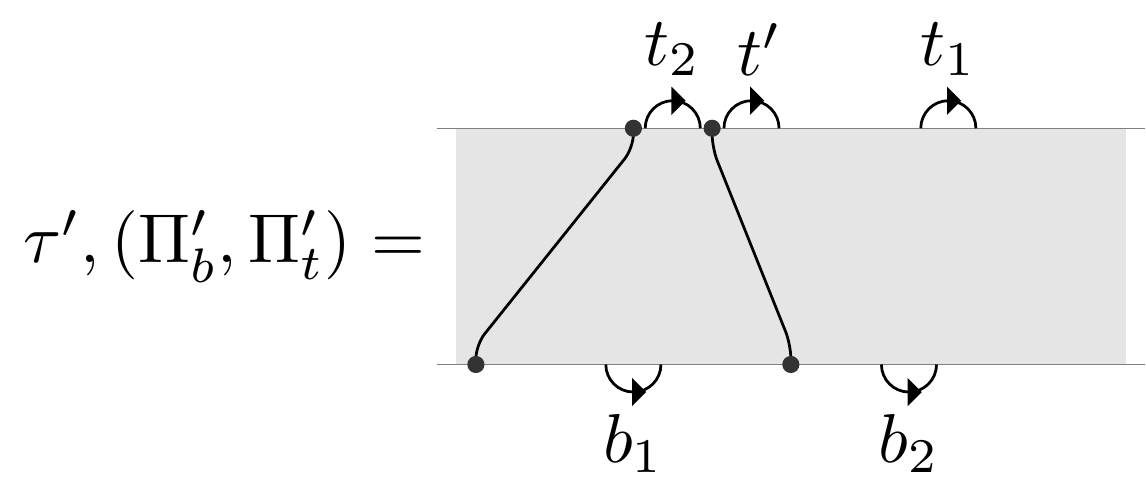}}\\[-.5cm]
\end{tikzcd}$
\end{center}
\caption{\label{fig_keep_track_label_for_L} Illustration that the edges, inserted within the labels of a reduced permutation, follow those labels when applying the boosted dynamic.}
\end{figure}
\qed
\end{pf}

The last main task of this section is to introduce the 2-point
monodromy theorem. As announced in the overview section,
Sec.\ \ref{ssec_pf_overview}, this theorem is the seventh (and last)
statement in the organisation of the main induction, in
Section~\ref{sec_induction}.

\begin{theorem}\label{thm_2monodromy}
Let $(\s,\Pi_b)$ be a permutation with invariant
$(\lam=\{\lam_1^{m_1},\ldots,\lam_{k}^{m_{k}}\},r)$,
equipped of a consistent labelling.
\begin{description}
\item{Cycle 1-shift:} Suppose $\lam$ does not contain any even cycle,
  or it has 2 or more even cycles. Let $i$ be a cycle of length
  $\lam_i$ in $\lam$. Then there exists a loop $S$ in the dynamics
  such that $\Pi'_b=S(\Pi_b)$ is a consistent labeling, and
  verifies 
  $\Pi'_b(\beta)= b_{\ell+1\mod \lam_i,\lam_i,i} \Leftrightarrow
  \Pi_b(\beta)=b_{\ell,\lam_i,i}$.
%% $\forall \beta$:
%% \[\Pi'_b(\beta)= \begin{cases} 
%% b_{\ell+1\mod \lam_i,\lam_i,i} & \text{if }\exists \ell\, / \Pi_b(\beta)=b_{\ell,\lam_i,i}  \\
%% \Pi'_b(\beta) & \text{otherwise}
%% \end{cases}\]
  In other words, there exists a loop in the dynamics that shifts the
  labels of the cycle $i$ by one (and thus, by taking powers, by any
  integer $m$). The positions of the other labels are, in principle,
  unknown, nonetheless they are constrained by the fact that $\Pi'_b$
  is also a consistent labelling.

\item{Cycle 2-shift:} Suppose now that $\lam$ has exactly one even
  cycle, of length $\lam_i$. Then there exists a loop $S$ such that
  $\Pi'_b=S(\Pi_b)$ verifies 
  $\Pi'_b(\beta)= b_{\ell+2\mod \lam_i,\lam_i,1} \Leftrightarrow
  \Pi_b(\beta)=b_{\ell,\lam_i,1}$.
%% $\forall \beta$:
%% \[\Pi'_b(\beta)= \begin{cases} 
%% b_{\ell+2\mod \lam_i,\lam_i,i} & \text{if }\exists \ell\, / \Pi_b(\beta)=b_{\ell,\lam_i,i}  \\
%% \Pi'_b(\beta) & \text{otherwise}
%% \end{cases}\]
In other words, there exists a loop that shifts the labels of the
unique even cycle by two (and thus by any even integer $2m$). This is
consistent with the fact that, if and only if $\lam_i$ is odd, then iterated
shifts by $2$ ultimately produce a shift by $1$.
%there exists a loop that shifts the labels by any $m$ (since 2 and $\lam_i$ are co-prime).

\item{Cycle jump:} For any two cycles $j_1$, $j_2$ of the same length
  $\lam_i$, there exists a loop $S$ such that
  $\Pi'_b=S(\Pi_b)$ verifies
  $\Pi'_b(\beta)= b_{\ell,\lam_i,j_2} \Leftrightarrow
  \Pi_b(\beta)=b_{\ell,\lam_i,j_1}$.
%%  $\forall \beta$:
%% \[\Pi'_b(\beta)= \begin{cases} 
%% b_{\ell+m\mod \lam_i,\lam_i,j_2} & \text{if }\exists \ell\, / \Pi_b(\beta)=b_{\ell,\lam_i,j_1}  \\
%% \Pi'_b(\beta) & \text{otherwise}
%% \end{cases}\]
In other words, there exists a loop that sends the labels of the cycle
$j_2$ on the positions of the labels of the cycle $j_1$, preserving
their ordering.
% with a shift of length $m$. 
\end{description}
\end{theorem}

It might not be clear a priori why such theorem implies the 2-point
monodromy. Indeed, in our proof we will choose the label of the top
arc to be $t^{\rm rk}_0$, thus, since it is fixed by Theorem
\ref{thm_consistent_lab}, we only need to consider the label of the
bottom arc. In this case we have the immediate corollary:

\begin{corollary}[2-point monodromy theorem]
\label{cor_2monodromy}
Let $(\s,\Pi_b)$ be a permutation with invariant
$(\lam=\{\lam_1^{m_1},\ldots,\lam_{k}^{m_{k}}\},r)$ and with a
consistent labelling. Let $(t,b)\in \Sigma_t \times \Sigma_b$ such
that $t=t^{\rm rk}_0$.
\begin{itemize}
\item If $b=b^{\rm rk}_i$ for some $i$, there exists a loop $S$ such
  that $\Pi'_b=S(\Pi_b)$ verifies $\Pi'_b(\beta)=b$ if and only if
  $\Pi_b(\beta)=b$.
\item Let us now distinguish two cases: \\
(1) $\lam$ has no even cycles or has at least two of them. Then:
\begin{itemize}
\item  If $b=b_{\ell,\lam_i,j_1}$ for some triple $(\ell,\lam_i,j_1)$,
  there exists a loop $S$ such that $\Pi'_b=S(\Pi_b)$ verifies
  $\Pi'_b(\beta)=b$ if and only if $\Pi_b(\beta)=b_{\ell',\lam_i,j_2}$
  for some $(\ell',j_2)$. 
\end{itemize}
(2) $\lam$ has exactly one even cycle, the cycle $i$ of length
$\lam_i$. Then
\begin{itemize}
\item  If $b=b_{\ell,\lam_j,j_1}$ for some $(\ell,\lam_j,j_1)$ with
  $j\neq i$, there exists a loop $S$ such that $\Pi'_b=S(\Pi_b)$
  verifies $\Pi'_b(\beta)=b$ if and only if
  $\Pi_b(\beta)=b_{\ell',\lam_j,j_2}$ for some $(\ell',j_2)$. 
\item  If $b=b_{\ell,\lam_i,1}$ for some $(\ell,\lam_i,1)$, there
  exists a loop $S$ such that $\Pi'_b=S(\Pi_b)$ verifies
  $\Pi'_b(\beta)=b$ if $\Pi_b(\beta)=b_{\ell',\lam_i,1}$ with
  $\ell'-\ell \equiv 0 \pmod 2$.
\end{itemize}
\end{itemize}
In both cases, $\Pi_t(0)=\Pi'_t(0)=t^{\rm rk}_0$ since the labels of
the rank are fixed in a consistent labelling.
\end{corollary}
Let us rephrase the content of this corollary. If we choose the label
of the top arc to be $t^{\rm rk}_0$, then if the bottom label is in
the rank it is fixed, otherwise, if it is in a cycle, we can find a
loop $S$ such that the label can move to any bottom label of a cycle
of the same length, with one exception: if the bottom label is in the
unique cycle of even length. In this case, we can only find a loop $S$
such that the label can move to the bottom labels of the cycle which
has the same parity of the index.

\proof
By Theorem \ref{thm_consistent_lab}, if $b=b^{\rm rk}_i$ then, since
$\Pi'_b$ is a consistent labelling and the labels of the rank are
fixed, $\Pi'_b(\beta)=b$ if and only if $\Pi_b(\beta)=b$.

If $b=b_{\ell,\lam_i,j_1}$, by Theorem \ref{thm_consistent_lab}, if
$\Pi'_b=S(\Pi_b)$ for some loop $S$, then it is a consistent
labelling, thus we must have $\Pi'_b(\beta)=b$ if and only if
$\Pi'_b(\beta)=b_{\ell',\lam_i,j_2}$ for some $\ell'$ and $j_2$, this
by Corollary~\ref{cor_consistent_shift_exchange}.

The rest of the statement is a straightforward application of
Theorem~\ref{thm_2monodromy}.
\qed

We can actually completely characterise the set $L(\s,(\Pi_b,\Pi_t))$. The following theorem answers this problem:

\begin{theorem}[Characterisation of $L(\s,(\Pi_b,\Pi_t))$ ]\label{thm_monodromy_total}
Let $\s$ be a permutation with invariant $(\lam=\{\lam_1^{m_1},\ldots,\lam_k^{m_k}\},r,s)$, and let $(\Pi_b,\Pi_t)$ be a consistent labelling. Then $L(\s,(\Pi_b,\Pi_t))$ consists of the labellings generated by: 
\begin{itemize}
\item Any exchange operators $Ex_{\lam_i,j_1,j_2}$.
\item Any 1-shift operators for cycles of odd length $\lam_i$ $Sh^1_{\lam_i,j}$.
\item Any 2-shift operators for cycles of even length $\lam_i$ $Sh^{2}_{\lam_i,j}$.
\item Any 1-shift operators on a pair of cycles of even lengths $\lam_{i_1},\lam_{i_2}$ $Sh^1_{\lam_{i_1},j_1}\circ Sh^1_{\lam_{i_2},j_2}$.
\end{itemize}
\end{theorem}

This theorem was the content of \cite{Boi13} (however we cannot use it in this article as the proof require the classification theorem). We can derive an alternative proof using the techniques of section \ref{ssec_induction_monodromy} (see remark \ref{rk_monodromy_thm_total}) but we will not do so in this paper.

\section{Cycle invariant and edge addition\label{sec_edge_addition}}

This preliminary section study the change of the cycle invariant when inserting a few consecutives edges  in a permutation. The results of this section will be used in sections \ref{sec_arf},\ref{sec_IIXshaped},\ref{sec_induction}.

\begin{figure}[b!]
\begin{center}
\includegraphics{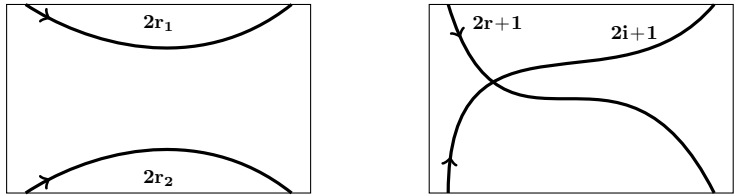}
\caption{\label{fig.typeXH}Left: a schematic representation of a
  permutation of type $H(r_1,r_2)$. Right: a representation of a
  permutation of type $X(r,i)$. These configurations have rank
  $r_1+r_2-1$ and $r$, respectively.}
  \end{center}
\end{figure}

First we reintroduce a notation from \cite{DS17} (it was also defined as $m|a$ and $m_l \cdot m_r$ in \cite{Del13}).
\begin{definition}
\label{def.XHtype}
A permutation $\s$ is \emph{of type $H$} if the rank path goes through
the $-1$ mark, and \emph{of type $X$} otherwise. In the case of a
type-$X$ permutation, we call \emph{principal cycle} the cycle going
through the $-1$ mark.
\end{definition}
\noindent See also figure \ref{fig.typeXH}.

\begin{figure}[b!]
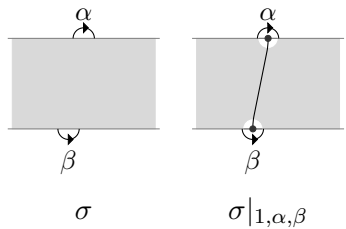

\begin{center}
\begin{tabular}{cc}
\includegraphics[scale=.5]{figure/fig_perm_add_one_edge}
&
\includegraphics[scale=.5]{figure/fig_perm_add_one_edge_1}\\
$\s$ 
&
$\s|_{1,\alpha,\beta}$
\end{tabular}
\end{center}
\caption{\label{fig_one_edge_example}The insertion of one edge within the arcs $\alpha$ and $\beta$.}
\end{figure}

\begin{notation}\label{not_add_edge}
let $\s$ be a permutation and let $\alpha$ be a top arc and $\beta$ be a bottom arc, we define $\s|_{i,\alpha,\beta}$ to be the permutation obtained from $\s$ by inserting $i\in \mathbb{N}$ consecutive and parallel edges within $\alpha$ and $\beta$. (see figure \ref{fig_one_edge_example} for an example with $i=1$).

For the special case of $i=2$ we call the two added edges a \emph{double-edge}.
\end{notation}

\begin{proposition}[One edge insertion]\label{pro_one_edge}
Let $\s$ be a permutation with cycle invariant $(\lambda,r)$. 

 Let $\s|_{1,\alpha,\beta}$ be the permutation resulting from the insertion of an edge within two arcs of two differents cycles of respective length $\ell$ and $\ell'$. Then the cycle invariant of  $\s|_{i,\alpha,\beta}$ is $\left(\lambda\setminus\{\ell,\ell'\}\bigcup\{\ell+\ell'+1\},r\right)$.

Let $\s|_{1,\alpha,\beta}$ be the permutation resulting from the insertion of an edge within an arc of the rank path and an arc of a cycle  (principal cycle or not) of length $\ell$. Then the cycle invariant of  $\s|_{i,\alpha,\beta}$ is $\left(\lambda\setminus\{\ell\},r+\ell+1\right)$.
\end{proposition}

\begin{pf}
See figure \ref{fig_add_one_edge}.
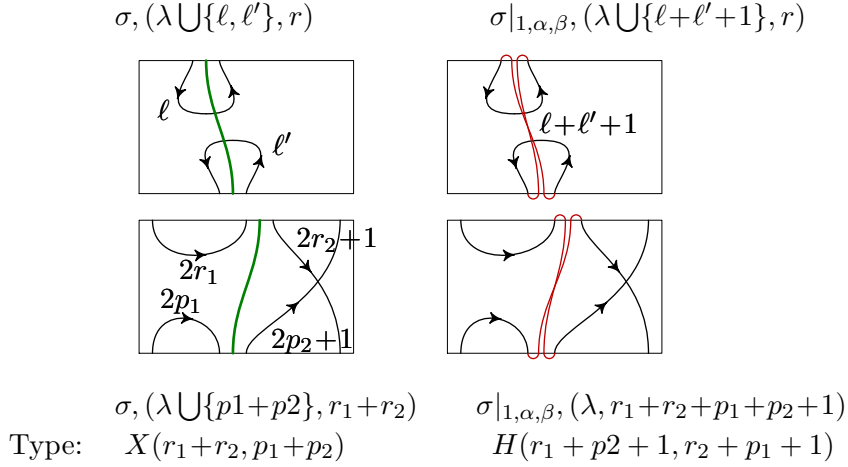
\begin{figure}[h!]
\begin{center}
\begin{tabular}{ll}
&$\s, (\lambda\bigcup\{\ell,\ell'\},r) \qquad \qquad\qquad \s|_{1,\alpha,\beta}, (\lambda\bigcup\{\ell\!+\!\ell'\!+\!1\},r)$\\[-2mm]
&\put(22,103){$\ell$}\put(22,103){$\ell$}
\put(65,90){$\ell'$}\put(65,90){$\ell'$}
\put(165,98){$\ell\!+\!\ell'\!+\!1$}\put(165,98){$\ell\!+\!\ell'\!+\!1$}
\put(30,43){$2r_1$}\put(30,43){$2r_1$}
\put(23,32){$2p_1$}\put(23,32){$2p_1$}
\put(74,55){$2r_2\!+\!1$}\put(74,55){$2r_2\!+\!1$}
\put(64.5,17){$2p_2\!+\!1$}\put(64.5,17){$2p_2\!+\!1$}
\includegraphics[scale=1]{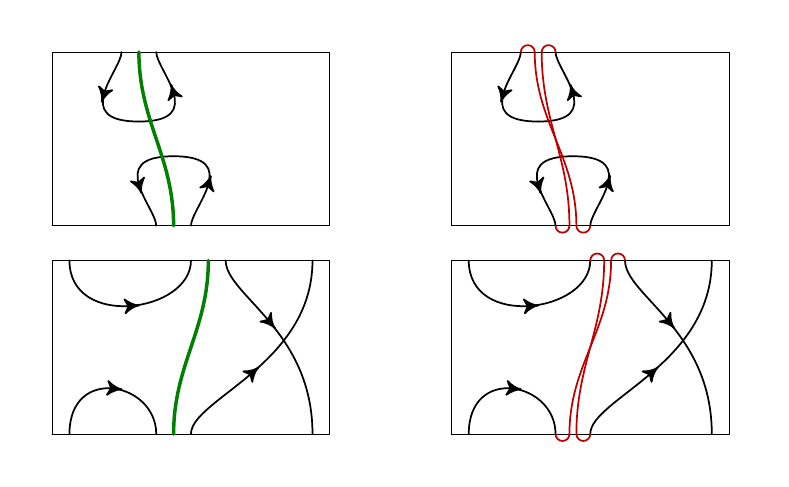}\\[-2mm]
&$\s, (\lambda\bigcup\{p1\!+\!p2\},r_1\!+\!r_2) \qquad  \s|_{1,\alpha,\beta}, (\lambda,r_1\!+\!r_2\!+\!p_1\!+\!p_2\!+\!1)$\\[.5mm]
Type: & $\ X(r_1\!+\!r_2,p_1\!+\!p_2) \qquad\qquad  \quad H(r_1+p2+1,r_2+p_1+1)$
\end{tabular}
\end{center}
\caption{\label{fig_add_one_edge}The first line represents the case: Top arc : any cycle. Bottom arc: any cycle. The second line represents the case: Top arc : rank path. Bottom arc: principal cycle.}
\end{figure}

 Some cases are not represented in the figure.
The missing cases are: \begin{itemize}
 \item Top arc: principal cycle. Bottom arc: any cycle. 
 \item Top arc: any cycle. Bottom arc: principal cycle. 
 \item Top arc: rank path. Bottom arc: any cycle. 
 \item Top arc: principal cycle. Bottom arc: rank path. 
 \item Top arc: any cycle. Bottom arc: rank path. 
\end{itemize} 
Their proof is nearly identical to the ones represented in the figure and are thus omitted. 
\qed
 \end{pf}

\begin{proposition}[double-edge insertion]\label{pro_double-edge}
Let $\s$ be a permutation with cycle invariant $(\lambda,r)$. 

Inserting a double-edge within two arcs of the same cycle (or rank path) increases the length of the cycle (or rank) by two.

 Likewise inserting a double-edge within two arcs of two differents cycles (or rank path) increases the length of each by 1. 
\end{proposition}

\begin{pf}
See figure \ref{fig_double-edge}.
\begin{figure}[h!]
\begin{center}
\begin{tabular}{ccc}
\includegraphics[scale=.5]{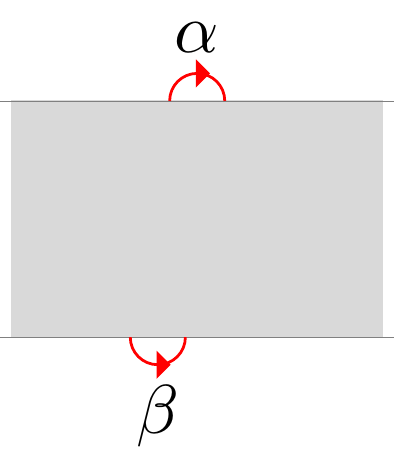}
&
\includegraphics[scale=.5]{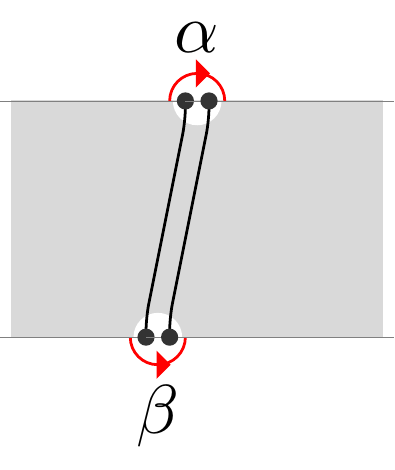}
&
\raisebox{14.5pt}{\includegraphics[scale=.5]{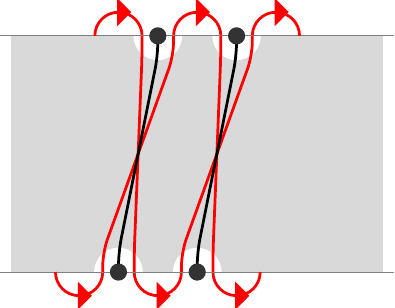}}\\
$\s,(\lambda'\bigcup\{{\color{ red}\ell}\},r)$ 
&
$\s'=\s|_{2,\alpha,\beta}$
&
$\s',(\lambda'\bigcup\{{\color{red}\ell+2}\},r)$\\
&\\
 \includegraphics[scale=.5]{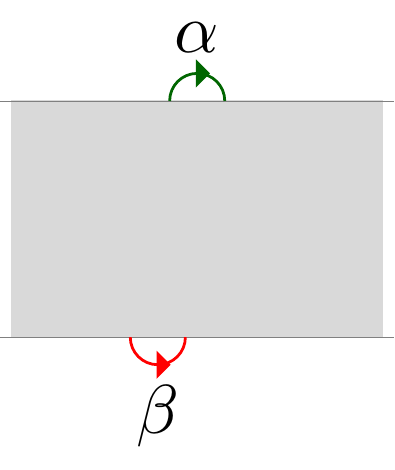}
&
 \includegraphics[scale=.5]{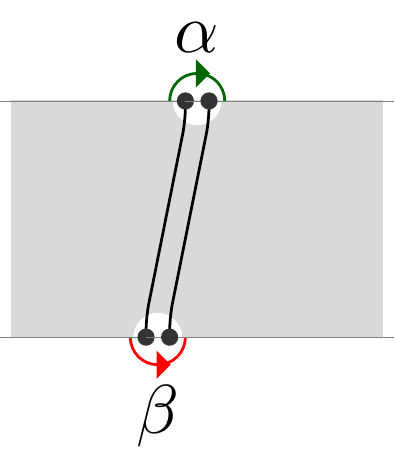}
&
\raisebox{14.5pt}{\includegraphics[scale=.5]{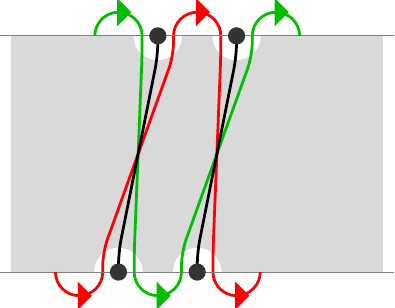}}\\
$\s,(\lambda'\bigcup\{{\color{ red}\ell},{\color{ green!50!black}\ell'}\},r)$ 
&
$\s'=\s|_{2,\alpha,\beta}$
& 
$\s',(\lambda'\bigcup\{{\color{red}\ell+1},{\color{ green!50!black}\ell'+1}\},r)$ 
\end{tabular}
\end{center}
\caption{\label{fig_double-edge}The first column represents the permutation $\s$ and its the cycle invariant. The second column represents the permutation $\s'$ obtained from the insertion of a double-edge within the arcs $\alpha$ and $\beta$ of $\s$ and the third column displays and demonstrates the new invariant of $\s'$.}
\end{figure}\qed

\end{pf}

\section{Shift-irreducible standard family\label{sec_shift_irreducible}}

In section 3.3 of \cite{DS17} we introduced (or reintroduced since it was already well-know in the litterature) the notions of standard permutation and standard family. We reproduce them here without proof for ease of reference. 

This section then concerns itselft with the more specific notion of \emph{shift-irreducible standard family}. Briefly a \emph{shift-irreducible standard} family is a standard family $(\s^{i})$ for which all but two permutations are irreducible after removing the edge $((\s^{i})^{-1}(1),1)$. They are extremly useful for the main induction as they resolve the irreduciblity issue of proposition \ref{pro_labelling_first_step} in the proof overwiew section \ref{ssec_pf_overview}. They are also our set $T_A$.

The permutation $\s$ is \emph{standard} if $\s(1)=1$. It is well know that every irreducible class contains a standard permutation and we can define the standard family of a standard permutation $\s$ as the collection of $n-1$ permutations $\{ \s^{(i)} := L^{i} \s \}_{0 \leq i \leq n-2}$.

Those notions allow one to prove the important following proposition (the original proposition 39 can be found page 43 of \cite{DS17}.)

\begin{proposition}[Properties of the standard family]
\label{pro_std_family}
Let $\s$ be a standard permutation with cycle invariant $(\lam,r)$, and 
$S=\{\s^{(i)}\}_i = \{L^i(\s)\}_i$ be its standard family. The latter has
the following properties:
\begin{enumerate}
  \item Every $\tau \in S$ has $\tau(1)=1$;
  \item The $n-1$ elements of $S$ are all distinct;
  \item There is a unique $\tau \in S$ such that $\tau(n)=n$;
  \item Let $m_i$ be the multiplicity of the integer $i$ in $\lambda$
    (i.e.\ the number of cycles of length $i$), and $r$ the rank.
    There are $i\, m_i$ permutations of $S$ which are of type
    $X(r,i)$, and $1$ permutation of type $H(r-j+1,j)$, for each $1
    \leq j\leq r$.\footnote{Note that, as $\sum_i i\, m_i+r=n-1$ by
      the dimension formula (\ref{eq.size_inv_cycle}), this list
      exhausts all the permutations of the family.}
 \end{enumerate}
\end{proposition}

Let us introduce a convenient notation.
\begin{notation}
Let $\s$ be a permutation, we define $d(\s)$ to be the permutation obtained from $\s$ by discarding the edge $(\s^{-1}(1),1)$, thus if $\s=\s(1),\ldots ,\s(\s^{-1}(1)),\ldots\s(n)$ then $d(\s)= \s(1)-1,\ldots,\widehat{\s(\s^{-1}(1))},\ldots,\s(n)-1$.
\end{notation}

Thus $d(\s)$ is the reduction of $(\s,c)$ where $c$ is the $(2n-2,2)$-coloring where the edge $(\s^{-1}(1),1))$ is gray.\\

Finally we describe the cycle invariant $(\lam',r')$ of $d(\s^i)$ in function of the cycle invariant $(\lam,r)$ of $\s^i$.

\begin{proposition}\label{pro_std_d_cycle_inv}
Let $(\s^i)_i$ be a standard family with cycle invariant $(\lam,r)$ then $\forall i$
\begin{itemize}
	\item If $\s^i$ has type $X(r,j)$, $d(\s^i)$ has type $H(j,r)$ and cycle invariant $(\lam\setminus\{j\},r+j-1)$. 
	\item If $\s^i$ has type $H(r_1,r_2)$, $d(\s^i)$ has type $X(r_2-1,r_1-1)$ and cycle invariant $(\lam\bigcup\{r_1-1\},r_2-1)$.
\end{itemize}
\end{proposition} 

\begin{pf}
The first case of this proposition was proven in lemma 5.16 in \cite{DS17}. The second case is proven likewise: This is the reverse implication of case \ding{177} in table \ref{table.addi} in \cite{DS17} (reproduced here), specialised to $s=0$. \qed
\end{pf}

\begin{table}[p!]
% \rule{390pt}{620pt}
\setlength{\unitlength}{72pt}
\begin{picture}(5.45,8.33)
\put(0,0){\includegraphics[scale=1.2]{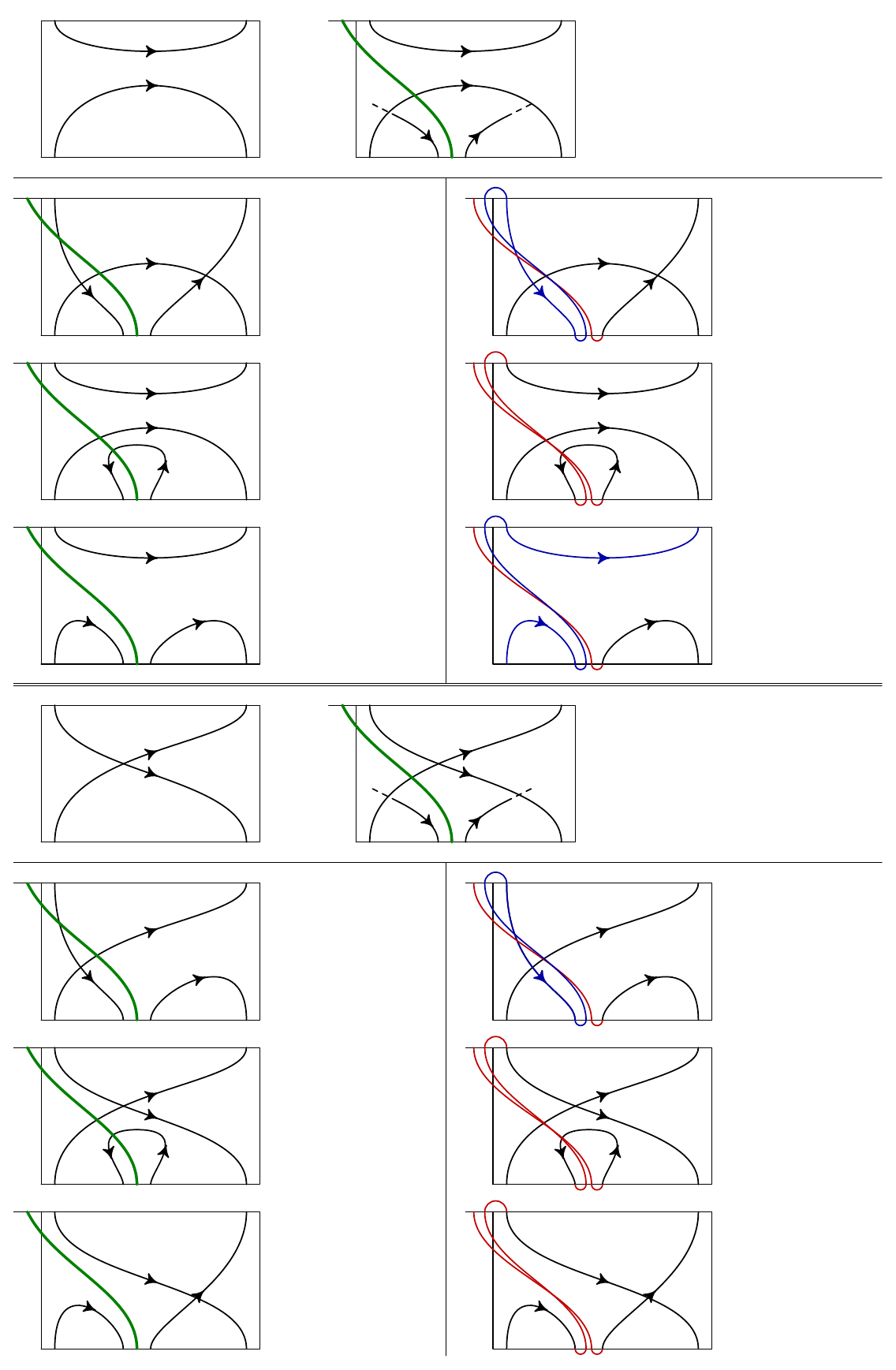}}
% --------------
\put(0.24,7.4){\begin{picture}(1.34,.83)
% \put(0,0){\rule{5pt}{5pt}}
% \put(1.34,.83){\rule{5pt}{5pt}}
\put(.6,.72){$2r_1$}
\put(.6,.3){$2r_2$}
\end{picture}}
% ---
\put(1,7.5){\begin{picture}(1.34,.83)
\put(1.25,0.27){?}
\put(2.25,0.27){?}
\put(2.8,.6){$H$-type\quad $H(r_1,r_2)$}
\end{picture}}
% ---
\put(0.24,6.315){\begin{picture}(1.34,.83)
% \put(0,0){\rule{5pt}{5pt}}
% \put(1.34,.83){\rule{5pt}{5pt}}
\put(.13,.7){$2s-1$}
\put(.42,.54){$2r_1-2s+1$}
\put(.6,.3){$2r_2$}
\put(2.25,.75){\ding{172}}
\put(1.42,.45){\begin{minipage}{2.5cm}
\begin{center}
Once per\\
$1 \leq s \leq r_1$
\end{center}
\end{minipage}}
\end{picture}}
% --
\put(3,6.315){\begin{picture}(1.34,.83)
% \put(0,0){\rule{5pt}{5pt}}
% \put(1.34,.83){\rule{5pt}{5pt}}
\put(.13,.7){$2s-1$}
\put(.42,.54){$2r_1-2s+1$}
\put(.6,.3){$2r_2$}
\put(1.42,.45){\begin{minipage}{2.5cm}
\begin{align*}
&H(r'_1,r'_2)\\
r_1'&=r_1\!-\!s\!+\!1 \\
r_2'&=r_2 \\
\lambda'&=\lambda \!\cup\! \{s\}
\end{align*}
\end{minipage}}
\end{picture}}
% --
\put(0.24,5.315){\begin{picture}(1.34,.83)
% \put(0,0){\rule{5pt}{5pt}}
% \put(1.34,.83){\rule{5pt}{5pt}}
\put(.6,.72){$2r_1$}
\put(.92,.43){$2r_2$}
\put(.82,.19){$2\ell$}
\put(2.25,.75){\ding{173}}
\put(1.42,.4){\begin{minipage}{2.5cm}
\begin{center}
$\ell$ times per\\
each cycle of $\lam$\\
of length $\ell$
\end{center}
\end{minipage}}
\end{picture}}
% --
\put(3,5.315){\begin{picture}(1.34,.83)
% \put(0,0){\rule{5pt}{5pt}}
% \put(1.34,.83){\rule{5pt}{5pt}}
\put(.6,.72){$2r_1$}
\put(.92,.43){$2r_2$}
\put(.82,.19){$2\ell$}
\put(1.42,.45){\begin{minipage}{2.5cm}
\begin{align*}
&H(r'_1,r'_2)\\
r_1'&=r_1\!+\!\ell+1 \\
r_2'&=r_2 \\
\lambda'&=\lambda \!\setminus\! \{\ell\}
\end{align*}
\end{minipage}}
\end{picture}}
% --
\put(0.24,4.315){\begin{picture}(1.34,.83)
% \put(0,0){\rule{5pt}{5pt}}
% \put(1.34,.83){\rule{5pt}{5pt}}
\put(.6,.72){$2r_1$}
\put(.17,.12){$2s$}
\put(.77,.31){$2r_2-2s$}
\put(2.25,.75){\ding{174}}
\put(1.42,.45){\begin{minipage}{2.5cm}
\begin{center}
Once per\\
$0 \leq s \leq r_2$
\end{center}
\end{minipage}}
\end{picture}}
% --
\put(3,4.315){\begin{picture}(1.34,.83)
% \put(0,0){\rule{5pt}{5pt}}
% \put(1.34,.83){\rule{5pt}{5pt}}
\put(.6,.72){$2r_1$}
\put(.17,.12){$2s$}
\put(.77,.31){$2r_2-2s$}
\put(1.42,.45){\begin{minipage}{2.5cm}
\begin{align*}
&X(r',i')\\
r'&=r_2\!-\!s \\
i'&=r_1\!+\!s \\
\lambda'&=\lambda \!\cup\! \{i'\}
\end{align*}
\end{minipage}}
\end{picture}}
% ----------------------------------------------
\put(0.24,3.21){\begin{picture}(1.34,.83)
% \put(0,0){\rule{5pt}{5pt}}
% \put(1.34,.83){\rule{5pt}{5pt}}
\put(.73,.71){$2i+1$}
\put(.9,.36){$2r+1$}
\end{picture}}
% ---
\put(1,3.31){\begin{picture}(1.34,.83)
\put(1.25,0.27){?}
\put(2.25,0.27){?}
\put(2.8,.6){$X$-type\quad $X(r,i)$}
\end{picture}}
% ---
\put(0.24,2.125){\begin{picture}(1.34,.83)
% \put(0,0){\rule{5pt}{5pt}}
% \put(1.34,.83){\rule{5pt}{5pt}}
\put(.73,.71){$2i+1$}
\put(.85,.32){$2r-2s$}
\put(.13,.7){$2s-1$}
\put(2.25,.75){\ding{175}}
\put(1.42,.45){\begin{minipage}{2.5cm}
\begin{center}
Once per\\
$1 \leq s \leq r$
\end{center}
\end{minipage}}
\end{picture}}
% --
\put(3,2.125){\begin{picture}(1.34,.83)
% \put(0,0){\rule{5pt}{5pt}}
% \put(1.34,.83){\rule{5pt}{5pt}}
\put(.73,.71){$2i+1$}
\put(.85,.32){$2r-2s$}
\put(.13,.7){$2s-1$}
\put(1.42,.45){\begin{minipage}{2.5cm}
\begin{align*}
&X(r',i')\\
r'&=r\!-\!s \\
i'&=i \\
\lambda'&=\lambda \!\cup\! \{s\}
\end{align*}
\end{minipage}}
\end{picture}}
% --
\put(0.24,1.125){\begin{picture}(1.34,.83)
% \put(0,0){\rule{5pt}{5pt}}
% \put(1.34,.83){\rule{5pt}{5pt}}
\put(.73,.71){$2i+1$}
\put(.9,.36){$2r+1$}
\put(.81,.13){$2\ell$}
\put(2.25,.75){\ding{176}}
\put(1.42,.4){\begin{minipage}{2.5cm}
\begin{center}
$\ell$ times per\\
each cycle of $\lam$\\
of length $\ell$
\end{center}
\end{minipage}}
\end{picture}}
% --
\put(3,1.125){\begin{picture}(1.34,.83)
% \put(0,0){\rule{5pt}{5pt}}
% \put(1.34,.83){\rule{5pt}{5pt}}
\put(.73,.71){$2i+1$}
\put(.9,.36){$2r+1$}
\put(.81,.13){$2\ell$}
\put(1.42,.45){\begin{minipage}{2.5cm}
\begin{align*}
&X(r',i')\\
r'&=r\!+\!\ell\!+\!1 \\
i'&=i \\
\lambda'&=\lambda \!\setminus\! \{\ell\}
\end{align*}
\end{minipage}}
\end{picture}}
% --
\put(0.24,0.125){\begin{picture}(1.34,.83)
% \put(0,0){\rule{5pt}{5pt}}
% \put(1.34,.83){\rule{5pt}{5pt}}
\put(.58,.71){$2i-2s+1$}
\put(.4,.54){$2r+1$}
\put(.17,.12){$2s$}
\put(2.25,.75){\ding{177}}
\put(1.42,.45){\begin{minipage}{2.5cm}
\begin{center}
Once per\\
$0 \leq s \leq i$
\end{center}
\end{minipage}}
\end{picture}}
% --
\put(3,0.125){\begin{picture}(1.34,.83)
% \put(0,0){\rule{5pt}{5pt}}
% \put(1.34,.83){\rule{5pt}{5pt}}
\put(.58,.71){$2i-2s+1$}
\put(.4,.54){$2r+1$}
\put(.17,.12){$2s$}
\put(1.42,.45){\begin{minipage}{2.5cm}
\begin{align*}
&H(r'_1,r'_2)\\
r'_1&=i\!-\!s\!+\!1 \\
r'_2&=r\!+\!s\!+\!1 \\
\lambda'&=\lambda \!\setminus\! \{i\}
\end{align*}
\end{minipage}}
\end{picture}}
% ---------------------------
\end{picture}
\caption{\label{table.addi}Modification to the cycle invariant between $\s$ and $\s|_{1,0,i}$, for every possible $i$. In green, the newly-added edge.  In red, parts which get added to the rank path. In blue, parts which are singled out to form a new cycle.}
\end{table}

 This proposition is very useful for the induction since it implies the following lemma

\begin{corollary}\label{cor_XH_d_remove}
Let $\s_1$ and $\s_2$ be two standard permutations/ If $\s_1$ and $\s_2$ have invariant $(\lambda,r,s)$ and same type $X(r,i)$ or $H(r_1,r_2)$ then the reducted permutations $\tau_1=d(\s_1)$ and $\tau_2=d(\s_2)$ have same cycle invariant.
\end{corollary}

This property (i.e same cycle invariant for $\tau_1$ and $\tau_2$) was the first of the three properties we where looking for in our proposition \ref{pro_labelling_first_step} of the proof overview.\\

For the purpose of guarantying the irreducibility of $\tau_1$ and $\tau_2$, we now define a more precise notion than just the standard family, we call it a \emph{shift-irreducible standard family}.

\begin{definition}\label{def_shift_irre}
Let $\s$ be a standard permutation and $S=(\s^i= L^i(\s))_{0\leq i \leq n-2}$ be its standard family. We say that $S$ is shift-irreducible if \[\forall i \in \{0,\ldots,n-2\}\setminus \left\{n-\sigma(2),n-\sigma(n)+1 \right\} \ d(\s^i) \text{ is irreducible.}\]
\end{definition}

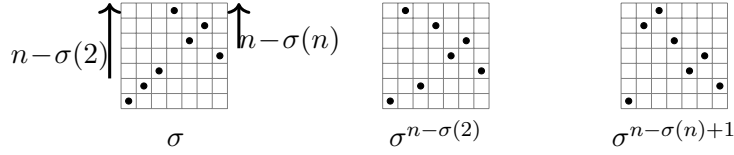
\begin{figure}[h!]
\begin{center}
\begin{tabular}{cccc}
\begin{tikzpicture}[scale=.2]
\permutation{1,2,3,7,5,6,4}
\draw[->,line width=.4mm] (.25,3)--(.25,8);
\draw[->,line width=.4mm] (8.75,5)--(8.75,8);
\node at (-3,4.5) {$n\!-\!\sigma(2)$};
\node at (12.1,5.5) {$n\!-\!\sigma(n)$};
\end{tikzpicture}
&
\begin{tikzpicture}[scale=.2]
\permutation{1,7,2,6,4,5,3}
\end{tikzpicture}
&
$\qquad$&
\begin{tikzpicture}[scale=.2]
\permutation{1,6,7,5,3,4,2}
\end{tikzpicture}
\\
$\s$ &
$\s^{n-\sigma(2)}$
&&
$\s^{n-\sigma(n)+1}$
\end{tabular}
\end{center}
\caption{\label{fig_shift_irreducible_ex} An example of a shift-irreducible family $S=(\s^i)_{0\leq i \leq n-2}$ with its two unavoidably reducible permutations $d(\s^{n-\sigma(2)})$ and $d(\s^{n-\sigma(n)+1}).$ For a proof of the shift-irreducibility, the characterisation proposition \ref{pro_chara_shift} is helpful.}
\end{figure}

In other words, a standard family $S$ is shift-irreducible if every $d(\s^i)$ that can be irreducible is indeed irreducible.  $d(\s^{n-\sigma(2)})$ and $d(\s^{n-\sigma(n)+1})$ are both always reducible since $\s^{n-\sigma(2)}(2)=n$ thus $d(\s^{n-\sigma(2)})(1)=n-1$ and  $\s^{n-\sigma(n)+1}(n)=2$ thus $d(\s^{n-\sigma(n)+1})(n-1)=1.$ Figure \ref{fig_shift_irreducible_ex} provides an example of a shift-irreducible family.

\begin{proposition}\label{pro_shift_standard}
Let $\s$ be a standard permutation with cycle invariant $(\lambda,r)$, and $S=(\s^i= L^i(\s))_{0\leq i \leq n-2}$ its standard family. 
Then $S$ is shift irreducible if and only if 
\begin{itemize}
\item For every $i\in \lambda$ the image by $d$ of the $i\cdot m_i$ permutations of $S$ of type $X(r,i)$ are irreducible (where $m_i$ is the multiplicity of $i$ in $\lambda$). 
\item For all $1< j <r$ the image by $d$ of the permutation of $S$ of type $H(r-j+1,j)$ is irreducible.
\end{itemize}
\end{proposition}

\begin{pf}
By proposition \ref{pro_std_family}, the permutations of $S$ are exactly the $i m_i$ permutations of type $X(r,i)$ for every $i\in\lambda$ and the permutations  of type $H(r-j+1,j)$ for all $1\leq j \leq r$. Clearly the permutations of type $H(1,r)$ and $H(r,1)$ are respectively $\s^{n-\sigma(2)}$ and $\s^{n-\sigma(n)+1}$ (see figure \ref{fig_perm_shift_pf}). Thus if the image by $d$ of every permutation besides those two are irreducible the family is shift-irreducible and reciproquely if the family is shift-irreducible the image by $d$ of every permutation besides those two are irreducible.
\begin{figure}[h!]
\begin{center}
\begin{tabular}{cccc}
\raisebox{2.3mm}{\includegraphics[scale=.5]{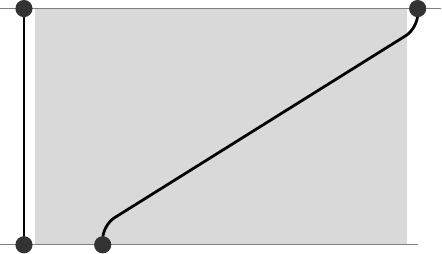}}
&
\raisebox{0mm}{\includegraphics[scale=.5]{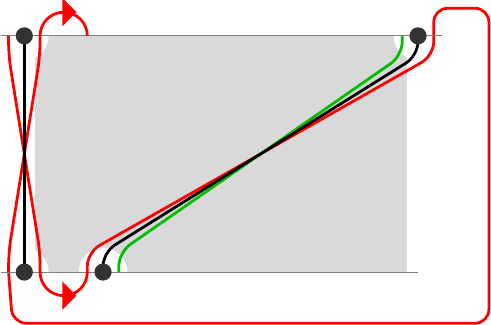}}
&
\raisebox{2.3mm}{\includegraphics[scale=.5]{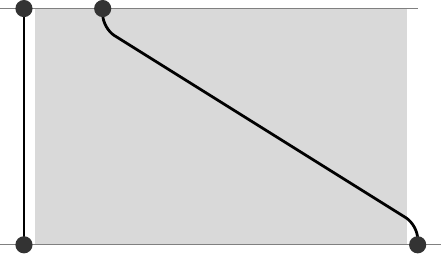}}
&
\raisebox{.8mm}{\includegraphics[scale=.5]{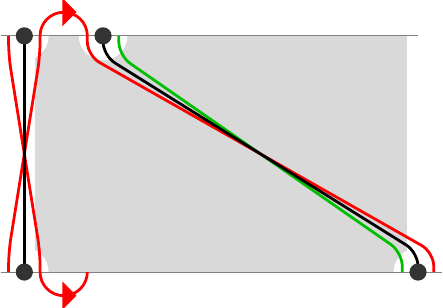}}\\
$\s^{n-\sigma(2)}$ &
has type $H(1,r)$
&
$\s^{n-\sigma(n)+1}$&
has type $H(r,1)$
\end{tabular}
\end{center}
\caption{\label{fig_perm_shift_pf} The two permutations $\s^{n-\sigma(2)}$ and $\s^{n-\sigma(n)+1}$ are of type $H(1,r)$ and $H(r,1)$ respectively.}
\end{figure} \qed
\end{pf}

As mentionned above the two permutations $\s^{n-\s(2)}$ and $\s^{n-\s(n)+1}$ -whose images by $d$ are reducible- have type $H(1,r)$ and $H(r,1)$. Thus $d(\s^{n-\s(2)})$ has cycle invariant $(\lam\bigcup\{0\},r-1)$ and $d(\s^{n-\s(n)+1})$ has cycle invariant $(\lam\bigcup\{r-1\},0)$ by proposition \ref{pro_std_d_cycle_inv}. Neither are cycle invariants that we allow for in the classification theorem (we never consider cycles or rank of length 0 since it always implies reducibility). 

In other term, every single one of the permutations $d(\s^i)$ that can be used for the induction are irreducible if $(\s^i)_i$ is an shift-irreducible standard family.\\

It is interesting that the fruitful notion of shift-irreducible family has a very simple characterisation in term of just the permutation $\s$ of the family with $\s(1)=1$ and $\s(2)=2.$

\begin{proposition}[Characterisation of shift-irreducible families]\label{pro_chara_shift}
Let $\s$ be a standard permutation with $\s(2)=2$ and $S=(\s^i= L^i(\s))_{0\leq i \leq n-2}$ its standard family. 
Then $S$ is shift-irreducible if and only if $\s$ does not have the following form : 
\begin{center}
\includegraphics{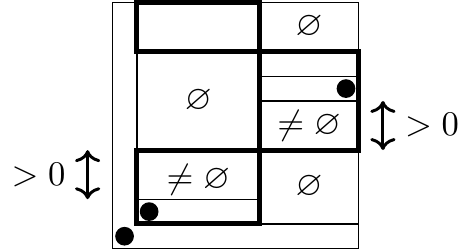}
\end{center}
\end{proposition}
\begin{pf}
We show that $S$ is not shift-irreducible if and only if $\s$ has the form described in the proposition.
If $S$ is not shift-irreducible then it means that there exists a $i \in \{0,\ldots,n-2\}\setminus \left\{n-\sigma(2),n-\sigma(n)+1 \right\}$ such that $d(\s^i)$ is reducible. 
Thus $d(\s^i)$ has the following form 
\begin{center}
\includegraphics{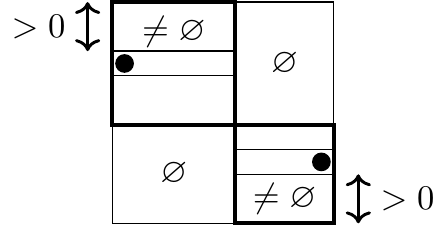}
\end{center}
The blocks with $\neq \varnothing$ must not be empty otherwise we would have either $\s^i= \s^{n-\sigma(2)}$ or $\s^{n-\sigma(n)+1}$. $\s^i$ has then the form
\begin{center}
\includegraphics{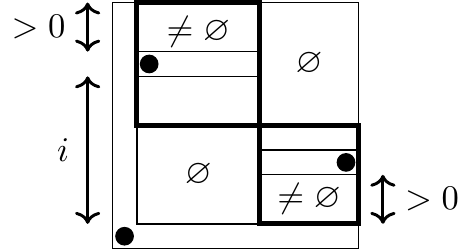}
\end{center}
and since by definition $\s^i=L^i(\s)$, $\s$ must have the form 
\begin{center}
\includegraphics{figure/fig_perm_shift_reducible_2.pdf}
\end{center}
as predicted. The reverse implication is obtained by reversing the steps of the proof. \qed
\end{pf}

In light of this proposition we will often say \emph{shift-irreducible permutation} to denote the permutation $\s$ of the shift-irreducible standard family with $\s(1)=1$ and $\s(2)=2$. Moreover to prove the existence of a shift-irreducible standard family we will just construct a shift-irreducible permutation.  

A particular type of shift-irreducible permutation is the \emph{$I_2X$-permutation}:

\begin{corollary}\label{cor_chara_shift}
Let $\s$ be a standard permutation with $\s(2)=2$
If $\s$ have the following form: 
\begin{center}
\includegraphics{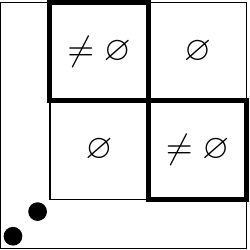} equivalently \includegraphics[scale=.5]{figure/fig_goal_induction.pdf}
\end{center}
Then $\s$ is a shift-irreducible permutation. We call $\s$ a \emph{$I_2X$-permutation}.
\end{corollary}
\begin{pf}
It is clear that a $I_2X$-permutation verifies the condition of the characterisation propositon \ref{pro_chara_shift}.
\qed
\end{pf}

%%%%%%%%%%%%%%%%%%%%%%%%%%%%%%%%%%%%%%%%%%%%%%%%%%%%%%%%
\section{The sign invariant}
\label{sec.signinv}
%%%%%%%%%%%%%%%%%%%%%%%%%%%%%%%%%%%%%%%%%%%%%%%%%%%%%%%

%-------------------------------------------------------
\subsection{Arf functions for permutations}
\label{ssec.arf_inv}
%--------------------------------------------------------

For $\s$ a permutation in $\kS_n$, let 
\be
\chi(\s)
= \# \{ 1\leq i<j \leq n \;|\; \s(i)<\s(j) \}
\ee
i.e.\ $\chi(\s)$ is the number of pairs of non-crossing edges in the
diagram representation of~$\s$. 
 
Let $E=E(\s)$ be the subset of $n$ edges in $\cK_{n,n}$ described by
$\s$.  For any $I \subseteq E$ of cardinality $k$, the permutation
$\s|_{I} \in \kS_k$ is defined in the obvious way, as the one
associated to the subgraph of $\cK_{n,n}$ with edge-set $I$, with
singletons dropped out, and the inherited total ordering of the two
vertex-sets.

Define the two functions
\begin{align}
A(\s)
&:=
\sum_{I \subseteq E(\s)} (-1)^{\chi(\s|_I)}
\ef;
&
\Abar(\s)
&:=
\sum_{I \subseteq E(\s)} (-1)^{|I|+\chi(\s|_I)}
\ef.
\end{align}
When $\s$ is understood, we will just write $\chi_I$ for
$\chi(\s|_I)$. The quantity $A$ is accessory in the forthcoming
analysis, while the crucial fact for our purpose is that the quantity
$\Abar$ is invariant in the $\perms_n$ dynamics.

In the following section, we define a technique to demonstrate identities of the arf invariant involving differents configurations.

%-------------------------------------------------------
\subsection{Automatic proofs of Arf identitites}
\label{ssec.arfcalcseasy}
%-------------------------------------------------------

We will \emph{not} try to evaluate Arf functions of
large configurations starting from scratch. We will rather compare the
Arf functions of two (or more) configurations, which differ by a
finite number of edges, and establish linear relations among their Arf
functions.  The method we develop here, gives an algorithm to find and check Arf identities.

In order to have the appropriate terminology for expressing this
strategy, let us define the following:
Given a permutation $\s$ define $\s_{k,\ell}$ to be a permutation with $k$ marks on its bottom line and $\ell$ marks on its top line. The marks are all at distinct positions and do not touch the corners of the permutation. These marks break the bottom (respectively top) line into $k+1$ open interval $P_{-,1},\ldots,P_{-,k+1}$ (respectively $\ell+1$ open interval $P_{+,1},\ldots,P_{+,\ell+1}$).

\noindent For example if $k=1,\ell=3$: 
\[
\put(50,-30){$P_{-,1}$}\put(100,-30){$P_{-,2}$}\put(26,25){$P_{+,1}$}\put(60,25){$P_{+,2}$}\put(90,25){$P_{+,3}$}\put(125,25){$P_{+,4}$}
\s_{k,\ell}=\raisebox{-20pt}{\includegraphics[scale=2.5]{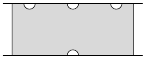}}
\]

Let $\s_{k,\ell,E'}$ be the permutation obtained by adding a set of edges $E'$ on the marks of permutation $\s_{k,\ell}$ with the following convention: an edge $e\in E'$ is a pair $(i.x,j.y)$. The edge connects the $i$th bottom mark and the $j$th top mark, and it is ordered as the $x$th edge within the bottom mark and the $y$th edge within the top mark. Note that if $i=0$ of $i=k+1$ (likewise of $j$) this implies that the edge is connected to a corner of the permutation.

\noindent For example if $k=1,\ell=3$ and $E'=\{(0.1,2.2),(1.1,3.1),(1.2,1.1),(1.3,2.1),(2,1.2)\}$: 
\[
\put(60,-30){$P_{-,1}$}\put(110,-30){$P_{-,2}$}\put(36,25){$P_{+,1}$}\put(70,25){$P_{+,2}$}\put(105,25){$P_{+,3}$}\put(135,25){$P_{+,4}$}
\s_{k,\ell,E'}=\raisebox{-20pt}{\includegraphics[scale=2.5]{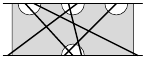}}
\]
We will define an algorithm that allows one to check if, for all $\s$, we have $\sum^n_{i=1} K_i\Abar(\s_{k,\ell,E^i})=0$ or  $\sum^n_{i=1} K_iA(\s_{k,\ell,E^i})=0$ for some $k,\ell,(E^i)_i,(K_i)_i,n.$

\begin{figure}[tb!]
\begin{center}
\includegraphics{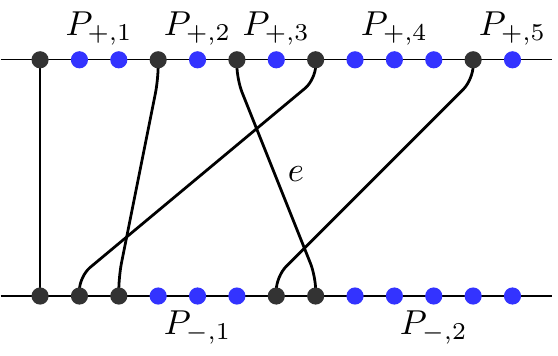}
\end{center}
\caption{\label{fig.arf_ex_def}The permutation $\s_{1,4,\{(0.1,0.1),(0.2,3.1),(0.3,1.1),(1.1,4.1),(1.2,2.1)\}}$. We cannot show the full matrix $Q$ for such a big
  example, but we can give one row, for the edge which has the label
  $e$ in the drawing. The row $Q_e$ reads
$(Q_e)_{11, 12, \ldots, 15, 21, \ldots, 25}=(1,1,0,0,0,\,0,0,1,1,1)$.}
\end{figure}

\begin{definition}
\label{def.637647}
Let $\s_{k,\ell,E'}$, $P_{-,1},\ldots,P_{-,k+1}$  and $P_{+,1},\ldots,P_{+,\ell+1}$ be as defined above.
Then define the $m \times (k+1)(\ell+1)$ matrix valued in $\gf_2$
\be
Q_{e,ij} :=
\left\{
\begin{array}{ll}
1 & \textrm{edge $e \in E'$ does not cross the segment connecting
  $P_{-,i}$ to $P_{+,j}$,}\\
0 & \textrm{otherwise.}
\end{array}
\right.
\ee
For $v \in \gf_2^{(k+1)(\ell+1)}$, let $|v|$ be the number of entries equal to
$1$. Similarly, identify $v$ with the corresponding subset of $[(k+1)(\ell+1)]$.
Given such a construction, introduce the following functions on
$(\gf_2)^{(k+1)(\ell+1)}$
\begin{align}
A_{k,\ell,E'}(v)
&:=
\sum_{u \in (\gf_2)^{E'}}
(-1)^{\chi_u + (u,Qv)}
\ef;
&
\Abar_{k,\ell,E'}(v)
&:=
\sum_{u \in (\gf_2)^{E'}}
(-1)^{|u|+\chi_u + (u,Qv)}
\ef.
\end{align}
\end{definition}
The construction is illustrated in
Figure~\ref{fig.arf_ex_def}.

Let us comment on the reasons for introducing such a definition. The
quantities $A_{k,\ell,E'}(v)$ (respectively $\Abar_{k,\ell,E'}(v)$)  do not depend on $\s$ and allows to sum together many contributions to the function $A(\s_{k,\ell,E'})$. Our goal is to have $E'$ of fixed size, while
$E$ (the edge set of $\s$) is arbitrary and of unbounded size, so that the
verification of our properties, as it is confined to the matrix $Q$,
involves a finite data structure. Thus the algorithm will be exponential in $|E'|$ which will not be a problem for small sizes.

Indeed, let us split in the natural way the sum over subsets $I$ of $E\bigcup E'$ the edge set of $\s_{k,\ell,E'}$ namely
\[
\sum_{I_0 \subseteq E\bigcup E'} f(I_0)
=
\sum_{I \subseteq E}
\sum_{I' \subseteq E'} 
 f(I \cup I')
\]
For $I$ and $J$ two disjoint sets of edges, call $\chi_{I,J}$ the number
of pairs $(i,j) \in I \times J$ which do not cross. Then clearly
\[
\chi_{I \cup J} = \chi_{I} + \chi_{J} + \chi_{I,J}
\]
Now let $u(I') \in \{0,1\}^{E'}$ be the vector with entries $u_e=1$ if
$e \in I'$ and $0$ otherwise. Let $m(I)=\{m_{ij}(I)\}$ be the 
$(k+1) \times (\ell+1)$ matrix describing the number of edges connecting the
intervals $P_{-,i}$ to $P_{+,j}$ in $\s$, and let
$v(I)=\{v_{ij}(I)\}$, $v_{ij} \in \{0,1\}$ be the parities of the
$m_{ij}$'s. Call $I_{ij}$ the restriction of $I$ to edges
connecting $P_{-,i}$ and $P_{+,j}$.  Clearly,
\[\chi_{I',I}=\sum_{ij} \chi_{I',I_{ij}} 
= \sum_{e,ij} u_e Q_{e,ij} m_{ij}=(u(I'),Qm(I)),\] 
which has the same parity as the analogous expression with $v$'s
instead of $m$'s. I.e. we have
\[(-1)^{\chi_{I',I}}=(-1)^{(u(I'),Qv(I))}.\]
Now, while the $m$'s are in $\bN$, the vector $v$ is
in a linear space of finite cardinality, which is crucial for allowing
a finite analysis of our expressions.

As a consequence,
\begin{align}\label{eq.arf_explain_def}
A(\s_{k,\ell,E^i})
&=
\sum_{I \subseteq E}(-1)^{\chi_I}A_{k,\ell,E'}(v(I))
\ef;
\\\label{eq.arf_explain_def_1}
\Abar(\s_{k,\ell,E^i})
&=\sum_{I \subseteq E}(-1)^{|I|+\chi_I}\Abar_{k,\ell,E'}(v(I))
\ef.
\end{align}

\noindent Thus we have the following proposition:

\begin{theorem}
\label{prop.ArfProp}\label{prop.Arf0}
Let $k,\ell \in \N$ and let $(E^i)_{1\leq i\leq n}$ be a family of edge set. Then the two following statements are equivalent:
\begin{enumerate}
\item For all $v \in GF_2^{(k+1)(\ell+1)},$ we have $\sum_{i=1}^n K_i A_{k,\ell,E^i}(v)=0$.
\item for all $\s,$ we have $\sum_{i=1}^n K_i A(\s_{k,\ell,E^i})=0.$ 
\end{enumerate}
The same statement holds for $\Abar$.
\end{theorem}

\begin{pf}
Statement 1 implies 2 due to equation (\ref{eq.arf_explain_def}).\\
Let us show the converse: If $\s$ has no edge then $\sum_{i=1}^n K_i A(\s_{k,\ell,E^i})=0$ is equivalent to $\sum_{i=1}^n K_i A_{k,\ell,E^i}(v)=0$ with $v$ being the zero vector of $(GF_2)^{(k+1)(\ell+1)}$. 

Then we choose the family of permutations $(\s^{a,b})_{a,b}$ with exactly one edge connecting $P_{-,a}$ to $P_{+,b}$. Then if we define $v_{a,b}$ to be the vector $(GF_2)^{(k+1)(\ell+1)}$ with exactly one 1 at position $ab$ we have
\begin{flalign*}
&\sum_{i=1}^n K_i A(\s^{a,b}_{k,\ell,E^i})=0\\
\iff& \sum_{i=1}^n K_i ( A_{k,\ell,E^i}(0) + A_{k,\ell,E^i}(v_{a,b})) =0\\
 \iff& \sum_{i=1}^n K_i  A_{k,\ell,E^i}(v_{a,b}) =0\\
\end{flalign*} in the last line we have used that $\sum_{i=1}^n K_i A_{k,\ell,E^i}(0)=0$.

Thus inductively we show that $\sum_{i=1}^n K_i A_{k,\ell,E^i}(v)=0$ for any $v\in (GF_2)^{(k+1)(\ell+1)}$ with at most $p\leq n$ ones. 
\qed
\end{pf}

This theorem is very important since it reduces the problem of calculating Arf identities for permutations of any size to a check of an Arf identity for $2^{(k+1)(\ell+1)}$ values. Thus in exponential time in $k\ell$ we can calculate $A_{k,\ell,E'}(v)$ for every $v \in GF_2^{(k+1)(\ell+1)}$ and decide if a given Arf identity is correct.

We can even do better:

\begin{proposition}
Let $k,\ell \in \N$ and let $(E^i)_{1\leq i\leq n}$ be a family of edge set. We can decide (in exponential time in $k\ell$) if there exists $x_1,\ldots,x_n$ such that  $\sum_{i=1}^n x_i A(\s_{k,\ell,E'})=0$.
\end{proposition} 

\begin{pf}
For every $v$ we have an equation $\sum_{i=1}^n x_i A_{k,\ell,E^i}(v) =0$ with the $n$ unknown variables. So there are $2^{(k+1)(\ell+1)}$ equations. We can find the subspace of solution in time exponential in $k\ell$.
\qed
\end{pf}

The previous proposition can be used when we suspect a relation between a few configurations without knowing the coefficients. The algorithm demands little more than the previous one for the verification so it remains usable for small $|E'|$.

Finally we can actually enumerate all the possible Arf identities:

\begin{proposition}
There is an algorithm that enumerate all the arf identities with at most $n$ terms and on an edge set $E'$ of size at most $h$.
\end{proposition} 

This algorithm is really not praticable. However it can be used in the following case:
 we have two terms and we want to find an arf identity relating them to one another but the previous algorithm failed (i.e there are no identity containing only those two terms). Then we use this algorithm to find a third term (or a fourth etc...) for which an identity exists.

We can even propose a generalisation of this framework: let us choose two permutations $\pi_{-}$ and $\pi_{+}$ of size $k+1$ and $\ell+1$ respectively then $\s_{k,\ell,E',\pi_{-},\pi_{+}}$ is obtained from $\s_{k,\ell,E'}$ by permuting the $P_{-,i}$ with $\pi_{-}$ and $P_{+,i}$ with $\pi_{+}$.

\noindent For example if $k=1,\ell=3$ and $E'=\{(0.1,2.2),(1.1,3.1),(1.2,1.1),(1.3,2.1),(2,1.2)\}$: 
\[
\put(60,-30){$P_{-,1}$}\put(110,-30){$P_{-,2}$}\put(36,25){$P_{+,1}$}\put(70,25){$P_{+,2}$}\put(105,25){$P_{+,3}$}\put(135,25){$P_{+,4}$}
\s_{k,\ell,E'}=\raisebox{-20pt}{\includegraphics[scale=2.5]{P05_Arf/figure/fig_arf_ex2.pdf}}
\put(130,-30){$P_{-,\pi_{-}1}$}\put(180,-30){$P_{-,\pi_{-}2}$}\put(100,25){$P_{+,\pi_{+}1}$}\put(140,25){$P_{+,\pi_{+}2}$}\put(175,25){$P_{+,\pi_{+}3}$}\put(210,25){$P_{+,\pi_{+}4}$}
\qquad \qquad\s_{k,\ell,E',\pi_{-},\pi_{+}}=\raisebox{-20pt}{\includegraphics[scale=2.5]{P05_Arf/figure/fig_arf_ex2.pdf}}.
\]
For an example with $\pi_{-}=id_2$ and $\pi_{+}=(2,1)$ (the reversing permutation $\omega$ at size 2) we have:
\[
\put(55,-30){$P_{-,1}$}\put(95,-30){$P_{-,2}$}\put(55,25){$P_{+,1}$}\put(85,25){$P_{+,2}$}
\s_{1,1,E'}=\raisebox{-20pt}{\includegraphics[scale=2.5]{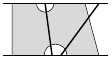}}
\put(125,-30){$P_{-,1}$}\put(165,-30){$P_{-,2}$}\put(120,50){$P_{+,2}$}\put(160,50){$P_{+,1}$}
\qquad \qquad\s_{1,1,E',\pi_{-},\pi_{+}}=\raisebox{-20pt}{\includegraphics[scale=2.5]{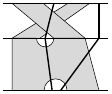}}
\]
It is easily checked that the previous theorems continue to hold for this generalisation once we introduce for $v \in \gf_2^{(k+1)(\ell+1)}$ the following function on
$(\gf_2)^{(k+1)(\ell+1)}$ (similar definition for $\Abar_{k,\ell,E',\pi_{-},\pi_{+}}$)
\begin{align}
A_{k,\ell,E',\pi_{-},\pi_{+}}(v)
&:=
\sum_{u \in (\gf_2)^{E'}}
(-1)^{\chi_u + (u,Q(P^{-1}_{\pi_{-}}vP_{\pi_{+}}))}
\ef;
\end{align}
Where in the expression $(P^{-1}_{\pi_{-}}vP_{\pi_{+}})$, $v$ is identified to the matrix of size $(k+1)\times (\ell+1)$ and $P_{\pi_{-}}$ and $P_{\pi_{+}}$ are the permutation matrices associated to $\pi_{-}$ and $\pi_{+}$.

The framework of automatic proofs of Arf identities we have developped is rather general. Most of the identities found in the litterature (see \cite{KZ03}, \cite{Boi13}, \cite{DS17}, \cite{Del13}, \cite{Gut17}) can be obtained in this setting.

Let us now apply the algorithm to find Arf identities.

It is convenient to introduce the notation
$\vec{A}(\s)=\begin{pmatrixsm} \Abar(\s) \\ A(\s) \end{pmatrixsm}$.
We have
\begin{proposition}
\label{prop.fingerred}
\begin{align}
\label{prop.signinvdyn}
\Abar\bigg(\,\tau=\raisebox{-12pt}{\includegraphics[scale=1.8]{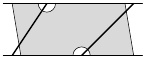}}\,\bigg)
&=
\Abar\bigg(\,\s=\raisebox{-12pt}{\includegraphics[scale=1.8]{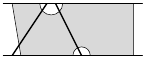}}\,\bigg)
\\
\label{eq.546455a}
\vec{A}\bigg(\,
%\tau=
\raisebox{-12pt}{\includegraphics[scale=1.8]{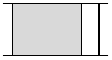}}\,\bigg)
&=
\begin{pmatrix}
1 & -1 \\ 1 & 1
\end{pmatrix}
\vec{A}\bigg(\,
%\s=
\raisebox{-12pt}{\includegraphics[scale=1.8]{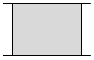}}\,\bigg)
\ef;
\\
\label{eq.546455b}
\vec{A}\bigg(\,
%\tau=
\raisebox{-12pt}{\includegraphics[scale=1.8]{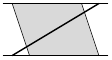}}\,\bigg)
&=
\begin{pmatrix}
0 & 0 \\ 0 & 2
\end{pmatrix}
\vec{A}\bigg(\,
%\s=
\raisebox{-12pt}{\includegraphics[scale=1.8]{FigFol/Figure4_fig_arfproof_s3.pdf}}\,\bigg)
\ef;
\end{align}
\end{proposition}

Clearly equation (\ref{prop.signinvdyn}) prove the invariance of the arf invariant for the dynamics since $\s=L(\tau)$ and the case $R$ is deduced by symmetry.

\subsection{Arf relation for the induction\label{ssec_arf_useful}\label{sec_arf}}

This section lists the arf identities required for the inductive proof using the labelling method.

This section contains four statements of note: propositions \ref{lem.certifopposign}, \ref{lem_q_i_sign} and \ref{cor.signsurjq2x} will be used in a technical way during the induction and proposition \ref{pro_opposite_sign} will be applied in section \ref{sec_IIXshaped} to construct pairs of permutations with same cycle invariant but opposite sign invariant.

We have a first proposition, involving the evaluation of $\Abar$ on
three distinct configurations
\begin{proposition}
\label{lem.certifopposign}
\be
\Abar\bigg(\,\s=
\raisebox{-12pt}{\makebox[0pt][l]{\raisebox{-8pt}{
\rule{15pt}{0pt}%
$e$$f$
\rule{15pt}{0pt}%
%\rule{11pt}{0}$\nu$\rule{11pt}{0}
$g$}}}%
\raisebox{-12pt}{\includegraphics[scale=1.8]{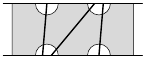}}\,\bigg)
+
\Abar\bigg(\,\tau=
\raisebox{-12pt}{\makebox[0pt][l]{\raisebox{-8pt}{
\rule{17pt}{0pt}%
$f$
\rule{15pt}{0pt}%
%\rule{11pt}{0}$\nu$\rule{11pt}{0}
$g$$e'$}}}%
\raisebox{-12pt}{\includegraphics[scale=1.8]{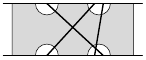}}\,\bigg)
=
2 \Abar\bigg(\,\rho=
\raisebox{-12pt}{\makebox[0pt][l]{\raisebox{-8pt}{
\rule{17pt}{0pt}%
$f$
\rule{17pt}{0pt}%
%\rule{11pt}{0}$\nu$\rule{11pt}{0}
$g$}}}%
\raisebox{-12pt}{\includegraphics[scale=1.8]{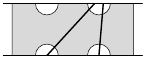}}\,\bigg)
\ee

\be
\Abar\bigg(\,\s=
\raisebox{-12pt}{\makebox[0pt][l]{\raisebox{-8pt}{
\rule{15pt}{0pt}%
$e$$f$
\rule{15pt}{0pt}%
%\rule{11pt}{0}$\nu$\rule{11pt}{0}
$g$}}}%
\raisebox{-12pt}{\includegraphics[scale=1.8]{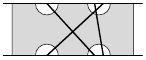}}\,\bigg)
+
\Abar\bigg(\,\tau=
\raisebox{-12pt}{\makebox[0pt][l]{\raisebox{-8pt}{
\rule{17pt}{0pt}%
$f$
\rule{15pt}{0pt}%
%\rule{11pt}{0}$\nu$\rule{11pt}{0}
$g$$e'$}}}%
\raisebox{-12pt}{\includegraphics[scale=1.8]{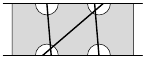}}\,\bigg)
=
2 \Abar\bigg(\,\rho=
\raisebox{-12pt}{\makebox[0pt][l]{\raisebox{-8pt}{
\rule{17pt}{0pt}%
$f$
\rule{17pt}{0pt}%
%\rule{11pt}{0}$\nu$\rule{11pt}{0}
$g$}}}%
\raisebox{-12pt}{\includegraphics[scale=1.8]{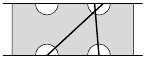}}\,\bigg)
\ee
\end{proposition}
\begin{pf}
Clearly we are in the framework developped in the previous section so the proof can be checked by proposition \ref{prop.ArfProp}.\qed
\end{pf}

Those identities are also found in \cite{Del13} lemmas 4.9 and 4.10 and \cite{Boi12} proposition 4.2. The correspondance is not exact due to the difference of language but the proof of those proposition/lemmas involve the use of those identities (once translated in their language).

Let $\s=$\raisebox{-4.5mm}{\includegraphics[scale=1.8]{figure/fig_arfproof_r10-eps-converted-to.pdf}}, we define $\s_{i,j} =$
 \makebox[0pt][l]{\raisebox{-4.5mm}{\includegraphics[scale=1.8]{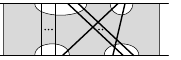}}}
\hspace{17pt}\raisebox{28pt-5mm}{$\rotatebox{90}{$\left. \rule{0pt}{2pt}\right\}\rotatebox{-90}{$\!\!i$} $}$}$
\hspace{4pt}\raisebox{28pt-5mm}{$\rotatebox{90}{$\left. \rule{0pt}{2.5pt}\right\}\rotatebox{-90}{$\!\!j$} $}$}$.\\

With this notation, the first equation of proposition \ref{lem.certifopposign} can be rewritten \[ \Abar(\s_{1,0})+\Abar(\s_{0,1})=2\Abar(\s_{0,0}).\]

\begin{lemma}\label{lem.certifopposign_2_3}
\be\label{eq_opposign_2}
\Abar(\s_{2,0})-\Abar(\s_{0,2})=2\left(\Abar(\s_{1,0})-\Abar(\s_{0,1})\right)
\ee
\be\label{eq_opposign_3}
\Abar(\s_{3,0})+\Abar(\s_{0,3})=3\Abar(\s_{2,0})+3\Abar(\s_{0,2})-4\Abar(\s_{0,0})
\ee
\end{lemma}
%\begin{pf}
%Let us start with equation \ref{eq_opposign_2}, by applying proposition \ref{lem.certifopposign} to $\s_{1,0}$ and $\s_{0,1}$ we have :
%\be\label{eq_opposign_4}
%\Abar(\s_{2,0})+\Abar(\s_{1,1})=2\Abar(\s_{1,0})
%\ee
%\be\label{eq_opposign_5}
%\Abar(\s_{0,2})+\Abar(\s_{1,1})=2\Abar(\s_{0,1})
%\ee
%Thus equation \ref{eq_opposign_4} minus equation \ref{eq_opposign_5} gives equation \ref{eq_opposign_2}.
%
%For equation \ref{eq_opposign_3}, we apply proposition \ref{lem.certifopposign} to $\s_{2,0}$, $\s_{1,1}$ and $\s_{0,2}$ : 
%\be\label{eq_opposign_6}
%\Abar(\s_{3,0})+\Abar(\s_{2,1})=2\Abar(\s_{2,0})
%\ee
%\be\label{eq_opposign_7}
%\Abar(\s_{2,1})+\Abar(\s_{1,2})=2\Abar(\s_{1,1})
%\ee
%\be\label{eq_opposign_8}
%\Abar(\s_{1,2})+\Abar(\s_{0,3})=2\Abar(\s_{0,2})
%\ee
%Now 
%\begin{align*}
%\Abar(\s_{3,0})+\Abar(\s_{0,3})&= \Abar(\s_{3,0})+\Abar(\s_{2,1})-\left(\Abar(\s_{2,1})+\Abar(\s_{1,2})\right)+\Abar(\s_{1,2})+\Abar(\s_{0,3}) &\\
%&= 2\Abar(\s_{2,0})+2\Abar(\s_{0,2})-2\Abar(\s_{1,1})&\text{by eqs \ref{eq_opposign_6},\ref{eq_opposign_7},\ref{eq_opposign_8}}\\
%&= 2\Abar(\s_{2,0})+2\Abar(\s_{0,2})- \left(2\Abar(\s_{1,0})-\Abar(\s_{2,0})\right)-\left(2\Abar(\s_{0,1})-\Abar(\s_{0,2})\right)&\text{by eqs \ref{eq_opposign_4},\ref{eq_opposign_5}}\\
%&=3\Abar(\s_{2,0})+3\Abar(\s_{0,2})-2\left(\Abar(\s_{1,0})+\Abar(\s_{0,1})\right)\\
%&=3\Abar(\s_{2,0})+3\Abar(\s_{0,2})-4\Abar(\s_{0,0}) &\text{by proposition \ref{lem.certifopposign}}
%\end{align*}
%\qed
%\end{pf}

Once again we are in the framework so the proposition can be verified by proposition \ref{prop.ArfProp}. \qed

The main use of those two lemmas is to compare the Arf invariant of
two configurations having the same cycle invariant and differing by just a few consecutives and parallel edges. The exact statement is the content of the next proposition.   

\begin{figure}[h!]
\begin{center}
\includegraphics[scale=1]{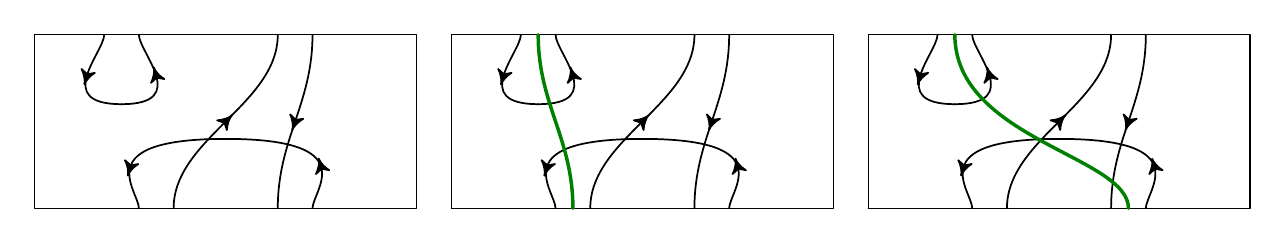}
\end{center}
\caption{\label{fig_opposite_sign}From left to right: the permutations $\s$, $\s|_{1,\alpha,\beta}$ and $\s|_{1,\alpha,\beta'}$. The proposition \ref{pro_opposite_sign} says that $\s|_{1,\alpha,\beta}$ and $\s|_{1,\alpha,\beta'}$ have same cycle invariant and opposite sign invariant.}
\end{figure}

\begin{proposition}[Opposite sign]\label{pro_opposite_sign}
Let $\s$ be a permutation with exactly two even cycles (one being possibly the rank). Let $\alpha$ be a top arc of the first cycle and $\beta$ and $\beta'$ be two consecutive bottom arcs of the second cycle. 

Define $\s|_{i,\alpha,\beta}$ and $\s|_{i,\alpha,\beta'}$ following notation \ref{not_add_edge} (as a reminder they are the two permutations obtained by adding $i$ parallel and consecutive edges within $\alpha, \beta$ and within $\alpha, \beta'$ respectively). See figure \ref{fig_opposite_sign} for the case $i=1$ and $\s$ has two even cycles (none being the rank path). 

 Then for $i\leq 3$, $\s|_{i,\alpha,\beta}$ and $\s|_{i,\alpha,\beta'}$ have invariant $(\lambda',r',s)$ and $(\lambda',r',-s)$ respectively.

To be more precise : \begin{enumerate}
\item If $\s$ has cycle invariant $(\lambda\bigcup\{2\ell,2\ell'\},r)$ we have \\[-9mm]
\begin{center}
$\begin{array}{|c|c|c|c|}
\cline{2-4}
\multicolumn{1}{c|}{}
&\raisebox{-5pt}{\rule{0pt}{16pt}}i=1&i=2&i=3\\ \hline
(\lam,r,s) \text{ of } \s|_{i,\alpha,\beta} 
&\raisebox{-5pt}{\rule{0pt}{16pt}}
(\lambda\bigcup\{2\ell\!+\!2\ell'\!+\!1\},r,s)
&(\lambda\bigcup\{2\ell\!+\!1,2\ell'\!+\!1\},r,s)
&(\lambda\bigcup\{2\ell\!+\!2\ell'\!+\!3\},r,s)
\\ 
\hline
(\lam,r,s) \text{ of }\s|_{i,\alpha,\beta'} 
&\raisebox{-5pt}{\rule{0pt}{16pt}}
(\lambda\bigcup\{2\ell\!+\!2\ell'\!+\!1\},r,\!-\!s)
&(\lambda\bigcup\{2\ell\!+\!1,2\ell'\!+\!1\},r,\!-\!s)
&(\lambda\bigcup\{2\ell\!+\!2\ell'\!+\!3\},r,\!-\!s)
\\ \hline
\end{array}$
\end{center} 

\item If $\s$ has cycle invariant $(\lambda\bigcup\{2\ell\},2r)$ we have \\[-5mm]
\begin{center}
$\begin{array}{|c|c|c|c|}
\cline{2-4}
\multicolumn{1}{c|}{}
&\raisebox{-5pt}{\rule{0pt}{16pt}}i=1&i=2&i=3\\ \hline
(\lam,r,s) \text{ of } \s|_{i,\alpha,\beta} 
&\raisebox{-5pt}{\rule{0pt}{16pt}}
(\lambda,2r\!+\!2\ell\!+\!1,s)
&(\lambda\bigcup\{2\ell\!+\!1\},2r\!+\!1,s)
&(\lambda,2r\!+\!2\ell\!+\!3,s)
\\ 
\hline
(\lam,r,s) \text{ of }\s|_{i,\alpha,\beta'} 
&\raisebox{-5pt}{\rule{0pt}{16pt}}
(\lambda,2r\!+\!2\ell\!+\!1,\!-\!s)
&(\lambda\bigcup\{2\ell\!+\!1\},2r\!+\!1,\!-\!s)
&(\lambda,2r\!+\!2\ell\!+\!3,\!-\!s)
\\ \hline
\end{array}$
\end{center} 
\end{enumerate}
\end{proposition}

\begin{pf}
Let us start with the cycle invariant. 
For the case $i=1$ and $i=2$ this is a strict application of propositions \ref{pro_one_edge} and \ref{pro_double-edge} respectively.

For $i=3$ note that $\s|_{3,\alpha,\beta}$ is obtained from $\s|_{1,\alpha,\beta}$ by adding a double-edge. In the first case $\s|_{1,\alpha,\beta}$ has cycle invariant $(\lambda\bigcup\{2\ell+2\ell'+1\},r)$ and the double-edge is inserted on two arcs of the cycle of length $2\ell+2\ell'+1$ thus by proposition \ref{pro_double-edge} the cycle invariant of $\s|_{3,\alpha,\beta}$ is $(\lambda\bigcup\{2\ell+2\ell'+3\},r)$.

 In the second case $\s|_{1,\alpha,\beta}$ has cycle invariant $(\lambda,2r+2\ell+1)$ and the double-edge is inserted on two arcs of the rank path of length $2r+2\ell+1$ thus by proposition \ref{pro_double-edge} the cycle invariant of $\s|_{3,\alpha,\beta}$ is $(\lambda,2r+2\ell+3)$.

 The same reasoning applies to $\s|_{3,\alpha,\beta'}.$

In the following the reasoning applies to both cases ($\s$ has cycle invariant $(\lambda\bigcup\{2\ell,2\ell'\},r)$ and $\s$ has cycle invariant $(\lambda\bigcup\{2\ell\},2r)$) so we will no longer differentiate.

By hypothesis on $\s$, $\lambda$ has no even part thus it can be verified on the two tables that the cycle invariant of $\s|_{i,\alpha,\beta}$ and $\s|_{i,\alpha,\beta'}$ has no even cycle for $i\leq3$. 

Consequently, $\Abar(\s|_{i,\alpha,\beta})=\pm 2^{\frac{n_i+\ell_i}{2}}$ and $\Abar(\s|_{i,\alpha,\beta'})=\pm 2^{\frac{n_i+\ell_i}{2}}$ where $n_i$ is the size of $\s|_{i,\alpha,\beta}$ and $\ell_i$ is the number of cycle (not including the rank) of $\s|_{i,\alpha,\beta}$ by theorem \ref{thm_arf_value}.

Let $n$ be the size of $\s$ and $\ell_0$ be the number of cycles (not including the rank) of $\s$. Then by inspecting the tables we know that 
\be\label{eq_sigma_1_0}
 \Abar(\s|_{1,\alpha,\beta})=\pm 2^{\frac{n+\ell_0}{2}}=\pm  \Abar(\s|_{1,\alpha,\beta'})
\ee  
\be\label{eq_sigma_2_0}
 \Abar(\s|_{2,\alpha,\beta})=\pm 2^{\frac{n+\ell_0}{2}+1}=\pm  \Abar(\s|_{2,\alpha,\beta'})
\ee  
 
and that 
\be\label{eq_sigma_0}
\Abar(\s)=0\ee
 by theorem \ref{thm_arf_value} since $\s$ has even cycles.

Moreover by applying proposition \ref{lem.certifopposign} to $\s|_{1,\alpha,\beta},\s|_{1,\alpha,\beta'}$ and $\s$  we have that
\begin{align*}
 \Abar(\s|_{1,\alpha,\beta})+ \Abar(\s|_{1,\alpha,\beta'})&=2\Abar(\s)&\\
&=0& \text{by eq \ref{eq_sigma_0}}
\end{align*}  
thus
\be\label{eq_sigma_1}
 \Abar(\s|_{1,\alpha,\beta})=- \Abar(\s|_{1,\alpha,\beta'})
\ee  

Now by applying lemma \ref{lem.certifopposign_2_3} (equation \ref{eq_opposign_2}) to $ \s|_{2,\alpha,\beta},\s|_{2,\alpha,\beta'}$ and $ \s|_{1,\alpha,\beta},\s|_{1,\alpha,\beta'}$ we have 
\begin{align*}
 \Abar(\s|_{2,\alpha,\beta})-\Abar(\s|_{2,\alpha,\beta'}) &=2(\Abar(\s|_{1,\alpha,\beta})-   \Abar(\s|_{1,\alpha,\beta'}))&\\
& =4\Abar(\s|_{1,\alpha,\beta}) &\text{by eq \ref{eq_sigma_1}} 
\end{align*}
Thus by equations \ref{eq_sigma_1_0} and \ref{eq_sigma_2_0} we must have 
\be\label{eq_sigma_2}
 \Abar(\s|_{2,\alpha,\beta})=- \Abar(\s|_{2,\alpha,\beta'})
\ee  

Finally by applying lemma \ref{lem.certifopposign_2_3} (equation \ref{eq_opposign_3}) to $ \s|_{3,\alpha,\beta},\s_{3,\alpha,\beta'}$ and $ \s|_{2,\alpha,\beta},\s|_{2,\alpha,\beta'}$ and $\s$ we have 
\begin{align*}
 \Abar(\s|_{3,\alpha,\beta})+\Abar(\s|_{3,\alpha,\beta'}) &=3(\Abar(\s|_{2,\alpha,\beta})+ \Abar(\s|_{2,\alpha,\beta'}))- 4\Abar(\s)  &\\
& =0 & \text{By eqs \ref{eq_sigma_2} and  \ref{eq_sigma_0}}
\end{align*}
thus \be \Abar(\s|_{3,\alpha,\beta})=-\Abar(\s_{3,\alpha,\beta'}).\ee
\qed
\end{pf}

The key point of the proposition is really that for $i\leq 3$, $\s|_{i,\alpha,\beta}$ and $\s|_{i,\alpha,\beta'}$ have same cycle invariant and opposite sign invariant. We will use this proposition to construct $I_2X$-permutations with same cycle invariant and opposite sign in the next section.

\begin{remark}\label{rk_dependencies_1}
The proof of proposition \ref{pro_opposite_sign} makes use of theorem \ref{thm_arf_value}. This theorem will only be proved during the induction, thus every time we use proposition \ref{pro_opposite_sign}, we must check that we have already proved theorem \ref{thm_arf_value} or indicate that the newly proven proposition is also dependent on theorem \ref{thm_arf_value}.
\end{remark}

\begin{proposition}\label{lem_q_i_sign}
 Let $\tau,(\Pi_b,\Pi_t)$ be a permutation with a consistent labelling and invariant $(\lam,r,s)$, then
\be\label{eq.546455b}
\Abar\bigg(\,\tau|_{1,t^{rk}_0,b^{rk}_{r}}=
%\t=
\raisebox{-12pt}{\reflectbox{\includegraphics[scale=1.8]{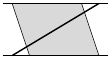}}}\,\bigg)
=0 
\ee
More generaly for $0\leq i\leq r,$
\be\label{eq.546455c}\Abar(\tau|_{1,t^{rk}_0,b^{rk}_{r-i}})=\begin{cases}
0, \text{ if } i \equiv 0 \mod 2\\
2\Abar(\tau) \text{ Otherwise.}
\end{cases}\ee
\end{proposition}
Note that the first equation corresponds to the case $i=0$ of the second equation.

\begin{pf}
The first equation is a straighforward application of proposition \ref{prop.Arf0}.
The second equation is derived from the first by induction on $i$.
As noted above the base case ($i=0$) of the induction is equation \ref{eq.546455b}. 
For the inductive step, we suppose equation \ref{eq.546455c} true for $i$. 

First case: if $i+1$ is even. Then $\Abar(\tau|_{1,t^{rk}_0,b^{rk}_{r-i}})=2\Abar(\tau)$ by induction, applying proposition \ref{lem.certifopposign}, we have that
 \[\Abar(\tau|_{1,t^{rk}_0,b^{rk}_{r-(i+1)}})+\Abar(\tau|_{1,t^{rk}_0,b^{rk}_{r-i}})=2\Abar(\tau)\] 
since the arcs $\beta=\Pi_b^{-1}(b^{rk}_{r-(i+1)})$ and $\beta'=\Pi_b^{-1}(b^{rk}_{r-i})$ are consecutive. Thus $\Abar(\tau|_{1,b^{rk}_0,t^{rk}_{r-(i+1)}})=0$.

Second case: if $i+1$ is odd is proved similarly.
\qed\end{pf}

Part of this identity can be found in the second half of proposition 5.11 of \cite{Del13}.\\

Now we have a relation that use the generalisation of the framework where the $(P_{\pm,i})$ can be permuted.
 
\begin{lemma}
\label{lem.signsurjq2}
\be 
\vec{A} \bigg(\,\tau_1=
\raisebox{-12pt}{\makebox[0pt][l]{\raisebox{32pt}{
\rule{3pt}{0pt}%
$A$
\rule{15pt}{0pt}%
$B$
}}}%
\raisebox{-12pt}{\makebox[0pt][l]{\raisebox{-8pt}{
\rule{5pt}{0pt}%
\phantom{$C$}
\rule{4pt}{0pt}%
\phantom{$e$$f$}
}}}%
\raisebox{-12pt}{\includegraphics[scale=1.8]{figure/fig_arfproof_s11-eps-converted-to.pdf}}\,\bigg)
=
\begin{pmatrix}
0 & -1 \\ -1 & 0
\end{pmatrix}
\vec{A} \bigg(\,\tau_2=
\raisebox{-12pt}{\makebox[0pt][l]{\raisebox{32pt}{
\rule{13pt}{0pt}%
$B$
\rule{15pt}{0pt}%
$A$
}}}%
\raisebox{-12pt}{\makebox[0pt][l]{\raisebox{-8pt}{
\rule{5pt}{0pt}%
\phantom{$C$}
\rule{4pt}{0pt}%
\phantom{$e$$f$}
}}}%
\raisebox{-12pt}{\includegraphics[scale=1.8]{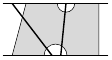}}\,\bigg)
\label{eq.546455x}
\ee
The letters $A$ and $B$ above the configurations denote that fact
that, contrarily to what was the case up to this point, the gray part
of the diagrams are permuted.
\end{lemma}

\proof We are in the extended framework: $\tau_1=\s_{1,1,E_1}$ and $\tau_2=\s_{1,1,E_2,\pi_{-},\pi_{+}}$ with $E_1=\{(1.1,1.1),(1.2,2.1)\}, E_2=\{(1.1,0.1),(1.1,1.1)\},\pi_{-}=id_2$ and $\pi_{+}=(2,1)$. \qed

We can now deduce a corollary, which is more relevant in the
following, and that involves three distinct configurations (of course could prove it directly with the theorem \ref{prop.Arf0} but since this was the first time we used the extended framework we choose to start with a simple case.)
\begin{proposition}
\label{cor.signsurjq2x}
\be
\Abar \bigg(\,\s=
\raisebox{-12pt}{\makebox[0pt][l]{\raisebox{32pt}{
\rule{3pt}{0pt}%
$A$
\rule{15pt}{0pt}%
$B$
}}}%
\raisebox{-12pt}{\makebox[0pt][l]{\raisebox{-8pt}{
\rule{5pt}{0pt}%
\phantom{$C$}
\rule{4pt}{0pt}%
\phantom{$e$$f$}
}}}%
\raisebox{-12pt}{\includegraphics[scale=1.8]{figure/fig_arfproof_s11-eps-converted-to.pdf}}\,\bigg)
+
\Abar \bigg(\,\tau=
\raisebox{-12pt}{\makebox[0pt][l]{\raisebox{32pt}{
\rule{13pt}{0pt}%
$B$
\rule{15pt}{0pt}%
$A$
}}}%
\raisebox{-12pt}{\makebox[0pt][l]{\raisebox{-8pt}{
\rule{5pt}{0pt}%
\phantom{$C$}
\rule{4pt}{0pt}%
\phantom{$e$$f$}
}}}%
\raisebox{-12pt}{\includegraphics[scale=1.8]{figure/fig_arfproof_t11-eps-converted-to.pdf}}\,\bigg)
=
\Abar \bigg(\,\rho=
\raisebox{-12pt}{\makebox[0pt][l]{\raisebox{32pt}{
\rule{13pt}{0pt}%
$A$
\rule{15pt}{0pt}%
$B$
}}}%
\raisebox{-12pt}{\makebox[0pt][l]{\raisebox{-8pt}{
\rule{5pt}{0pt}%
\phantom{$C$}
\rule{4pt}{0pt}%
\phantom{$e$$f$}
}}}%
\raisebox{-12pt}{\includegraphics[scale=1.8]{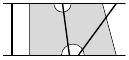}}\,\bigg)
\ee
\end{proposition}
\proof Combine equation (\ref{eq.546455a}) (or rather the mirror of it), and equation
(\ref{eq.546455x}) applied to $\rho$. \qed

\section{A $I_2X$-permutation for every $(\lambda,r,s)$\label{sec_IIXshaped}}
In this section, we will construct a non exceptional $I_2X$-permutation for every invariant $(\lambda,r,s)$.
\begin{center} \includegraphics[scale=.5]{figure/fig_goal_induction.pdf}\end{center}

 The tools we will use from the preceding sections are : \begin{itemize}
\item The double-edge insertion proposition \ref{pro_double-edge} to obtain the correct cycle invariant $(\lambda,r)$.
\item The opposite sign proposition \ref{pro_opposite_sign} to obtain the correct sign invariant $s.$
\end{itemize}

Let us give an overview of the section before detailling the actual construction of a $I_2X$-permutation for every a given $(\lam,r,s).$

Proposition \ref{pro_ins_cross_perm} and \ref{pro_ins_cross_perm_odd} allows to add a (or a pair of) cycle  of any given length into a permutation without breaking the property that is it $I_2X$. 

Propositions \ref{pro_X_family} and \ref{pro_X_family_odd} construct $I_2X$-permutations for some specific values of $(\lambda',r')$, more precisely for every $r'$ and for some $\lambda'$ of small cardinalities. We call those permutations \emph{base permutations}.

Finally Propositions \ref{pro_opposite_sign} and \ref{pro_two_oppo_sign} allows us to produce a $I_2X$-permutation with invariant $(\lam,r,-s)$ from a $I_2X$-permutation with invariant $(\lam,r,s)$.

The construction of a $I_2X$-permutation $\s$ for a given $(\lam,r,s)$, proceeds in three steps. Of course, at each step of the construction we garantee that the resulting permutations are still $I_2X$
\begin{enumerate}
\item We choose a base permutation $\s^0$ with invariant $(\lam',r)$ such that $\lam' \subseteq \lam$.
\item We add cycles on $\s^0$ until the new permutation $\s^1$ has invariant $(\lam,r)$. 
\item If $\s^1$ has no even cycle then the sign of $\s^1$ is $\pm 1$ (refer to theorem \ref{thm_arf_value} for a reminder of the relationship between the sign invariant $s$ and the cycle invariant $(\lam,r)$)). Since we have no control on whether it is $+1$ or $-1$, we need to construct another permutation $\s^2$ with opposite sign so as to insure that either $\s^1$ or $\s^2$ has invariant $(\lam,r,s)$.

This is done by using proposition \ref{pro_two_oppo_sign} in most cases and by proposition \ref{pro_opposite_sign} for a few remaining cases. The constuction is the subject of theorem \ref{thm_I_2X}. 

 If $\s^1$ has even cycles then the sign of $\s^1$ is necessarily 0 and we are done. The construction is the subject of theorem \ref{thm_I_2X_odd}
\item Finally we will verify (in a remark after the theorem) that the constructed permutations are not exceptional.
\end{enumerate}

\noindent Let us define $C_p$ for any $p\in \mathbb{N}$ the cross permutation
\begin{center}
 $\makebox[0pt][l]{\raisebox{14pt}{$C_{p}=\,$}  \raisebox{0mm}{\includegraphics[scale=.5]{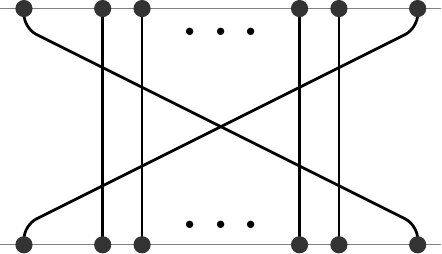}}}
\hspace{42pt}\raisebox{34pt}{$\rotatebox{90}{$\left. \rule{0pt}{22.5pt}\right\}\rotatebox{-90}{$\!\!p$} $}$}$
\end{center}

\begin{figure}[h!]
\begin{center}
\begin{tabular}{lll}
\makebox[0pt][l]{\raisebox{0mm}{\includegraphics[scale=.5]{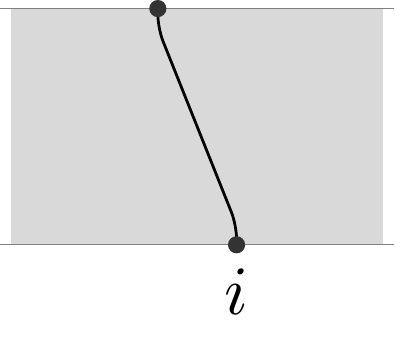}}}
&$\qquad\qquad\qquad$&
\makebox[9pt][l]{\raisebox{4.2mm}{\includegraphics[scale=.5]{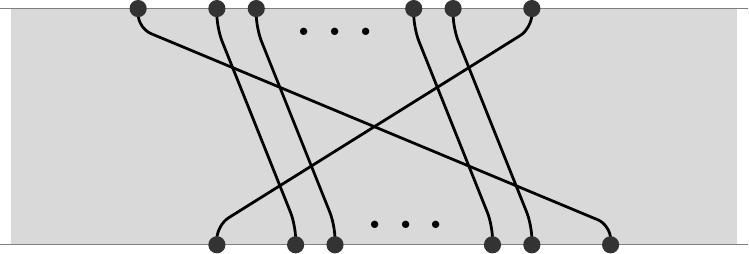}}}
\hspace{15pt}\raisebox{34pt+4.2mm}{$\rotatebox{90}{$\left. \rule{0pt}{22.5pt}\right\}\rotatebox{-90}{$\!\!p$} $}$}
\end{tabular}
\end{center}
\caption{\label{fig_cross_perm} Left: a permutation $\s$. Right : the permutation $\s_i(C_p)$. }
\end{figure}

\begin{proposition}[Adding cycles 1]\label{pro_ins_cross_perm}
let $\s$ be a permutation with cycle invariant $(\lambda,r)$, we call $\s_i(C_p)$ the permutation obtained by replacing the ith edge of $\s$ by a cross permutation $C_p$ (see figure \ref{fig_cross_perm}). The cycle invariant $(\lambda',r)$ of $\s_i(C_p)$ depends on $p$ in the following way :
\begin{itemize}
\item if $p=4k$ then $\lambda'=\lambda\bigcup\{p+1\}$.
\item if $p=4k+1$ then $\lambda'=\lambda\bigcup\{2k+1,2k+1\}$.
\item if $p=4k+2$ then $\lambda'=\lambda\bigcup\{p+1\}$.
\item if $p=4k+3$ then $\lambda'=\lambda\bigcup\{2k+2,2k+2\}$.
\end{itemize}
\end{proposition}

\begin{pf}
By induction on $p$. The base cases are for $p=0$ and $p=1$. 

For those cases, the permutation $\s$ with invariant $(\lambda,r)$ becomes respectively $\s_i(C_0)$  with invariant $(\lambda \bigcup\{1\},r)$ and $\s_i(C_1)$ with invariant $(\lambda \bigcup\{1,1\},r)$, as displayed in figure \ref{fig_pf_base_case} :\\
\begin{figure}
\begin{center}
\begin{tabular}{ccc}
\includegraphics[scale=.5]{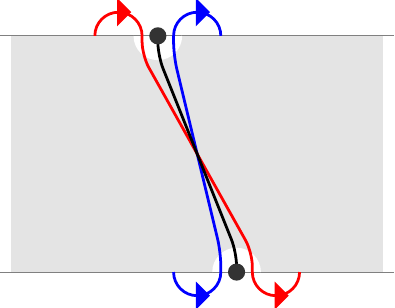}
&
\includegraphics[scale=.5]{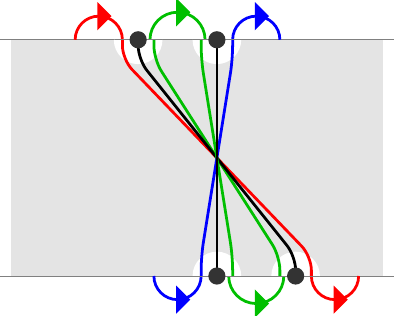}
&
\includegraphics[scale=.5]{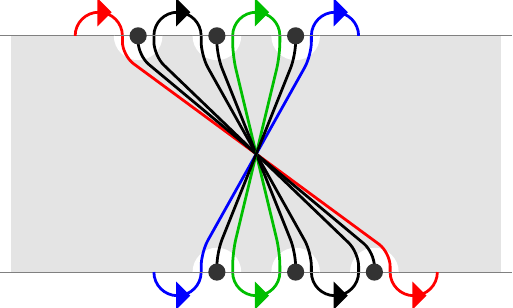}\\
$\s$, $(\lambda,r)$&
 $\s_i(C_0)$, $(\lambda \bigcup\{1\},r)$&
$\s_i(C_1)$,  $(\lambda \bigcup\{1,1\},r)$
\end{tabular}
\end{center}
\caption{\label{fig_pf_base_case} The base cases of the induction of the proof of proposition \ref{pro_ins_cross_perm}}
\end{figure}
Then the statement follows by induction from the insertion of a double-edge (the resulting change of the cycle invariant are described in proposition \ref{pro_double-edge}).\qed
\end{pf}

Note that if $p=4k+3$ the two cycles of the $C_p$ structure attached on $\s'=\s_i(C_p)$ are even. Let $\alpha$ be the first top arc of the $C_p$ structure of $\s'$ and $\beta$ and $\beta'$ its first and third bottom arcs (refer to figure \ref{fig_change_sign} left). Then by proposition \ref{pro_opposite_sign} for $j=1,2,3$ the sign invariant of $\s'|_{j,\alpha,\beta} $ and  $\s'|_{j,\alpha,\beta'}$ are opposite while their cycle invariants remain egal (figure \ref{fig_change_sign} middle and right). 

For clarity, let us call $C_{p,j}$ for any $p,j$ the permutation: 
\begin{center}
\raisebox{14pt}{$C_{p,j}=\,$}  \raisebox{0mm}{\includegraphics[scale=.5]{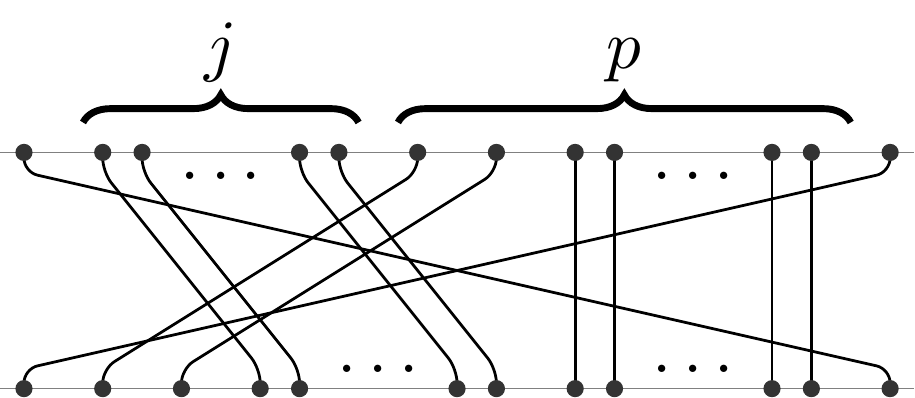}}
\end{center}

Then we have $\s'|_{j,\alpha,\beta}=\s_i(C_{p+j})$ and $\s'|_{j,\alpha,\beta'}=\s_i(C_{p,j})$ and our discussion implies the following statement:

\begin{figure}[t!]
\begin{center}
\begin{tabular}{ccc}
\includegraphics[scale=.5]{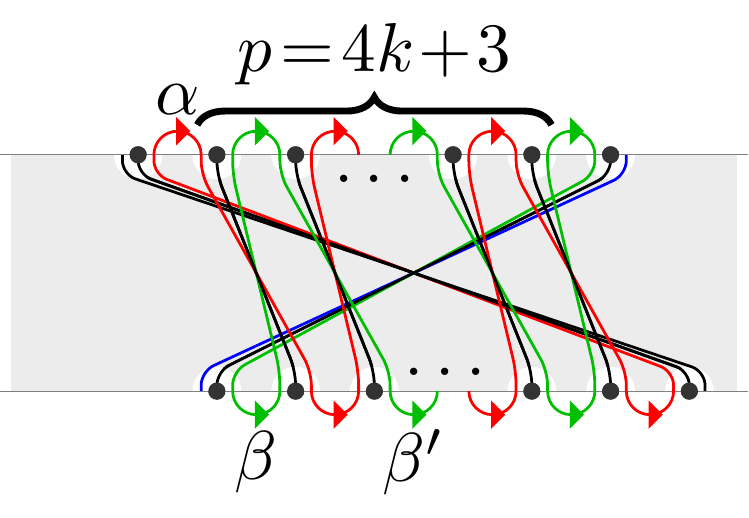}
&
\includegraphics[scale=.5]{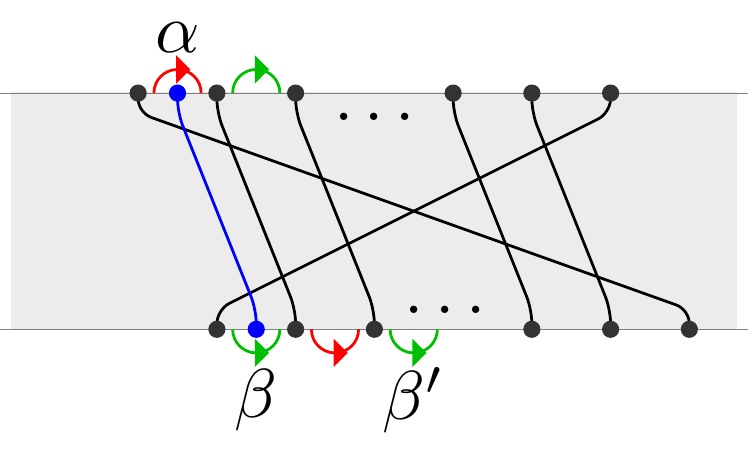}
&
\includegraphics[scale=.5]{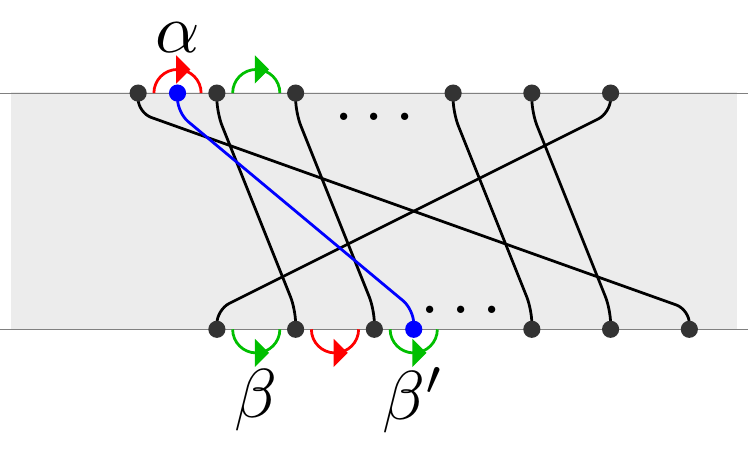}\\
$\s'=\s_i(C_p)$ & $\s'|_{1,\alpha,\beta}=\s_i(C_{p+1})$ & $\s'|_{1,\alpha,\beta'}=\s_i(C_{p,1})$ 
\end{tabular}
\end{center}
\caption{\label{fig_change_sign} Left: a permutation with a $C_p$ structure containing two even cycles of even length. Middle and Right: By proposition \ref{pro_opposite_sign} the two permutations have same cycle invariant and opposite sign.}
\end{figure}

\begin{proposition}[Two opposite signs]\label{pro_two_oppo_sign}
Let $\s$ be a permutation and let $p=4k+j$ with $0\leq j<3$ and $k>0$. Then $\s_i(C_p)$ and $\s_i(C_{p-(j+1),j+1})$ have invariant $(\lam,r,s)$ and $(\lam,r,-s)$ respectively for some $(\lam,r,s)$.
\end{proposition}

This statement will be crucial for our construction in theorem \ref{thm_I_2X}. Indeed, as outlined in the beginning of the section, we will construct for every $(\lam,r,s)$ a permutation with invariant $(\lam,r,s)$ by adding cycles on a base permutation with invariant $(\lam'\subseteq \lam,r)$. This construction will be performed through the means of proposition \ref{pro_ins_cross_perm}) .

Thus in the final step of this procedure, we construct $\s_{1}=\s_i(C_p)$ with invariant $(\lam,r)$. In order to insure that our constructed permutation has sign invariant $s$, we also consider $\s'_{1}=\s_i(C_{p-(j+1),j+1})$ for a correct $j$. Then by proposition \ref{pro_two_oppo_sign} either $\s_{1}$  or $\s'_{1}$ will have invariant $(\lam,r,s)$.\\

\noindent The last ingredient of our proof of theorem \ref{thm_I_2X} are the base permutations. 

The following proposition provides the base permutations for the case where $(\lambda,r,s)$ has no even cycle (first line of figure \ref{fig_perm_x_family}). It will also allow us to apply proposition \ref{pro_opposite_sign} to obtain two permutations with opposite sign (second line of figure \ref{fig_perm_x_family}) for the few cases not covered by proposition \ref{pro_two_oppo_sign}.

Let us define the permutations $X_{p,p'}$, $X_{p,p',p''}$ for any $p,p',p''$ to be 
\[\raisebox{14pt}{$X_{p,p'}=\,$} \makebox[0pt][l]{\raisebox{0mm}{\includegraphics[scale=.5]{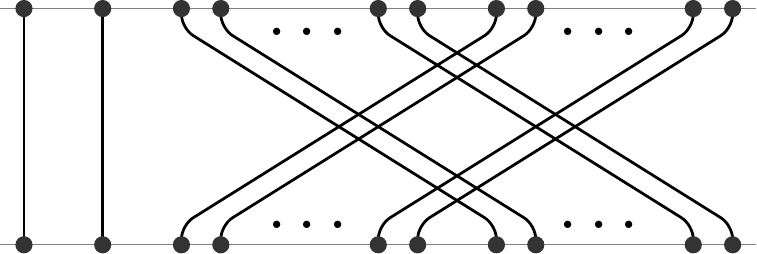}}}
\hspace{21pt}\raisebox{34pt}{$\rotatebox{90}{$\left. \rule{0pt}{22.5pt}\right\}\rotatebox{-90}{$\!\!\! p$} $}$}
\hspace{7pt}\raisebox{34pt}{$\rotatebox{90}{$\left. \rule{0pt}{22.5pt}\right\}\rotatebox{-90}{$\!\!\!p'$} $}$}
 \qquad\raisebox{14pt}{$X_{p,p',p''}=\,$} \makebox[0pt][l]{\raisebox{0mm}{\includegraphics[scale=.5]{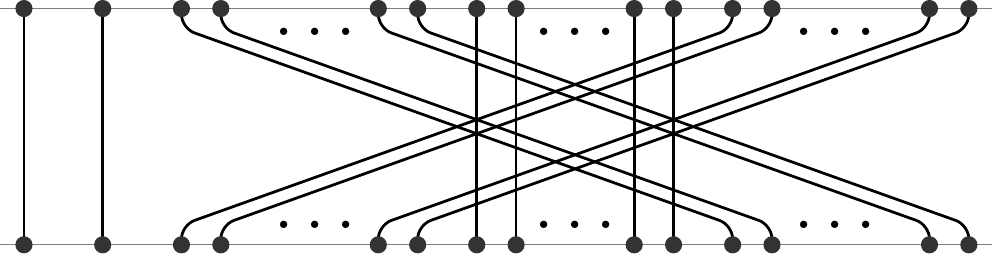}}}
\hspace{21pt}\raisebox{34pt}{$\rotatebox{90}{$\left. \rule{0pt}{22.5pt}\right\}\rotatebox{-90}{$\!\!\! p$} $}$}
\hspace{6pt}\raisebox{34pt}{$\rotatebox{90}{$\left. \rule{0pt}{19.5pt}\right\}\rotatebox{-90}{$\!\!\!p'$} $}$}
\hspace{4pt}\raisebox{34pt}{$\rotatebox{90}{$\left. \rule{0pt}{22.5pt}\right\}\rotatebox{-90}{$\!\!\!p''$} $}$}\]
Then we have:

\begin{proposition}[Base permutations, no even cycle] \label{pro_X_family}
For every $k \geq 1.$

The permutations $X_{1,2k}$, $X_{2,1,2k}$ and $X_{2,2k}$ (described in the first line of figure \ref{fig_perm_x_family} have cycle invariant $(\lambda=\{2k\!+\!1\},r=1)$, $(\lambda=\{2k\!+\!1\},r=3)$ and $(\lambda=\emptyset,r=2k\!+\!3)$ respectively.

The permutations $X_{1,4k+3}$, $X_{2,1,4k+3}$ and $X_{2,4k+3}$ (described in the second line of figure \ref{fig_perm_x_family}) have cycle invariant $(\lambda=\{{\color{red}2k\!+\!2},{\color{green!50!black}2k\!+\!2}\},r=1)$, $(\lambda=\{{\color{red}2k\!+\!2},{\color{green!50!black}2k\!+\!2}\},r=3)$ and $(\lambda=\{{\color{red}2k\!+\!2}\},r={\color{green!50!black}2k\!+\!4})$ respectively.
\begin{figure}[h!]
\begin{center}
\begin{tabular}{lll}
\makebox[0pt][l]{\raisebox{0mm}{\includegraphics[scale=.46]{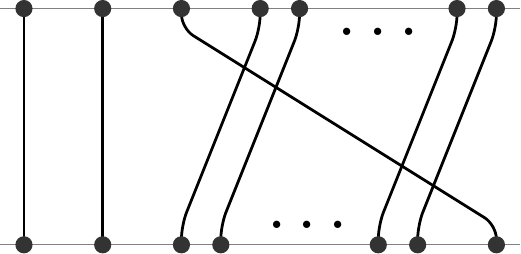}}}
\hspace{28pt}\raisebox{31pt}{$\rotatebox{90}{$\left. \rule{0pt}{21.5pt}\right\}\rotatebox{-90}{$\!\!\!\! 2k$} $}$}
&
\makebox[0pt][l]{\raisebox{0mm}{\includegraphics[scale=.46]{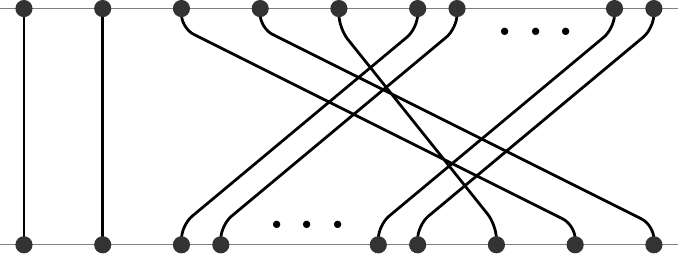}}}
\hspace{49pt}\raisebox{31pt}{$\rotatebox{90}{$\left. \rule{0pt}{21.5pt}\right\}\rotatebox{-90}{$\! \! \! \! 2k$} $}$}&
\makebox[0pt][l]{\raisebox{0mm}{\includegraphics[scale=.46]{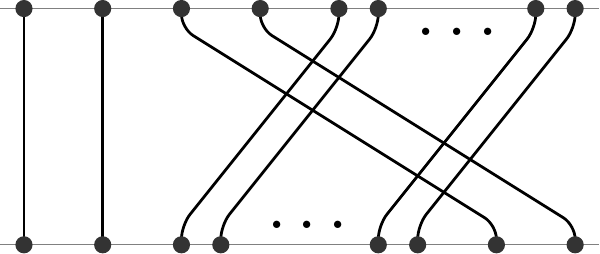}}}
\hspace{38.5pt}\raisebox{31pt}{$\rotatebox{90}{$\left. \rule{0pt}{21.5pt}\right\}\rotatebox{-90}{$\! \! \! \! 2k$} $}$}\\
$\lambda=\{2k\!+\!1\},r=1$ & $\lambda=\{2k\!+\!1\},r=3$ & $\lambda=\emptyset,r=2k\!+\!3$\\
&&\\

\makebox[0pt][l]{\raisebox{0mm}{\includegraphics[scale=.46]{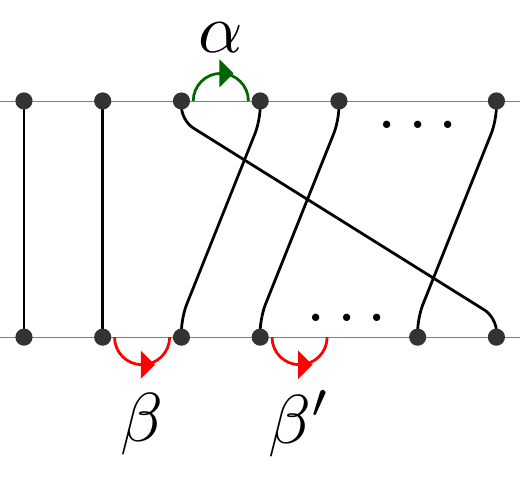}}}
\hspace{29pt}\raisebox{49pt}{$\rotatebox{90}{$\left. \rule{0pt}{21.5pt}\right\}\rotatebox{-90}{$\!\!\!\!\!\!\!\!\! 4k\!+\!3$} $}$}
&
\makebox[0pt][l]{\raisebox{0mm}{\includegraphics[scale=.46]{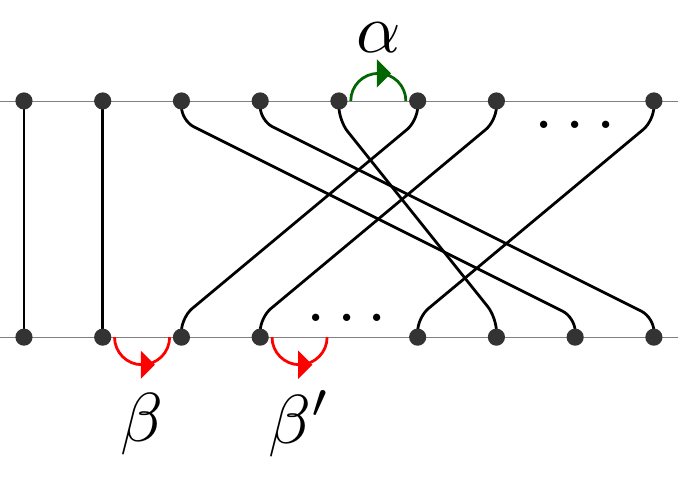}}}
\hspace{50pt}\raisebox{49pt}{$\rotatebox{90}{$\left. \rule{0pt}{21.5pt}\right\}\rotatebox{-90}{$\!\!\!\!\!\! \! \! \! 4k\!+\!3$} $}$}&
\makebox[0pt][l]{\raisebox{0mm}{\includegraphics[scale=.46]{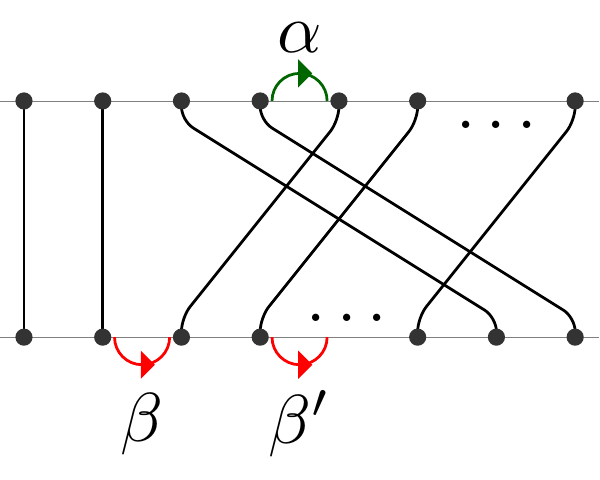}}}
\hspace{39.5pt}\raisebox{49pt}{$\rotatebox{90}{$\left. \rule{0pt}{21.5pt}\right\}\rotatebox{-90}{$\!\!\!\!\!\! \! \! \! 4k\!+\!3$} $}$}\\
$\lambda=\{{\color{red}2k\!+\!2},{\color{green!50!black}2k\!+\!2}\},r=1$ & $\lambda=\{{\color{red}2k\!+\!2},{\color{green!50!black}2k\!+\!2}\},r=3$ & $\lambda=\{{\color{red}2k\!+\!2}\},r={\color{green!50!black}2k\!+\!4}$
\end{tabular}
\end{center}
\caption{\label{fig_perm_x_family} Two families of base permutations with their respective cycle invariant.}
\end{figure}
\end{proposition}

\begin{pf}
By induction on $k$. The base cases for the first family are respectively :\begin{center}
 \scalebox{.55}{\begin{tikzpicture}[scale=.8]\permutationmatching{1,2,4,5,3}\end{tikzpicture}},
 \scalebox{.55}{\begin{tikzpicture}[scale=.8]\permutationmatching{1,2,6,7,5,3,4}\end{tikzpicture}} and
\scalebox{.55}{\begin{tikzpicture}[scale=.8]\permutationmatching{1,2,5,6,3,4}\end{tikzpicture}}
.\end{center}

And they have cycle invariant $(\lambda,r)$ :  $(\{3\},1)$,  $(\{3\},3)$ $(\emptyset,5)$ as shown just below :
\begin{center}\includegraphics[scale=1]{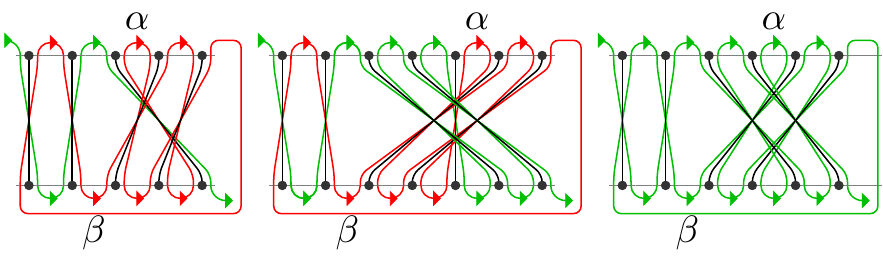}\end{center}

Then the statement follows by induction from the insertion a double-edge within $\alpha$ and $\beta$ (the resulting change of the cycle invariant are described in proposition \ref{pro_double-edge}).

The base cases for the second family are respectively: 
\begin{center}
 \scalebox{.55}{\begin{tikzpicture}[scale=.8]\permutationmatching{1,2,4,5,6,3}\end{tikzpicture}},
 \scalebox{.55}{\begin{tikzpicture}[scale=.8]\permutationmatching{1,2,6,7,8,5,3,4}\end{tikzpicture}} and
\scalebox{.55}{\begin{tikzpicture}[scale=.8]\permutationmatching{1,2,5,6,7,3,4}\end{tikzpicture}}
.\end{center}

And they have cycle invariant $(\lambda,r)$ :  $\lambda=\{{\color{red}2},{\color{green!50!black}2}\},r=1$, $\lambda=\{{\color{red}2},{\color{green!50!black}2}\},r=3$, $\lambda=\{{\color{red}2}\},r={\color{green!50!black}4}$ as shown just below :
\begin{center}\includegraphics[scale=1]{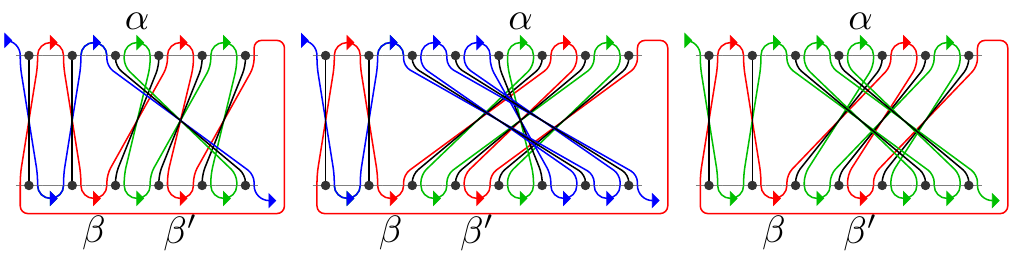}\end{center}

Then the statement follows by induction from the insertion of two double-edges within $\alpha$ and $\beta$ (the resulting change of the cycle invariant are described in proposition \ref{pro_double-edge}). \qed
\end{pf}

We can finally state and prove the first theorem of this section.

\begin{theorem}\label{thm_I_2X}
Let $(\lambda,r,s)$ be an invariant with no even cycle, then there exists a $I_2X$-permutation with invariant $(\lambda,r,s)$.
\end{theorem}

\begin{pf}
We can always consider that the size of the permutation is at least 10, for smaller size the result can be obtained by automatic search.

Let $(\lam,r,s)$ be an invariant with no even cycles. Following the proof sketch of the introduction, we fist construct a base permutation $\s^0$ with invariant $(\lam'\subseteq\lam,r)$ then add cycles to obtain a $I_2X$-permutation $\s_1$ with invariant $(\lam,r)$ and finally we use either proposition \ref{pro_two_oppo_sign} or proposition \ref{pro_opposite_sign} to obtain two $I_2X$-permutations with opposite sign. 

\begin{itemize}
 \item If the rank $r=1$, then the base permutation $\s^0$ with invariant $(\lambda'=\{2\ell+1\}, 1)$ for any $\ell\geq1$ is exactly $X_{2\ell}$ according to proposition \ref{pro_X_family} (first line, first case of figure \ref{fig_perm_x_family}).
 \item If the rank $r=3$, then the base permutation $\s^0$ with invariant $(\lambda'=\{2\ell+1\}, 3)$ for any $\ell\geq1$ is exactly $X_{2,1,2\ell}$ according to proposition \ref{pro_X_family} (first line, second case of figure \ref{fig_perm_x_family})).
\item If the rank $r>5$, then the base permutation $\s^0$ with invariant $(\lambda'=\emptyset, r)$ is exactly $X_{2,r-3}$ according to proposition \ref{pro_X_family} (first line, third case of figure \ref{fig_perm_x_family} indeed $r$ odd implies $r-3$ is even) ).
\end{itemize}

Next we add (by the means of proposition \ref{pro_ins_cross_perm}) cycles one by one and then by pair on the last edge of the current permutation $\s^i$ so as to make sure that the last $C_p$ attached is not a $C_2$. This is always possible unless $\lam=\lam'$ (in which case $\s_1=\s^0$) or $\lam'=\lam\bigcup\{3\}$. In this case we only have to add a cycle a length 3 to the permutation and the procedure begins and ends with the attachement of a $C_2$: $\s_1=\s^0_n(C_2)$). 

The procedure must finish with the addition of a $C_p$ with $p>2$ in order to make the application of proposition \ref{pro_two_oppo_sign} possible.

 Note that since we always attach $C_p$ on the last edge $(|\s^i|,\s^i(|\s^i|))$ of the current permutation $\s^i$, the successive permutations from $\s^0$ to $\s_1$ are all $I_2X$. See figure \ref{fig_add_cycle_to_IIX}.

\begin{figure}[tb]
\begin{center}
\hspace*{-2cm}\makebox[90pt][l]{\raisebox{0mm}{\includegraphics[scale=.5]{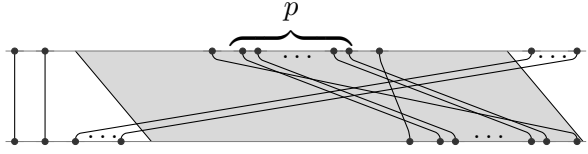}}}
\hspace{-10pt}\raisebox{35pt}{$\rotatebox{90}{$\left. \rule{0pt}{26pt}\right\}\rotatebox{-90}{$\!\!\!p$} $}$}
\end{center}
\caption{\label{fig_add_cycle_to_IIX} The permutation resulting from adding a $C_p$ on the last edge of $\s^0$ or one of its successor. It is clearly also $I_2X$.}
\end{figure}

Let $2\ell_1+1<\ldots<2\ell_k+1$ be the length and $(m_i)_{1\leq i\leq k}$ the multiplicity of the cycles to be added on $\s^0$ (i.e. the cycles of $\lam\setminus \lam'$). The procedure is divided into two steps:

$\bullet$ First we look at the parity of $m_i$ from $i=1$ to $i=k$ and if $m_i$ is odd we attach a $C_{2\ell_i}$ on the last edge of the current permutation (which adds a cycle of length $2\ell_i+1$ to the cycle invariant). 

More precisely, let $i_1,\ldots,i_m$ be the indices such that the $(m_{i_j})_j$ are odd, we define:
\[ \s^{j}=\s^{j-1}_{|\s^{j-1}|}(C_{2\ell_j})\text{ for }1 \leq j\leq m.\] 
Let $(\lam^j,r)$ be the cycle invariant of $\s^j$, by proposition \ref{pro_ins_cross_perm}, we have $\lam^j=\lam^{j-1}\bigcup \{2\ell_j+1\}.$
Thus the multiplicities $(2m'_i)_{1\leq i\leq k}$ of the cycles of length $(2\ell_i+1)_{1\leq i\leq k}$ to be added on $\s^{m}$ are all even.

$\bullet$ In the second step, we attach a $C_{4\ell_i+1}$ on the last edge of the current permutation (which add two cycles of length $2\ell_i+1$ to the cycle invariant) consecutively $m_i/2$ time for $i=1$ to $i=k$.

More specifically, we define $\s^{j,i}$ and its cycle invariant $(\lam^{j,i},r)$ by \[ \begin{cases}
 \s^{j,i}= \s^{j-1,i}_{|\s^{j-1,i}|}(C_{4\ell_i+1}) \text{ and } \lam^{j,i}= \lam^{j-1,i}\bigcup \{2\ell_i\!+\!1,2\ell_i\!+\!1\} \!\! \!\!& \text{ for } 1\!\leq \! j\leq\! m'_i/2 \text{ and } 1\!\leq\! i \!\leq\! k\\[1mm]
 \s^{0,i}=\s^{m'_{i-1}/2,i-1} \text{ and } \lam^{0,i}= \lam^{m'_{i-1}/2,i-1} & \text{ for  } 1\!<\! i \!\leq\! k,\\
 \s^{0,1}=\s^m \text{ and } \lam^{0,1}= \lam^{m}.
\end{cases}\] Once again, the values of the cycle invariants are justified by proposition \ref{pro_ins_cross_perm}

Let us call the permutation obtained by the procedure $\s_1$. 

By construction $\s_1$ has cycle invariant $(\lam,r)$ and the last added $C_p$ by the procedure is $C_2$ if and only if there was just one cycle of length 3 to be added.  \\

Let us now deal with the sign invariant.
 
We can divide the problem in two cases : The last $C_p$ attached has $p>2$ or not (the latter case corresponds to attaching exactly one $C_2$ or not attaching anything, thus remaining with the base permutation). Moreover, in accordance with the procedure, $p=4k+j$, $k>0$ and $0\geq j<3$. Thus $p$ and k verify the condition

$\ast$ If the last $C_p$ attached has $p>2$. We have $\s_1=\s_{|\s|}(C_p)$ for some $\s$ (more specifically $\s=\s^{m'_k/2-1, k}$ the second to last permutation of the procedure) and $\s_1$ has the form shown in figure \ref{fig_add_cycle_to_IIX}. 

Moreover, in accordance with the procedure, $p=4k+j$, $k>0$ and $0\geq j<3$. Thus $p$ verify the condition of proposition \ref{pro_two_oppo_sign} and $\s'_1=\s_{|\s|}(C_{p-(j+1),j+1})$ has invariant cycle $(\lam,r)$ and sign invariant opposite to $\s_1$, thus either one of them has invariant $(\lam,r,s)$. By construction $\s'_1$ is also $I_2X$. (See figure \ref{fig_I_2X_j2} for a representation of $\s_1$ and $\s'_1$ for $j=2$).

\begin{changemargin}{-1.15cm}{0cm} 
\begin{figure}[t]
\begin{center}
\begin{tabular}{cc}
\includegraphics[scale=.46]{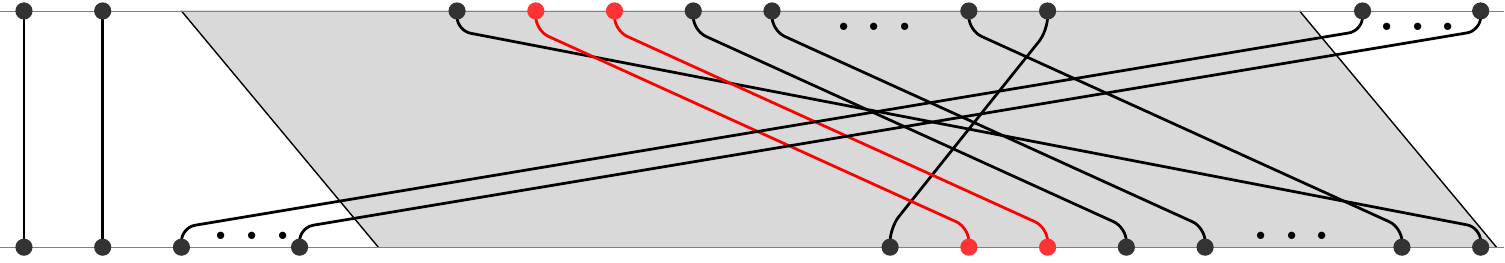}
&
\includegraphics[scale=.46]{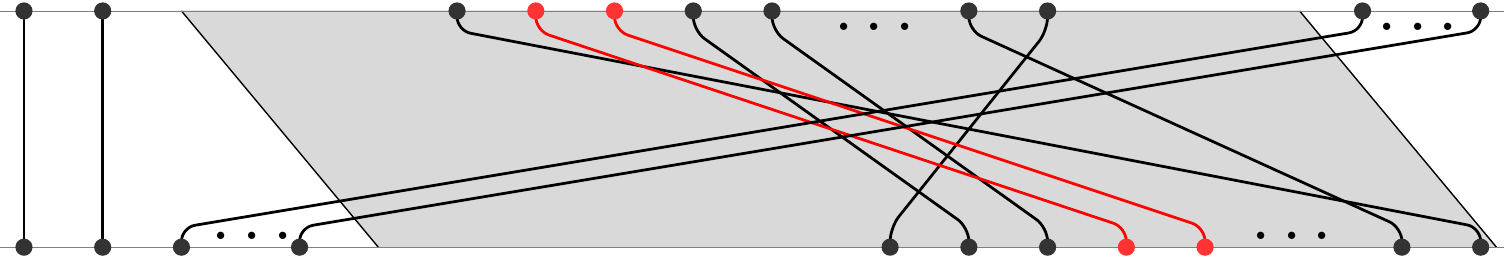}

\\
 $\s_1$ & $ \s'_1$
\end{tabular}
\end{center}\caption{\label{fig_I_2X_j2} The two constructed $I_2X$-permutations $\s_1=\s_{|\s|}(C_p)$ and $\s'_1=\s_{|\s|}(C_{p-(2),2})$ with invariant $(\lam,r,\pm s)$. In this case $p=4k+1$.}
\end{figure}
\end{changemargin}

$\ast$ If no $C_p$ or just one $C_2$ are attached. We will only consider the case no $C_p$ are added, if a $C_2$ is added on the last edge, the reasoning is strictly identical since adding a $C_p$ does not change the already existing cycles and arcs of the permutation (it only adds cycles and consecutive arcs).
Futhermore among the three cases $X_{1,2\ell}$, $X_{2,1,2\ell}$ and $X_{2,r-3}$ we will only handle the first one, the other two being similar. 

Let $\s^0=\s_1=X_{1,2\ell}$ for some $\ell\geq2$ (this is always the case since otherwise $X_{1,2\ell}$ has size 5)
\begin{itemize}
\item If $2\ell=4k$ then removing the first parallel edge of the $2\ell$ consecutives and parallel edges of $\s_1=X_{1,2\ell}$ we get $\tau=X_{1,4(k-1)+1}$ by proposition \ref{pro_opposite_sign} and proposition \ref{pro_X_family} (second line first case of figure \ref{fig_perm_x_family}) the permutations $ \tau|_{1,\alpha,\beta}=\s^1$ and $ \tau|_{1,\alpha,\beta'}$ have same cycle invariant and opposite sign.

\item If $2\ell=4k+2$ then removing the first three parallel edge of the $2\ell$ consecutives and parallel edges of $\s_1=X_{1,2\ell}$ we get $\tau=X_{1,4(k-1)+1}$ by proposition \ref{pro_opposite_sign} and proposition \ref{pro_X_family}(second line first case of fig \ref{fig_perm_x_family}) the permutations $ \tau|_{3,\alpha,\beta}=\s^1$ and $ \tau|_{3,\alpha,\beta'}$ have same cycle invariant and opposite sign.
\end{itemize}
Schematically the two cases (for $X_{1,2\ell}$, $X_{2,1,2\ell}$ and $X_{2,r-3}$ with or without a $C_2$ attached on the last edge) are :\begin{center}
\begin{tabular}{cc}
\includegraphics[scale=.46]{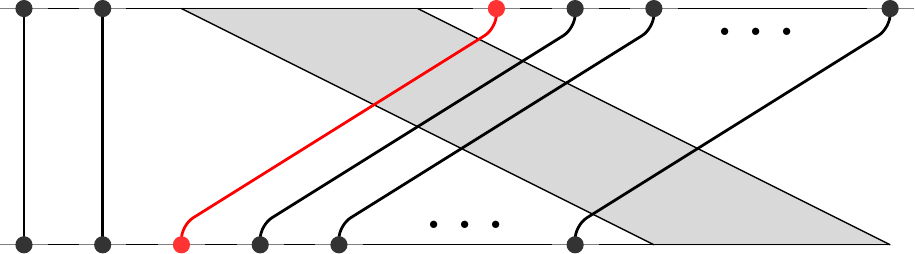}
&
\includegraphics[scale=.46]{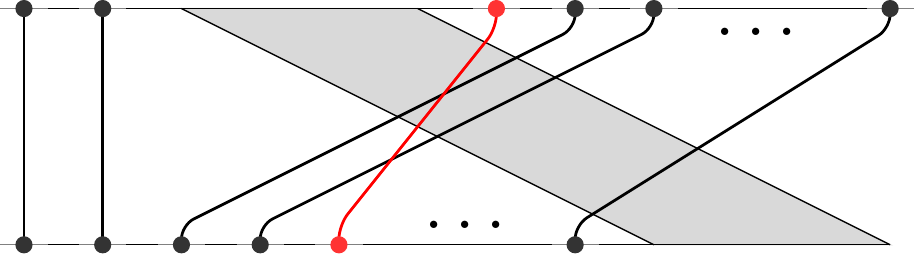}

\\
 $ \tau|_{1,\alpha,\beta}=\s^1$ & $ \tau|_{1,\alpha,\beta'}=\s^2$
\end{tabular} 
\end{center}
and 
\begin{center}
\begin{tabular}{cc}
\includegraphics[scale=.46]{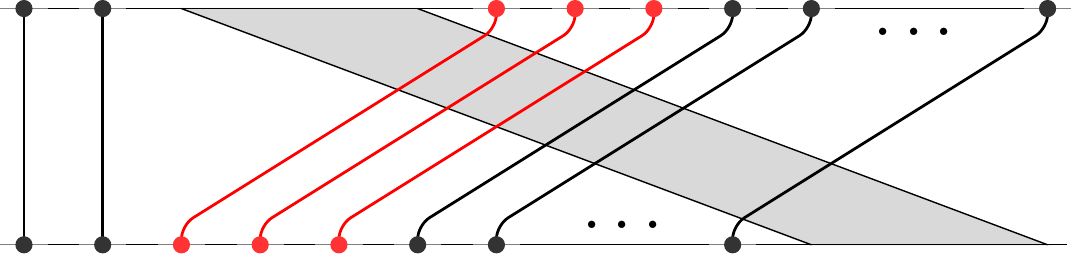}
&
\includegraphics[scale=.46]{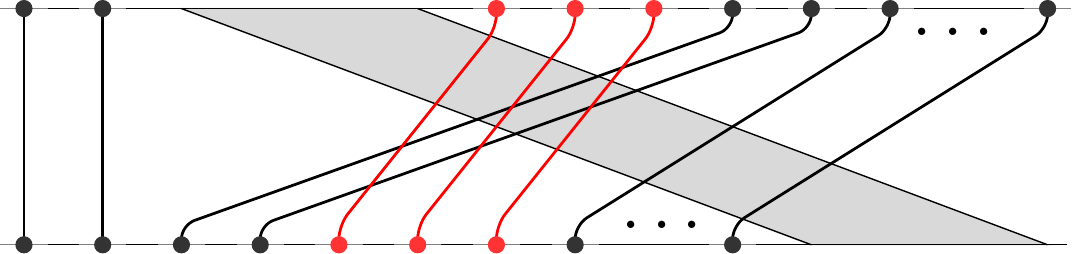}

\\
 $ \tau|_{3,\alpha,\beta}=\s^1$ & $ \tau|_{3,\alpha,\beta'}=\s^2$
\end{tabular}
\end{center}

Finally the permutations obtained are clearly $I_2X$.
\qed
\end{pf}

\begin{remark}\label{rk_dependencies_2}
Theorem \ref{thm_I_2X} makes use of proposition \ref{pro_opposite_sign} directly and indirectly in its use of proposition \ref{pro_two_oppo_sign}. Therefore theorem \ref{thm_I_2X} is also dependent on theorem \ref{thm_arf_value} and we must check that we have proved theorem \ref{thm_arf_value} before using theorem \ref{thm_I_2X}.
\end{remark}

\begin{remark}
It is unfortunate that in the case $r=1$ and $\lam=\{2\ell+1\}$ one of the two permutations with invariant $(\lam,r,\pm s)$ produced by theorem \ref{thm_I_2X} is exactly $X_{1,2\ell}=\id'_{2\ell+3}$ i.e. a permutation from an exceptional class. We rectify this by constructing a third permutation with the same invariant as $X_{1,2\ell}$.

\noindent $\bullet$ if $2\ell=4k$
\begin{center}
\begin{tabular}{lll}
\makebox[0pt][l]{\raisebox{0mm}{\includegraphics[scale=.46]{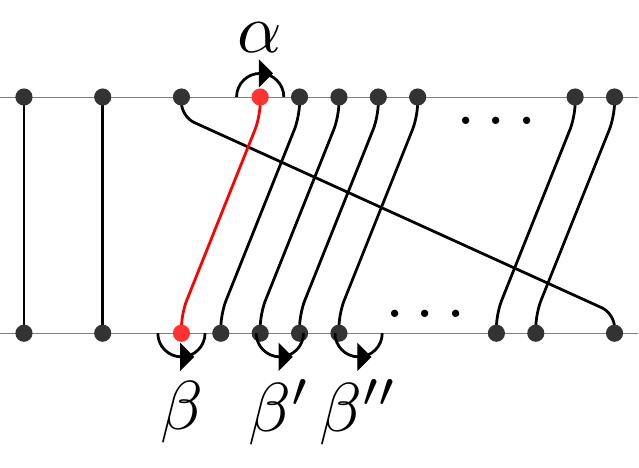}}}
\hspace{36pt}\raisebox{48pt}{$\rotatebox{90}{$\left. \rule{0pt}{26pt}\right\}\rotatebox{-90}{$\!\!\!\!\!\!\!\!\!\!\!\!\!\! 4(k\!-\!1)\!+\!3$} $}$}
&
\makebox[0pt][l]{\raisebox{0mm}{\includegraphics[scale=.46]{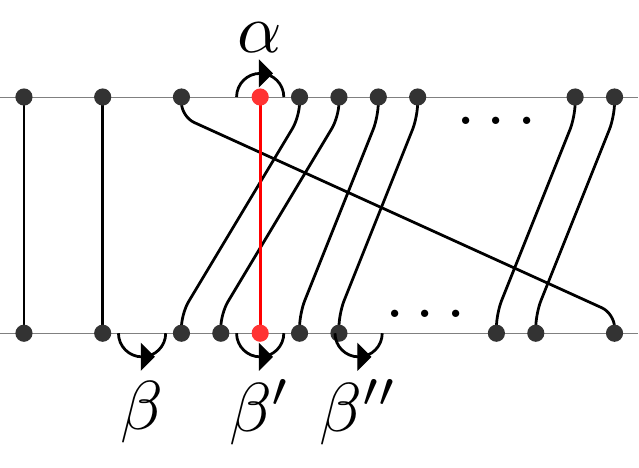}}}
\hspace{36pt}\raisebox{48pt}{$\rotatebox{90}{$\left. \rule{0pt}{26pt}\right\}\rotatebox{-90}{$\!\!\!\!\!\!\!\!\!\!\!\!\!\! 4(k\!-\!1)\!+\!3$} $}$}&
\makebox[0pt][l]{\raisebox{0mm}{\includegraphics[scale=.46]{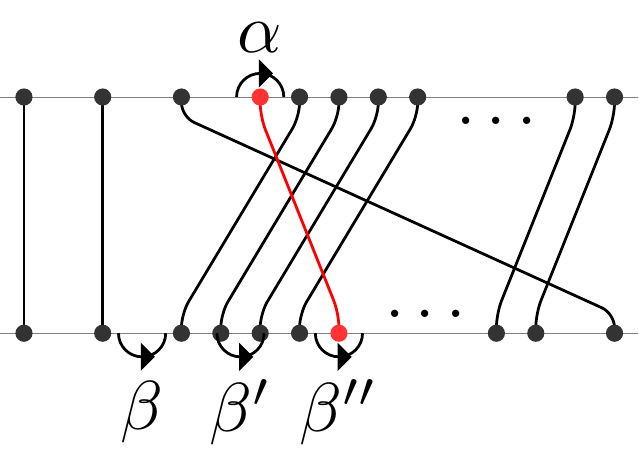}}}
\hspace{36pt}\raisebox{48pt}{$\rotatebox{90}{$\left. \rule{0pt}{26pt}\right\}\rotatebox{-90}{$\!\!\!\!\!\!\!\!\!\!\!\!\!\! 4(k\!-\!1)\!+\!3$} $}$}\\
$(\{2\ell+1\},1,s)$ & $(\{2\ell+1\},1,-s)$ & $(\{2\ell+1\},s)$

\end{tabular}
\end{center}
The first two permutations are the ones with invariant $(\{2\ell+1\},1,\pm s)$ constructed by theorem \ref{thm_I_2X}. The invariant of the third one are justified by applying proposition \ref{pro_opposite_sign} on $\alpha,\beta'$ and $\beta''$ since $\alpha$ is part of a even cycle of length $2k$ and $\beta''$ and $\beta'$' are consecutive arc of another even cycle of length $2k$ by proposition \ref{pro_X_family} (second line third case of figure \ref{fig_perm_x_family}).

\noindent $\bullet$ if $2\ell=4k+2$:
\begin{center}
\begin{tabular}{lll}
\makebox[0pt][l]{\raisebox{0mm}{\includegraphics[scale=.46]{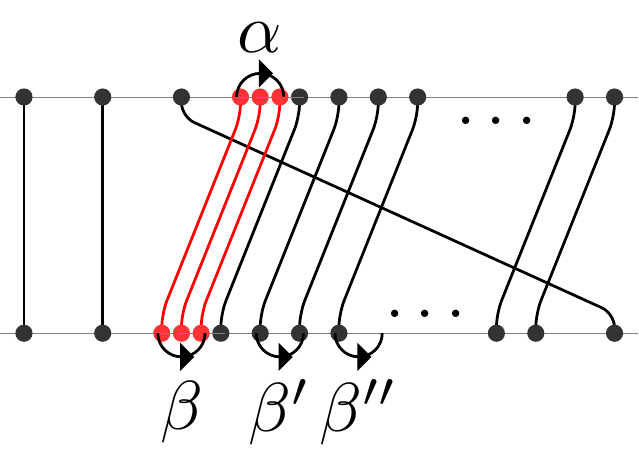}}}
\hspace{36pt}\raisebox{48pt}{$\rotatebox{90}{$\left. \rule{0pt}{26pt}\right\}\rotatebox{-90}{$\!\!\!\!\!\!\!\!\!\!\!\!\!\! 4(k\!-\!1)\!+\!3$} $}$}
&
\makebox[0pt][l]{\raisebox{0mm}{\includegraphics[scale=.46]{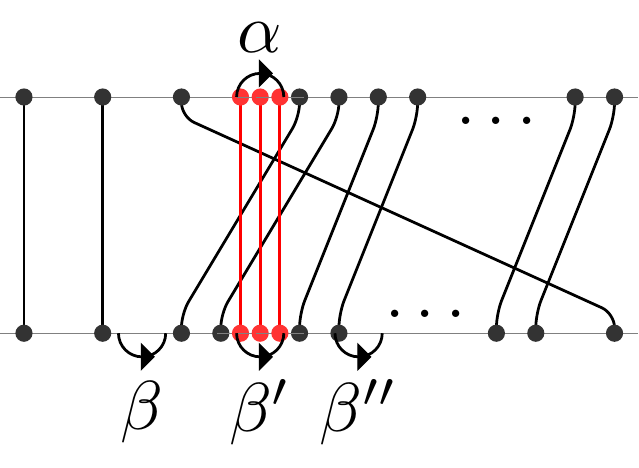}}}
\hspace{36pt}\raisebox{48pt}{$\rotatebox{90}{$\left. \rule{0pt}{26pt}\right\}\rotatebox{-90}{$\!\!\!\!\!\!\!\!\!\!\!\!\!\! 4(k\!-\!1)\!+\!3$} $}$}&
\makebox[0pt][l]{\raisebox{0mm}{\includegraphics[scale=.46]{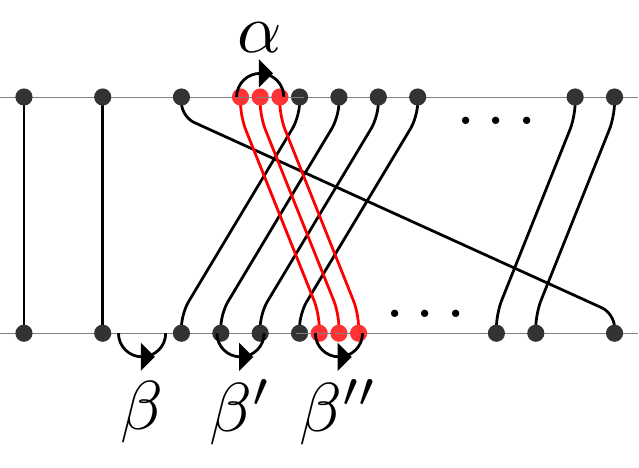}}}
\hspace{36pt}\raisebox{48pt}{$\rotatebox{90}{$\left. \rule{0pt}{26pt}\right\}\rotatebox{-90}{$\!\!\!\!\!\!\!\!\!\!\!\!\!\! 4(k\!-\!1)\!+\!3$} $}$}\\
$(\{2\ell+1\},1,s)$ & $(\{2\ell+1\},1,-s)$ & $(\{2\ell+1\},s)$
\end{tabular}
\end{center}
Once again the first two permutations are the one with invariant $(\lam,r,\pm s)$ constructed by theorem \ref{thm_I_2X} . The invariant of the third one are justified by applying proposition \ref{pro_opposite_sign} on $\alpha,\beta'$ and $\beta''$.\\

No other permutation produced by theorem \ref{thm_I_2X} are exceptional since by appendix C in \cite{DS17} (see also appendix A in this article), there are only two exceptional permutations starting with $\s(1)=1$ and $\s(2)=2$ those are $\id_n$ and $\id'_n$. We just solved the case of $\id'_n$ and $\id_n$ is not $I_2X$.
\end{remark}
~\\

We now consider the case $(\lam,r)$ contains even cycles. For that purpose, proposition \ref{pro_ins_cross_perm} is not sufficient since it does not allow one to add pairs of even cycles of differing lengths. Recall that by theorem \ref{thm_arf_value} (and more generaly the classification theorem) that even cycles must be in even number in a permutation, thus they can only be added in at least pairs.

The following proposition complements proposition \ref{pro_ins_cross_perm} and makes it possible to add pair of even cycles of differing lengths.

\begin{figure}[h!]
\begin{center}
\begin{center}
\begin{tikzcd}
\raisebox{2.2mm}{\includegraphics[scale=.5]{figure/fig_perm_block_add1.pdf}} \arrow[d, "1"]
&   \raisebox{0mm}{\includegraphics[scale=.5]{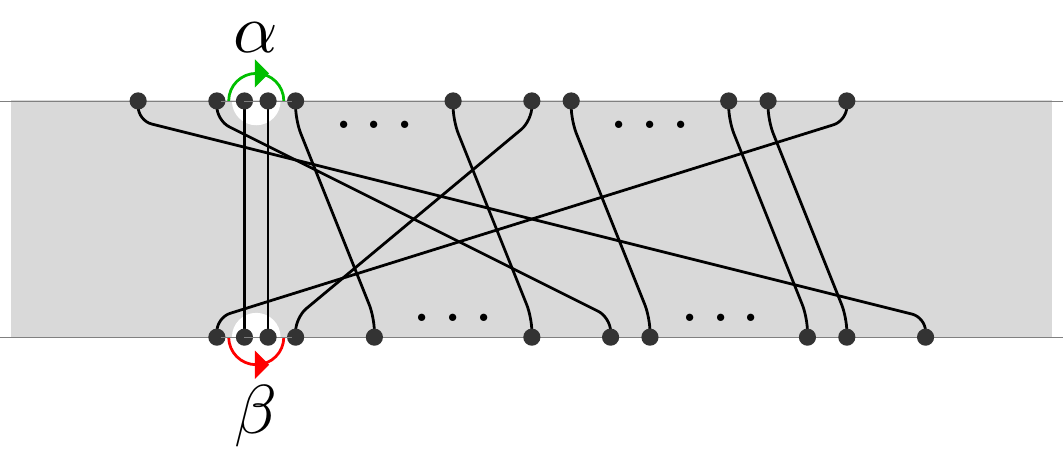}}
 \\
\raisebox{-10mm}{\includegraphics[scale=.5]{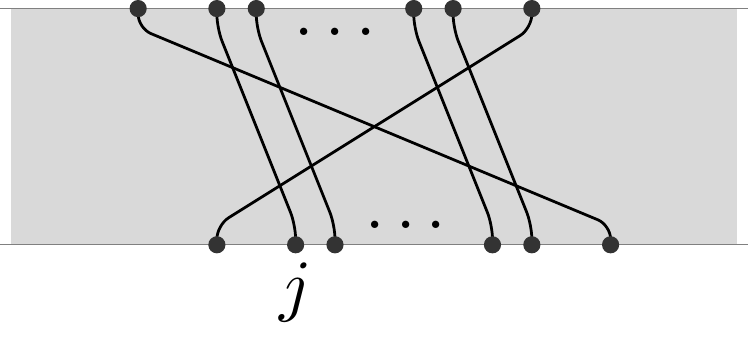}}
\hspace{-81pt}\raisebox{34pt-5.4mm}{$\rotatebox{90}{$\left. \rule{0pt}{22.5pt}\right\}\rotatebox{-90}{$\!\!\!\!2p$} $}$} \hspace{35pt}
\arrow[r, "2"]
&  \raisebox{-11.6mm}{\includegraphics[scale=.5]{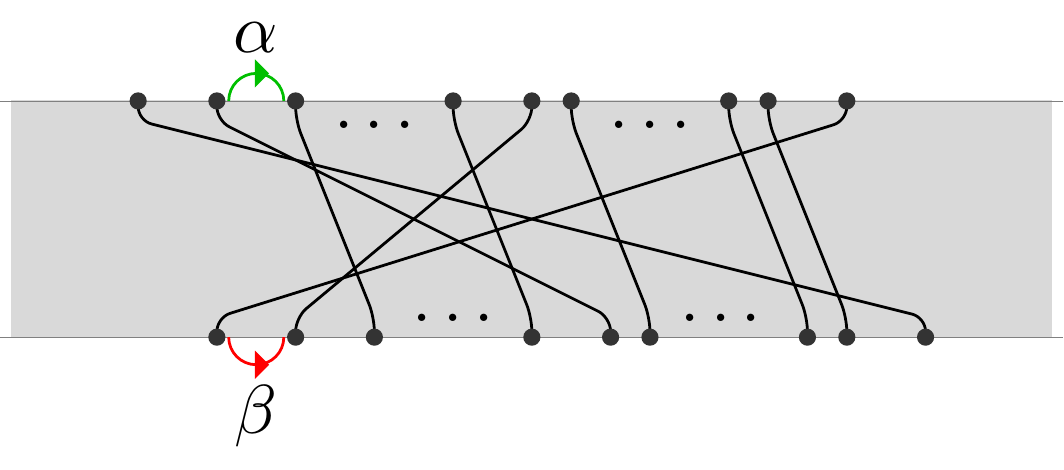}}
\hspace{-116pt}\raisebox{34pt-5.4mm}{$\rotatebox{90}{$\left. \rule{0pt}{19pt}\right\}\rotatebox{-90}{$\!\!\!\!2p'$} $}$}\hspace{10pt}\raisebox{34pt-5.4mm}{$\rotatebox{90}{$\left. \rule{0pt}{19.5pt}\right\}\rotatebox{-90}{$\!\!\!\!\!\!\!\!\!2p\!-\!1$} $}$} \hspace{35pt} \arrow[u, "3"]
\end{tikzcd}
\end{center}
\end{center}
\caption{\label{fig_double_cross_perm}The scheme to add two different even cycles to a permutation. }
\end{figure}

\begin{proposition}[Adding cycles 2]\label{pro_ins_cross_perm_odd}
let $\s$ be a permutation with cycle invariant $(\lambda,r)$ and let $p>p'$, we call $\s_i(C_{2p}\bigcup C_{2p'})$ the permutation obtained by : (see also figure \ref{fig_double_cross_perm})
\begin{enumerate}
\item replacing the $i$th edge of $\s$ by the cross permutation $C_{2p}$
\item replacing the first parallel edge of $C_{2p}$ by the cross permutation $C_{2p'}$ 
\item and finally inserting a double-edge within $\alpha$ the leftmost top arc of $C_{2p'}$ and $\beta$ the leftmost bottom arc of $C_{2p}$. 
\end{enumerate}
The cycle invariant $(\lambda',r)$ of $\s_i(C_{2p}\bigcup C_{2p'})$ verifies $\lambda'=\lambda\bigcup \{ 2(p+1),2(p'+1)\}$ :
\end{proposition}

\begin{pf}
The proof follows from the diagram below.
\begin{center}
\begin{tikzcd}
\s,\ (\lambda,r)
\arrow[d, "1"]
&  \s_i(C_{2p}\bigcup C_{2p'}),\ (\lambda\bigcup\{{\color{red}2p+2},{\color{green!50!black}2p'+2}\},r)
 \\
\s'\!=\!\s_i(C_{2p}),\ (\lambda\bigcup\{2p+1\},r)
\arrow[r, "2"]
&s''\!=\!\s'_j(C_{2p'}),\ (\lambda\bigcup\{{\color{red}2p+1},{\color{green!50!black}2p'+1}\},r) \arrow[u, "3"]
\end{tikzcd}
\end{center}
The modifications of the cycle invariants are justified in step 1 and 2 by the proposition \ref{pro_ins_cross_perm} and in step 3 by the double-edge insertion proposition \ref{pro_double-edge} since $\s_i(C_{2p}\bigcup C_{2p'})=\s''|_{2,\alpha,\beta}$ and $\alpha$ and $\beta$ are two arcs of two different cycles of length $2p+1$ and $2p'+1$ respectively.
\qed
\end{pf}

The second theorem handles the case $\s$ has even cycles. Fortunately since the sign is 0 in this case there are no need two produce two permutations with opposite sign. However there are more base permutations to consider so the proof is not much shorter.

\begin{figure}[tb!]
\begin{center}
\begin{tabular}{lll}
\makebox[0pt][l]{\raisebox{0mm}{\includegraphics[scale=.46]{figure/fig_perm_base.pdf}}}
\hspace{18pt}\raisebox{31pt}{$\rotatebox{90}{$\left. \rule{0pt}{21.5pt}\right\}\rotatebox{-90}{$\!\!\!\!\!\!\!\!\! 2k\!+\!1$} $}$}
\hspace{0pt}\raisebox{31pt}{$\rotatebox{90}{$\left. \rule{0pt}{21.5pt}\right\}\rotatebox{-90}{$\! \! \! \! 2k'$} $}$}&
\makebox[0pt][l]{\raisebox{0mm}{\includegraphics[scale=.46]{figure/fig_perm_base.pdf}}}
\hspace{18pt}\raisebox{31pt}{$\rotatebox{90}{$\left. \rule{0pt}{21.5pt}\right\}\rotatebox{-90}{$\! \! \! \! 2k$} $}$}
\hspace{0pt}\raisebox{31pt}{$\rotatebox{90}{$\left. \rule{0pt}{21.5pt}\right\}\rotatebox{-90}{$\!\!\!\!\!\!\!\!\! 2k'\!+\!1$} $}$}&\\
$\lambda=\{k\!+\!2k'\!+\!1\},r=k\!+\!1$ & $\lambda=\{k'\!+\!1\},r=2k\!+\!k'\!+\!1$\\
&&\\

\makebox[15pt][l]{\raisebox{0mm}{\includegraphics[scale=.46]{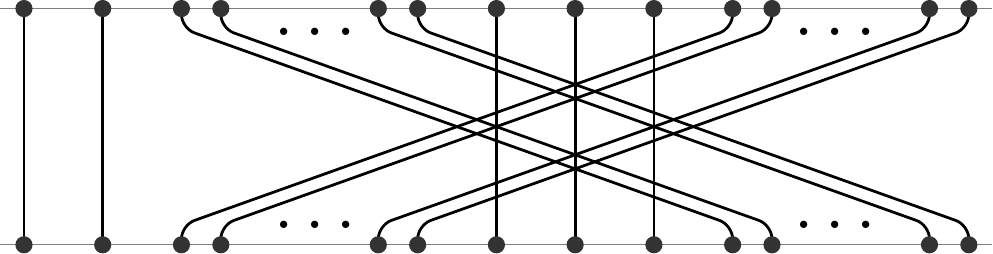}}}
\hspace{2pt}\raisebox{31pt}{$\rotatebox{90}{$\left. \rule{0pt}{21.5pt}\right\}\rotatebox{-90}{$\!\!\!\! 2k$} $}$}
\hspace{34pt}\raisebox{31pt}{$\rotatebox{90}{$\left. \rule{0pt}{21.5pt}\right\}\rotatebox{-90}{$\!\!\!\! 2k$} $}$}
&
\makebox[15pt][l]{\raisebox{0mm}{\includegraphics[scale=.46]{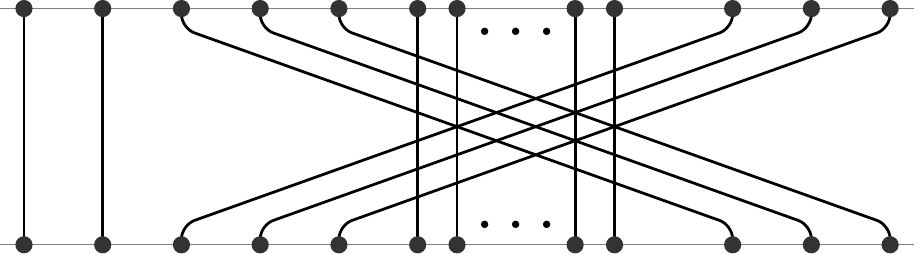}}}
\hspace{34pt}\raisebox{31pt}{$\rotatebox{90}{$\left. \rule{0pt}{18.5pt}\right\}\rotatebox{-90}{$\!\!\!\! 2k $} $}$}

&
\scalebox{.46}{\begin{tikzpicture}[scale=.8]\permutationmatching{1,2,7,8,9,6,3,4,5}\end{tikzpicture}}

\\
$\lambda=\{2k+2\},r=2k\!+\!2$ & $\lambda=\{2,2k\!+3\},r=2$ & $\lambda=\{2,2,2\},r=2$
\end{tabular}
\end{center}
\caption{\label{fig_perm_x_family_odd} Two families of $I_2X$-permutations with their respective cycle invariant.}
\end{figure}

\begin{proposition}[Base permutations 2]\label{pro_X_family_odd}
For every $k,k' \geq 0,$ 

The permutations $X_{2k+1,2k'}$ and, $X_{2k,2k'+1}$ (described in the first line of figure \ref{fig_perm_x_family} have cycle invariant $(\lambda=\{k\!+\!2k'\!+\!1\},r=k\!+\!1)$ and $(\lambda=\{k'\!+\!1\},r=2k\!+\!k'\!+\!1)$ respectively.

The permutations $X_{2k,3,2k}$, $X_{3,2k,3}$ and $X_{2,1,3}$ (described in the second line of figure \ref{fig_perm_x_family}) have cycle invariant $(\lambda=\{2k+2\},r=2k\!+\!2)$, $(\lambda=\{2,2k\!+3\},r=2)$ and $(\lambda=\{2,2,2\},r=2)$ respectively.
\end{proposition}

\begin{pf}
By induction on $k,k'$. In order to highlight the structure of the cycle invariant we start the base cases $k,k'\geq1$. The base cases for the first line are respectively :\begin{center}
 \scalebox{.5}{\begin{tikzpicture}[scale=.8]\permutationmatching{1,2,6,7,3,4,5}\end{tikzpicture}} and \scalebox{.5}{\begin{tikzpicture}[scale=.8]\permutationmatching{1,2,5,6,7,3,4}\end{tikzpicture}}.\end{center}

And they have cycle invariant $(\lambda,r)$ : $(\{{\color{blue}4}\},{\color{red}2})$,  $(\{{\color{blue}2}\},{\color{red}4})$ as shown just below :\begin{center}
\includegraphics{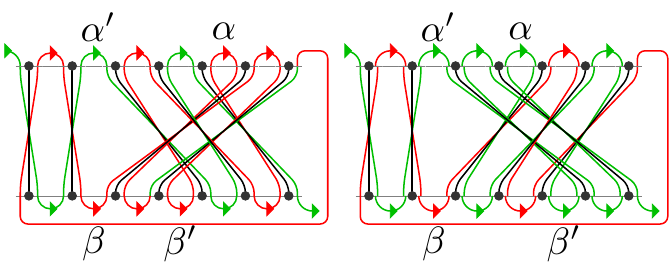}
\end{center}
Then the statement follows by induction from the insertion a double-edge within $\alpha$ and $\beta$ or within $\alpha'$ and $\beta'$ (the resulting change of the cycle invariant are described in proposition \ref{pro_double-edge})

The base cases for the second family are respectively (for the second we make $k$ start at 0 since the strcuture of the cycle invariant is just as explicit here) :\begin{center}
\scalebox{.5}{\begin{tikzpicture}[scale=.8]\permutationmatching{1,2,8,9,5,6,7,3,4}\end{tikzpicture}} \scalebox{.5}{\begin{tikzpicture}[scale=.8]\permutationmatching{1,2,6,7,8,3,4,5}\end{tikzpicture}} and \scalebox{.5}{\begin{tikzpicture}[scale=.8]\permutationmatching{1,2,7,8,9,6,3,4,5}\end{tikzpicture}}.\end{center}

And they have cycle invariant $(\lambda,r)$ : $(\{{\color{blue}4}\},{\color{red}4})$, $(\{{\color{blue}2},{\color{green!50!black}3}\},{\color{red}2})$,  $(\{{\color{blue}2},{\color{green!50!black}2},{\color{black}2}\},{\color{red}2})$ as shown below :
\begin{center}
\includegraphics[scale=1]{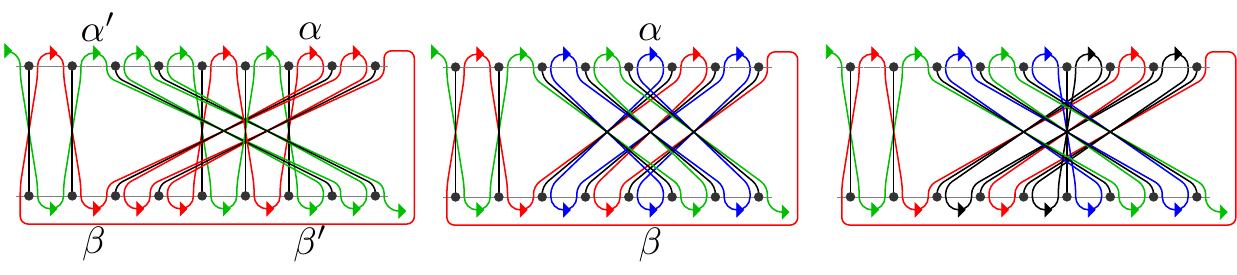}
 \end{center}

Then the statement follows by induction from the insertion of a pair of double-edges within $\alpha$ and $\beta$ and within $\alpha'$ and $\beta'$ for the first permutation. For the second permutation the double-edge must be inserted within $\alpha$ and $\beta$ (the resulting change of the cycle invariant are described in proposition \ref{pro_double-edge}). \qed
\end{pf}

\begin{theorem}\label{thm_I_2X_odd}
Let $(\lambda,r,0)$ be an invariant with even cycle, then there exists a shift-irreducible family with invariant $(\lambda,r,0)$.
\end{theorem}

\begin{pf}
We can always consider that the size of the permutation is at least 10, for smaller size the result can be obtained by automatic search.

Let $(\lam,r,s)$ be an invariant with even cycles. Following the proof sketch of the introduction, we fist construct a base permutation $\s^0$ with invariant $(\lam'\subseteq\lam,r)$ then add cycles to obtain a $I_2X$-permutation $\s_1$ with invariant $(\lam,r,0)$. 

\begin{itemize}
\item If the rank is odd and there is at least one odd cycle $2\ell+1$, then the base permutation is either $X_{1,2\ell}$, $X_{2,1,2\ell}$ ot $X_{2,r-3}$ as in theorem \ref{thm_I_2X}.
\item if the rank is odd and there are no odd cycles :
\begin{itemize}
	 \item If the rank $r=1$, then there are two base cases: 
If $\lam$ has 2 even cycles of the same length $2\ell$, the base permutation $\s^0$ with invariant $(\lambda'=\{2\ell,2\ell\}, 1)$ for any $\ell\geq1$ is $X_{1,4(\ell-1)+3}$ according to lemma \ref{pro_X_family} (second line, first case).\\

If $\lam$ does not have two even cycles of the same length, the base permutation $\s^0$ with invariant
$(\lambda'=\{2\ell,2\ell'+2\ell\}, 1)$ $\ell',\ell\geq1$ is obtained in two steps. First we take, as above, $\s'^0=X_{1,4(\ell-1)+3}$,  it has invariant $(\lambda'=\{2\ell,2\ell\}, 1)$. Then we choose the two arcs $\alpha$ and $\beta$ as below:
\begin{center}
\raisebox{11.5mm}{$X_{1,4(\ell-1)+3}=$ }\raisebox{0mm}{\includegraphics[scale=.46]{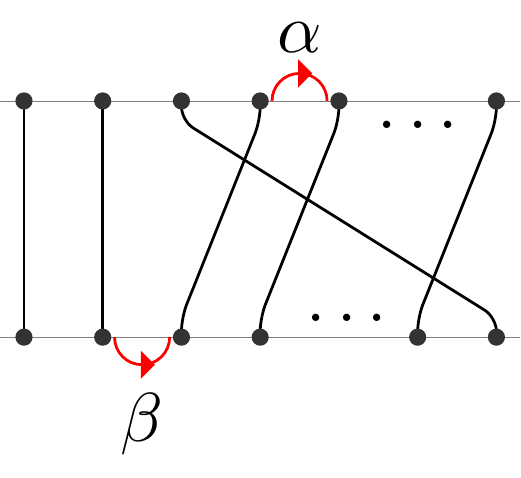}}
\end{center}
$\alpha$ and $\beta$ are in the same cycle of length $2\ell$, thus by the double-edge insertion proposition  \ref{pro_double-edge} $\s'^0|_{2\ell',\alpha,\beta}$ the permutation resulting from the insertion of $\ell$ double-edges within $\alpha$ and $\beta$ has invariant  $(\lambda'=\{2\ell,2\ell'+2\ell\}, 1)$.

	 \item If the rank $r=3$, then there are two base permutations:

If $\lam$ has 2 even cycle of the same length $2\ell$, the base permutation $\s^0$ with invariant $(\lambda'=\{2\ell,2\ell\}, 3)$ for any $\ell\geq1$ is $X_{2,1,4(\ell-1)+3}$ according to lemma \ref{pro_X_family} (second line, second case).

If $\lam$ does not have two even cycles of the same length, the base permutation $\s^0$ with invariant
$(\lambda'=\{2\ell,2\ell'+2\ell\}, 1)$ $\ell',\ell\geq1$ is obtained in two steps. First we take, as above, $\s'^0=X_{2,1,4(\ell-1)+3}$,  it has invariant $(\lambda'=\{2\ell,2\ell\}, 3)$. Then we choose the two arcs $\alpha$ and $\beta$ as below:

\begin{center}
\raisebox{11.5mm}{$X_{2,1,4(\ell-1)+3}=$ }\raisebox{0mm}{\includegraphics[scale=.46]{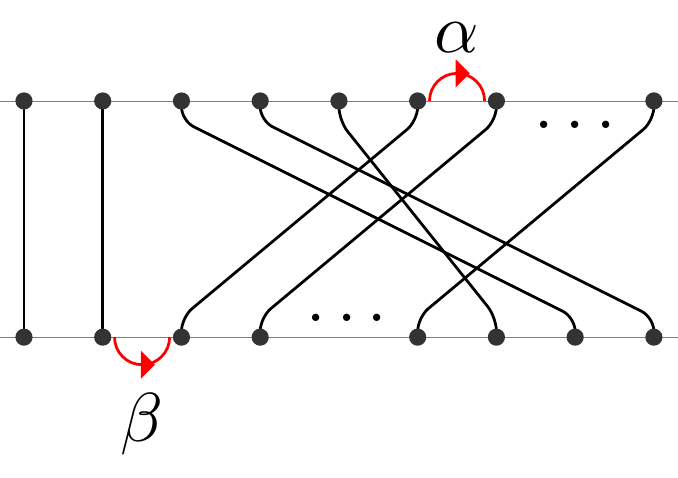}}
\end{center}
$\alpha$ and $\beta$ are in the same cycle of length $2\ell$, thus by the double-edge insertion proposition  \ref{pro_double-edge} $\s'^0|_{2\ell',\alpha,\beta}$ the permutation resulting from the insertion of $\ell$ double-edges within $\alpha$ and $\beta$ has invariant  $(\lambda'=\{2\ell,2\ell'+2\ell\}, 3)$.

\item If the rank $r>5$, then the base permutation $\s^0$ with invariant $(\lambda'=\emptyset, r)$ is exactly $X_{2,r-3}$ according to proposition \ref{pro_X_family} (first line, third case of figure \ref{fig_perm_x_family} since $r$ is odd implies $r-3$ is even  ).

\end{itemize}
\item If the rank $r$ is even and the longest even cycle of $\lam$ has length $2\ell+r$, $\ell\geq 1$ then $\s^0$ is $X_{2r-1,2\ell}$ and has invariant $(\lambda=\{2\ell+r\},r)$ by proposition \ref{pro_X_family_odd} (first line, first case).
\item If the rank $r$ is even and the longest even cycle of $\lam$ has length $2\ell$, $r>2\ell\geq 1$ then $\s^0$ is $X_{r-2\ell,4\ell-1}$ and has invariant $(\lambda=\{2\ell\},r)$ by proposition \ref{pro_X_family_odd} (first line, second case).
\item If the rank $r$ is even and the longest even cycle of $\lam$ has also length $r$ :
\begin{itemize}
	\item  If $r>2$ then $\s^0$ is $X_{r-2,3,r-2}$ and has invariant $(\lambda=\{r\},r)$ by proposition \ref{pro_X_family_odd} (second line, first case).
	\item If $r=2$ and there is a odd cycle of length $2\ell'+1,$ $\ell'\geq1$ in $\lam$ then $\s^0$ is $X_{3,2\ell'-1,3}$ and has invariant $(\lambda=\{2,2\ell'+1\},2)$ by proposition \ref{pro_X_family_odd} (second line, second case).
	\item If $r=2$ and there are no odd cycle (thus every cycle has length two) in $\lam$ then $\s^0$ is $X_{2,1,3}$ and has invariant $(\lambda=\{2,2,2\},2)$ by proposition \ref{pro_X_family_odd} (second line, third case).
\end{itemize}
\end{itemize}
\qed
\end{pf}

\section{The induction\label{sec_induction}}

Let us list a few lemma before beginning the induction.

\begin{lemma}\label{lem_top_bot_label_fixed}
Let $\s$ be a permutation and let $c$ be a $(2k,2r)$-coloring such that either the edge $e_1=(\s^{-1}(1),1)$ or $e_2=(n,\s(n))$ are grayed. Let $\tau$ be the corresponding reduction and let $S$ be a sequence. Define $\tau'=S(\tau)$. 

Then $(\s',c')=B(S)(\s,c)$ has the following property: the gray edge $e_1$ is $(\s'^{-1}(1),1)$ or the gray edge $e_2$ is $(n,\s'(n))$ in $(\s',c')$.
\end{lemma}

\begin{pf}
Let us do the case for $e_1$, the case for $e_2$ is identical

Let $(\Pi_b,\Pi_t)$ be a consistent labelling for $\tau$ then the gray edge $e_1$ of $\s,c$ is inserted within the arcs with label $r^{rk}_0$ and $b\in \Sigma_b$ (since $e_1=(\s'^{-1}(1),1)$. By theorem \ref{thm_consistent_lab} the image of a consistent labelling is a consistent labelling thus $\Pi'_t=S(\Pi_t)$ verifies $\Pi'^{-1}_t(r^{rk}_{0})=1$ by definition. Since the labelling is compatible with the boosted dynamics (cf proposition \ref{thm_keep_track_label}) the top endpoint of the gray edge $e_1$ of $(\s',c')$ is still within $r^{rk}_{0}$ and thus inserted within the top arc with position $1$ in $(\tau',\Pi'_b)$.
\qed
\end{pf}

In other words, the leftmost top endpoint of the edge $(\s^{-1}(1),1)$ and the rightmost bottom endpoint of the edge $(n,\s(n))$ are fixed by the dynamics. It had already been proved many times but this proof is a good illustration of how we will employ the labelling and the boosted dynamics.\\

\begin{lemma}[$d(\s)$ for $\s$ of type $X$]\label{lem_Xtype_s1}
Let $\s$ be a standard permutation with invariant $(\lam,r,s)$ of type $X(r,i)$ and let $\tau=d(\s)$. 
Then $\tau$ has cycle invariant $(\lambda\setminus\{i\},r+i-1)$ and type $H(i,r)$. 
Moreover for any consistent labelling $(\Pi_b,\Pi_t)$ of $\tau$ we have $\s_1=R(\s)=\tau|_{1,t^{rk}_0,b^{rk}_{i-1}}$.
\end{lemma}

\begin{pf}
Let $\s$ and $\tau$ be as in the lemma. Then $\tau$ has indeed cycle invariant  $(\lambda\setminus\{i\},r+i-1)$ and type $H(i,r)$, by proposition \ref{pro_std_d_cycle_inv}.

\begin{figure}[tb]
\begin{center}
$\begin{tikzcd}[%
    ,row sep = 3ex
    ,/tikz/column 1/.append style={anchor=base east}
    ,/tikz/column 2/.append style={anchor=base east}
    ,/tikz/column 3/.append style={anchor=base east}
    ]
 \raisebox{-32pt}{\includegraphics[scale=.5]{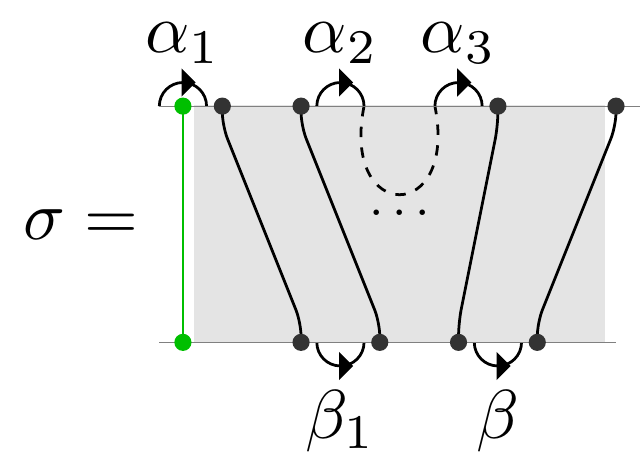}}\arrow{dr}{d} \arrow{r}{R} & 
 \raisebox{-32pt}{\includegraphics[scale=.5]{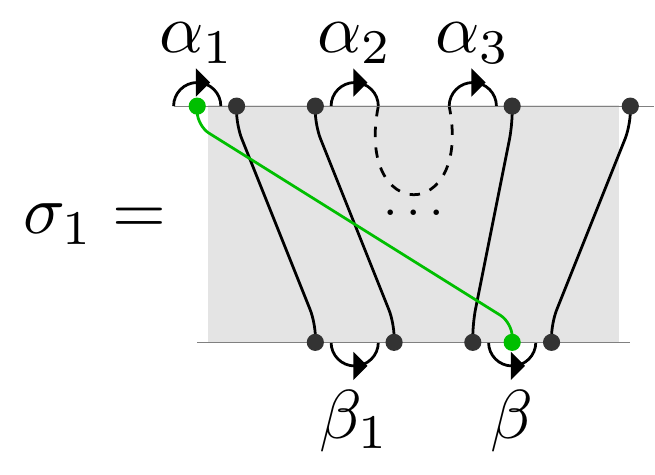}}\arrow{d}{red} \\&
  \raisebox{-42pt}{\includegraphics[scale=.5]{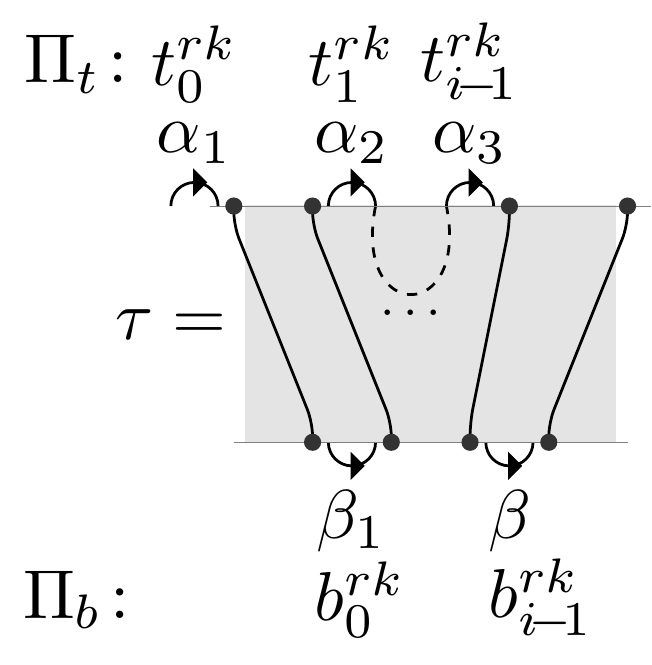}}\\[-.5cm]
\end{tikzcd}$
\end{center}
\caption[caption]{\label{fig_pf_lemma_1_intro} The case $\s$ has type $X(r,i)$ of lemma \ref{lem_Xtype_s1}. We have $\s_1=\tau_{1,\alpha_1=1,\beta}=\tau_{1,t^{rk}_0,b^{rk}_{i\!-\!1}}$ }
\end{figure}

 let $(\Pi_b,\Pi_t)$ be a consistent labelling for $\tau$ and define $\beta=\tau^{-1}(n)-1$ the bottom arc to the left of the edge $(\tau^{-1}(n),n)$. Then clearly $\s_1=R(\s)$ is $\tau_{1,1,\beta}$, moreover $\Pi_b(\beta)= b^{rk}_{i\!-\!1}$ since $\tau$ has type $H(i,r)$ and therefore the top part of the rank (connecting the top left corner to the top right corner) has length $i$. 

Thus we also have $\s_1=\tau_{1,t^{rk}_0,b^{rk}_{i-1}}$ for $\tau,(\Pi_b,\Pi_t).$ See figure \ref{fig_pf_lemma_1_intro} 
\qed
\end{pf}

\begin{lemma}[$d(\s)$ for $\s$ of type $H$]\label{lem_Htype_s1}
Let $\s$ be a standard permutation with invariant $(\lam,r,s)$ of type $H(i,r_2)$ and let $\tau=d(\s)$. 
Then $\tau$ has cycle invariant $(\lambda\bigcup\{i-1\},r_2-1)$ and type $X(r_2-1,i-1)$. 
Moreover there exists a consistent labelling $(\Pi_b,\Pi_t)$ of $\tau$ such that we have $\s_1=R(\s)=\tau|_{1,t^{rk}_0,b_{i-2,i-1,1}}$.
\end{lemma}

\begin{pf}
Let $\s$ and $\tau$ be as in the lemma. Then $\tau$ has indeed cycle invariant  $(\lambda\bigcup\{i-1\},r_2-1)$ and type $H(i-1,r_2-1)$, by proposition \ref{pro_std_d_cycle_inv}.

\begin{figure}[tb]
\begin{center}
$\begin{tikzcd}[%
    ,row sep = 3ex
    ,/tikz/column 1/.append style={anchor=base east}
    ,/tikz/column 2/.append style={anchor=base east}
    ,/tikz/column 3/.append style={anchor=base east}
    ]
 \raisebox{-32pt}{\includegraphics[scale=.5]{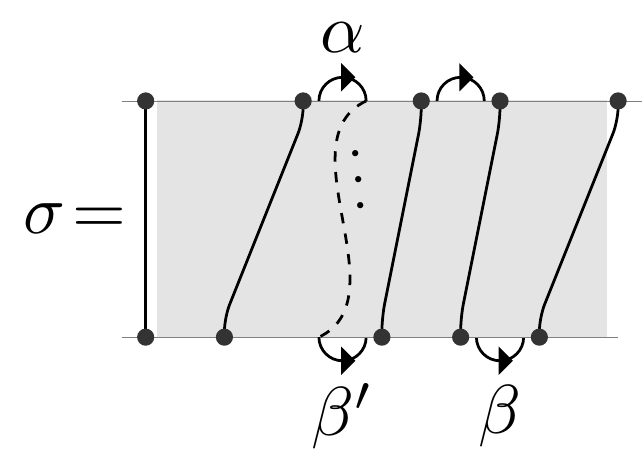}}\arrow{dr}{d} \arrow{r}{R} & 
 \raisebox{-32pt}{\includegraphics[scale=.5]{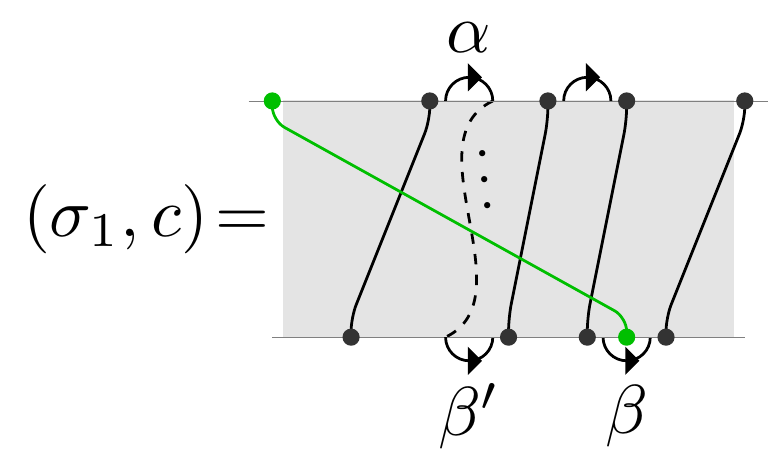}}\arrow{d}{red} \\&
  \raisebox{-42pt}{\includegraphics[scale=.5]{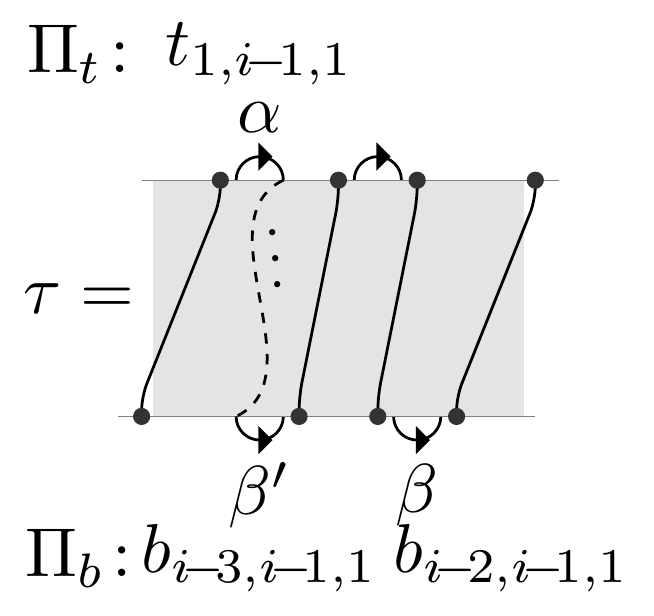}}\\[-.5cm]
\end{tikzcd}$
\end{center}
\caption[caption]{\label{fig_pf_lemma_2_intro} The case $\s$ has type $H(i,r_2)$ of lemma \ref{lem_Htype_s1}. We have $\s_1=\tau_{1,\alpha_1=1,\beta}=\tau_{1,t^{rk}_0,b_{i-1,i\!-\!1,1}}$ }
\end{figure}

Let $(\Pi_b,\Pi_t)$ be a consistent labelling for $\tau$ such that the label of the first top arc of the principal cycle (i.e the arc with position $\alpha=\s(1)+1$) is $t_{1,i-1,1}$. Such a labelling exists since the principal cycle has length $i-1$.

 Then consider $\beta'=\Pi^{-1}(b_{2,3,1})$ and $\beta=\Pi^{-1}(b_{3,3,1})$ the last two consecutives bottom arcs (in this order) of the principal cycle. Thus $\beta=\tau^{-1}(n-1)-1$ the bottom arc to the left of the edge $(\tau^{-1}(n-1),n-1)$ ($\tau$ has size $n-1$). See figure \ref{fig_pf_lemma_2_intro} bottom.

Clearly $\s_1=R(\s)$ is $\tau_{1,1,\beta}=\tau_{1,t^{rk}_0,b_{i-2,i - 1,1}}$. See figure \ref{fig_pf_lemma_2_intro} top: right and middle.
\qed
\end{pf}
%
%\begin{remark}
%Note that if $\s$ has invariant $(\lam,r,s)$ and type $H(r_1,r_2)$ then $r_1+r_2=r+1$.
%\end{remark}

\noindent Let us proceed with the induction. The statements 1 to 7 are true at small size $<10$.

\noindent \textbf{Inductive case:} By induction we suppose that statements 1 to 7 are true at size up to $n-1$, let us prove that they are true at size $n$:

\subsection[Every non exceptional class has a shift-irreducible family]
{Statement 1: Every non exceptional class has a shift-irreducible family}

Let $C$ be a non exceptional class with invariant $(\lam,r,s)$ and $\s \in C$ a standard permutation with $\s(2)=2$. Let $c$ be the $(n-2,2)$-coloring in which the edge $e=(n,\s(n))$ of $\s$ is grayed and $\tau$ the corresponding reduction.
Clearly $\tau$ is irreducible since $\tau(1)=1$ (and $\tau(2)=2$). 

We distinguish two cases : $\tau$ is not in an exceptional class or it is.

\paragraph*{ \textbullet~ Suppose that $\tau$ is not in an exceptional class.} By induction (proposition \ref{pro_i2x} statement 6) there exists a $I_2X$ permutation $\tau'$ in the class of $\tau$. Thus there exists $S$ such that $S(\tau)=\tau'$.

Let $B(S)$ be the boosted sequence of $S$, we have $B(S)(\s,c)=(\s',c')$ with $c'$ having the following property: the gray edge $e$ is $(n,\s'(n))$ by lemma \ref{lem_top_bot_label_fixed}.

Thus the class $C$ contains a permutation $\s'$ such that removing the last edge $(n,\s'(n))$ gives a $I_2X$-permutation $\tau'$. By the characterisation proposition for shift-irreducible permutations (proposition \ref{pro_chara_shift}) a simple case study on the value of $\s'(n)$ shows that the standard family of $\s'$ is shift-irreducible.

\begin{remark}Let us stop here for a moment. We introduced a trick to finding normals forms in section \ref{sec_diff_lab} by defining two sets $T_A$ and $T_B$ and in the proof overview we specified that $T_A$ were the shift-irreducible standard families and $T_B$ the $I_2X$-permutations. We have finally used the trick here: by the assumption that $I_2X$-permutation exists at size $n'<n$ we prove that shift-irreducible standard families exists at size $n$. One can verify that if we had the weaker induction hypothesis that shift-irreducible standard families exists at size $n'<n$ then the case study on the value of $\s'(n)$ would show that the standard family of $\s'$ is not always shift-irreducible.

Of course to complete the trick we still need to prove that every class contains a $I_2X$-permutation at size $n$, but as will be shown in section \ref{sec_proof_i2x}, this is direct consequence of the construction of an $I_2X$-permutation for every valid invariant $(\lam,r,s)$ done in section \ref{sec_IIXshaped} and the classification theorem at size $n$.
\end{remark}

\noindent If $\tau$ is in an exceptional class then it is either in $\Id_n'$ or $\Id_n$. 

\paragraph*{\textbullet~ Suppose $\tau\in\Id'_{n-1}$.} By proposition \ref{pro_excep_std}) since $\tau(1)=1$ and $\tau(2)=2,$ we must have $\tau=\id'^\ell_{n-1}$.  But then, the standard family of $\s$ is shift-irreducible by the characterisation proposition \ref{pro_chara_shift} since $\s$ has the following form :

\[\s= \raisebox{-19pt}{\includegraphics[scale=.5]{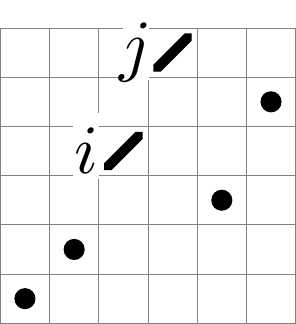}}\quad i,j\geq 0 \]

\paragraph*{\textbullet~ Suppose $\tau\in\Id_{n}$.} By proposition \ref{pro_excep_std})  since $\tau(1)=1$ and $\tau(2)=2,$ we must have $\tau=\id_{n-1}$, thus $\s$ has the form:

\[\s= \raisebox{-15pt}{\includegraphics[scale=.5]{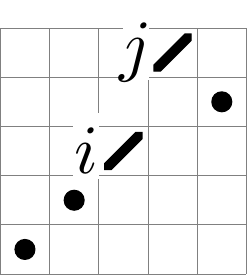}}\quad i>0, j> 0 \]

We must have $i>0$ since otherwise $\s=\id'^\ell_{n-1}$ and $j>0$ since otherwise  $\s=\id_{n}$, in both cases $\s$ would be in an exceptional class which is false by hypothesis. 

Now the permutation $\s'=R^{-1}L^iR^{-j}L^j(\s)$ has the form: 

\[\s'= \raisebox{-15pt}{\includegraphics[scale=.5]{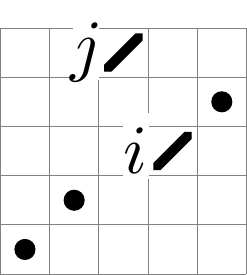}}\quad i,j> 0. \]
and the standard family of $\s'$ is shift-irreducible by the characterisation proposition \ref{pro_chara_shift}.

\subsection[Proof of theorem \ref{thm_arf_value}]
{Statement 2: Proof of theorem \ref{thm_arf_value}}

\paragraph*{ \textasteriskcentered} Let $C$ be a class with invariant $(\lambda,r)$, we show that the list $\lambda\bigcup\{r\}$ has an even number of even parts.

Let $\s$ be a standard permutation with $\s(2)=2$, then $\tau=d(\s)$ is irreducible for it is standard and has invariant $(\lam',r')$. 

By induction hypothesis the list $\lambda'\bigcup\{r'\}$ has an even number of even parts. Moreover, by inspection of table \ref{table.addi} it is easy to see that the parity of the number of even parts is preserved when adding an edge $(i,1)$ to a permutation for every $i$. Thus since $d(\s)=\tau$, $\s$ is obtained by adding the edge $(1,1)$ to $\tau$ and the list $\lambda\bigcup\{r\}$ must therefore contain an even number of even parts

 \paragraph \textasteriskcentered Let $\s\in C$ be a permutation with cycle invariant $(\lam,r)$ we must prove that : \[\Abar(\s)=\begin{cases}\pm2^{\frac{n+\ell}{2}}& \text{ If the number of even parts of the list $\lambda\bigcup\{r\}$ is 0.}\\
0 &\text{ Otherwise.} \end{cases}\]
Where $\ell$ is the number of parts in $\lambda$.\\

We distinguish two cases :
\paragraph*{\textbullet~ Suppose $\lambda=\emptyset$.} We must prove that $\Abar(C)=\pm2^{\frac{n}{2}}$ since there are no even cycles and $\ell=0$. \textbf{Proof idea:} We apply proposition \ref{lem.certifopposign} to two permutations $\s$ and $\s'$ of $C$ and conclude by induction. 
\begin{center}\begin{tabular}{ccc}
\includegraphics[scale=.5]{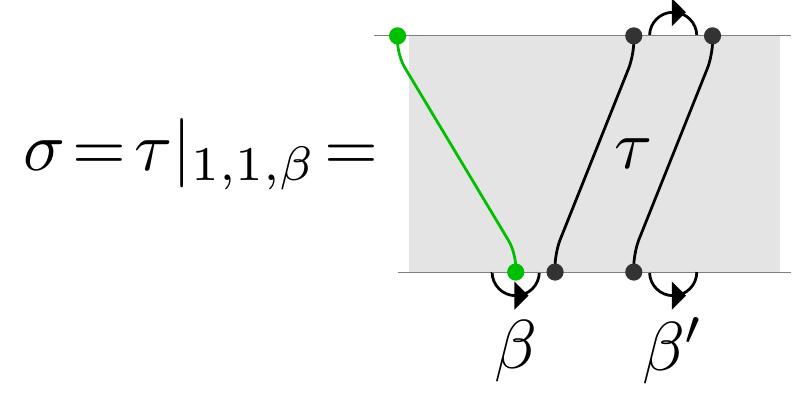} && \includegraphics[scale=.5]{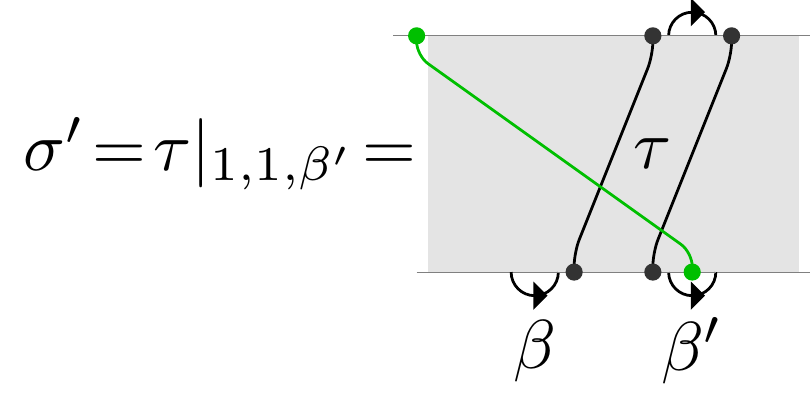} 
\end{tabular}
\end{center}

More formaly, the proof proceeds from the two following lemma :
\begin{lemma*} Let $C$ be as above, There exists $\tau$ with $\Abar(\tau)=\pm2^{\frac{n}{2}}$ and two consecutive bottom arcs $\beta,\beta'$ such that $\tau|_{1,1,\beta}$ and  $\tau|_{1,1,\beta'}$ are in $C$.
\end{lemma*}
and 
\begin{lemma*} Let $\tau$ be a permutation and let $\beta,\beta'$ be two consecutives bottom arcs such that $\tau|_{1,1,\beta}$ and  $\tau|_{1,1,\beta'}$ are in the same class $C$ then $\Abar(C)=\Abar(\tau).$ 
\end{lemma*}

\noindent Clearly the two lemmas put together prove the statement.
\begin{proof}[Proof of the first lemma]
 Let $St$ be a shift-irreducible family of $C$ and let $\s\in St$ the unique permutation with type $H(4,r-4+1)$ (it exists by proposition \ref{pro_std_family}). Let $\tau=d(\s)$, $(\Pi_b,\Pi_t)$ and $\s_1=R(\s)=\tau|_{1,t^{rk}_0,b_{2,3,1}}$ be as in lemma \ref{lem_Htype_s1}. Then $\tau$ has cycle invariant $(\{3\},r-4)$, type $X(r-4,3)$ and is irreducible since the family is shift-irreducible (proposition \ref{pro_shift_standard}). Moreover by induction it has $\Abar(d(\s))=\pm2^{\frac{n-1+1}{2}}$. Let us say wlog $\Abar(d(\s))=2^{\frac{n}{2}}$. 

\begin{figure}[tb]
\begin{center}
$\begin{tikzcd}[%
    ,row sep = 3ex
    ,/tikz/column 1/.append style={anchor=base east}
    ,/tikz/column 2/.append style={anchor=base east}
    ,/tikz/column 3/.append style={anchor=base east}
    ]
 \raisebox{-32pt}{\includegraphics[scale=.5]{figure/fig_perm_s_lemma1.pdf}}\arrow{dr}{d} \arrow{r}{R} & 
 \raisebox{-32pt}{\includegraphics[scale=.5]{figure/fig_perm_s1_lemma1.pdf}}\arrow{d}{red} \arrow{r}{B(S)} & 
 \raisebox{-32pt}{\includegraphics[scale=.5]{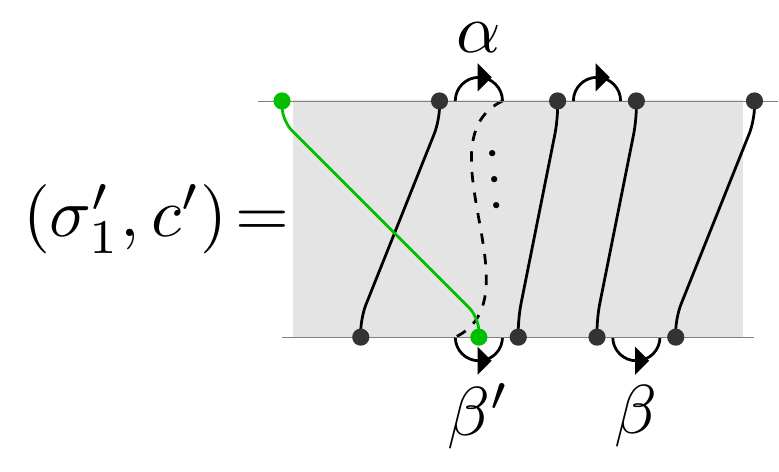}}\arrow{d}{red}   \\&
  \raisebox{-42pt}{\includegraphics[scale=.5]{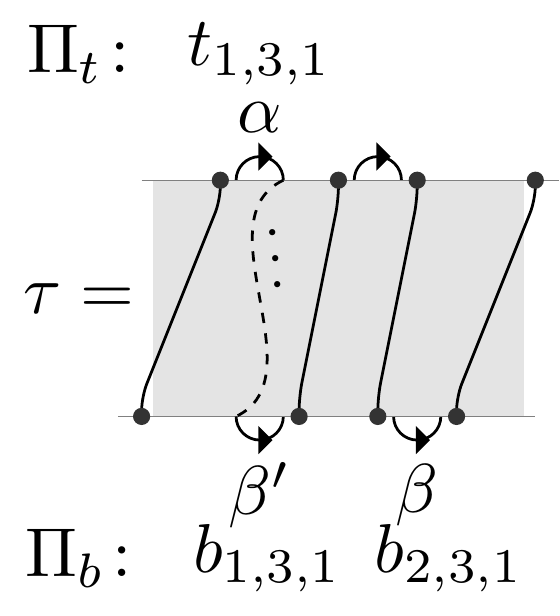}}\arrow{r}{S} &    \raisebox{-42pt}{\includegraphics[scale=.5]{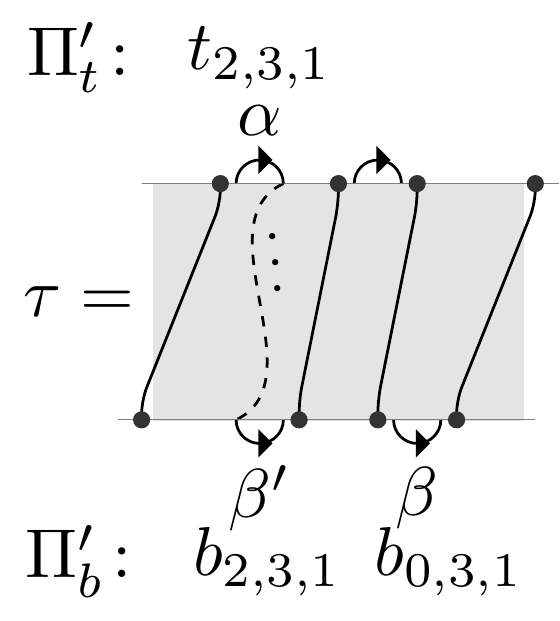}}\\[-.5cm]
\end{tikzcd}$
\end{center}
\caption[caption]{\label{fig_pf_lemma_1_induction} The construction of the proof of the lemma. We have $\s_1=\tau_{1,1,\beta'}$, $\s_1'=\tau_{1,1,\beta}$ and $\s_1=B(S)(\s_1)$ moreover $\Abar(\tau)=\pm 2^{\frac{n}{2}}$.}
\end{figure}

Consider $\beta'=\Pi^{-1}(b_{1,3,1})$ and $\beta=\Pi^{-1}(b_{2,3,1})$ the last two consecutives bottom arcs of the principal cycle. Thus $\beta=\tau^{-1}(n-1)-1$ is the bottom arc to the left of edge $(\tau^{-1}(n-1),n-1)$. See figure \ref{fig_pf_lemma_1_induction} bottom: middle.

By induction we apply the corollary \ref{cor_2monodromy} of the monodromy theorem \ref{thm_2monodromy}(which is true by induction hypothesis). Thus there exists a loop $S$ on $\tau$ such that $\Pi'_b=S(\Pi_b)$ verifies $\Pi'_b(\beta')=b_{2,3,1}$ since the cycle is odd (has length 3). 

Let $c$ be $(n-2,2)$-coloring of $\s_1$ in which the edge $e=(\s_1^{-1}(1),1)$ is grayed (by choice of $c$ the reduction is $\tau$). Thus $e\in (r^{rk}_0,b_{2,3,1})$ in $\tau,(\Pi_b,\Pi_t)$. See figure \ref{fig_pf_lemma_1_induction} top: right and middle.

 By lemma \ref{lem_top_bot_label_fixed} and since the labelling is compatible with the boosted dynamics (theorem \ref{thm_keep_track_label}) $(s_1',c')=B(S)(\s_1,c)$ verifies $\s_1'=\tau_{1,1,\beta'}$. Indeed the gray edge $e \in (r^{rk}_0,b_{2,3,1})$ in $S(\tau,(\Pi_b,\Pi_t))=\tau (\Pi'_b,\Pi'_t)$ i.e. the edge $e$ is inserted within the top arc with position 1 and bottom arc with position $\beta'=\Pi'^{-1}_b(b_{2,3,1})$. See figure \ref{fig_pf_lemma_1_induction} top and bottom right.

The permutations constructed $\s_1',\s_1$ and $\tau$ prove the lemma.
\end{proof}

\begin{proof}[Proof of the second lemma]
Let us apply proposition \ref{lem.certifopposign} to the permutations $\tau|_{1,1,\beta}$, $\tau|_{1,1,\beta'}$ and $\tau$. 
 \[\Abar(\tau|_{1,1,\beta})+ \Abar(\tau|_{1,1,\beta'})=2\Abar(\tau)\]

We know that \[\Abar(\tau|_{1,1,\beta})= \Abar(\tau|_{1,1,\beta'})\]
since they are in the same class.

Thus \[\Abar(C)=\Abar(\tau|_{1,1,\beta})=\Abar(\tau)\]

\end{proof}

\noindent Let us now consider the second case.

\paragraph*{\textbullet~ Suppose $\lambda\neq \emptyset.$} \textbf{Proof idea:} Since $\lambda \neq \emptyset$ there must be a  shift-irreducible standard permutation $\s$ of $C$ of type $X(r,i)$.
Then the result follows from an application of proposition \ref{cor.signsurjq2x} and lemma \ref{lem_Xtype_s1}.\\

Let $\s$ be a standard shift-irreducible permutation of type $X(r,i)$, and let $\tau=d(\s)$ and $\s_1=R(\s)=\tau|_{1,t^{rk}_0,b^{rk}_{i-1}}$ be as in lemma \ref{lem_Xtype_s1}. Then $\tau$ is irreducible since $\s$ is shift-irreducible (proposition \ref{pro_shift_standard}) with invariant $(\lambda\setminus\{i\},r'=r+i-1)$ and type $H(i,r)$.

Let us choose $j=r'-i+1=r$ then $r'-j=i-1$. According to lemma \ref{lem_q_i_sign} we have \be\label{eq.546455d}\Abar(\s_1= \tau_{1,t^{rk}_0,b^{rk}_{r-j}})=\begin{cases}
0, \text{ if } j \equiv 0 \mod 2\\
2\Abar(\tau) \text{ Otherwise.}
\end{cases}\ee

Clearly if $j=r$ is even $\Abar(\s_1)=0.$ However since $r$ is even and the list $\lam\bigcup\{r\}$ contains even parts.

If $j=r$ is odd and $\lambda$ does contain some even cycles then the list $\lambda\setminus \{i\}\bigcup \{r'\}$ also contains even parts and by induction $\Abar(\tau)=0$, thus $\Abar(R(\s))=2\Abar(\tau)=0$.

If $j=r$ is odd and $\lambda$ does not contain even cycles then the list $\lambda\setminus \{i\}\bigcup \{r'\}$ does not contains even parts either, thus by induction $\Abar(\tau)=\pm2^{\frac{n-1+\ell}{2}}$ (where $\ell$ is the number of parts of $\lambda\setminus \{i\}$) thus \[\Abar(\s_1)=2\cdot \pm2^{\frac{n-1+\ell}{2}}=\pm2^{\frac{n+\ell+1}{2}}=\pm2^{\frac{n+\ell'}{2}}\] where $\ell'=\ell+1$ is the number of parts of $\lambda$.

Note that we have proven the following lemma (that will be used in the section just below)

\begin{lemma}\label{lem_short_thing}
Let $\s$ be a shift-irreducible permutation with invariant $(\lam\neq\varnothing,r,s)$ of type $X(r,i)$, let $\tau=d(\s)$ and suppose the list $\lam\bigcup\{r\}$ contains no even cycle. Then $\tau$ is irreducible and has invariant $(\lam\setminus \{i\},r'=r+i-1,s).$
\end{lemma}

\subsection[Existence: for every valid invariant
  $(\lam,r,s)$ there exists a permutation with invariant
  $(\lam,r,s)$.]
{Statement 3, \textbf{Existence:} For every valid invariant $(\lam,r,s)$ there exists a permutation with invariant $(\lam,r,s)$.}

The statement is a direct consequence of theorems \ref{thm_I_2X} and \ref{thm_I_2X_odd} where we constructed a $I_2X$-permutation for every valid $(\lambda,r,s)$. 

Note that we can indeed apply theorem \ref{thm_I_2X} since we have just proved theorem \ref{thm_arf_value} up to $n$ (cf warning from dependency remarks \ref{rk_dependencies_2} and \ref{rk_dependencies_1}).

\subsection[Proposition \ref{pro_labelling_first_step}: First step of the labelling method]
{Statement 4, Proposition \ref{pro_labelling_first_step}: First step of the labelling method}

We already outlined the proof in the proof overview section. 

Let $\s_a$ and $\s_b$ with invariant $(\lam,r,s)$. We shall find $\s'_a$ and $\s'_b$ two standard shift-irreducible permutations of the same type ($X$ or $H$) connected to $\s_a$ and $\s_b$ respectively, such that $\tau_1=d(\s'_a)$ and $\tau_2=d(\s'_b)$ have the same invariant $(\lam',r',s')$. 

Then the proposition \ref{pro_labelling_first_step} is proven by taking $\s_1'=R(\s'_a)$ and $\s'_2=R(\s'_b)$ as in lemma \ref{lem_Xtype_s1} or lemma \ref{lem_Htype_s1}.
Indeed if $\s'_a$ has type $X(r_i)$ by lemma \ref{lem_Xtype_s1} $\s'_1=R(\s'_a)=\tau_1|_{1,t^{rk}_0,b^{rk}_{i-1}}$. Let $c_1$ be the $(n-2,2)$-coloring of $\s'_1$ where the edge $e=(\s^{\prime -1}_1(1),1)$ is grayed. Then the reduction of $(\s'_1,c_1)$ is $\tau_1$. The same reasoning apply to $\s'_b$ which gives a $(\s'_2,c_2)$ with reduction $\tau_2$.

The case $H(i,r_2)$ is similar with the use of lemma \ref{lem_Htype_s1} instead.

\begin{remark}
There is a reason for our normal forms not to be the standard permutation $(\s'_a,\s'_b)$ but the permutations $\s_1'=R(\s'_a)$ and $\s'_1=R(\s'_b)$ instead; despite the fact that the reductions are $\tau_1$ and $\tau_2$ in both cases. If we chose $(\s'_a,\s'_b)$ the gray edge would be the edge $(1,1)$ which is a pivot, thus $(\s'_a,c_a)$ and $(\s'_b,c_b)$ are not proper (cf section \ref{ssec.reddyn} on the boosted dynamics) while $(\s_1',c_1)$ and $(\s_2',c_2)$ are. 
\end{remark}

\paragraph*{Case 1: Suppose $\lam\neq \varnothing$ and $\lam\bigcup\{r\}$ has no even parts.} Let $i\in \lam$, let $\s'_a,\s'_b$ be two standard shift-irreducible permutation of type $X(r,i)$ in the class of $\s_a$ and $\s_b$ respectively (they exist by proposition \ref{pro_std_family}). 

By lemma \ref{lem_short_thing} $\tau_1=d(\s'_a)$ and $\tau_2=d(\s'_b)$ are irreducible and have the same invariant.\qed

\paragraph*{Case 2: Suppose $\lam\neq \varnothing$ and $\lam$ has at least two even cycles (equivalently $\lam\bigcup\{r\}$ has at least four even parts).} Let $i\in \lam$, let $\s'_a,\s'_b$ be two standard shift-irreducible permutations of type $X(r,i)$ in the class of $\s_a$ and $\s_b$ respectively (they exist by proposition \ref{pro_std_family}). 

By lemma \ref{lem_Xtype_s1} $\tau_1$ and $\tau_2$ have cycle invariant $(\lambda\setminus\{i\},r+i-1)$, thus the list $\lambda\setminus\{i\}\bigcup \{r+i-1\}$ must still contains even cycles and by theorem \ref{thm_arf_value} $\tau_1$ and $\tau_2$ have sign $0$. Thus they have the same invariant. Moreover they are irreducible since $\s_a'$ and $\s'_b$ are shift-irreducible (proposition \ref{pro_shift_irreducible}).\qed

\paragraph*{Case 3: Suppose $r\geq 4$.} Let $\s'_a,\s'_b$ be two standard shift-irreducible permutations of type $H(3,r-2)$ in the class of $\s_a$ and $\s_b$ respectively (they exist by proposition \ref{pro_std_family} since $r\geq 4$).

By lemma \ref{lem_Htype_s1} $\tau_1$ and $\tau_2$ have cycle invariant $(\lambda\bigcup\{2\},r-3)$. Thus the list $\lambda\bigcup\{2\}\bigcup \{r-3\}$ contains even cycles and by theorem \ref{thm_arf_value} $\tau_1$ and $\tau_2$ have sign $0$. Therefore they have the same invariant. Moreover they are irreducible since $\s'_a$ and $\s'_b$ are shift-irreducible (proposition \ref{pro_shift_irreducible}).\qed

\paragraph*{Case 4: Suppose $\lam={2k}$ and $r=2$.} We will make use of proposition \ref{cor.signsurjq2x}. 
We first prove the following lemma:
\begin{lemma}\label{lem_small_lemma}
Let $\s$ be a standard permutation with invariant $(\lam,r,0)$ of type $X(r,i)$. Let $e=(\s^{-1}(n),n)$ and $e'=(\s^{-1}(n)-1,\s(\s^{-1}(n)-1))$ and define $j=\s(\s^{-1}(n)-1)$  then $\s_{next}=L^{n-j+1}(\s)$ is of type $X(r,i)$. 

Moreover $\Abar(d(\s_{next}))=-\Abar(d(\s))$.
\end{lemma}

\begin{figure}
\begin{center}
$\begin{tikzcd}[%
    ,row sep = 3ex
    ,/tikz/column 1/.append style={anchor=base east}
    ,/tikz/column 2/.append style={anchor=base east}
    ,/tikz/column 3/.append style={anchor=base east}
    ]
 \raisebox{-20pt}{\includegraphics[scale=.5]{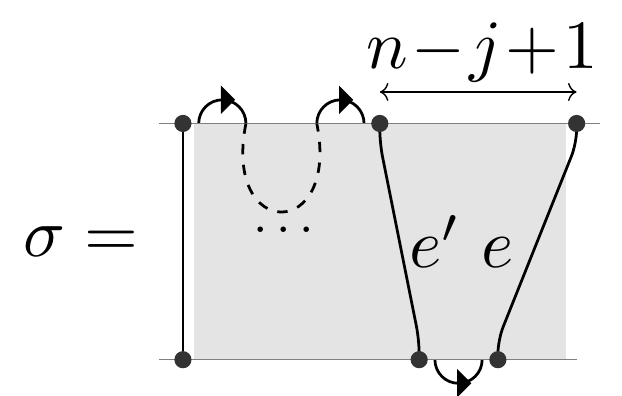}}\arrow{d}{d} \arrow{r}{L^{n-j+1}} & 
 \raisebox{-20pt}{\includegraphics[scale=.5]{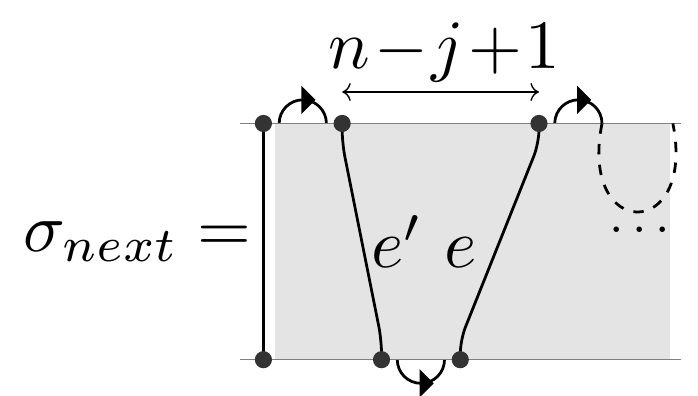}}\arrow{d}{d}   \\
  \raisebox{-42pt}{\includegraphics[scale=.5]{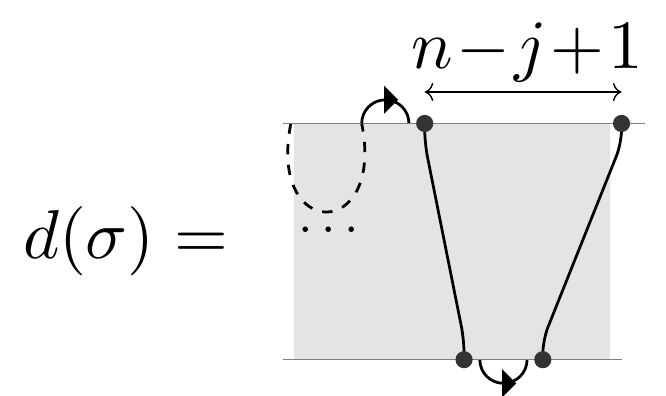}}&    \raisebox{-42pt}{\includegraphics[scale=.5]{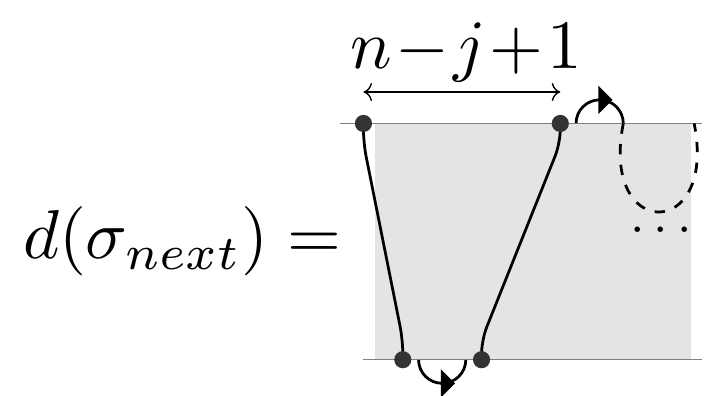}}\\[-.5cm]
\end{tikzcd}$
\end{center}
\caption[caption]{\label{fig_pf_lemma_2_induction} The construction of the proof of the lemma \ref{lem_small_lemma}. Both $\s$ and $\s_{next}$ have type $X(r,i)$ and $\Abar(d(\s))=-\Abar(d(\s_{next}))$.}
\end{figure}

\begin{pf}
Clearly $\s_{next}$ is of type $X(r,i)$ (see figure \ref{fig_pf_lemma_2_induction}).

For the second part of the statement we apply proposition \ref{cor.signsurjq2x} to $d(\s_{next})$, $d(\s)$ and $\s$. we have $\Abar(d(\s_{next}))+\Abar(d(\s))=\Abar(\s)$. Thus   $\Abar(d(\s_{next}))=-\Abar(d(\s))$ since $\Abar(\s)=0.$
\qed
\end{pf}

Let $\s'_a,\s'_b$ be two standard shift-irreducible permutation of type $X(r,2k)$ in the class of $\s_a$ and $\s_b$ respectively (they exist by proposition \ref{pro_std_family}).  Applying lemma \ref{lem_small_lemma} to $\s'_a$ and $\s'_b$ we obtain respectively $\s'_{a,next}$ and $\s'_{b,next}$ of type $X(r,2k)$ such that $d(\s'_a)$ and $d(\s'_{a,next})$ as well as $d(\s'_b)$ and $d(\s'_{b,next})$ have opposite sign. 

Thus choosing $\tau_1=d(\s'_a)$ and $\tau_2=d(\s'_b)$ or $\tau_2=d(\s'_{b,next})$, $\tau_1$ and $\tau_2$ have the same cycle invariant by lemma \ref{lem_Xtype_s1} and the same sign invariant (by choice of $\tau_2$). They are also irreducible since $\s'_a$ and $\s'_b$ or $\s'_{b,next}$ are shift-irreducible (proposition \ref{pro_shift_irreducible}). \qed

\paragraph*{} The four cases overlap somewhat, however they do cover all possibilities. Indeed 'no cycles' is covered by case 3, 'no even parts and some cycles' by case 1, 'at least four even parts' by case 2, 'exactly two even parts' by case 3 (for $r\geq 4$) and case 4 (for $r=2$). Recall that there are always an even number of even parts in $\lam\bigcup\{r\}$ by theorem \ref{thm_arf_value} so every single possible invariant is handled. 

\paragraph*{We must now justify that the $\tau_i$ obtained are not in an exceptional class.}
First note that in case 1,2 and 4 the rank of $\tau_i$ is strictly more than one and in case $3$ $\tau$ has a cycle of length 2 thus $\tau\notin \Id'_{n-1}$ since the cycle invariants are not compatible.
Likewise in case 3, $\tau\notin \Id_{n-1}$ since the cycle invariants are not compatible.

It thus remains to show that $\tau\notin \Id_{n-1}$ in the cases 1,2,4 all with a rank $<4$. 
Consider the following lemma (proven in the appendix A). 
\begin{lemma}
\label{lem.StFamNotmanyId_tobeproveninAppC}
Let $C$ be a non-exceptional class and let $St$ be a standard family. Then at most one $\s \in St$ has $d(\s)\in\Id_{n-1}.$
\end{lemma}
Since the rank is small, there must be either many cycles or a cycle of large length (since $n$ is at least 10 in the induction). If we are unlucky and the permutation $\s'_a$ or $\s'_b$ of type $X(r,i)$ has $d(\s'_a)\in Id_{n-1}$ or $d(\s'_b)\in Id_{n-1}$ we know that it is the only such one in the standard family and we can choose another of type $X(r,i)$ (since there are $im_i$ permutation of type $X(r,i)$ by proposition \ref{pro_std_family}).

%-------------------------------------------------------
\subsection[Completeness: Every pair of permutations $(\s,\s')$ with invariant $(\lam,r,s)$ are connected]
{Statement 5, \textbf{Completeness:} Every pair of permutations $(\s,\s')$ with invariant $(\lam,r,s)$ are connected. }
\label{ssec.rau2sta5}

The completeness statement is demonstrated by the labelling method. (Refer to section \ref{sec_labelling_method})

\begin{itemize}
\item Let $\s_a$ and $\s_b$ be two irreducible permutations with invariant $(\lam,r,s)$. By proposition \ref{pro_labelling_first_step} there exists $(\s'_1,c_1)$ and $(\s'_2,c_2)$ connected to $\s_a$ and $\s_b$ respectively with the following property: 

$c_1$ and $c_2$ are the $(2n-2,2)$-coloring of $\s'_1$ and $\s'_2$ where the edge $e_1=(\s_1^{\prime-1}(1),1)$ and $e_2=(\s_2'^{-1}(1),1)$ are grayed respectively and $\tau_1$ and $\tau_2$ the reductions of $(\s'_1,c_1)$ and $(\s'_2,c_2)$ are irreducible, have the same invariant $(\lam',r',s')$ and are not in exceptional classes.

Moreover by the proof of proposition \ref{pro_labelling_first_step} $\s'_a=R^{-1}(\s'_1)$ and $\s'_b=R^{-1}(\s'_2)$ are standard and have the same type $X(r,i)$ or $H(i,r_2)$. 

If they have type $X(r,i)$ then by lemma \ref{lem_Xtype_s1} there exists a consistent labelling $(\Pi_b,\Pi_t)$ of $\tau_1$ and $(\Pi'_b,\Pi'_t)$ of $\tau_2$ such that $\s'_1=\tau_1|_{1,t^{rk}_0,b^{rk}_{i-1}}$ and $\s'_2=\tau_2|_{1,t^{rk}_0,b^{rk}_{i-1}}$. Define $t=t^{rk}_0$ and $b=b^{rk}_{i-1}$, then $e_1$ is within $(t,b)$ in $(\tau_1,(\Pi_b,\Pi_t))$ and $e_2$ is within $(t,b)$ in $(\tau_2,(\Pi'_b,\Pi'_t))$. 

If they have type $H(i,r_2)$ then by lemma \ref{lem_Htype_s1} there exists a consistent labelling $(\Pi_b,\Pi_t)$ of $\tau_1$ and $(\Pi'_b,\Pi'_t)$ of $\tau_2$ such that $\s'_1=\tau_1|_{1,t^{rk}_0,b_{i-2,i-1,1}}$ and $\s'_2=\tau_2|_{1,t^{rk}_0,b_{i-2,i-1,1}}$. Define $t=t^{rk}_0$ and $b=b_{i-2,i-1,1}$ then $e_1$ is within $(t,b)$ in $(\tau_1,(\Pi_b,\Pi_t))$ and $e_2$ is within $(t,b)$ in $(\tau_2,(\Pi'_b,\Pi'_t))$. Moreover $\lam'=\lam\bigcup{i-1}$ (lemma \ref{lem_Htype_s1} again).

\item Since $\tau_1$ and $\tau_2$ are irreducible and have the same invariant they are in the same class by the classification theorem (which is true by induction hypothesis). Therefore there exists $S$ such that $\tau_2=S(\tau_1)$.

\item Since there exists a boosted dynamics and the labelling is compatible with the boosted dynamics (theorem \ref{thm_keep_track_label}) there exists $B(S)$ such that the reduction of $(\s'_3,c_3)=B(S)(\s'_1,c_1)$ is $\tau_2$ and the gray edge $e_1$ of $(\s'_3,c_3)$ is inserted within the arcs with labels $t$ and $b$ of $\tau_2,(\Pi''_b,\Pi''_t)=S(\tau_1,(\Pi_b,\Pi_t))$. Moreover by theorem \ref{thm_consistent_lab} $(\Pi''_b,\Pi''_t)$ is a consistent labelling.

\begin{remark}
The next point is the last difficulty we emphasized in section \ref{sec_diff_lab}: we need to prove that $P((\tau_2,(\Pi''_b.\Pi''_t),t,b,\alpha,\beta)$ is true. Where $P(\tau_2,(\Pi''_b,\Pi''_t),t,b,\alpha,\beta)$ is true if and only if there exists a loop $S_*$ such that $\tau_2,(\Pi'''_b,\Pi'''_t)=S_*(\tau_2,(\Pi''_b,\Pi''_t))$ verifies  $\Pi^{\prime\prime\prime-1}_b(b)=\beta$ and $\Pi^{\prime\prime\prime-1}_t(t)=\alpha.$

Contrarily to the involution dynamics (see \cite{DS18}) the propriety $P$ can be verified. Essentially the proof is based on the fact that the labellings are consistent: if one endpoint of a gray edge is within an arc of a cycle of length $\lam_i$ in $x,(\Pi_b,\Pi_t)$ then this endpoint will still be within an arc of a cycle of length $\lam_i$ in $S(x,(\Pi_b,\Pi_t))$ for any $S.$ This remark combined with the 2-point monodromy theorem (corollary \ref{cor_2monodromy}) is enough to deduce that the sequence $S_*$ always exists.
\end{remark}

\item \textbf{First case: $b=b^{rk}_{i-1}$.} Both $(\Pi''_b,\Pi''_t)$ and $(\Pi'_b,\Pi'_t)$ are consistent labellings of $\tau_2$ and the labels of the rank are fixed by definition, thus we have $\beta=\Pi''^{-1}_b(b^{rk}_{i-1})=\Pi'^{-1}_b(b^{rk}_{i-1}$) (and $\alpha=1=\Pi''^{-1}_t(t^{rk}_{0})=\Pi'^{-1}_t(t^{rk}_{0}$)). Therefore the gray edge $e_1$ of $(\s'_3,c_3)$ and $e_2$ of $(\s_2',c_2)$ are both inserted within the arcs with position $\alpha=1$ and $\beta$. Thus 
$(\s'_2,c_2)=(\s_3',c_3).$

\textbf{Second case: $b=b_{i-2,i-1,1}$}. Since $(\Pi'_b,\Pi'_t)$ is a consistent labelling of $\tau_2$, the arc with position $\beta=\Pi'^{-1}_b(b_{i-2,i-1,1})$ is an arc of a cycle of length $i-1$ (More precisely it is the last bottom arc of the principal cycle by lemma \ref{lem_Htype_s1}). Since $(\Pi''_b,\Pi''_t)$ is also a consistent labelling of $\tau_2$, $\Pi''_b(\beta)=b_{c,i-1,d}$ for some $0\leq c< i-1$ and $1\leq d\leq m_{i-1}$ where $m_{i-1}$ is the multiplicity of $i-1$ in $\lam'$. 

\textbf{If $\lam'$ has no even cycle or at least two even cycles,} then by corollary \ref{cor_2monodromy} of theorem \ref{thm_2monodromy} (which is true by induction hypothesis) there exists a loop $S_*$ of $\tau_2$ such that $\Pi'''_b=S_*(\Pi''_b)$ verifies $\Pi'''_b(\beta)=b_{i-2,i-1,1}$. Let $B(S_*)$ be the boosted sequence of $S_*$ then the gray edge $e_1$ of $(s'_4,c_4)=B(S_*)(s'_3,c_3)$ is inserted within the arcs with position $\alpha=1$ and $\beta$ therefore $(s'_4,c_4)=(s_2',c_2).$

\textbf{If $\lam'$ has exactly one even cycle: the cycle of length $i-1$.} In this case $d=1$ and $\Pi''_b(\beta)=b_{c,i-1,1}$ then by corollary \ref{cor_2monodromy} of theorem \ref{thm_2monodromy} (which is true by induction hypothesis) there exists a loop $S_*$ of $\tau_2$ such that $\Pi'''_b=S_*(\Pi''_b)$ verifies $\Pi'''_b(\beta)=b_{i-2,i-1,1}$ if $i-2-c \equiv 0 \mod (2) $  and $\Pi'''_b(\beta)=b_{0,i-1,1}$ otherwise.
In the first case we conclude as in the previous paragraph. 

Let us show that $\Pi'''_b(\beta)=b_{0,i-1,1}$ cannot happen. Refer also to figure \ref{fig_pf_stat5}.

We constructed $\s'_2=\tau_2|_{1,t^{rk}_0,b_{i-2,i-1,1}}$ and $(\Pi'_b,\Pi'_t)$ by lemma \ref{lem_Htype_s1} therefore, as said above, $\beta$ is the last arc of the principal cycle of $\tau_2$ (to the left of the edge $(\tau_2^{-1}(n-1),n-1)$) and $\lam=\lam'\setminus\{i-1\}$, thus the list $\lam\bigcup\{r\}$ has no even parts and $s\neq0$.

If $\Pi'''_b(\beta)=b_{0,i-1,1}$, let $\beta'=\Pi'''^{-1}_b(b_{i-2,i-1,1})$ then $\beta'$ and $\beta$ are consecutive (in this order) since $\Pi'''_b$ is a consistent labelling and $i-2$ and 0 are consecutive  (modulo $i-1$) indices. Moreover $(\s'_4,c_4)=B(S_*)(\s'_3,c_3)$ verifies $\s'_4=\tau_2|_{1,1,\beta'}$ since the gray edge $e_1$ of $(\s'_4,c_4)$ is inserted within the arcs with labels $t^{rk_0}$ and $b_{i-2,i-1,1}$ of $\tau_2,(\Pi'''_b,\Pi'''_t)$. 

\begin{figure}[tb]
\begin{center}
$\begin{tikzcd}[%
    ,row sep = 3ex
    ,/tikz/column 2/.append style={anchor=base east}
    ,/tikz/column 3/.append style={anchor=base east}
    ]
 \raisebox{-32pt}{\includegraphics[scale=.5]{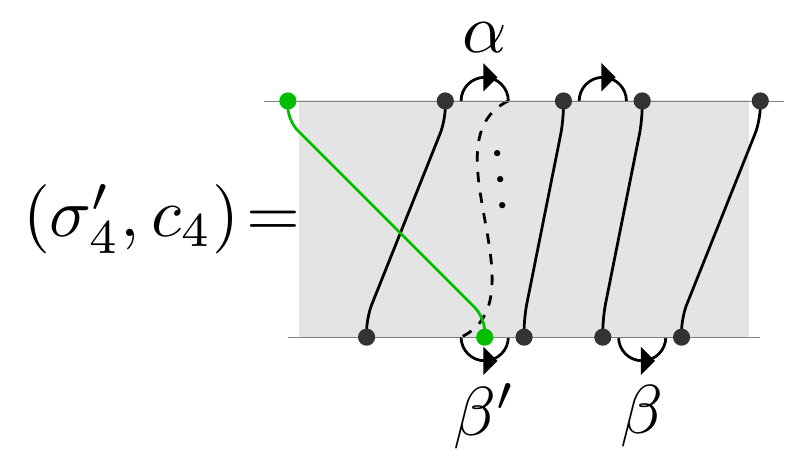}}\arrow{d}{red} & 
 \raisebox{-32pt}{\includegraphics[scale=.5]{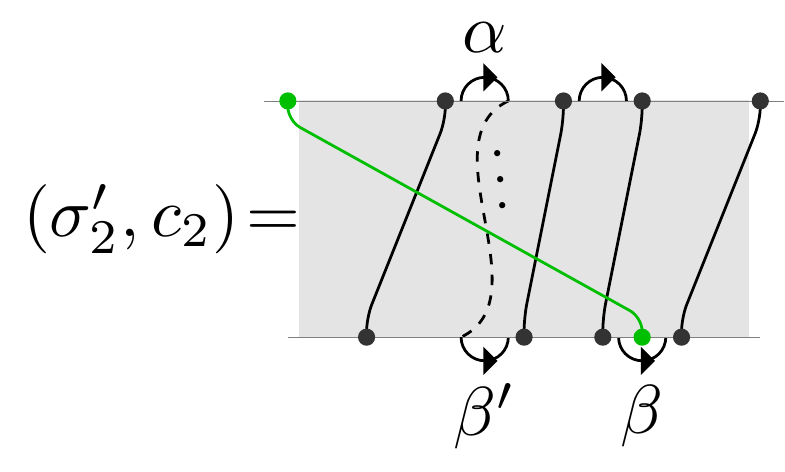}}\arrow{d}{red} \\   \raisebox{-42pt}{\includegraphics[scale=.5]{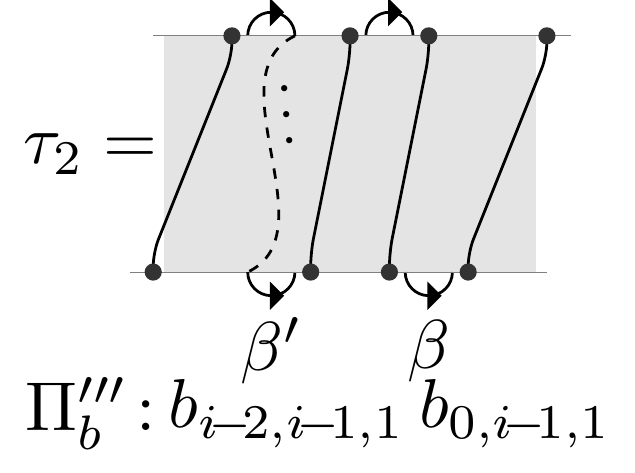}}&
  \raisebox{-42pt}{\includegraphics[scale=.5]{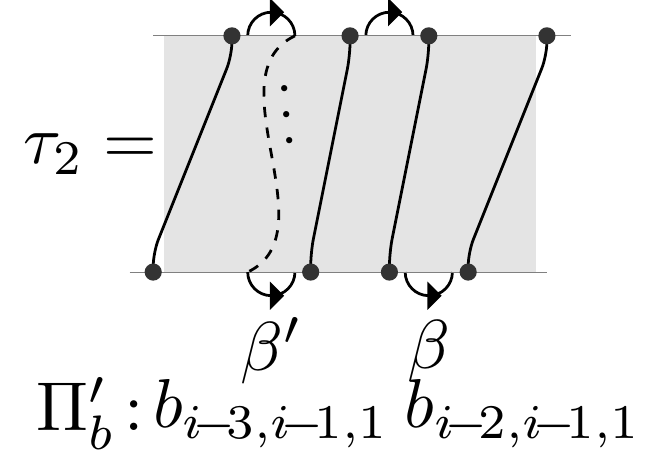}}\\[-.5cm]
\end{tikzcd}$
\end{center}
\caption[caption]{\label{fig_pf_stat5}$(\s'_4,\s'_2)$ are both permutations with invariant $(\lam,r,s\neq 0)$ and $\tau_2$ has sign invariant 0. However proposition \ref{lem.certifopposign} applied to  $\s'_4,\s'_2$ and $\tau_2$ implies that $\s'_4$ should have invariant $-s$, This is contradictory and thus the case $\Pi'''_b(\beta)=b_{0,i-1,1}$ cannot happen.}
\end{figure}

However, we have $\Abar(\s'_4)= \Abar(\s'_2)=\pm 2^{\frac{n+\ell}{2}}$ since both have invariant $(\lam,r,s)$ and $\Abar(\s'_2)=\pm 2^{\frac{n+\ell}{2}}$ by theorem \ref{thm_arf_value}. Moreover $\Abar(\tau_2)=0$ since $\lam'$ has an even cycle. Let us apply proposition \ref{lem.certifopposign} to $\s'_4,\s'_2$ and $\tau'$.
\[\Abar(\s'_4)+ \Abar(\s'_1)=2\cdot \pm 2^{\frac{n+\ell}{2}}=0\] this is contradictory thus $\Pi'''_b(\beta)=b_{0,i-1,1}$ cannot happen.

\end{itemize}

The classification theorem (theorem \ref{thm.Main_theorem}) is a consequence of statements 3 and 5.
 
\subsection[Proposition \ref{pro_i2x}: 
Every non exceptional class contains a $I_2X$-permutation]
{Statement 6, Proposition \ref{pro_i2x}: 
Every non exceptional class contains a $I_2X$-permutation.\label{sec_proof_i2x}}

The statement is a direct consequence of theorems \ref{thm_I_2X} and \ref{thm_I_2X} where we constructed a $I_2X$ permutation for every valid $(\lambda,r,s)$ and of the classification theorem \ref{thm.Main_theorem} which says that there is only one non-exceptional class per invariant.

\subsection[The 2-point monodromy theorem \ref{thm_2monodromy}]
{Statement 7: the 2-point monodromy theorem \ref{thm_2monodromy}}
\label{ssec_induction_monodromy}

Before tackling the monodromy theorem, let us consider the following problem:  we have two standard permutations $\s$ and $\s'$ of type $X(r,i)$ and a labelling $(\Pi_b,\Pi_t)$ of $\s$. The labels of the principal cycle of $\s$ are $t_{c,i,j}, \ldots,t_{c+i-1\mod i,i,j}$ in that order. 

We wish to find a sequence $S'$ such that $S'((\Pi_b,\Pi_t)$) is a labelling of $\s'$ and the labels of the principal cycle of $\s'$ are also $t_{c,i,j}, \ldots,t_{c+i-1\mod i,i,j}$ in that order. 

The following proposition tells us that it is possible if $\tau=d(\s)$ and $\tau'=d(\s')$ are in the same class.

\begin{proposition}\label{pro_label_moves_well}
Let $\s$ and $\s'$ be two standard permutations with invariant $(\lam,r,s)$ and type $X(r,i)$.
Let $(\Pi_b,\Pi_t)$ be a consistent labelling of $\s$, and let $(\alpha_1,\ldots,\alpha_i)$ and $(\alpha'_1,\ldots,\alpha'_i)$ be the arcs  (in that order) of the principal cycle of $\s$ and $\s'$ respectively. Then \[\Pi_t(\alpha_1)=t_{c,i,j}, \ldots,\Pi_t(\alpha_i)=t_{c+i-1\mod i,i,j} \text{ for some $c$ and $j$ }.\]

 Let $\tau=d(\s)$ and $\tau'=d(\s')$ and suppose there exists $S$ such that $\tau'=S(\tau)$.

Finally let $\s',(\Pi'_b,\Pi'_t)=R^{-1}B(S)R(\s,(\Pi_b,\Pi_t))$. We have \[\Pi'_t(\alpha'_1)=t_{c,i,j}, \ldots,\Pi'_t(\alpha'_i)=t_{c+i-1\mod i,i,j}.\]
\end{proposition}

\begin{pf}
Let $\tau,(\Pi^{\ell}_b,\Pi^{\ell}_t)$, $\s_1=R(\s)=\tau|_{1,t^{rk}_0,b^{rk}_{i-1}}$ and $\tau',(\Pi^{\ell_1}_b,\Pi^{\ell_1}_t)$, $\s_1=R(\s')=\tau'|_{1,t^{rk}_0,b^{rk}_{i-1}}$ be as in lemma \ref{lem_Xtype_s1}. 

Let $c$ and $c'$ be the $(n-2,2)$-coloring of $\s_1$ and $\s'_1$ such that the corresponding reductions are $\tau$ and $\tau'$ (i.e $e\in(t^{rk}_0,b^{rk}_{i-1})$ in $\tau,(\Pi^{\ell}_b,\Pi^{\ell}_t)$ likewise for $e').$

Let $B(S)$ be the boosted sequence of $S$ then $B(S)(\s_1,c)=(\s'_1,c')$ since both $e$ and $e'$ are inserted within $t^{rk}_0,b^{rk}_{i-1}$ in $\tau',(\Pi'^{\ell}_b,\Pi'^{\ell}_t)=S(\tau,(\Pi^{\ell}_b,\Pi^{\ell}_t))$ and $\tau',(\Pi^{\ell_1}_b,\Pi'^{\ell_1}_t)$ respectively, and the labels of the rank are fixed. (This part of the proof follows that of statement 5). 

%
%Let $\s_1,(\Pi'_b,\Pi'_t)=R(\s,(\Pi_b,\Pi_t)). In $\s_1$, $\Pi'_t(\alpha_1)=t_{c,i,j}, \ldots,\Pi'_t(\alpha_i)=t_{c+i-1\mod i,i,j}$ however in $\tau$, $\Pi^{\ell}_t(\alpha_1)=t^{rk}_{0}, \ldots,\Pi^{\ell}_t(\alpha_i)=t^{rk}_{i-1}$. Thus the labels of the cycle of $\s_1$ are attached to the labels of the rank of $\tau$. 

For the following, refer to figure \ref{fig_pf_monodromy_1}

\begin{figure}[tb]
\begin{center}
$\begin{tikzcd}[%
    ,row sep = 3.5ex
    ,column sep = 13ex
    ,/tikz/column 1/.append style={anchor=base east}
    ,/tikz/column 2/.append style={anchor=base east}
    ,/tikz/column 3/.append style={anchor=base east}
    ]
 \raisebox{-32pt}{\includegraphics[scale=.5]{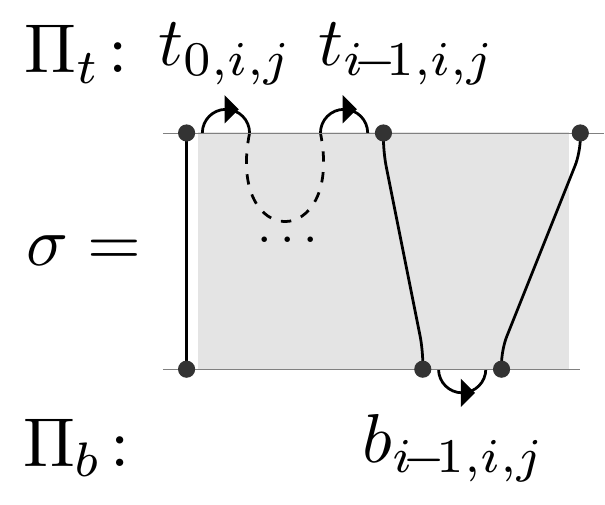}}\arrow{d}{R} \arrow{r}{R^{-1}B(S)R} & \raisebox{-32pt}{\includegraphics[scale=.5]{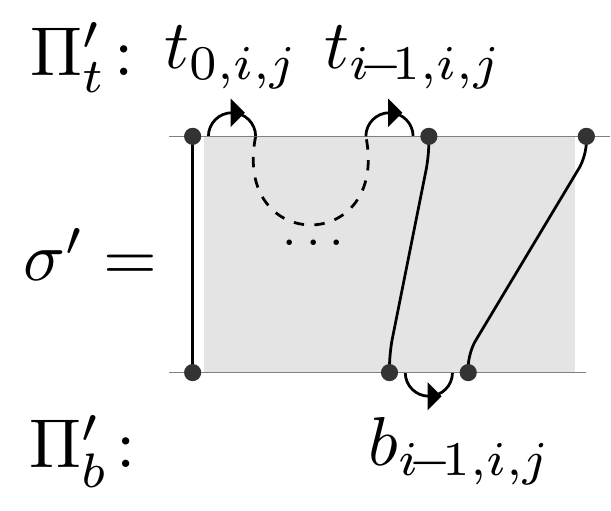}}\\
 \raisebox{-32pt}{\includegraphics[scale=.5]{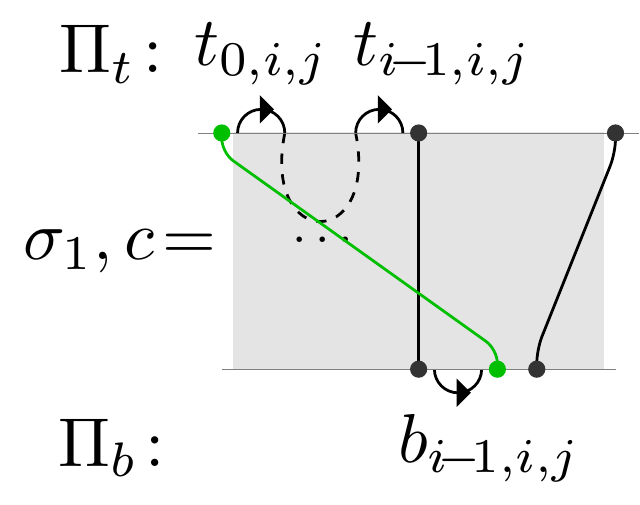}}\arrow{d}{red} \arrow{r}{B(S)} & \raisebox{-32pt}{\includegraphics[scale=.5]{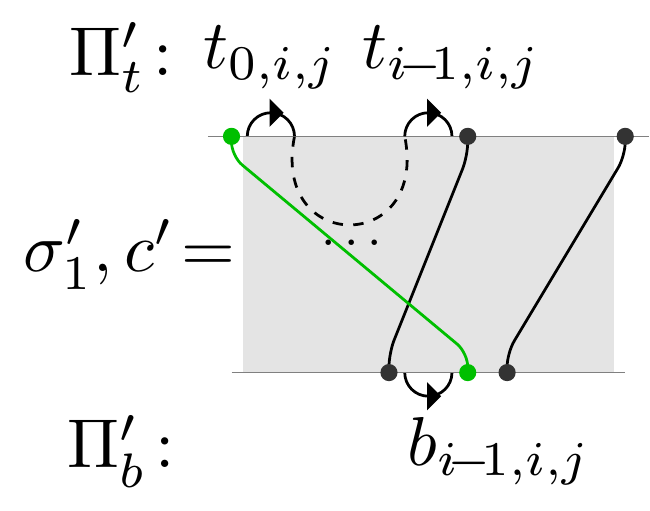}}\arrow{d}{red} \arrow{u}{R^{-1}} \\
  \raisebox{-42pt}{\includegraphics[scale=.5]{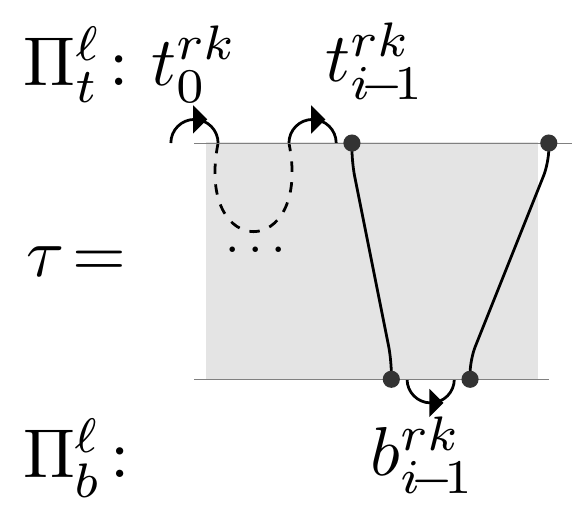}}\arrow{r}{S}&
  \raisebox{-42pt}{\includegraphics[scale=.5]{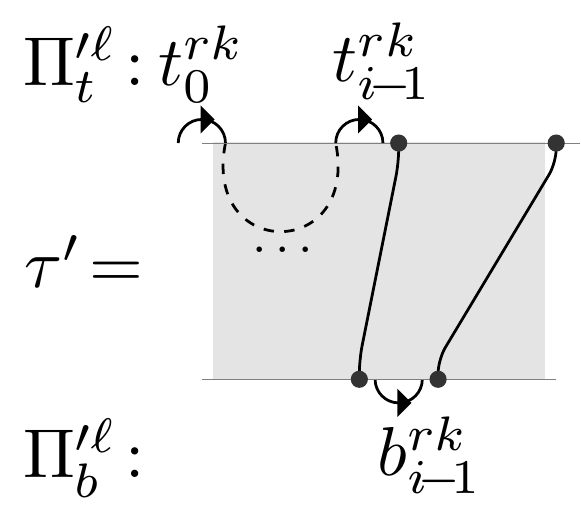}}\\[-.5cm]
\end{tikzcd}$
\end{center}
\caption[caption]{\label{fig_pf_monodromy_1} $\s$ and $\s'$ are standard of type $X(r,i)$. The labels of the principal cycle of $\s$ are send to the arcs of the principal cycle of $\s'$ by the sequence $R^{-1}B(S)R$. Indeed they are attached to the $i$th first labels of the rank of $\tau$ and the labels of the rank are fixed. Thus they are also attached to $i$th first the labels of the rank of $\tau'$ which are the arcs of the principal cycle of $\s'$.\\ In the figure, we choose $\Pi_t(\alpha)=t_{0,i,j}$ instead of $t_{c,i,j}$ for some $c$ for space-saving purpose.}
\end{figure}

In $\s$, $\Pi_t(\alpha_1)=t_{c,i,j}, \ldots,\Pi_t(\alpha_i)=t_{c+i-1\mod i,i,j}$ however in $\tau$, $\Pi^{\ell}_t(\alpha_1)=t^{rk}_{0}, \ldots,$ $\Pi^{\ell}_t(\alpha_i)=t^{rk}_{i-1}$. Thus the labels of the principal cycle of $\s$ are attached to the $i$th first labels of the rank of $\tau$ and they will move along with them in the boosted sequence. 

In $\tau',(\Pi'^{\ell}_b,\Pi'^{\ell}_t)$ we have $\Pi'^{\ell}_t(\alpha'_1)=t^{rk}_{0}, \ldots,\Pi'^{\ell}_t(\alpha'_i)=t^{rk}_{i-1}$ since the labels of the rank are fixed, thus for $\s',(\Pi'_b,\Pi'_t)=R^{-1}B(S)R(\,s(\Pi_b,\Pi_t))$, we must also have $\Pi'_t(\alpha'_1)=t_{c,i,j}, \ldots,\Pi'_t(\alpha'_i)=t_{c+i-1\mod i,i,j}$.

\qed
\end{pf}

We now prove the monodromy theorem \ref{thm_2monodromy}. Let $C$ be a class with invariant $(\lam,r,s)$, if we establish the theorem for a given permutation with a consistent labelling $\s,(\Pi_b,\Pi_t)$ then by conjugation the theorem is valid for any $\s'\in C$ and any consistent labelling of $\s'$. We had already noted that in the labelling method section \ref{sec_labelling_method} in the discussion below definition \ref{def_monomprob} and in the 'difficulties' section (sec. \ref{sec_diff_lab}) in the part concerning the $r$-point monodormy.\\

The proof of every case will follow the same pattern: we start with a standard permutation $\s$ of type $X(r,i)$ with a consistent labelling $(\Pi_b,\Pi_t)$, the top labels of the principal cycle are $t_1,\ldots,t_i\in \Sigma_t$ in that order. The objective is to find a loop $S$ such that the top labels of the principal cycle becomes $t'_1,\ldots,t'_i$ in that order. 

For that purpose, we construct a standard permutation $\s'$ of type $X(r,i)$ with the following properties:
\begin{itemize}
\item $\s',(\Pi'_b,\Pi'_b)=S_1(\s,(\Pi_b,\Pi_t))$ for some $S_1$ and the labels of the principal cycle of $\s'$ are $t'_1,\ldots,t'_i$ in that order.
\item $\tau'=d(\s)$ and $\tau=d(\s)$ are in the same class, thus there is a sequence $S_2$ such that $\tau=S_2(\tau')$.
\end{itemize}
Then we conclude by applying proposition \ref{pro_label_moves_well} to $\s'$ and $\s$, and $S=R^{-1}B(S_2)RS_1$. 

In our application, $\s'$ is always another permutation of the standard family of $\s$ so $S_1=L^k$ for some $k$.

\paragraph*{\textbullet~  Let us demonstrate by the statement on the cycle 1-shift. In this case $\lam$ has no even cycles or at least 2 even cycles} Let $i\in \lam$, let $\s$ a standard shift-irreducible permutation of type $X(r,i)$ and let $(\Pi_b,\Pi_t)$ be a choosen consistent labelling of $\s$. Let $\alpha_1,\ldots,\alpha_i$ be the top arcs of the principal cycle, we have $\Pi_t(\alpha_1)=t_{0,i,j}$ and $\Pi_t(\alpha_2)=t_{1,i,j}, \ldots, \Pi_t(\alpha_i)=t_{i-1,i,j}$. 

We construct a loop $S$ of $\s$ such that $(\Pi''_b,\Pi''_t)=S(\Pi_b,\Pi_t)$ verifies $\Pi''_t(\alpha_1)=t_{1,i,j}$ and $\Pi''_t(\alpha_2)=t_{2,i,j} \ldots \Pi''_t(\alpha_i)=t_{0,i,j}$.

\begin{figure}[tb]
\begin{center}
$\begin{tikzcd}[%
    ,row sep = 3.5ex
    ,column sep = 10ex
    ,/tikz/column 1/.append style={anchor=base east}
    ,/tikz/column 2/.append style={anchor=base east}
    ,/tikz/column 3/.append style={anchor=base east}
    ]
 \raisebox{-32pt}{\includegraphics[scale=.5]{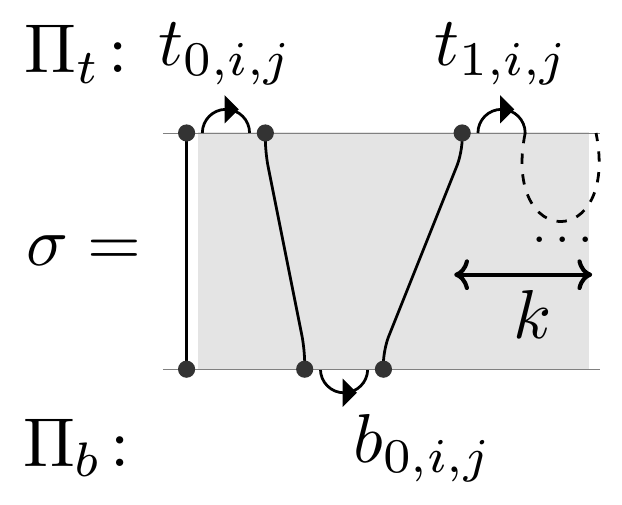}}\arrow{r}{L^{-k}}  & \raisebox{-32pt}{\includegraphics[scale=.5]{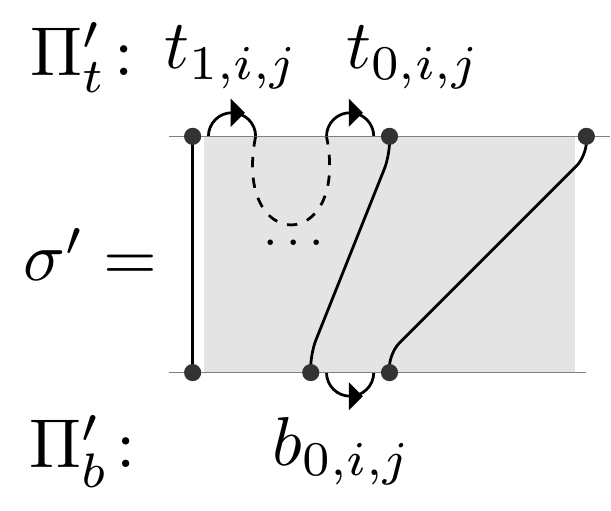}}\arrow{r}{R^{-1}B(S)R} &  \raisebox{-32pt}{\includegraphics[scale=.5]{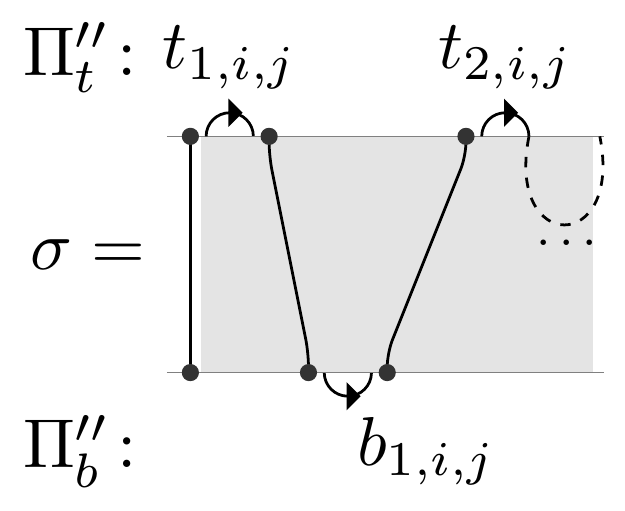}}\\[-.5cm]
\end{tikzcd}$
\end{center}
\caption[caption]{\label{fig_pf_monodromy_2} The case of the cycle 1-shift. We apply proposition \ref{pro_label_moves_well} on $\s'$ and $\s$.}
\end{figure}

Let $k=\alpha_2-2$ then $\s',(\Pi'_b,\Pi'_t)=L^{-k}(\s,(\Pi_b,\Pi_t))$ is a standard shift-irreducible permutation of type $X(r,i)$. Let $\alpha'_1,\ldots,\alpha'_i$ be the arc of the princicpal cycle, we have $\Pi'_t(\alpha'_1)=t_{1,i,j}, \Pi'_t(\alpha'_2)=t_{2,i,j}, \ldots, \Pi'_t(\alpha'_i)=t_{0,i,j}$. See figure \ref{fig_pf_monodromy_2} left and middle.

By the proof of proposition \ref{pro_labelling_first_step} case 1 and 2, $\tau=d(\s)$ and $\tau'=d(\s')$ are irreducible and have the same invariant, thus there is $S$ such that $S(\tau')=\tau$ by the classification theorem. Therefore by proposition \ref{pro_label_moves_well} let us define $\s,(\Pi'_b,\Pi'_t)=R^{-1}B(S)R(\s',(\Pi'_b,\Pi'_t))$, we have 
$\Pi''_t(\alpha_1)=t_{1,i,j}, \Pi''_t(\alpha_2)=t_{2,i,j}, \ldots, \Pi''_t(\alpha_i)=t_{0,i,j}$ as expected. See figure \ref{fig_pf_monodromy_2} right.

\paragraph*{\textbullet~ Let us demonstrate the statement on the cycle jump. In this case $\lam$ has at least two cycles of length $i$} Let $i\in \lam$, let $\s$ a standard shift-irreducible permutation of type $X(r,i)$ and let $(\Pi_b,\Pi_t)$ be a choosen consistent labelling of $\s$. Let $\alpha_1,\ldots,\alpha_i$ be the top arcs of the principal cycle, we have $\Pi_t(\alpha_1)=t_{0,i,j}$ and $\Pi_t(\alpha_2)=t_{1,i,j} \ldots \alpha_i=t_{i-1,i,j}$. Let $\alpha'$ be an arc of another cycle of length $i$, $\Pi_t(\alpha')=t_{c,i,j'}$ for some $c,j'$.

 We construct a loop $S$ of $\s$ such that $(\Pi''_b,\Pi''_t)=S(\Pi_b,\Pi_t)$ verifies $\Pi''_t(\alpha_1)=t_{c,i,j'}$ and $\Pi''_t(\alpha_2)=t_{c+1\mod i ,i,j}, \ldots,\Pi''_t(\alpha_i)=t_{c+i-1 \mod i,i,j}$

\begin{figure}[tb]
\begin{center}
$\begin{tikzcd}[%
    ,row sep = 3.5ex
    ,column sep = 10ex
    ,/tikz/column 1/.append style={anchor=base east}
    ,/tikz/column 2/.append style={anchor=base east}
    ,/tikz/column 3/.append style={anchor=base east}
    ]
 \raisebox{-32pt}{\includegraphics[scale=.5]{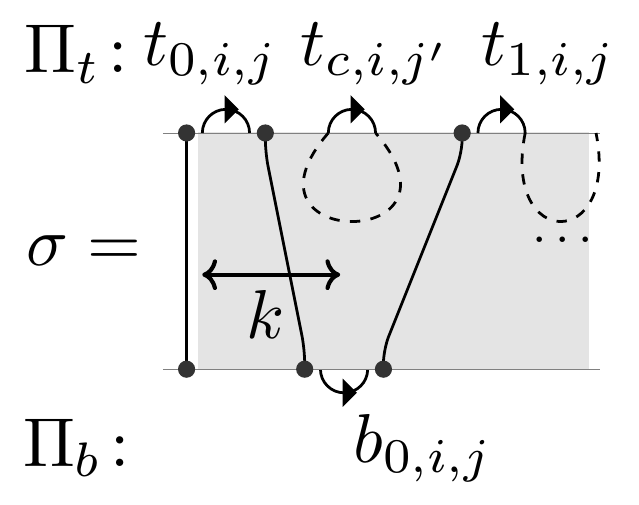}}\arrow{r}{L^{-k}}  & \raisebox{-32pt}{\includegraphics[scale=.5]{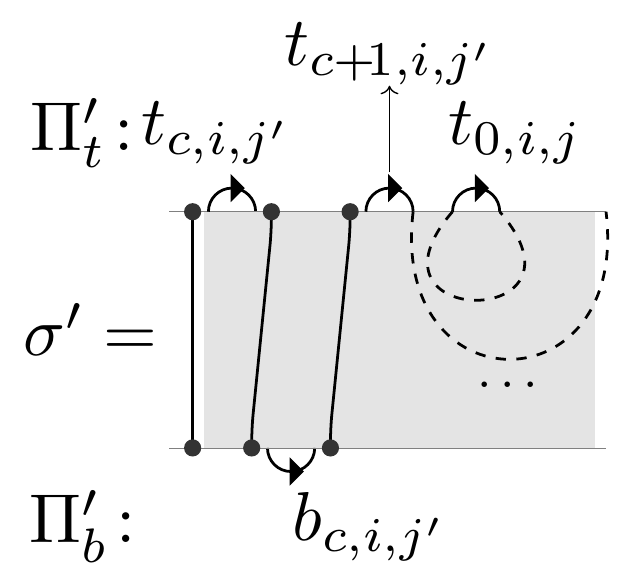}}\arrow{r}{R^{-1}B(S)R} &  \raisebox{-32pt}{\includegraphics[scale=.5]{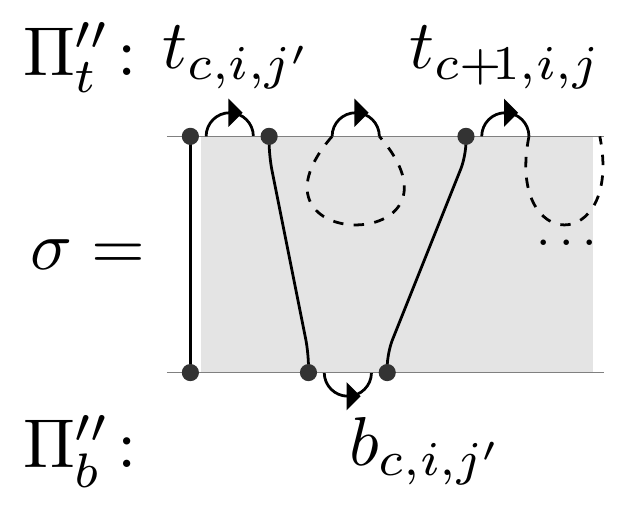}}\\[-.5cm]
\end{tikzcd}$
\end{center}
\caption[caption]{\label{fig_pf_monodromy_3} The case of the cycle jump. We apply proposition \ref{pro_label_moves_well} on $\s'$ and $\s$.}
\end{figure}

Let $k=\alpha'-2$ then $\s',(\Pi'_b,\Pi'_t)=L^{-k}(\s,(\Pi_b,\Pi_t))$ is a standard shift-irreducible permutation of type $X(r,i)$. let $\alpha'_1,\ldots,\alpha'_i$ be the arc of the princicpal cycle, we have $\Pi'_t(\alpha'_1)=t_{c,i,j'}, \Pi'_t(\alpha'_2)=t_{c+1\mod i,i,j'}, \ldots, \Pi'_t(\alpha'_i)=t_{c+i-1 \mod i,i,j'}$. See figure \ref{fig_pf_monodromy_3} left and middle.

By the proof of proposition \ref{pro_labelling_first_step} case 1 and 2, $\tau=d(\s)$ and $\tau'=d(\s')$ are irreducible and have the same invariant, thus there is $S$ such that $S(\tau')=\tau$ by the classification theorem. Therefore, by proposition \ref{pro_label_moves_well} let us define $\s,(\Pi'_b,\Pi'_t)=R^{-1}B(S)R(\s',(\Pi'_b,\Pi'_t))$, we have 
$\Pi''_t(\alpha_1)=t_{c,i,j'}, \Pi''_t(\alpha_2)=t_{c+1 \mod i,i,j'}, \ldots, \Pi''_t(\alpha_i)=t_{c+i-1 mod i,i,j'}$ as expected. See figure \ref{fig_pf_monodromy_3} right.

\paragraph*{\textbullet~ Let us demonstrate the statement on the cycle 2-shift. In this case $\lam$ has at exactly one cycle of even length, the cycle of length $i$.}Let $\s$ a standard shift-irreducible permutation of type $X(r,i)$ and let $(\Pi_b,\Pi_t)$ be a choosen consistent labelling of $\s$. Let $\alpha_1,\ldots,\alpha_i$ be the top arcs of the principal cycle, we have $\Pi''_t(\alpha_1)=t_{0,i,j}$ and $\Pi''_t(\alpha_2)=t_{1,i,j}, \ldots, \Pi''_t(\alpha_i)=t_{i-1,i,j}$. 

We construct a loop $S$ of $\s$ such that $(\Pi''_b,\Pi''_t)=S(\Pi_b,\Pi_t)$ verifies $\Pi_t(\alpha_1)=t_{2,i,j}$ and $\Pi_t(\alpha_2)=t_{3,i,j} \ldots \alpha_i=t_{1 \mod i,i,j}$.

\begin{figure}[tb]
\begin{center}
$\begin{tikzcd}[%
    ,row sep = 3.5ex
    ,column sep = 10ex
    ,/tikz/column 1/.append style={anchor=base east}
    ,/tikz/column 2/.append style={anchor=base east}
    ,/tikz/column 3/.append style={anchor=base east}
    ]
 \raisebox{-32pt}{\includegraphics[scale=.5]{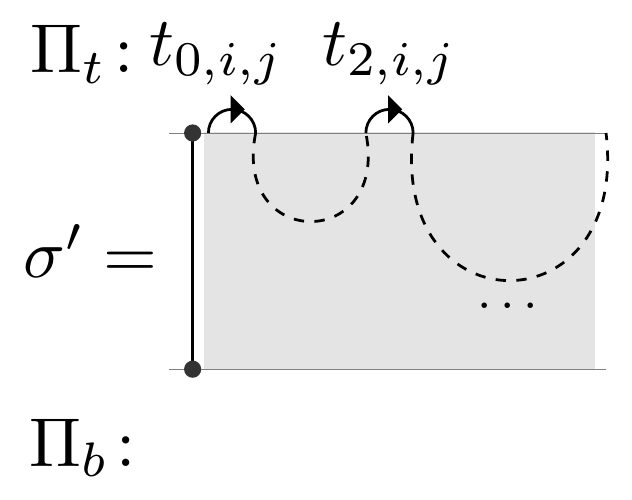}}\arrow{r}{L^{-k}}  & \raisebox{-32pt}{\includegraphics[scale=.5]{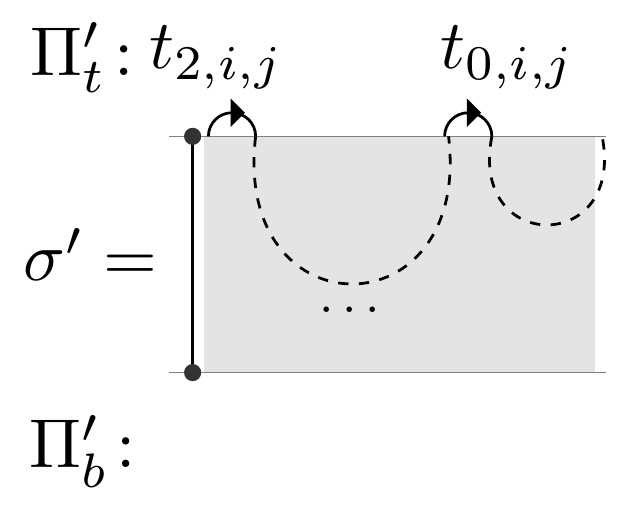}}\arrow{r}{R^{-1}B(S)R} &  \raisebox{-32pt}{\includegraphics[scale=.5]{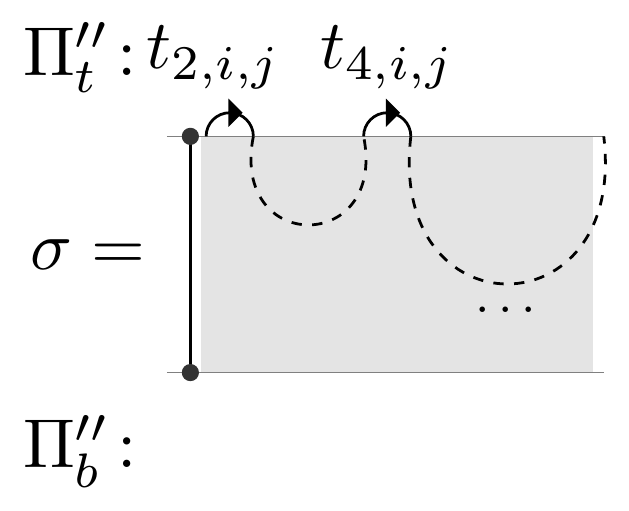}}\\[-.5cm]
\end{tikzcd}$
\end{center}
\caption[caption]{\label{fig_pf_monodromy_4} The case of the cycle 2-shift. We apply proposition \ref{pro_label_moves_well} on $\s'$ and $\s$.}
\end{figure}

Let $k=\alpha_3-2$ then $\s',(\Pi'_b,\Pi'_t)=L^{-k}(\s,(\Pi_b,\Pi_t))$ is a standard shift-irreducible permutation of type $X(r,i)$. Let $\alpha'_1,\ldots,\alpha'_i$ be the arc of the princicpal cycle, we have $\Pi'_t(\alpha'_1)=t_{2,i,j}, \Pi'_t(\alpha'_2)=t_{3,i,j}, \ldots, \Pi'_t(\alpha'_i)=t_{1,i,j}$. See figure \ref{fig_pf_monodromy_4} left and middle.

We must establish that $\tau=d(\s)$ and $\tau'=d(\s')$ are irreducible and have the same invariant. For the irreducibility and the cycle invariant, this follows from lemma \ref{lem_Xtype_s1} and the fact that $\s$ and $\s'$ are shift irreducible.

The sign invariant is slightly more complicated, we employ lemma \ref{lem_small_lemma} as we already did in case 4 of the proof of proposition \ref{pro_labelling_first_step}. In a few words,
define $\s''=\s'_{next}$ as in the lemma, then we have $\s''_{next}=\s$ thus $\Abar(\s)=-\Abar(\s'')=\Abar(\s')$.  

Thus there is $S$ such that $S(\tau')=\tau$ by the classification theorem. Then by proposition \ref{pro_label_moves_well} let $\s,(\Pi'_b,\Pi'_t)=R^{-1}B(S)R(\s',(\Pi'_b,\Pi'_t))$, we have 
$\Pi''_t(\alpha_1)=t_{3,i,j}, \Pi''_t(\alpha'_2)=t_{4,i,j}, \ldots, \Pi''_t(\alpha_i)=t_{1,i,j}$ as expected. See figure \ref{fig_pf_monodromy_4} right.
\qed

This complete the proof of theorem \ref{thm_2monodromy}.

\begin{remark}\label{rk_monodromy_thm_total}
Once the classification theorem is proven we can prove theorem \ref{thm_monodromy_total}.
Indeed If instead of using the induction hypothesis of theorem \ref{thm_2monodromy} we use the stronger induction hypothesis of that of theorem \ref{thm_monodromy_total} then the proof above (with some added technicallities) transforms into a proof of theorem \ref{thm_monodromy_total}.
\end{remark}

\appendix

\section{Exceptional classes}

In this appendix, when using a matrix representation of
configurations, it is useful to adopt the following notation:
The symbol $\epsilon$ denotes the $0 \times 0$ empty matrix.
The symbol
\raisebox{-2.5pt}{\includegraphics[width=4.8mm]{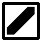}}
% $\includegraphics[width=4.8mm]{Imgs_matr_fig_matr_id.pdf}$
denotes a square block in a matrix (of any size $\geq 0$), filled with
an identity matrix.  A diagram, containing these special symbols and
the ordinary bullets used through the rest of the appendix, describes the
set of all configurations that could be obtained by specifying the
sizes of the identity blocks.
% the set $\id = \cup_{n \geq 0}\; \id_n$, empty space stand for
% zeroes, and bullets stand for ones. 
In such a syntax, we can write
equations of the like
% syntax:
% \raisebox{-(2n-1)mm}{\includegraphics[width=(4n).8mm]{Imgs_matr/(file nxn).pdf}}
\begin{align}
\id 
&:=
\raisebox{-1mm}{\includegraphics[width=4.8mm]{FigFol/FigureA2_fig_matr_id.pdf}}
=
\epsilon \cup
\raisebox{-3mm}{\includegraphics[width=8.8mm]{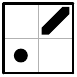}}
=
\epsilon \cup
\raisebox{5.8mm}{\includegraphics[width=8.8mm,
    angle=180]{FigFol/FigureA2_fig_matr_Bid.pdf}}
\ef;
% \\
% \intertext{and}
&
\id'
&:=
\raisebox{-7mm}{\includegraphics[width=16.8mm]{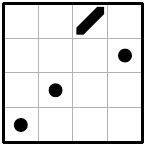}}
\ef.
\end{align}
The sets $\id$ and $\id'$ contain one element per size, $\id_n$ and
$\id'_n$, for $n \geq 0$ and $n \geq 3$ respectively.

The two exceptional classes $\Id_n$ and $\tree_n$ contain the
configurations $\id_n$ and $\id'_n$, respectively.

\noindent We have the following proposition:
\begin{proposition}\label{pro_excep_std}
The permutation $\s=\id_n$ (respectively $\s=\id'_n)$ is the only permutation of $\Id_n$ (respectively $\Id'_n$) with $\s(1)=1$ and $\s(2)=2$. 
\end{proposition}

\begin{figure}[tb!]
\begin{center}
\begin{align*}
X_{LL}^{(k)}
&=
\quad
\makebox[0pt][l]{\raisebox{-16mm}{\includegraphics[width=33mm]{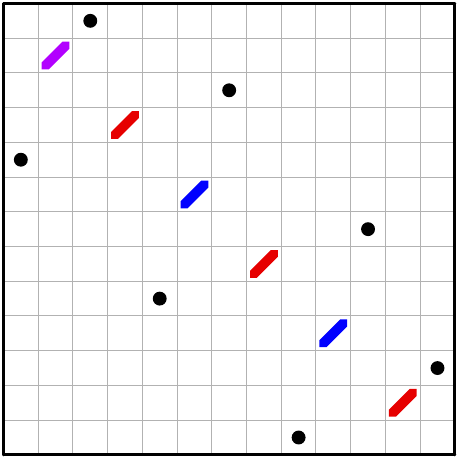}}}
\hspace{-16pt}
\raisebox{-30pt}{$\rotatebox{45}{$k \left\{\rule{0pt}{50pt}\right. $}$}
&
X_{RR}^{(k)}
&=
\makebox[0pt][l]{\raisebox{-16mm}{\includegraphics[width=33mm]{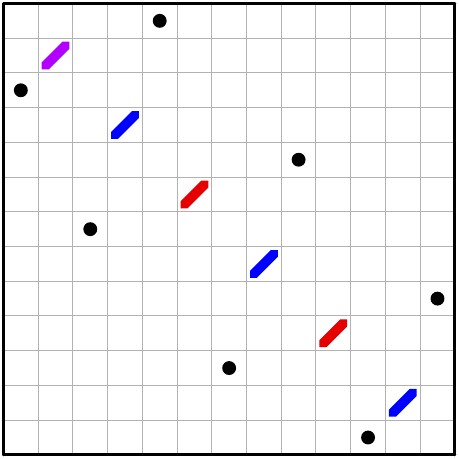}}}
\hspace{28pt}\raisebox{14pt}{$\rotatebox{45}{$\left. \rule{0pt}{50pt}\right\} k $}$}
\\
X_{RL}^{(k)}
&=
\makebox[0pt][l]{\raisebox{-16mm}{\includegraphics[width=38mm]{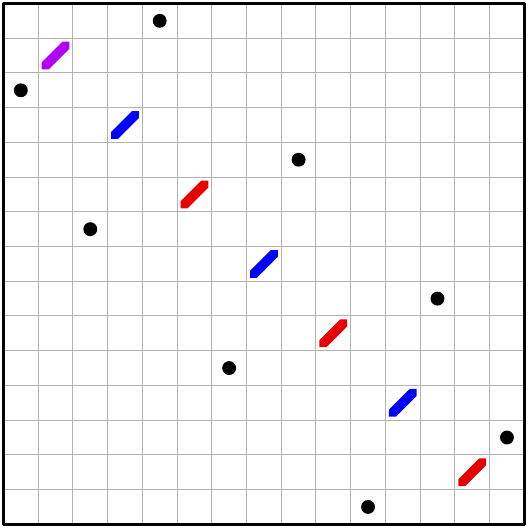}}}
\hspace{-2pt}
\raisebox{-30pt}{$\rotatebox{45}{$k \left\{\rule{0pt}{50pt}\right. $}$}
&
X_{LR}^{(k)}
&=
\makebox[0pt][l]{\raisebox{-16mm}{\includegraphics[width=38mm]{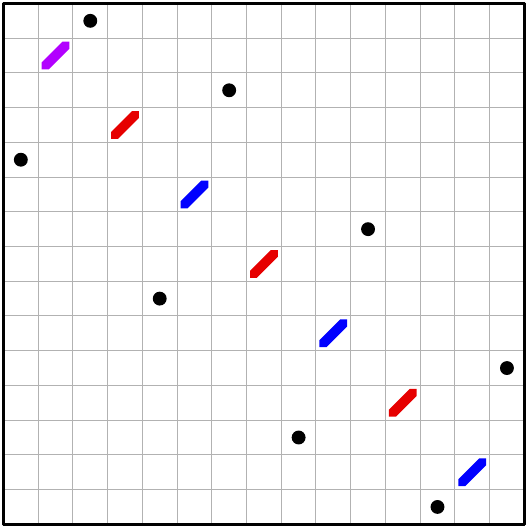}}}
\hspace{42pt}\raisebox{14pt}{$\rotatebox{45}{$\left. \rule{0pt}{50pt}\right\} k $}$}
\end{align*}
\end{center}
\caption{Description of the structure of configurations in
the identity class $\Id_n$, besides the configuration $\id_n$.
\label{fig.struct_id}}
\end{figure}

The structure of the classes $\Id_n$ is summarised by the following
relation:
\be
% \begin{split}
\Id := \bigcup_n \, \Id_n
=
\bigg(
\bigcup_{k \geq 1}
(X_{RL}^{k} \cup X_{LR}^{k} \cup X_{LL}^{k} \cup X_{RR}^{k})
\bigg)
\cup
\id
% \raisebox{-1mm}{\includegraphics[width=4.8mm]{FigFol/FigureA2_fig_matr_id.pdf}}
\ee
where the configurations $X_{\cdot \cdot}^{k}$ are defined as in
figure~\ref{fig.struct_id} (discard colours for the moment).\\

We can now prove the lemma \ref{lem.StFamNotmanyId_tobeproveninAppC} that we introduced in section 10.4 !!!!!!!!!.\\

\pfof{Lemma \ref{lem.StFamNotmanyId_tobeproveninAppC}}
This is equivalent to say that there are no pairs of
permutations $\s_1, \s_2 \in \Id_n$ which allow for a block
decomposition
\begin{align}
\s_1 &= \begin{array}{|c|}
\hline
\rule{5pt}{0pt} 
A 
\rule{5pt}{0pt} 
\\ 
\hline
B \\
\hline
\end{array}
&
\s_2 &= \begin{array}{|c|}
\hline
\rule{5pt}{0pt} 
B
\rule{5pt}{0pt} 
\\ 
\hline
A \\
\hline
\end{array}
\end{align}
If the block $A$ has $\ell$ rows, we say that $\s_2$ is the result of
shifting $\s_1$ by $\ell$.

Clearly, at the light of the structure of configurations that we have
presented (refer in particular to Figure \ref{fig.struct_id}), this
pattern is incompatible with $\s_1$ or $\s_2$ being $\id_n$ (as a
non-trivial shift produces a configuration which is not even
irreducible), so we have excluded the cases in which, still with
reference to the figure, we have only one violet block, and the number
of black points is at least 3, for $X^{(k)}_{LL}$ and $X^{(k)}_{RR}$,
and at least 4, for $X^{(k)}_{LR}$ and $X^{(k)}_{RL}$. Note that the
black points are the positions in the grid which are south-west or
north-east extremal (i.e., positions $(i,j) \in \s$ such that there is
no $(i',j') \in \s$ with $i'<i$ and $j'<j$, or the analogous statement
with $i'>i$ and $j'>j$). Let us call \emph{number of records},
$\rho(\s)$, this parameter. Thus we have that configurations in
$X^{(k)}_{LL}$ and $X^{(k)}_{RR}$ have $\rho=2k+1$, and configurations
in $X^{(k)}_{LR}$ and $X^{(k)}_{RL}$ have $\rho=2k+2$.

Now, if we perform a shift within one block of consecutive ascents, it
is easily seen, by investigation of the sub-configuration at the right
of the entry of the new configuration in the bottom-most row, or the
one at the left of the entry of the new configuration in the top-most
row, that the resulting structure is incompatible with the structure
of $\Id$. The same reasoning apply if we perform the shift at the
beginning/end of a non-empty diagonal block, which is not the one at
the bottom-right/top-left.  On the other side, if we perform a shift
in any other configuration, we have a new configuration in which
$\rho$ has strictly decreased. As $\s_2$ is a non-trivial shift of
$\s_1$, and $\s_1$ is a non-trivial shift of $\s_2$, we can thus
conclude.
\qed

%%%%%%%%%%%%%%%%%%%%%%%%%%%%%%%%%%%%%%%%%%%%%%%%%%%%%%%%
%
%\input{ppA1_rauzy1appends.tex}
%
%\input{ppA2_involutappends.tex}
%
%%%%%%%%%%%%%%%%%%%%%%%%%%%%%%%%%%%%%%%%%%%%%%%%%%%%%%%%

\end{document}